\documentclass{amsproc}
\usepackage{amsmath}
\usepackage{graphicx}
\usepackage{amsfonts}
\usepackage{amssymb}
\newtheorem{theorem}{Theorem}[section]
\newtheorem{lemma}[theorem]{Lemma}
\newtheorem{corollary}[theorem]{Corollary}
\newtheorem{proposition}[theorem]{Proposition}
\theoremstyle{definition}
\newtheorem{definition}[theorem]{Definition}
\newtheorem{example}[theorem]{Example}
\theoremstyle{remark}
\newtheorem{remark}[theorem]{Remark}

\newtheorem{thm}{Theorem}[section]

\newtheorem{rem}[thm]{Remark}
\theoremstyle{remark}
\newtheorem{note}{Note}

\numberwithin{equation}{section}

\begin{document}
\title[Spectral Invariants and Partitioned Manifolds]{Spectral Invariants of Operators of Dirac Type \\on Partitioned Manifolds}
\author{David Bleecker}
\address{David Bleecker\\
Department of Mathematics\\
University of Hawaii\\
Honolulu, HI 96822\\
USA}
\email{bleecker@math.hawaii.edu}
\author{Bernhelm Booss--Bavnbek}
\address{Bernhelm Booss--Bavnbek\\
Institut for matematik og fysik\\
Roskilde Universitetscenter\\
Postboks 260\\
DK-4000 Roskilde\\
Denmark}
\curraddr{}
\email{booss@ruc.dk}
\thanks{Printed \today}
\date{November 2002}

\begin{abstract}
We review the concepts of the \textit{index} of a Fredholm operator, the
\textit{spectral flow} of a curve of self-adjoint Fredholm operators, the
\textit{Maslov index} of a curve of Lagrangian subspaces in symplectic Hilbert
space, and the \textit{eta invariant} of operators of Dirac type on closed
manifolds and manifolds with boundary. We emphasize various (occasionally
overlooked) aspects of rigorous definitions and explain the quite different
stability properties. Moreover, we utilize the heat equation approach in
various settings and show how these topological and spectral invariants are
mutually related in the study of additivity and nonadditivity properties on
partitioned manifolds.
\end{abstract}\maketitle
\tableofcontents

\bigskip

%
%
%
%
%
%
%
%
%
%
%
%
%
%

%
%
%
%
%
%
%
%
%
%
%
%
%
%

%
%
%
%
%
%
%
%
%
%
%
%
%
%

%
%
%
%
%
%
%
%
%
%
%
%
%
%

%
%
%
%
%
%
%
%
%
%
%
%
%
%

%
%
%
%
%
%
%
%
%
%
%
%
%
%

%
%
%
%
%
%
%
%
%
%
%
%
%
%

%
%
%
%
%
%
%
%
%
%
%
%
%
%

%
%
%
%
%
%
%
%
%
%
%
%
%
%

\bigskip

\begin{center}
\textbf{Introduction}

\bigskip
\end{center}

For many decades now, workers in differential geometry and mathematical
physics have been increasingly concerned with differential operators (exterior
differentiation, connections, Laplacians, Dirac operators, etc.) associated to
underlying Riemannian or space--time manifolds. Of particular interest is the
interplay between the spectral decomposition of such operators and the
geometry/topology of the underlying manifold. This has become a large, diverse
field involving index theory, the distribution of eigenvalues, zero sets of
eigenfunctions, Green functions, heat and wave kernels, families of elliptic
operators and their determinants, canonical sections, etc.. Moreover,
Donaldson's analysis of moduli of solutions of the nonlinear Yang--Mills
equations and Seiberg--Witten theory have led to profound insights into the
classification of four--manifolds, which were not accessible by techniques
that are effective in higher dimensions.

We shall touch upon many of these topics, but we focus on three invariants
characterizing the asymmetry of the spectrum of operators of Dirac type: the
\textit{index} which gives the chiral asymmetry of the kernel or null space
(i.e., the difference in the number of independent left and right zero modes)
of a total Dirac operator; the \textit{spectral flow} of a curve of Dirac
operators that counts the net number of eigenvalues moving from the negative
half line over to the positive; and the \textit{eta invariant} which describes
the overall asymmetry of the spectrum of a Dirac operator. The three
invariants appear both for Dirac operators and curves of Dirac operators on a
closed manifold or on a smooth compact manifold with boundary subject to
suitable boundary conditions.

\smallskip

Index and spectral flow can be described in general functional analytic terms,
namely for bounded and unbounded, closed Fredholm operators and
\textit{curves} of bounded and unbounded self--adjoint Fredholm operators.
Correspondingly, stability properties of index and spectral flow and their
topological and geometric meaning are relatively well understood. There is,
however, not an easily identifiable operator class for eta invariant, and its
behavior under perturbations is rather delicate.

\smallskip

Most significant differences between the behavior of these three invariants
are met when we address splitting properties on partitioned manifolds (with
product metrics assumed near the partitioning hypersurface). It is well known
that the index can be described in local terms. This is one aspect of the
Atiyah--Singer Index Theorem. So, splitting formulas for the index are
relatively easily obtainable, once one specifies and understands boundary and
transmission conditions for the gluing of the parts of the underlying
manifold. In this process, precise additivity is obtained only for a few
invariants, Euler characteristic and signature, which actually can be
characterized by cutting and pasting invariance. In general, an error term
appears which can be expressed either by spectral flow or, equivalently, by
the index of a boundary value problem on a cylinder over the separating
hypersurface or, alternatively, by the index of a suitable Fredholm pair.

Surprisingly, simple splitting formulas can be obtained also for the spectral
flow and eta invariant. This is particularly surprising for the eta invariant
where we only have a local formula for the first derivative. Now the integer
error terms are expressed by the Maslov index for curves of Cauchy data
spaces. These formulas relate the symmetric category of self--adjoint Dirac
operators over closed partitioned manifolds (and self--adjoint boundary value
problems over compact manifolds with boundary) to the symplectic analysis of
Lagrangian subspaces (the Cauchy data spaces).

One message of this review is, that index, spectral flow, and eta invariant,
in spite of their quite different appearance, share various features which
become most visible on partitioned manifolds. Roughly speaking, one reason for
that is that the spectral flow of a path of Dirac type operators (say, on a
closed manifold) with unitarily equivalent ends $A_{1}=gA_{0}g^{-1}$ equals
the index of the induced suspension operator $\partial_{t}+A_{t}$ on the
underlying mapping torus (see \cite[Theorems 17.3, 17.17, and Proposition
25.1]{BoWo93}). Another reason is that the integer part of the derivative of
the eta invariant along a path of Dirac operators or boundary problems can be
expressed by the spectral flow (see, e.g., \cite{DoWo91}, \cite[p. 39]%
{LeWo96}, and \cite[Section 3]{KiLe00}).

We also summarize recent discussion on additivity and non--additivity of the
zeta regularized determinant. It should be mentioned that the various
splitting formulas all depend decisively on the well--established unique
continuation property (UCP) for operators of Dirac type. We give a full proof
of weak UCP below in Section \ref{WeakUCP}. Moreover (modulo some
technicalities for the computation of the index density form, to appear in our
forthcoming book \cite{BlBo03}), we provide proofs of the Atiyah--Singer and
Atiyah--Patodi--Singer Index Theorems in important special cases using heat
equation methods. That we explain concrete calculations only for the case of
the Atiyah--Patodi--Singer (spectral -- ``APS'') boundary condition is no big
loss of generality since each admissible boundary condition for an operator of
Dirac type $\mathcal{D}$ with tangential part $\mathcal{B}$ can be written as
the APS projection of a perturbed operator $\mathcal{B}^{\prime}$ (see below
Lemma \ref{l:grubb}, following a recent result, Gerd Grubb \cite{Gr02}).

\bigskip

The second author wants to thank K.P. Wojciechowski (Indiana University --
Purdue University Indianapolis) for many discussions about the subject(s) of
this review and the Erwin Schr\"{o}dinger International Institute for
Mathematical Physics at Vienna for generous hospitality during the finalizing
phase of this review. We both thank the referees for their corrections,
thoughtful comments, and helpful suggestions which led to many improvements.
They clearly went beyond the call of duty, and we are in their debt.

\bigskip

\section{Basic Notations and Results\label{s:basic}}

\subsection{Index of Fredholm Operators and Spectral Flow of Curves of
Self--Adjoint Fredholm Operators}

\quad

\smallskip

%

\subsubsection{Notation}

\label{sss:notations} Let $H$ be a separable complex Hilbert space. First let
us introduce some notation for various spaces of operators in $H$:
\[%
\begin{array}
[c]{rl}%
\mathcal{C}(H):= & \text{closed, densely defined operators on $H$},\\
\mathcal{B}(H):= & \text{bounded linear operators $H\rightarrow H$},\\
\mathcal{U}(H):= & \text{unitary operators $H\rightarrow H$},\\
\mathcal{K}(H):= & \text{compact linear operators $H\rightarrow H$},\\
\mathcal{F}(H):= & \text{bounded Fredholm operators $H\rightarrow H$},\\
\mathcal{CF}(H):= & \text{closed, densely defined Fredholm operators on $H$}.
\end{array}
\]

If no confusion is possible we will omit ``$(H)$'' and write $\mathcal{C},$
$\mathcal{B},\mathcal{K},$ etc.. By $\mathcal{C}^{\operatorname{sa}%
},\mathcal{B}^{\operatorname{sa}}$ etc., we denote the set of self--adjoint
elements in $\mathcal{C}, $ $\mathcal{B},$ etc..

\smallskip

\subsubsection{Operators With Index -- Fredholm Operators}

\label{sss:fredholm} The topology of the operator spaces $\mathcal{U}(H)$,
$\mathcal{F}(H)$, and $\mathcal{F}^{\operatorname{sa}}(H)$ is quite well
understood. The key results are

\begin{enumerate}
\item  the \textit{Kuiper Theorem} which states that $\mathcal{U}(H)$ is contractible;

\item  the \textit{Atiyah--J\"{a}nich Theorem} which states that
$\mathcal{F}(H)$ is a classifying space for the functor $K$. Explicitly, the
construction of the index bundle of a continuous family of bounded Fredholm
operators parametrized over a compact topological space $X$ yields a
homomorphism of semigroups, namely $\operatorname{index}:\left[
X,\mathcal{F}\right]  \rightarrow K\left(  X\right)  .$ In particular, the
index is a homotopy invariant, and it provides a one--to--one correspondence
of the connected components $\left[  point,\mathcal{F}\right]  $ of
$\mathcal{F}$ with $K\left(  point\right)  =\mathbb{Z}.$

\item  the corresponding observation that $\mathcal{F}^{\operatorname{sa}}(H)$
consists of three connected components, the contractible subsets
$\mathcal{F}^{\operatorname{sa}}(H)_{+}$ and $\mathcal{F}^{\operatorname{sa}%
}(H)_{-}$ of essentially positive, respectively essentially negative, Fredholm
operators, and the topologically nontrivial component $\mathcal{F}_{\ast
}^{\operatorname{sa}}(H)$ which is a classifying space for the functor
$K^{-1}$\thinspace. In particular, the spectral flow gives an isomorphism of
the fundamental group $\pi_{1}\left(  \mathcal{F}_{\ast}^{\operatorname{sa}%
}(H)\right)  $ onto the integers.
\end{enumerate}

Full proofs can be found in \cite{BlBo03} for items 1 and 2, and in
\cite{AtSi69}, \cite{BoWo93}, and \cite{Ph96} for item 3.

So much for the bounded case. Of course, the Dirac operators of interest
to\ us are not bounded in $L^{2}$. On closed manifolds, however, they can be
considered as bounded operators from the first Sobolev space $H^{1}$ into
$L^{2}$, and by identifying these two Hilbert spaces we can consider a Dirac
operator as a bounded operator of a Hilbert space in itself. The same
philosophy can be applied to Dirac operators on compact manifolds with
boundary when we consider the domain not as dense subspace in $L^{2}$ but as
Hilbert space (with the graph inner product) and then identify.

Strictly speaking, the concept of unbounded operator is dispensable here.
However, for varying boundary conditions it is necessary to keep the
distinction and to follow the variation of the domain as a variation of
subspaces in $L^{2}$.

Following \cite{CorLab}, we generalize the concept of \emph{Fredholm
operators} to the unbounded case.

\smallskip

\begin{definition}
\label{d:fredholm} Let $H$ be a complex separable Hilbert space. A linear (not
necessarily bounded) operator $F$ with domain $\operatorname{Dom}(F)$,
null--space $\operatorname{Ker}(F)$, and range $\operatorname{Im}(F)$ is
called \emph{Fredholm} if the following conditions are satisfied.

\begin{enumerate}
\item [(i)]$\operatorname{Dom}(F)$ is dense in $H$.

\item[(ii)] $F$ is closed.

\item[(iii)] The range $\operatorname{Im}(F)$ of $F$ is a closed subspace of
$H$.

\item[(iv)] Both $\dim\operatorname{Ker}(F)$ and $\operatorname{codim}%
\operatorname{Im}(F)=\dim\operatorname{Im}(F)^{\perp}$ are finite. The
difference of the dimensions is called $\operatorname{index}(F)$.
\end{enumerate}
\end{definition}

\smallskip

So, a closed operator $F$ is characterized as a Fredholm operator by the same
properties as in the bounded case. Moreover, as in the bounded case, $F$ is
Fredholm if and only if $F^{\ast}$ is Fredholm (proving the closedness of
$\operatorname{Im}(F^{\ast})$ is delicate: see \cite[Lemma 1.4]{CorLab}) and
(clearly) we have
\[
\operatorname{index}F=\dim\operatorname{Ker}F-\dim\operatorname{Ker}F^{\ast
}=-\operatorname{index}F^{\ast}.
\]
In particular, $\operatorname{index}F=0$ in case $F$ is self--adjoint.\bigskip

The composition of (not necessarily bounded) Fredholm operators yields again a
Fredholm operator. More precisely, we have the following composition rule. For
the proof, which is considerably more involved than that in the bounded case,
we refer to \cite{GoKr57}, \cite[Lemma 2.3 and Theorem 2.1]{CorLab}.

\begin{theorem}
\label{t:composition} \emph{(Gohberg, Krein)} If $F$ and $G$ are (not
necessarily bounded) Fredholm operators then their product $GF$ is densely
defined with
\[
\operatorname{Dom}(GF)=\operatorname{Dom}(F)\cap F^{-1}\operatorname{Dom}(G)
\]
and is a Fredholm operator. Moreover,
\[
\operatorname{index}GF=\operatorname{index}F+\operatorname{index}G\,.
\]
\end{theorem}

\subsubsection{Metrics on the Space of Closed Operators\label{sss:metrics}}

For $S,\,T\in\mathcal{C}(H)$ the orthogonal projections $P_{\frak{G}(S)}$,
$P_{\frak{G}(T)}$ onto the graphs of $S,\,T$ in $H\oplus H$ are bounded
operators and
\[
\gamma(S,T):=\Vert P_{\frak{G}(T)}-P_{\frak{G}(S)}\Vert
\]
defines a metric for $\mathcal{C}(H)$, the \textit{projection metric}.

It is also called the \textit{gap metric} and it is (uniformly) equivalent
with the metric given by measuring the distance between the (closed) graphs,
namely
\[
d(\frak{G}(S),\frak{G}(T)):=\sup_{\left\{  x\in\frak{G}(S):\left\|  x\right\|
=1\right\}  }d\left(  x,\frak{G}(T)\right)  +\sup_{\left\{  x\in
\frak{G}(T):\left\|  x\right\|  =1\right\}  }d\left(  x,\frak{G}(S)\right)  .
\]
For details and the proof of the following Lemma and Theorem, we refer to
\cite[Section 3]{CorLab}.

\begin{lemma}
\label{l:projection} For $T\in\mathcal{C}(H)$ the orthogonal projection onto
the graph of $T$ in $H\oplus H$ can be written (where\ \ $R_{T}%
:=(\operatorname{I}+T^{\ast}T)^{-1}$) as
\[
P_{\frak{G}(T)}=
\begin{pmatrix}
R_{T} & R_{T}T^{\ast}\\
TR_{T} & TR_{T}T^{\ast}%
\end{pmatrix}
\,=
\begin{pmatrix}
R_{T} & T^{\ast}R_{T^{\ast}}\\
TR_{T} & TT^{\ast}R_{T^{\ast}}%
\end{pmatrix}
=
\begin{pmatrix}
R_{T} & T^{\ast}R_{T^{\ast}}\\
TR_{T} & \operatorname{I}-R_{T^{\ast}}%
\end{pmatrix}
.
\]
\end{lemma}

\begin{theorem}
\label{t:maincordes}\emph{(Cordes, Labrousse)} a) The space $\mathcal{B}(H)$
of bounded operators on $H$ is dense in the space $\mathcal{C}(H)$ of all
closed operators in $H$. The topology induced by the projection ($\cong$ gap)
metric on $\mathcal{B}(H)$ is equivalent to that given by the operator norm.

\noindent b) Let $\mathcal{CF}(H)$ denote the space of closed (not necessarily
bounded) Fredholm operators. Then the index is constant on the connected
components of $\mathcal{CF}(H)$ and yields a bijection between the integers
and the connected components.
\end{theorem}

\begin{example}
\label{ex:fuglede} Let $H$ be a Hilbert space with $e_{1},e_{2},\dots$ a
complete orthonormal system (i.e., an orthonormal basis for $H$). Consider the
multiplication operator $M_{\operatorname{id}}$, given by the domain
\[
D:=\operatorname{Dom}(M_{\operatorname{id}}):=\left\{  \sum\nolimits_{j=1}%
^{\infty}c_{j}e_{j}\ |\ \sum\nolimits_{j=1}^{\infty}j^{2}|c_{j}|^{2}%
<+\infty\right\}
\]
and the operation
\[
D\ni u=\sum c_{j}e_{j}\quad\mapsto\quad M_{\operatorname{id}}(u):=\sum{j}%
c_{j}e_{j}\,.
\]
It is a densely defined closed operator which is injective and surjective. Let
$P_{n}$ denote the orthogonal projection of $H$ onto the linear span of the
$n$-th orthonormal basis element $e_{n}$\thinspace. Clearly the sequence
$(P_{n})$ does not converge in $\mathcal{B}(H)$ in the operator norm. However,
the sequence $M_{\operatorname{id}}-2nP_{n}$ of self--adjoint Fredholm
operators converges in $\mathcal{C}(H)$ with the projection metric to
$M_{\operatorname{id}}$.

This can be seen by the following argument: On the subset of self--adjoint
(not necessarily bounded) operators in the space $\mathcal{C}%
^{\operatorname{sa}}(H)$ the projection metric is uniformly equivalent to the
metric $\gamma$ given by
\[
\gamma(A_{1},A_{2}):=\Vert(A_{1}+i)^{-1}-(A_{2}+i)^{-1}\Vert\,,
\]
(see below Theorem \ref{S1.1}). Then for $T_{n}:=M_{\operatorname{id}}%
-2nP_{n}$,
\[
\left\|  \left(  T_{n}+i\right)  ^{-1}-\left(  M_{\operatorname{id}}+i\right)
^{-1}\right\|  =\left\|  \left(  -n+i\right)  ^{-1}e_{n}-\left(  n+i\right)
^{-1}e_{n}\right\|  =\tfrac{2n}{n^{2}+1}\rightarrow0.
\]
\end{example}

\begin{remark}
\label{r:cordes} The results by Heinz Cordes and Jean--Philippe Labrousse may
appear to be rather counter--intuitive. For (a), it is worth mentioning that
the operator--norm distance and the projection metric on the set of bounded
operators are equivalent but not uniformly equivalent since the operator norm
is complete while the projection metric is not complete on the set of bounded
operators. Actually, this is the point of the first part of (a); see also the
preceding example.

Assertion (b) says two things: (i) that the index is a homotopy invariant,
i.e., two Fredholm operators have the same index if they can be connected by a
continuous curve in $\mathcal{CF}(H)$; (ii) that two Fredholm operators having
the same index always can be connected by a continuous curve in $\mathcal{CF}%
(H)$. Note that the topological results are not as far reaching as for bounded
Fredholm operators.\bigskip
\end{remark}

\subsubsection{Self--Adjoint Fredholm Operators and Spectral
Flow\label{sss:self-adjoint fred ops}}

We investigate the topology of the subspace of self--adjoint (not necessarily
bounded) Fredholm operators. Many users of the notion of spectral flow feel
that the definition and basic properties are too trivial to bother with.
However, there are some difficulties both with extending the definition of
spectral flow from loops to paths and from curves of bounded self--adjoint
Fredholm operators to curves of not necessarily bounded self--adjoint Fredholm operators.

To overcome the second difficulty, the usual way is to apply the \textit{Riesz
transformation} which yields a bijection
\begin{equation}%
\begin{array}
[c]{l}%
\mathcal{R}:\mathcal{C}^{\operatorname{sa}}\longrightarrow\{S\in
\mathcal{B}^{\operatorname{sa}}\mid\Vert S\Vert\leq1\text{ and $S\pm I$ both
injective}\},\text{ where}\\
\mathcal{R}\left(  T\right)  :=T(I+T^{2})^{-1/2}%
\end{array}
\label{G1.1}%
\end{equation}

In \cite{BoFu98} the following theorem was proved:

\begin{theorem}
\label{t:nest} Let $S$ be a self--adjoint operator with compact resolvent in a
real separable Hilbert space $\mathcal{H}$ and let $C$ be a bounded
self--adjoint operator. Then the sum $S+C$ also has compact resolvent and is a
closed Fredholm operator. We have
\[
\Vert\mathcal{R}({S+C})-\mathcal{R}(S)\Vert\leq c\Vert C\Vert\,,
\]
where the constant $c$ does not depend on $S$ or on $C$.
\end{theorem}

The preceding theorem is applied in the following form:

\begin{corollary}
Curves of self--adjoint (unbounded) Fredholm operators in a separable real
Hilbert space of the form $\{{S}+C_{t}\}_{t\in I}$ are mapped into continuous
curves in $\mathcal{F}^{\operatorname{sa}}$ by the transformation
$\mathcal{R}$ when ${S}$ is a self--adjoint operator with compact resolvent
and $\{C_{t}\}_{t\in I}$ is a continuous curve of bounded self--adjoint operators.\bigskip
\end{corollary}

\begin{remark}
\label{r:gap}Define
\[
\mathcal{T}_{{S}}:\mathcal{B}^{\operatorname{sa}}\longrightarrow
\mathcal{CF}^{\operatorname{sa}}\text{ by \ }\mathcal{T}_{{S}}\left(
C\right)  :={S}+C
\]
This is translation by ${S}$, mapping bounded self--adjoint operators on
$\mathcal{H}$ into self--adjoint Fredholm operators in $\mathcal{H}$. On
$\mathcal{CF}^{\operatorname{sa}}$, the gap topology is defined by the metric
\[
g(A_{1},A_{2}):=\sqrt{\Vert R_{A_{1}}-R_{A_{2}}\Vert^{2}+\Vert A_{1}R_{A_{1}%
}-A_{2}R_{A_{2}}\Vert^{2}}\,,
\]
where $R_{A_{.}}:=(\operatorname{I}+A_{.}^{2})^{-1}$ as before (see Cordes and
Labrousse, \cite{CorLab} and also Kato, \cite{Kat}). Theorem \ref{t:nest} says
that the composition $\mathcal{R}\circ\mathcal{T}_{{S}}$ is continuous.
\[%
\begin{array}
[c]{rll}%
\mathcal{B}^{\operatorname{sa}} & \overset{\mathcal{T}_{{S}}}{\longrightarrow}%
& \mathcal{CF}^{\operatorname{sa}}\\
& \!\!\!\!\!\!\underset{\mathcal{R}\circ\mathcal{T}_{{S}}}{\searrow} &
{\LARGE \downarrow}_{^{%
\genfrac{}{}{0pt}{}{\mathcal{R}}{{}}%
}}\\
&  & \mathcal{F}^{\operatorname{sa}}%
\end{array}
\]
\end{remark}

Further, we can prove that the translation operator $\mathcal{T}_{S}$ is a
continuous operator from $\mathcal{B}^{\operatorname{sa}}$ onto the subspace
$\mathcal{B}^{\operatorname{sa}}+S\subset\mathcal{CF}^{\operatorname{sa}}$.\bigskip

The preceding arguments permit to treat continuous curves of Dirac operators
in the same way as continuous curves of self--adjoint bounded Fredholm
operators under the precondition that the domain is fixed and the perturbation
is only by bounded self--adjoint operators. That precondition is satisfied
when we have a curve of Dirac operators on a closed manifold which differ only
by the underlying connection. It is also satisfied for curves of Dirac
operators on a manifold with boundary as long the perturbation is bounded. In
particular, this demands that the domain remains fixed.

\smallskip

A closer look at Example \ref{ex:fuglede} shows that the preceding argument
cannot be generalized: The sequence of the Riesz transforms $\mathcal{R}%
(T_{n})$ of $T_{n}:=M_{\operatorname{id}}-2nP_{n}$ does not converge, since
$\mathcal{R}(T_{n})$ is not a Cauchy sequence:
\begin{align*}
&  \left\|  \mathcal{R}(T_{n})-\mathcal{R}\left(  T_{n+1}\right)  \right\|
\geq\Vert\mathcal{R}(T_{n})e_{n}-\mathcal{R}\left(  M_{\operatorname{id}%
}\right)  e_{n}\Vert\\
&  =\left\|  \tfrac{n-2n}{\sqrt{1+\left(  n-2n\right)  ^{2}}}e_{n}-\tfrac
{n}{\sqrt{1+n^{2}}}e_{n}\right\|  =\tfrac{2n}{\sqrt{1+n^{2}}}\rightarrow
2\text{\ \ as\ $\;n\rightarrow\infty$}.
\end{align*}
In particular, the Riesz transformation is not continuous on the whole space
$\mathcal{C}^{\operatorname{sa}}$, nor on $\mathcal{CF}^{\operatorname{sa}}$,
neither on the whole space of self--adjoint operators with compact resolvent.
Other methods are needed for working with varying domains.\bigskip

Of course, we can define a different metric in $\mathcal{C}^{\operatorname{sa}%
}$\thinspace; e.g. the metric which makes the Riesz transformation a
homeomorphism. That approach was chosen by L. Nicolaescu in \cite{Ni00}. He
shows that quite a large class of naturally arising curves of Dirac operators
with varying domain are continuous under this `Riesz metric', as opposed to
the aforementioned example.

\medskip

Here, we choose a different approach and adopt the gap metric. Note that
continuity in the gap metric is much easier to establish than continuity in
the Riesz metric. We follow \cite{BoLePh01} where the proofs can be found. A
feature of this approach is the use of the Cayley Transform:

\smallskip

\begin{theorem}
\label{S1.1}
%
%
%
%
%
%
%
%
%
%
%
%
(a) On $\mathcal{C}^{\operatorname{sa}}$ the gap metric is (uniformly)
equivalent to the metric $\gamma$ given by
\[
\gamma(T_{1},T_{2})=\Vert(T_{1}+i)^{-1}-(T_{2}+i)^{-1}\Vert.
\]
\noindent(b) Let $\kappa:\mathbb{R}\rightarrow S^{1}\setminus\{1\},$
$x\mapsto\frac{x-i}{x+i}$ denote the Cayley transform. Then $\kappa$ induces a
homeomorphism
\begin{equation}%
\begin{split}
\mathbf{\kappa}  &  :\mathcal{C}^{\operatorname{sa}}(H)\longrightarrow
\{U\in\mathcal{U}(H)\mid U-I\text{ is injective }\}=:\mathcal{U}%
_{\operatorname{inj}}\\
&  T\mapsto\mathbf{\kappa}(T)=(T-i)(T+i)^{-1}.
\end{split}
\label{G1.2}%
\end{equation}
More precisely, the gap metric is (uniformly) equivalent to the metric
$\tilde{\delta}$ defined by $\tilde{\delta}(T_{1},T_{2})=\Vert\mathbf{\kappa
}(T_{1})-\mathbf{\kappa}(T_{2})\Vert$.
\end{theorem}

\smallskip

We note some immediate consequences of the Cayley picture:

\smallskip

\begin{corollary}
\label{c:fred-open} (a) With respect to the gap metric the set $\mathcal{B}%
^{\operatorname{sa}}(H)$ is dense in $\mathcal{C}^{\operatorname{sa}}%
(H)$.\newline \noindent(b) For $\lambda\in\mathbb{R}$ the sets
\[
\{T\in\mathcal{C}^{\operatorname{sa}}(H)\mid\lambda\notin\operatorname{spec}%
T\}\qquad\text{and}\qquad\{T\in\mathcal{C}^{\operatorname{sa}}(H)\mid
\lambda\not \in\operatorname{spec}_{\operatorname{ess}}T\}
\]
are open in the gap topology.\newline \noindent(c) The set $\mathcal{CF}%
^{\operatorname{sa}}=\{T\in\mathcal{C}^{\operatorname{sa}}\mid0\notin
\operatorname{spec}_{\operatorname{ess}}T\}=\mathbf{\kappa}^{-1}%
(_{\mathcal{F}}\mathcal{U})$, where $_{\mathcal{F}}\mathcal{U}:=\{U\in
\mathcal{U}\mid-1\notin\operatorname{spec}_{\operatorname{ess}}U\}=\{U\in
\mathcal{U}\mid U+I$ \ Fredholm operator$\}$, of (not necessarily bounded)
self--adjoint Fredholm operators is open in $\mathcal{C}^{\operatorname{sa}}$.
\end{corollary}

\smallskip

The preceding Corollary implies that the set $\mathcal{F}^{\operatorname{sa}}$
is dense in $\mathcal{CF}^{\operatorname{sa}}$ with respect to the gap
metric.\newline Contrary to the bounded case, we have the following somewhat
surprising result in the unbounded case: In particular, it shows that not
every gap continuous path in $\mathcal{CF}^{\operatorname{sa}}$ with endpoints
in $\mathcal{F}^{\operatorname{sa}}$ can be continuously deformed into an
operator norm continuous path in $\mathcal{F}^{\operatorname{sa}}$, in spite
of the density of $\mathcal{F}^{\operatorname{sa}}$ in $\mathcal{CF}%
^{\operatorname{sa}}$.

\smallskip

\begin{theorem}
\label{t:connected} (a) $\mathcal{CF}^{\operatorname{sa}}$is path connected
with respect to the gap metric.

\noindent(b) Moreover, its Cayley image
\[
{_{\mathcal{F}}}\mathcal{U}_{\mathrm{\operatorname{inj}}}:=\{U\in
\mathcal{U}\mid U+I\text{ Fredholm and $U-I$ injective}\}=\mathbf{\kappa
}(\mathcal{CF}^{\operatorname{sa}})
\]
is dense in ${_{\mathcal{F}}}\mathcal{U}$.
\end{theorem}

\smallskip

\begin{proof}
(a) Once again we look at the Cayley transform picture. Note that so far we
have introduced three different subsets of unitary operators
\begin{align}
\mathcal{U}_{\mathrm{\operatorname{inj}}}  &  :=\left\{  U\in\mathcal{U}\mid
U-I\text{ \ injective}\right\}  =\mathbf{\kappa}(\mathcal{C}%
^{\operatorname{sa}})\\
{_{\mathcal{F}}}\mathcal{U}  &  :=\left\{  U\in\mathcal{U}\mid U+I\text{
\ Fredholm}\right\}  ,\text{ and}\\
{_{\mathcal{F}}}\mathcal{U}_{\mathrm{\operatorname{inj}}}  &  :={_{\mathcal{F}%
}}\mathcal{U\cap U}_{\mathrm{\operatorname{inj}}}=\mathbf{\kappa}%
(\mathcal{CF}^{\operatorname{sa}}).
\end{align}
Let $U\in{_{\mathcal{F}}}\mathcal{U}_{\mathrm{\operatorname{inj}}}$. Then $H$
is the direct sum of the spectral subspaces $H_{\pm}$ of $U$ corresponding to
$[0,\pi)$ and $[\pi,2\pi]$ respectively and we may decompose $U=U_{+}\oplus
U_{-}$. More precisely, we have
\[
\operatorname{spec}(U_{+})\subset\{e^{it}\mid t\in\lbrack0,\pi)\}\text{ and
}\operatorname{spec}(U_{-})\subset\{e^{it}\mid t\in\lbrack\pi,2\pi]\}\,.
\]
Note that there is no intersection between the spectral spaces in the
endpoints: if $-1$ belongs to $\operatorname{spec}(U)$, it is an isolated
eigenvalue by our assumption and hence belongs only to $\operatorname{spec}%
(U_{-})$; if $1$ belongs to $\operatorname{spec}(U)$, it can belong both to
$\operatorname{spec}(U_{+})$ and $\operatorname{spec}(U_{-})$, but in any
case, it does not contribute to the decomposition of $U$ since, by our
assumption, $1$ is not an eigenvalue at all.

By spectral deformation (``squeezing the spectrum down to $+i$ and $-i$'') we
contract $U_{+}$ to $iI_{+}$ and $U_{-}$ to $-iI_{-}$\thinspace, where
$I_{\pm}$ denotes the identity on $H_{\pm}$\thinspace. We do this on the upper
half arc and the lower half arc, respectively, in such a way that 1 does not
become an eigenvalue under the course of the deformation: actually it will no
longer belong to the spectrum; neither will $-1$ belong to the spectrum. That
is, we have connected $U$ and $iI_{+}\oplus-iI_{-}$ within $\mathbf{\kappa
}\left(  \mathcal{CF}^{\operatorname{sa}}\right)  $.%
\begin{figure}
[ptb]
\begin{center}
\includegraphics[
height=4.7928in,
width=3.3209in
]%
{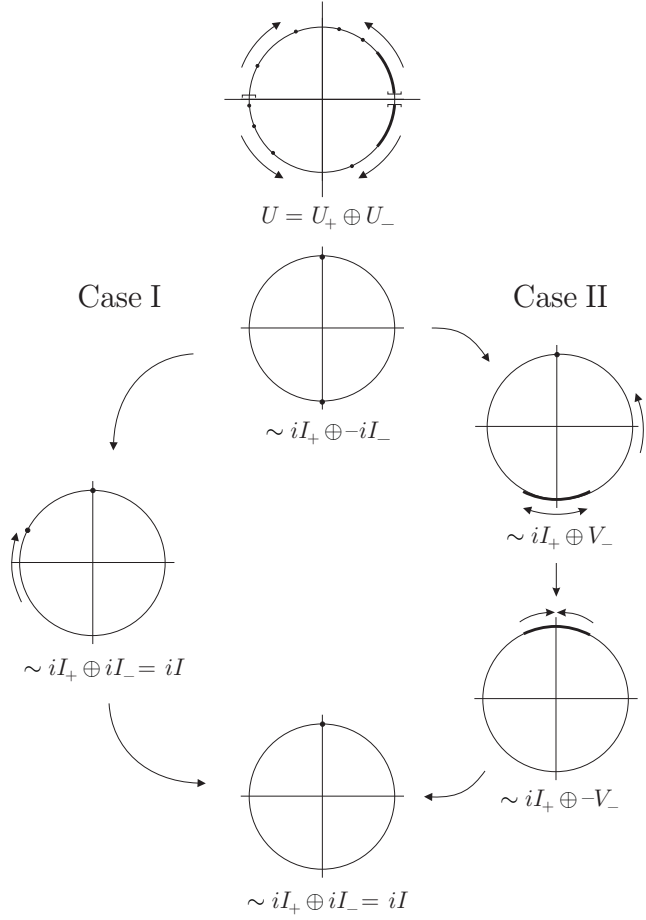}%
\caption{Connecting a fixed $U$ in $_{\mathcal{F}}\mathcal{U}_{\mathrm{inj}}$
to $iI$. Case I (finite rank $U_{-}$) and Case II (infinite rank $U_{-}$)}%
\label{f:deformation}%
\end{center}
\end{figure}

We distinguish two cases: If $H_{-}$ is \emph{finite-dimensional}, we now
rotate $-iI_{-}$ up through $-1$ into $iI_{-}$\.{.} More precisely, we
consider $\{iI_{+}\oplus e^{i(\pi/2+(1-t)\pi)}I_{-}\}_{t\in\lbrack0,1]}%
$\thinspace. This proves that we can connect $U$ with $iI_{+}\oplus iI_{-}=iI$
within $\mathbf{\kappa}\left(  \mathcal{CF}^{\operatorname{sa}}\right)  $ in
this first case.

If $H_{-}$ is \emph{infinite-dimensional}, we ``dilate'' $-iI_{-}$ in such a
way that no eigenvalues remain. To do this, we identify $H_{-}$ with
$L^{2}([0,1])$. Now multiplication by $-i$ on $L^{2}([0,1])$ can be connected
to multiplication by a function whose values are a short arc centered on $-i$
and so that the resulting operator $V_{-}$ on $H_{-}$ has no eigenvalues. This
will at no time introduce spectrum near $+1$ or $-1$. We then rotate this arc
up through $+1$ (which keeps us in the right space) until it is centered on
$+i$. Then we contract the spectrum on $H_{-}$ to be $+i$. That is, also in
this case we have connected our original operator $U$ to $+iI$. To sum up this
second case (see also Figure \ref{f:deformation}):
\begin{multline*}
U\sim iI_{+}\oplus-iI_{-}\sim iI_{+}\oplus V_{-}\sim iI_{+}\oplus e^{it\pi
}V_{-}\text{ for $t\in\lbrack0,1]$}\\
\sim iI_{+}\oplus-V_{-}\sim iI_{+}\oplus-(-iI_{-})\sim iI\,.
\end{multline*}
To prove (b), we just decompose any $V\in\,_{\mathcal{F}}\mathcal{U}$ into
$V=U\oplus I_{1}$, where $U\in{_{\mathcal{F}}}\mathcal{U}%
_{\mathrm{\operatorname{inj}}}(H_{0})$ and $I_{1}$ denotes the identity on the
1-eigenspace $H_{1}=\operatorname{Ker}(V-I)$ of $V$ with $H=H_{0}\oplus H_{1}$
an orthogonal decomposition. Then for $\varepsilon>0$, $U\oplus
e^{i\varepsilon}I_{1}\in{_{\mathcal{F}}}\mathcal{U}%
_{\mathrm{\operatorname{inj}}}$ approaches $U$ for $\varepsilon\rightarrow0$.
\end{proof}

\bigskip

\begin{remark}
\label{r:connected} The preceding proof shows also that the two subsets of
$\mathcal{C}\mathcal{F}^{\operatorname{sa}}$
\[
\mathcal{C}\mathcal{F}_{\pm}^{\operatorname{sa}}=\{T\in\mathcal{C}%
\mathcal{F}^{\operatorname{sa}}\mid\operatorname{spec}_{\operatorname{ess}%
}(T)\subset\mathbb{C}_{\mp}\},
\]
the spaces of all essentially positive, resp. all essentially negative,
self-adjoint Fredholm operators, are no longer open. The third of the three
complementary subsets
\[
\mathcal{C}\mathcal{F}_{\ast}^{\operatorname{sa}}=\mathcal{C}\mathcal{F}%
^{\operatorname{sa}}\setminus\left(  \mathcal{C}\mathcal{F}_{+}%
^{\operatorname{sa}}\cup\mathcal{C}\mathcal{F}_{-}^{\operatorname{sa}%
}\right)
\]
is also not open. We do not know whether the two ``trivial'' components are
contractible as in the bounded case nor whether the whole space is a
classifying space for $K^{1}$ as the nontrivial component in the bounded case.
Independently of Example \ref{ex:fuglede}, the connectedness of $\mathcal{C}%
\mathcal{F}^{\operatorname{sa}}$ and the disconnectedness of $\mathcal{F}%
^{\operatorname{sa}}$ show that the Riesz map is not continuous on
$\mathcal{C}\mathcal{F}^{\operatorname{sa}}$ in the gap topology.
\end{remark}

\bigskip

In analogy to \cite{Ph96}, we can give an explicit description of the winding
number (spectral flow across $-1$) $\operatorname{wind}(f)$ of a curve $f$ in
$_{\mathcal{F}}\mathcal{U}$. Alternatively, it can be used as a definition of
$\operatorname{wind}$:

\bigskip

\begin{proposition}
\label{p:wind-def} Let $f:[0,1]\rightarrow\,_{\mathcal{F}}\mathcal{U}$ be a
continuous path.\newline \noindent(a) There is a partition $\{0=t_{0}%
<t_{1}<\dots<t_{n}=1\}$ of the interval and positive real numbers
$0<\varepsilon_{j}<\pi$, $j=1,\dots,n$, such that $\operatorname{Ker}%
(f(t)-e^{i(\pi\pm\varepsilon_{j})})=\{0\}$ for $t_{j-1}\leq t\leq t_{j}%
$\thinspace.\newline \noindent(b) Then
\[
\operatorname{wind}(f)=\sum_{j=1}^{n}k(t_{j},\varepsilon_{j})-k(t_{j-1}%
,\varepsilon_{j}),
\]
where
\[
k(t,\varepsilon_{j}):=\sum_{0\leq\theta<\varepsilon_{j}}\dim\operatorname{Ker}%
(f(t)-e^{i(\pi+\theta)}).
\]
\noindent(c) In particular, this calculation of $\operatorname{wind}(f)$ is
independent of the choice of the partition of the interval and of the choice
of the barriers.
\end{proposition}

\begin{proof}
In (a) we use the continuity of $f$ and the fact that $f\left(  t\right)  \in$
$_{\mathcal{F}}\mathcal{U}$. (b) follows from the path additivity of
$\operatorname{wind}$. (c) is immediate from (b).
\end{proof}

\bigskip

This idea of a \emph{spectral flow across $-1$} was introduced first in
\cite[Sec. 1.3]{BoFu98}, where it was used to give a definition of the Maslov
index in an infinite dimensional context (see also Definition \ref{d:maslov}
further below).

\medskip

After these explanations the definition of spectral flow for paths in
$\mathcal{CF}^{\operatorname{sa}}$ is straightforward:

\medskip

\begin{definition}
\label{S2.2} Let $f:[0,1]\rightarrow\mathcal{CF}^{\operatorname{sa}}(H)$ be a
continuous path. Then the \emph{spectral flow} $\operatorname{SF}(f)$ is
defined by
\[
\operatorname{SF}(f):=\operatorname{wind}(\mathbf{\kappa}\circ f).
\]
\end{definition}

\medskip

From the properties of $\mathbf{\kappa}$ and of the winding number we infer immediately:

\medskip

\begin{proposition}
\label{S2.3} $\operatorname{SF}$ is path additive and homotopy invariant in
the following sense: let $f_{1},f_{2}:[0,1]\rightarrow\mathcal{CF}%
^{\operatorname{sa}}$ be continuous paths and let $h:[0,1]\times
\lbrack0,1]\rightarrow\mathcal{CF}^{\operatorname{sa}}$ be a homotopy such
that $h(0,t)=f_{1}(t),h(1,t)=f_{2}(t)$ and such that $\dim\operatorname{Ker}%
h(s,0),\dim\operatorname{Ker}h(s,1)$ are independent of $s$. Then
$\operatorname{SF}(f_{1})=\operatorname{SF}(f_{2})$. In particular,
$\operatorname{SF}$ is invariant under homotopies leaving the endpoints fixed.
\end{proposition}

\medskip

From Proposition \ref{p:wind-def} we get

\medskip

\begin{proposition}
\label{S2.4} For a continuous path $f:[0,1]\rightarrow\mathcal{F}%
^{\operatorname{sa}} $ our definition of spectral flow coincides with the
definition in \cite{Ph96}.
\end{proposition}

\medskip

Note that also the conventions coincide for $0\in\operatorname{spec}f(0)$ or
$0\in\operatorname{spec}f(1)$.

\medskip

\begin{corollary}
\label{c:spectral_flow}For any $S\in\mathcal{CF}^{\operatorname{sa}}$ with
compact resolvent and any continuous path $C:[0,1]\rightarrow\mathcal{B}%
^{\operatorname{sa}} $ we have $\operatorname{SF}(S+C)=\operatorname{SF}%
(\mathcal{R}\left(  S+C\right)  )$ where $\mathcal{R}$ denotes the Riesz
transformation of (\ref{G1.1}).
\end{corollary}

\medskip

Note that the curve $S+C$ is in $\mathcal{CF}^{\operatorname{sa}}$, so that
$\operatorname{SF}(S+C)$ is defined via Cayley transformation, whereas the
curve $\mathcal{R}\left(  S+C\right)  $ of the Riesz transforms is in
$\mathcal{F}^{\operatorname{sa}}$.\bigskip

\begin{remark}
The spectral flow induces a surjection of $\pi_{1}\left(  \mathcal{CF}%
^{\operatorname{sa}}\right)  $ onto $\mathbb{Z}$. Because $\mathbb{Z}$ is
free, there is a right inverse of $\operatorname{SF}$ and a normal subgroup
$G$ of $\pi_{1}\left(  \mathcal{CF}^{\operatorname{sa}}\right)  $ such that we
have a split short exact sequence
\[
0\rightarrow G\rightarrow\pi_{1}\left(  \mathcal{CF}^{\operatorname{sa}%
}\right)  \rightarrow\mathbb{Z}\rightarrow0.
\]

\smallskip

For now, an open question is whether $G$ is trivial: does the spectral flow
distinguish the homotopy classes? That is, the question is whether each loop
with spectral flow 0 can be contracted to a constant point, or equivalently,
whether two continuous paths in $\mathcal{CF}^{\operatorname{sa}}$ with same
endpoints and with same spectral flow can be deformed into each other? Or is
$\pi_{1}\left(  \mathcal{CF}^{\operatorname{sa}}\right)  \cong\mathbb{Z}%
\times\!\!|$ $G$ the semi-direct product of a nontrivial factor $G$ with
$\mathbb{Z}$? In that case, homotopy invariants of a curve in $\mathcal{CF}%
^{\operatorname{sa}}$ are not solely determined by the spectral flow (contrary
to the folklore behind parts of the topology and physics literature). For now,
we can only speculate about the existence of an additional invariant and its
possible definition. For example, one can try to define a \emph{spectral flow
at infinity}. Then, continuity of the Riesz transformation $\mathcal{R}$ on a
subclass $\mathcal{S}\subset\mathcal{CF}^{\operatorname{sa}}$ would imply
vanishing spectral flow at infinity. Non-vanishing spectral flow at infinity
will typically appear with families of the type discussed in Example
\ref{ex:fuglede} (and after Remark \ref{r:gap}). One may also expect it with
curves of differential operators of second order. However, the results of
\cite{Ni00}, though only obtained under quite restrictive conditions, may
indicate that perhaps spectral flow at infinity will not be exhibited for
continuous curves of Dirac operators. If this is true, it will also explain
why the mentioned unfounded folklore has not yet led to clear contradictions.
\end{remark}

\bigskip

\subsection{Symmetric Operators and Symplectic Analysis\label{ss:Symmetric
Operators}}

In an interview with Victor M. Buchstaber in the Newsletter of the European
Mathematical Society \cite[p. 20]{No01}, Sergej P. Novikov recalls his idea of
the late 1960's and the early 1970's, which were radically different of the
main stream in topology at that time, ``that the explanation of higher
signatures and of other deep properties of multiply connected manifolds had a
symplectic origin... In 1971 I. Gelfand went into my algebraic ideas: they
impressed him greatly. In particular, he told me of his observation that the
so-called von Neumann theory of self-adjoint extensions of symmetric operators
is simply the choice of a Lagrangian subspace in a Hilbert space with
symplectic structure.'' Closely following \cite{BoFu98}, \cite{BoFu99}, and
\cite{BoFuOt01} (see also the Krein--Vishik--Birman theory summarized in
\cite{AlSi80} and \cite{LaSnTu75}, for first pointing to symplectic aspects),
we will elaborate on that thought.

\medskip

\subsubsection{Symplectic Hilbert Space, Fredholm Lagrangian Grassmannian, and
Mas\-lov Index}

We fix the following notation. Let $(\mathcal{H},\left\langle .,.\right\rangle
)$ be a separable real Hilbert space with a fixed symplectic form $\omega$,
i.e., a skew-symmetric bounded bilinear form on $\mathcal{H}\times\mathcal{H}$
which is nondegenerate. We assume that $\omega$ is compatible with
$\left\langle .,.\right\rangle $ in the sense that there is a corresponding
almost complex structure $J:\mathcal{H}\rightarrow\mathcal{H}$ defined by
\begin{equation}
\omega(x,y)=\left\langle Jx,y\right\rangle \label{e:jot-omega}%
\end{equation}
with $J^{2}=-\operatorname{I},\ ^{t}\!J=-J$, and $\left\langle
Jx,Jy\right\rangle =\left\langle x,y\right\rangle $. Here $^{t}\!J$ denotes
the transpose of $J$ with regard to the (real) inner product $\left\langle
x,y\right\rangle $. Let $\mathcal{L}=\mathcal{L}(\mathcal{H})$ denote the set
of all Lagrangian subspaces $\lambda$ of $\mathcal{H}$ (i.e., $\lambda
=(J\lambda)^{\perp}$, or equivalently, let $\lambda$ coincide with its
annihilator $\lambda^{0} $ with respect to $\omega$). The topology of
$\mathcal{L}$ is defined by the operator norm of the orthogonal projections
onto the Lagrangian subspaces.

\smallskip

\begin{figure}
[ptb]
\begin{center}
\includegraphics[
height=2.2857in,
width=4.9052in
]%
{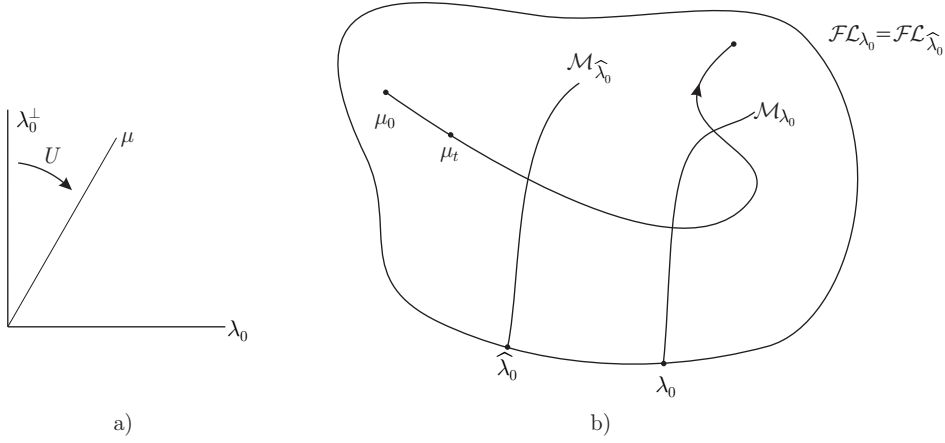}%
\caption{a) The generation of $\mathcal{L}=\left\{  \mu=U(\lambda_{0}^{\bot
})|U\in\mathcal{U}\left(  \mathcal{H}\right)  \right\}  $\newline b) One curve
and two Maslov cycles in $\mathcal{FL}\lambda_{0}=\mathcal{FL}\widehat
{\lambda}_{0}$}%
\label{f:symplectic}%
\end{center}
\end{figure}

Let $\lambda_{0}\in\mathcal{L}$ be fixed. Then any $\mu\in\mathcal{L}$ can be
obtained as the image of $\lambda_{0}^{\perp}$ under a suitable unitary
transformation
\[
\mu=U(\lambda_{0}^{\perp})
\]
(see also Figure \ref{f:symplectic}a). Here we consider the real symplectic
Hilbert space $\mathcal{H}$ as a complex Hilbert space via $J$. The group
$\mathcal{U}(\mathcal{H})$ of unitary operators of $\mathcal{H}$ acts
transitively on $\mathcal{L}$; i.e., the mapping
\begin{equation}%
\begin{matrix}
\rho:\  & \mathcal{U}(\mathcal{H}) & \longrightarrow & \mathcal{L}\\
\  & U & \mapsto &  U({\lambda_{0}}^{\perp})
\end{matrix}
\label{e:rho}%
\end{equation}
is surjective and defines a principal fibre bundle with the group of
orthogonal operators $\mathcal{O}({\lambda_{0}})$ as structure group.

\medskip

\begin{example}
(\emph{a}) In finite dimensions one considers the space $\mathcal{H}%
:=\mathbb{R}^{n}\oplus\mathbb{R}^{n}$ with the symplectic form
\[
\omega\left(  (x,\xi),(y,\eta)\right)  :=-\left\langle x,\eta\right\rangle
+\left\langle \xi,y\right\rangle \quad\text{ for $(x,\xi),\,(y,\eta
)\in\mathcal{H}$}.
\]
To emphasize the finiteness of the dimension we write $\operatorname*{Lag}%
(\mathbb{R}^{2n}):=\mathcal{L}(\mathcal{H})$. For linear subspaces of
$\mathbb{R}^{2n} $ one has
\[
l\in\operatorname*{Lag}(\mathbb{R}^{2n})\iff\dim l=n\text{\ and\ }l\subset
l^{0}:=\omega\text{-annihilator of }l\text{;}%
\]
i.e., Lagrangian subspaces are true half-spaces which are maximally isotropic
(`isotropic' means $l\subset l^{0}$).

Note that $\mathbb{R}^{2n}=\mathbb{R}^{n}\otimes\mathbb{C}$ with the Hermitian
product
\[
\left(  (x,\xi),(y,\eta)\right)  _{\mathbb{C}}=(x+i\xi)(y-i\eta):=\left\langle
x,y\right\rangle +\left\langle \xi,\eta\right\rangle +i\left\langle
\xi,y\right\rangle -i\left\langle x,\eta\right\rangle .
\]
Then every $U\in\mathrm{U}(n)$ can be written in the form $U=A+iB=\left(
\begin{array}
[c]{cc}%
A & -B\\
B & A
\end{array}
\right)  $ with $^{t}\!AB=\ ^{t}\!BA$, $A^{t}\!B=B^{t}\!A$, and $A^{t}%
\!A+B^{t}\!B=I$, and
\[
\mathrm{O}(n)\ni A\rightarrow\left(
\begin{array}
[c]{cc}%
A & 0\\
0 & A
\end{array}
\right)  =A+0i\in\mathrm{U}\left(  n\right)
\]
gives the embedding of $\mathrm{O}(n)$ in $\mathrm{U}(n)$. One finds
$\operatorname*{Lag}(\mathbb{R}^{2n})\cong\mathrm{U}(n)/\mathrm{O}(n)$ with
the fundamental group
\[
\pi_{1}(\operatorname*{Lag}(\mathbb{R}^{2n})),\lambda_{0})\cong\mathbb{Z}.
\]
The mapping is given by the `Maslov index' of loops of Lagrangian subspaces
which can be described as an intersection index with the `Maslov cycle'. There
is a rich literature on the subject, see e.g. the seminal paper \cite{Ar67},
the systematic review \cite{CaLeMi94}, or the cohomological presentation
\cite{Go97}.

\medskip

\noindent(\emph{b}) Let $\{\varphi_{k}\}_{k\in\mathbb{Z}\setminus\{0\}}$ be a
complete orthonormal system for $\mathcal{H}$. We define an almost complex
structure, and so by \eqref{e:jot-omega} a symplectic form, by
\[
J\varphi_{k}:=\operatorname*{sign}(k)\,\varphi_{-k}\,.
\]
Then the spaces $\mathcal{H}_{-}:=\operatorname{span}\{\varphi_{k}\}_{k<0}$
and $\mathcal{H}_{+}:=\operatorname{span}\{\varphi_{k}\}_{k>0}$ are
complementary Lagrangian subspaces of $\mathcal{H}$.
\end{example}

\medskip

In infinite dimensions, the space $\mathcal{L}$ is contractible due to
Kuiper's Theorem (see \cite{BlBo03}, Part I) and therefore topologically not
interesting. Also, we need some restrictions to avoid infinite-dimensional
intersection spaces when counting intersection indices. Therefore we replace
$\mathcal{L}$ by a smaller space. This problem can be solved, as first
suggested in Swanson \cite{Sw78}, by relating symplectic functional analysis
with the space $\operatorname{Fred}(\mathcal{H})$ of Fredholm operators, we
obtain finite dimensions for suitable intersection spaces and at the same time
topologically highly nontrivial objects.%

\bigskip

\begin{definition}
\label{d:flg_mas} (\emph{a}) The space of \emph{Fredholm pairs} of closed
infinite-dimensional subspaces of $\mathcal{H}$ is defined by
\begin{multline*}
\operatorname{Fred}^{2}(\mathcal{H}):=\{({\lambda},{\mu})\mid\dim\left(
{\lambda}\cap{\mu}\right)  <+\infty\text{ and }{\lambda}\text{$+{\mu}%
\subset\mathcal{H}$ closed}\\
\text{and $\dim\mathcal{H}/($}{\lambda}\text{$+{\mu})<+\infty$}\}.
\end{multline*}
Then
\[
\operatorname{index}({\lambda},{\mu}):=\dim\left(  {\lambda}\cap{\mu}\right)
-\text{$\dim\mathcal{H}/($}{\lambda}\text{$+{\mu})$}%
\]
is called the \emph{Fredholm intersection index of }$({\lambda},{\mu})$.

\noindent(\emph{b}) The \emph{Fredholm-Lagrangian Grassmannian} of a real
symplectic Hilbert space $\mathcal{H}$ at a fixed Lagrangian subspace
${\lambda}_{0}$ is defined as
\[
\mathcal{FL}_{\lambda_{0}}:=\{{\mu}\in\mathcal{L}\mid({\mu},{\lambda_{0}}%
)\in\operatorname{Fred}^{2}(\mathcal{H})\}.
\]

\noindent(\emph{c}) The \emph{Maslov cycle} of ${\lambda}_{0}$ in
$\mathcal{H}$ is defined as
\[
\mathcal{M}_{\lambda_{0}}:=\mathcal{FL}_{\lambda_{0}}\setminus\mathcal{FL}%
_{\lambda_{0}}^{(0)}\,,
\]
where $\mathcal{FL}_{\lambda_{0}}^{(0)}$ denotes the subset of Lagrangians
intersecting ${\lambda}{_{0}}$ transversally; i.e., ${\mu}\cap{\lambda}{_{0}%
}=\{0\}$.
\end{definition}

Recall the following algebraic and topological characterization of general and
Lagrangian Fredholm pairs (see \cite{BoFu98} and \cite{BoFu99}, inspired by
\cite{BoWo85}, Part 2, Lemma 2.6).

\begin{proposition}
\label{p:characterization} (\emph{a}) Let ${\lambda},\mu\in\mathcal{L}$ and
let ${\pi}_{\lambda},{\pi}_{\mu}$ denote the orthogonal projections of
$\mathcal{H}$ onto ${\lambda}$ respectively $\mu$. Then
\[
({\lambda},\mu)\in\operatorname{Fred}^{2}(\mathcal{H})\Longleftrightarrow{\pi
}_{{\lambda}}+{\pi}_{\mu}\in\operatorname{Fred}(\mathcal{H}%
)\Longleftrightarrow{\pi}_{{\lambda}}-{\pi}_{\mu}\in\operatorname{Fred}%
(\mathcal{H}).
\]

\noindent(\emph{b}) The fundamental group of $\mathcal{FL}_{\lambda_{0}}$ is
$\mathbb{Z}$, and the mapping of the loops in $\mathcal{FL}_{\lambda_{0}%
}(\mathcal{H})$ onto $\mathbb{Z}$ is given by the \textit{Maslov index}
\[
\mathbf{mas}:\pi_{1}(\mathcal{FL}_{\lambda_{0}}(\mathcal{H})))\overset{\cong
}{\longrightarrow}\mathbb{Z}.
\]
\end{proposition}

To define the Maslov index, one needs a systematic way of counting, adding and
subtracting the dimensions of the intersections $\mu_{s}\cap{\lambda}_{0} $ of
the curve $\{\mu_{s}\}$ with the Maslov cycle $\mathcal{M}_{\lambda_{0}} $. We
recall from \cite{BoFu98}, inspired by \cite{Ph96}, a functional analytical
definition for continuous curves without additional assumptions.

First we recall from (\ref{e:rho}) that any $\mu\in\mathcal{FL}_{\lambda_{0}%
}(\mathcal{H})$ can be obtained as the image of $\lambda_{0}^{\bot}$ under a
suitable unitary transformation $\mu=U(\lambda_{0}^{\bot})$. As noted before,
the operator $U$ is not uniquely determined by $\mu$. Actually, from
$\beta\cong\lambda_{0}\otimes\mathbb{C}$ we obtain a complex conjugation so
that we can define the transpose by the following formula
\[
^{T}\!U:=\overline{U^{\ast}}%
\]
and obtain a unitary operator $W(\mu):=U^{T}\!U$ which can be defined
invariantly as the \textit{complex generator} of the Lagrangian space $\mu$
relative to $\lambda_{0}$.

\begin{proposition}
\label{p:complex_generator}The composed mapping $W$%
\[
\mathcal{FL}_{\lambda_{0}}(\mathcal{H})\ni\mu\longmapsto U\longmapsto
W(\mu):=U^{T}\!U\in\,_{\mathcal{F}}\mathcal{U}\left(  \mathcal{H}\right)
\]
is well defined. In particular, we have
\[
\operatorname{Ker}\left(  I+W\right)  =\left(  \mu\cap\lambda_{0}\right)
\otimes\mathbb{C}=\left(  \mu\cap\lambda_{0}\right)  +J\left(  \mu\cap
\lambda_{0}\right)
\]
\end{proposition}

\begin{proof}
The operator $I+W$ is a Fredholm operator because $\left(  \mu,\lambda
_{0}\right)  $ is a Fredholm pair.
\end{proof}

\bigskip

\begin{definition}
\label{d:maslov}Let $\mathcal{H}$ be a symplectic Hilbert space and
$\lambda_{0}\in\mathcal{L}(\mathcal{H})$. Let
\[
\left[  0,1\right]  \ni s\longmapsto\mu_{s}\in\mathcal{FL}_{\lambda_{0}%
}(\mathcal{H})
\]
be a continuous curve. Then $W\circ\mu$ is a continuous curve in
$_{\mathcal{F}}\mathcal{U}\left(  \mathcal{H}\right)  $, and the Maslov index
can be defined by
\[
\mathbf{mas}\left(  \mu,\lambda_{0}\right)  :=\operatorname{wind}\left(
W\circ\mu\right)  ,
\]
where $\operatorname{wind}$ is defined as in Proposition \ref{p:wind-def}.
\end{definition}

\begin{remark}
\label{r:maslov_def} To define the Maslov (intersection) index $\mathbf{mas}%
\left(  \left\{  \mu_{s}\right\}  ,\lambda_{0}\right)  $, we count the change
of the eigenvalues of $W_{s}$ near $-1$ little by little. For example, between
$s=0$ and $s=s^{\prime}$ we plot the spectrum of the complex generator $W_{s}$
close to $e^{i\pi}$. In general, there will be no parametrization available of
the spectrum near $-1$. For sufficiently small $s^{\prime}$, however, we can
find barriers $e^{i(\pi+\theta)}$ and $e^{i(\pi-\theta)}$ such that no
eigenvalues are lost through the barriers on the interval $[0,s^{\prime}]$.
Then we count the number of eigenvalues (with multiplicity) of $W_{s}$ between
$e^{i\pi}$ and $e^{i(\pi+\theta)}$ at the right and left end of the interval
$[0,s^{\prime}]$ and subtract. Repeating this procedure over the length of the
whole $s$-interval $[0,1]$ gives the Maslov intersection index $\mathbf{mas}%
\left(  \left\{  \mu_{s}\right\}  ,\lambda_{0}\right)  $ without any
assumptions about smoothness of the curve, `normal crossings', or
non-invertible endpoints. \newline It is worth mentioning that the
construction can be simplified for a \emph{complex} symplectic Hilbert space
$\mathcal{H}$. Then each Lagrangian subspace of $\mathcal{H}$ is the graph of
a uniquely determined unitary operator from $\operatorname{Ker}\left(
J-iI\right)  $ to $\operatorname{Ker}\left(  J+iI\right)  $. Moreover, a pair
$\left(  \lambda,\mu\right)  $ of Lagrangian subspaces is a Fredholm pair if
and only if $U^{-1}V$ is Fredholm, where $\lambda=\frak{G}\left(  U\right)  $
and $\mu=\frak{G}\left(  V\right)  $. We have $\mathbf{mas}\left(  \left\{
\mu_{t}\right\}  ,\lambda_{0}\right)  =\operatorname{SF}\left(  U_{t}%
^{-1}V,1\right)  $ where the $1$ indicates that the spectral flow is taken at
the eigenvalue 1.
\end{remark}

\begin{remark}
\label{r:horm} (\textit{a}) By identifying $\mathcal{H}\cong\lambda_{0}%
\otimes\mathbb{C}\cong\lambda_{0}\oplus\sqrt{-1}\,\lambda_{0}$\thinspace, we
split in \cite{BoFuOt01} any $U\in\mathcal{U}(\mathcal{H})$ into a real and
imaginary part
\[
U=X+\sqrt{-1}\,Y
\]
with $X,Y:\lambda_{0}\rightarrow\lambda_{0}$\thinspace. Let $\mathcal{U}%
(\mathcal{H})^{{\text{\textrm{Fred}}}}$ denote the subspace of unitary
operators which have a Fredholm operator as real part. This is the total space
of a principal fibre bundle over the Fredholm Lagrangian Grassmannian
$\mathcal{FL}_{\lambda_{0}}$ as base space and with the orthogonal group
$\mathcal{O}(\lambda_{0})$ as structure group. The projection is given by the
restriction of the trivial bundle $\rho:\mathcal{U}(\mathcal{H})\rightarrow
\mathcal{L}$ of \eqref{e:rho} . This bundle
\[
\mathcal{U}(\mathcal{H})^{{\text{\textrm{Fred}}}}\overset{\rho}%
{\longrightarrow}\mathcal{FL}_{\lambda_{0}}%
\]
may be considered as the infinite-dimensional generalization of the familiar
bundle $\mathrm{U}(n)\rightarrow\operatorname*{Lag}(\mathbb{R}^{2n})$ for
finite $n$ and provides an alternative proof of the homotopy type of
$\mathcal{FL}_{\lambda_{0}}$.\newline \noindent(\textit{b}) The Maslov index
for curves depends on the specified Maslov cycle $\mathcal{M}_{\lambda_{0}}$.
It is worth emphasizing that two \textit{equivalent} Lagrangian subspaces
$\lambda_{0}$ and $\widehat{\lambda}_{0}$ (i.e., $\dim\lambda_{0}/(\lambda
_{0}\cap\widehat{\lambda}_{0})<+\infty$) always define the same Fredholm
Lagrangian Grassmannian $\mathcal{FL}_{\lambda_{0}}=\mathcal{FL}%
_{\widehat{\lambda}_{0}}$ but may define different Maslov cycles
$\mathcal{M}_{\lambda_{0}}\neq\mathcal{M}_{\widehat{\lambda}_{0}}$\thinspace.
The induced Maslov indices may also become different
\begin{equation}
\mathbf{mas}(\{\mu_{s}\}_{s\in\lbrack0,1]},\lambda_{0})-\mathbf{mas}(\{\mu
_{s}\}_{s\in\lbrack0,1]},\widehat{\lambda}_{0})\overset{\text{in general}%
}{\neq}0 \label{e:horm}%
\end{equation}
(see \cite{BoFu99}, Proposition 3.1 and Section 5). However, if the curve is a
loop, then the Maslov index does not depend on the choice of the Maslov cycle.
From this property it follows that the difference in \eqref{e:horm}, beyond
the dependence on $\lambda_{0}$ and $\widehat{\lambda}_{0}$, depends only on
the initial and end points of the path $\{\mu_{s}\}$ and may be considered as
the infinite-dimensional generalization $\sigma_{\text{\textbf{H\"{o}r}}}%
(\mu_{0},\mu_{1};\lambda_{0},\widehat{\lambda}_{0})$ of the H{\"{o}}rmander
index. It plays a part as the transition function of the universal covering of
the Fredholm Lagrangian Grassmannian (see also Figure \ref{f:symplectic}b).
\end{remark}

\bigskip

\subsubsection{Symmetric Operators and Symplectic Analysis\label{sss:Symmetric
Operators}}

Let $H$ be a real separable Hilbert space and $A$ an (unbounded) closed
symmetric operator defined on the domain $D_{\min}$ which is supposed to be
dense in $H$. Let $A^{\ast}$ denote its adjoint operator with domain $D_{\max
}$. We have that $A^{\ast}|_{D_{\min}}=A$ and that $A^{\ast}$ is the maximal
closed extension of $A$ in $H$. Note that $D_{\max}$ is a Hilbert space with
the graph scalar product
\[
\left(  x,y\right)  _{\mathcal{G}}:=\left(  x,y\right)  +\left(  A^{\ast
}x,A^{\ast}y\right)  ,
\]
and $D_{\min}$ is a closed subspace of this $D_{\max}$ since $A$ is closed
(each sequence in $D_{\min}$ which is Cauchy relative to the graph norm
defines a sequence in the graph $\frak{G}\left(  A\right)  $ which is Cauchy
relative to the simple norm in $H\times H$).

We form the space $\mathbf{\beta}$ of \emph{natural boundary values} with the
\emph{natural trace map} $\gamma$ in the following way:
\begin{align*}
&  D_{\max}\overset{\gamma}{\longrightarrow}D_{\max}/D_{\min}=:\mathbf{\beta
}\\
&  \;\;x\;\;\longmapsto\gamma\left(  x\right)  =\left[  x\right]  :=x+D_{\min
}.
\end{align*}
The space $\mathbf{\beta}$ becomes a symplectic Hilbert space with the induced
scalar product and the symplectic form given by Green's form
\[
\omega\left(  \left[  x\right]  ,\left[  y\right]  \right)  :=\left(  A^{\ast
}x,y\right)  -\left(  x,A^{\ast}y\right)  \;\;\;\text{for }\left[  x\right]
,\left[  y\right]  \in\mathbf{\beta}.
\]

We define the \emph{natural Cauchy data space} $\Lambda:=\gamma\left(
\operatorname{Ker} A^{\ast}\right)  $. It is a Lagrangian subspace of
$\mathbf{\beta}$ under the assumption that $A$ admits at least one
self-adjoint Fredholm extension $A_{D}$. Actually, we shall assume that $A$
has a self-adjoint extension $A_{D}$ with compact resolvent. Then $\left(
\Lambda,\gamma\left(  D\right)  \right)  $ is a Fredholm pair of subspaces of
$\mathbf{\beta}$; i.e., $\Lambda\in\mathcal{FL}_{\gamma\left(  D\right)
}\left(  \mathbf{\beta}\right)  $.

We consider a continuous curve $\left\{  C_{s}\right\}  _{s\in\lbrack0,1]}$ in
the space of bounded self-adjoint operators on $H$. We assume that the
operators $A^{\ast}+C_{s}-r$ have no `inner solutions'; i.e., they satisfy the
\emph{weak inner unique continuation property} (UCP)
\begin{equation}
\operatorname{Ker}\left(  A^{\ast}+C_{s}-r\right)  \cap D_{\min}=\left\{
0\right\}  \label{e:1.9-new}%
\end{equation}
for $s\in\left[  0,1\right]  $ and $\left|  r\right|  <\varepsilon_{0}$ with
$\varepsilon_{0}>0$. For a discussion of UCP see Section \ref{WeakUCP} below.

Clearly, the domains $D_{\max}$ and $D_{\min}$ are unchanged by the
perturbation $C_{s}$ for any $s$. So, $\mathbf{\beta}$ does not depend on the
parameter $s$. Moreover, the symplectic form $\omega$ is invariantly defined
on $\mathbf{\beta}$ and so also independent of $s$. It follows (see
\cite[Theorem 3.9]{BoFu98}) that the curve $\left\{  \Lambda_{s}%
:=\gamma\left(  \operatorname{Ker}\left(  A^{\ast}+C_{s}\right)  \right)
\right\}  $ is continuous in $\mathcal{FL}_{\gamma\left(  D\right)  }\left(
\mathbf{\beta}\right)  $.

We summarize the basic findings:

\begin{proposition}
\label{p:lagr} \emph{(a)} Assume that there exists a self--adjoint Fredholm
extension $A_{D}$ of $A$ with domain $D$. Then the Cauchy data space
$\Lambda(A)$ is a closed Lagrangian subspace of $\mathbf{\beta}$ and belongs
to the Fredholm--Lagrangian Grassmannian $\mathcal{FL}_{\gamma\left(
D\right)  }\left(  \mathbf{\beta}\right)  $.

\noindent\emph{(b)} For arbitrary domains $D$ with $D_{\min}\subset D \subset
D_{\max}$ and $\gamma(D)$ Lagrangian, the extension $A_{D}:= A_{\max}|_{D}$ is
self--adjoint. It becomes a Fredholm operator, if and only if the pair
$(\gamma(D),\Lambda(A))$ of Lagrangian subspaces of $\mathbf{\beta}$ becomes a
Fredholm pair.

\noindent\emph{(c)} Let $\{C_{t}\}_{t\in I}$ be a continuous family (with
respect to the operator norm) of bounded self--adjoint operators. Here the
parameter $t$ runs within the interval $I=[0,1]$. Assume the weak inner UCP
for all operators $A^{\ast}+C_{t}$\thinspace. Then the spaces $\gamma
(\operatorname{Ker}(A^{\ast}+C_{t},0))$ of Cauchy data vary continuously in
$\mathbf{\beta}$.
\end{proposition}

Given this, the family $\left\{  A_{D}+C_{s}\right\}  $ can be considered at
the same time in the spectral theory of self--adjoint Fredholm operators,
defining a spectral flow, and in the symplectic category, defining a Maslov
index. Under the preceding assumptions, the main result obtainable at that
level is the following general spectral flow formula (proved in \cite[Theorem
5.1]{BoFu98} and inviting to generalizations for varying domains $D_{s}$
instead of fixed domain $D$):

\begin{theorem}
\label{t:sff-old}Let $A_{D}$ be a self-adjoint extension of $A$ with compact
resolvent and let $\left\{  A_{D}+C_{s}\right\}  _{s\in\left[  0,1\right]  }$
be a family satisfying the weak inner UCP assumption. Let $\Lambda_{s}$ denote
the Cauchy data space $\gamma\left(  \operatorname*{Ker}\left(  A^{\ast}%
+C_{s}\right)  \right)  $ of $A^{\ast}+C_{s}$. Then
\[
\operatorname{SF}\left\{  A_{D}+C_{s}\right\}  =\mathbf{mas}\left(  \left\{
\Lambda_{s}\right\}  ,\gamma\left(  D\right)  \right)  .
\]
\end{theorem}

\subsection{Operators of Dirac Type and their Ellipticity\label{ss: Ops of
Dirac Type}}

\bigskip

There are different notions of operators of Dirac type. We shall not discuss
the original hyperbolic Dirac operator (in the Minkowski metric), but restrict
ourselves to the elliptic case related to Riemannian metrics.

Recall that, if $(M,g)$ is a compact smooth Riemannian manifold (with or
without boundary) with $\dim M=m$, we denote by $\frak{Cl}(M)=\{\frak{Cl}%
(TM_{x},g_{x})\}_{x\in M}$ the bundle of Clifford algebras of the tangent
spaces. For $S\rightarrow M$ a smooth complex vector bundle of Clifford
modules, the \textit{Clifford multiplication} is a bundle map $\mathbf{c}%
:\frak{Cl}(M)\rightarrow\operatorname*{Hom}(S,S)$ which yields a
representation $\mathbf{c}:\frak{Cl}(TM_{x},g_{x})\rightarrow
\operatorname*{Hom}_{\mathbb{C}}(S_{x},S_{x})$ in each fiber. We may assume
that the bundle $S$ is equipped with a Hermitian metric which makes the
Clifford multiplication skew-symmetric
\[
\langle\mathbf{c}(v)s,s^{\prime}\rangle=-\langle s,\mathbf{c}(v)s^{\prime
}\rangle\qquad\text{ for $v\in TM_{x}\text{ and }s\in S_{x}$}.
\]
We note that it is not necessary to assume that $(M,g)$ admits a spin
structure in order that $\frak{Cl}(M)$ and $S$ exist. Indeed, one special case
is obtained by taking $S=\Lambda^{\ast}\left(  TM\right)  $ and letting
$\mathbf{c}$ be the extension of $\mathbf{c}\left(  v\right)  \left(
\alpha\right)  =v\wedge\alpha-v\llcorner\alpha$ for $v\in TM_{x}%
\subset\frak{Cl}(TM_{x},g_{x})$ and $\alpha\in\Lambda^{\ast}\left(
TM_{x}\right)  ,$ where ``$\llcorner$'' denotes interior product (i.e., the
dual of the exterior product $\wedge$). The extension is guaranteed by the
fact that $\mathbf{c}\left(  v\right)  ^{2}=-g_{x}\left(  v,v\right)
\operatorname{I}$. However, we do need a spin structure in the case where $S$
is a bundle of spinors.

Any choice of a smooth connection
\[
\nabla:\mathrm{C}^{\infty}(M;S)\rightarrow\mathrm{C}^{\infty}(M;T^{\ast
}M\otimes S)
\]
defines an \emph{operator of Dirac type} $\ \mathcal{D}:=\mathbf{c}\circ
\nabla$ under the Riemannian identification of the bundles $TM$ and $T^{\ast
}M$. In local coordinates we have ${}\mathcal{D}:=\sum_{j=1}^{m}%
\mathbf{c}(e_{j})\nabla_{e_{j}}$ for any orthonormal base $\{e_{1},\dots
,e_{m}\}$ of $TM_{x}$. Actually, we may choose a local frame in such a way
that
\[
\nabla_{e_{j}}=\frac{\partial}{\partial x_{j}}+\text{ zero order terms}%
\]
for all $1\leq j\leq m$. So, locally, we have
\begin{equation}
{}\mathcal{D}:=\sum\nolimits_{j=1}^{m}\mathbf{c}(e_{j})\frac{\partial
}{\partial x_{j}}+\text{ zero order terms }. \label{e:dirac}%
\end{equation}
It follows at once that the principal symbol $\sigma_{1}(\mathcal{D})(x,\xi)$
is given by Clifford multiplication with $i\xi$, so that any operator of Dirac
type is elliptic with symmetric principal symbol. If the connection $\nabla$
is \textit{compatible} with Clifford multiplication (i.e. $\nabla\mathbf{c}%
=0$), then the operator $\mathcal{D}\;$itself becomes symmetric. We shall,
however, admit incompatible metrics. Moreover, the \textit{Dirac Laplacian}
$\mathcal{D}^{2}$ has principal symbol $\sigma_{2}(\mathcal{D}^{2})(x,\xi)$
given by scalar multiplication by $\Vert\xi\Vert^{2}$ using the Riemannian
metric. So, the principal symbol of $\mathcal{D}^{2}$ is a real multiple of
the identity, and $\mathcal{D}^{2}$ is elliptic. In the special case above
where $S=\Lambda^{\ast}\left(  TM\right)  $ and we identify $\Lambda^{\ast
}\left(  TM\right)  $ with $\Lambda^{\ast}\left(  T^{\ast}M\right)  $ by means
of the metric $g$, $\mathcal{D}$ becomes $d+\delta:\Omega^{\ast}\left(
M,\mathbb{C}\right)  \rightarrow\Omega^{\ast}\left(  M,\mathbb{C}\right)  $
and $\mathcal{D}^{2}=d\delta+\delta d$ is the Hodge Laplacian on the space
$\Omega^{\ast}\left(  M,\mathbb{C}\right)  $ of complex-valued forms on $M$.

On a closed manifold $M$, a key result for any operator $\mathcal{D}$ of Dirac
type (actually, for any linear elliptic operator of first order) is the
\textit{a priori} estimate (\emph{G\aa rding's inequality})
\begin{equation}
\left\|  \psi\right\|  _{1}\leq C\left(  \left\|  \mathcal{D}\psi\right\|
_{0}+\left\|  \psi\right\|  _{0}\right)  \text{ for all }\psi\in H^{1}\left(
M;S\right)  . \label{e:gaarding}%
\end{equation}
Here $\left\|  \,\cdot\,\right\|  _{0}$ denotes the $L^{2}$ norm and
$H^{1}\left(  M;S\right)  $ denotes the first Sobolev space with the norm
$\left\|  \,\cdot\,\right\|  _{1}$. Note that the same symbol $\mathcal{D}$ is
used for the original operator (defined on smooth sections) and its closed
$L^{2}$ extension with domain $H^{1}\left(  M;S\right)  $.

Combined with the simple continuity relation $\left\|  \mathcal{D}%
\psi\right\|  _{0}\leq C^{\prime}\left\|  \psi\right\|  _{1}$, inequality
(\ref{e:gaarding}) shows that the first Sobolev norm $\left\|  \,\cdot
\,\right\|  _{1}$ and the graph norm coincide on $H^{1}\left(  M;S\right)  $.\bigskip

\subsection{Weak Unique Continuation Property\label{WeakUCP}}

\bigskip

A linear or non-linear operator $\frak{D}$, acting on functions or sections of
a bundle over a compact or non-compact manifold $M$ has the \emph{weak Unique
Continuation Property (UCP)} if any solution $\psi$ of the equation
$\frak{D}\psi=0$ has the following property: if $\psi$ vanishes on a nonempty
open subset $\Omega$ of $M$, then it vanishes on the whole connected component
of $M$ containing $\Omega$. Note that weak inner UCP, as defined in
(\ref{e:1.9-new}), follows from weak UCP, but not vice versa.

\medskip

There is also a notion of \emph{strong UCP}, where, instead of assuming that a
solution $\psi$ vanishes on an open subset, one assumes only that $\psi$
vanishes `of high order' at a point. The concepts of weak and strong UCP
extend a fundamental property of analytic functions to \textit{some} elliptic
equations other than the Cauchy-Riemann equation.

\medskip

Up to now, (almost) all work on UCP goes back to two seminal papers
\cite{Ca33}, \cite{Ca39} by Torsten Carleman, establishing an inequality of
Carleman type (see~ our inequality \ref{e:8.1} below). In this approach, the
difference between weak and strong UCP and the possible presence of more
delicate nonlinear perturbations are related to different choices of the
weight function in the inequality, and to whether $L^{2}$ estimates suffice or
$L^{p}$ and $L^{q}$ estimates are required.

\medskip

The weak UCP is one of the basic properties of an operator of Dirac type
$\mathcal{D}$. Contrary to common belief, UCP is \textit{not} a general fact
of life for elliptic operators. See \cite{Pl61} where counter-examples are
given with smooth coefficients. Lack of UCP invalidates the continuity of the
Cauchy data spaces and of the Calder{\'{o}}n projection (Proposition
\ref{p:lagr}c and Theorem \ref{t:dependence}b) and of the main continuity
lemma (Lemma \ref{l:grass}). It corrupts the invertible double construction
(Section \ref{sss:InvertibleExtension}) and threatens Bojarski type theorems
(like Proposition \ref{p:boj} and Theorem \ref{t:sff-mwb}). For partitioned
manifolds $M=M_{1}\cup_{\Sigma}M_{2}$\thinspace(see Subsection \ref{ss:Index
of Elliptic Boundary Value Problems}), it guarantees that there are no
\textit{ghost} solutions of $\mathcal{D}\psi=0$; that is, there are no
solutions which vanish on $M_{1}$ and have nontrivial support in the interior
of $M_{2}$. This property is also called UCP \textit{from open subsets} or
\textit{across any hypersurface}. For Euclidean (classical) Dirac operators
(i.e., Dirac operators on $\mathbb{R}^{m}$ with constant coefficients and
without perturbation), the property follows by squaring directly from the
well-established UCP for the classical (constant coefficients and no
potential) Laplacian.

\medskip

From \cite[Chapter 8]{BoWo93} we recall a very simple proof of the weak UCP
for operators of Dirac type, inspired by \cite[Sections 6-7, in particular the
proof of inequality (7.11)]{Ni73} and \cite[Section II.3]{Tr80}. We refer to
\cite{Boo00} for a further slight simplification and a broader perspective,
and to \cite{BoMaWa02} for perturbed equations.

The proof does not use advanced arguments of the Aronszajn/Cordes type (see
\cite{Ar57} and \cite{Co56}) regarding the diagonal and real form of the
principal symbol of the Dirac Laplacian nor any other reduction to operators
of second order (like \cite{We82}), but only the following product property of
Dirac type operators (besides G{\aa}rding's inequality).

\medskip

\begin{lemma}
\label{l:product} Let $\Sigma$ be a closed hypersurface of $M$ with orientable
normal bundle. Let $u$ denote a normal variable with fixed orientation such
that a bicollar neighborhood ${N}$ of $\Sigma$ is parameterized by
$[-\varepsilon,+\varepsilon]\times\Sigma$. Then any operator of Dirac type can
be rewritten in the form
\begin{equation}
\mathcal{D}|_{{N}}=\mathbf{c}(du)\left(  \frac{\partial}{\partial u}%
+B_{u}+C_{u}\right)  , \label{e-product}%
\end{equation}
where $B_{u}$ is a self-adjoint elliptic operator on the parallel hypersurface
$\Sigma_{u}$, and $C_{u}:S|_{\Sigma_{u}}\rightarrow S|_{\Sigma_{u}}$ a
skew-symmetric operator of $0$-th order, actually a skew-symmetric bundle homomorphism.\bigskip
\end{lemma}

\medskip

\begin{proof}
Let $(u,y)$ denote the coordinates in a tubular neighborhood of $\Sigma$.
Locally, we have $y=(y_{1},\dots,y_{m-1})$. Let $\mathbf{c}_{u},\mathbf{c}%
_{1},\dots,\mathbf{c}_{m-1}$ denote Clifford multiplication by the unit
tangent vectors in normal, resp. tangential, directions. By \eqref{e:dirac},
we have
\begin{align*}
&  \mathcal{D}=\mathbf{c}_{u}\frac{\partial}{\partial u}+\sum\nolimits_{k=1}%
^{m-1}\mathbf{c}_{k}\frac{\partial}{\partial y_{k}}+\text{ zero order
terms}=\mathbf{c}_{u}\left(  \frac{\partial}{\partial u}+\mathcal{B}%
_{u}\right)  ,\text{ where}\\
&  \mathcal{B}_{u}:=\sum\nolimits_{k=1}^{m-1}-\mathbf{c}_{u}\mathbf{c}%
_{k}\frac{\partial}{\partial y_{k}}+\text{ zero order terms.}%
\end{align*}
We shall call $\mathcal{B}_{u}$ the \textit{tangential} operator component of
the operator $A$. Clearly it is an elliptic differential operator of first
order over $\Sigma_{u}$. From the skew-hermicity of $\mathbf{c}_{u}$ and
$\mathbf{c}_{k}$, we have
\begin{align*}
\left(  \mathbf{c}_{u}\mathbf{c}_{k}\frac{\partial}{\partial y_{k}}\right)
^{\ast}  &  =\left(  -\frac{\partial}{\partial y_{k}}\right)  (-\mathbf{c}%
_{k})(-\mathbf{c}_{u})=-\mathbf{c}_{k}\,\mathbf{c}_{u}\frac{\partial}{\partial
y_{k}}\,+\text{ zero order terms}\\
&  =\mathbf{c}_{u}\mathbf{c}_{k}\,\frac{\partial}{\partial y_{k}}+\text{ zero
order terms}.
\end{align*}
So,
\[
\mathcal{B}_{u}^{\ast}=\mathcal{B}_{u}+\text{ zero order terms}.
\]
Hence, the principal symbol of $\mathcal{B}_{u}$ is self-adjoint. Then the
assertion of the lemma is proved by setting
\begin{equation}
B_{u}:=\frac{1}{2}\left(  \mathcal{B}_{u}+\mathcal{B}_{u}^{\ast}\right)
\qquad\text{and}\qquad C_{u}:=\frac{1}{2}\left(  \mathcal{B}_{u}%
-\mathcal{B}_{u}^{\ast}\right)  \label{e:b+}%
\end{equation}
\end{proof}

\bigskip

\begin{remark}
\label{r-homo-perturbation} (a) It is worth mentioning that the product form
\eqref{e-product} is invariant under perturbation by a bundle homomorphism.
More precisely: Let $\frak{D}$ be an operator on $M$ which can be written in
the form \eqref{e-product} close to any closed hypersurface $\Sigma$%
\thinspace, with $B_{u}$ and $C_{u}$ as explained in the preceding Lemma. Let
$R$ be a bundle homomorphism. Then
\[
\left(  \frak{D}+R\right)  |_{N}=\mathbf{c}(du)\left(  \frac{\partial
}{\partial u}+B_{u}+C_{u}\right)  +\mathbf{c}(du)S|_{N}%
\]
with $T|_{N}:=\mathbf{c}(du)^{\ast}R|_{N}$. Splitting $T=\frac{1}{2}%
(T+T^{\ast})+\frac{1}{2}(T-T^{\ast})$ into a symmetric and a skew-symmetric
part and adding these parts to $B_{u}$ and $C_{u}$, respectively, yields the
desired form of $\left(  \frak{D}+R\right)  |_{N}$.\newline \qquad(b) For
operators of Dirac type, it is well known that a perturbation by a bundle
homomorphism is equivalent to modifying the underlying connection of the
operator. This gives an alternative argument for the invariance of the form
\eqref{e-product} for operators of Dirac type under perturbation by a bundle
homomorphism.\newline \qquad(c) By the preceding arguments (a), respectively
(b), establishing weak UCP for sections belonging to the kernel of a Dirac
type operator, respectively an operator which can be written in the form
\eqref{e-product}, implies weak UCP for all eigensections. Warning: for
general linear elliptic differential operators, weak UCP for ``zero-modes''
does not imply weak UCP for all eigensections.
\end{remark}

\bigskip

To prove the weak UCP, in combination with the preceding lemma, the standard
lines of the UCP literature can be radically simplified, namely with regard to
the weight functions and the integration order of estimates. These
simplifications make it also very easy to generalize the weak UCP to the
perturbed case.

We replace the equation $\mathcal{D}\psi=0$ by
\begin{equation}
\widetilde{\mathcal{D}}\psi:=\mathcal{D}\psi+\frak{P}_{A}(\psi)=0\,,
\label{e-dirac_pert}%
\end{equation}
where $\frak{P}_{A}$ is an \emph{admissible} perturbation, in the following sense.

\medskip

\begin{definition}
\label{d-pert} A perturbation is \emph{admissible} if it satisfies the
following estimate:
\begin{equation}
\left|  \frak{P}_{A}(\psi)|_{x}\right|  \leq P(\psi,x)|\psi(x)|\qquad\text{for
$x\in M$} \label{e-pert}%
\end{equation}
for a real-valued function $P(\psi,\cdot)$ which is locally bounded on $M$ for
each fixed $\psi$.
\end{definition}

\medskip

\begin{example}
\label{ex-pert} Some typical examples of perturbations satisfying the
admissibility condition of Definition \ref{d-pert} are:\newline \smallskip1.
Consider a nonlinear perturbation
\[
\frak{P}_{A}(\psi)|_{x}:=\omega(\psi(x))\cdot\psi(x)\,,
\]
where $\omega(\psi(x))|_{x\in M}$ is a (bounded) function which depends
continuously on $\psi(x)$, for instance, for a fixed (bounded) spinor section
$a(x)$ we can take
\[
\omega(\psi(x)):=\langle\psi(x),a(x)\rangle
\]
with $\langle\cdot,\cdot\rangle$ denoting the Hermitian product in the fiber
of the spinor bundle over the base point $x\in M$. This satisfies
\eqref{e-pert} .\newline \smallskip2. Another interesting example is provided
by (linear) nonlocal perturbations with
\[
\omega(\psi,x)=\left|  \int k(x,z)\psi(z)dz\right|
\]
with suitable integration domain and integrability of the kernel $k$. These
also satisfy \eqref{e-pert} .\newline \smallskip3. Clearly, an unbounded
perturbation may be both nonlinear and global at the same time. This will, in
fact, be the case in our main application. In all these cases the only
requirement is the estimate \eqref{e-pert} with bounded $\omega(\psi(\cdot))$.
\end{example}

\medskip

We now show that (admissible) perturbed Dirac operators always satisfy the
weak Unique Continuation Property. In particular, we show that this is true
for unperturbed Dirac operators.

\bigskip

\begin{theorem}
\label{t-pert} Let $\mathcal{D}$\ be an operator of Dirac type and
$\frak{P}_{A}$ an admissible perturbation. Then any solution $\psi$ of the
perturbed equation \eqref{e-dirac_pert} vanishes identically on any connected
component of the underlying manifold $M$ if it vanishes on a nonempty open
subset of the connected component.
\end{theorem}

\begin{proof}
Without loss of generality, we assume that $M$ is connected. Let $\psi$ be a
solution of the perturbed (or, in particular) unperturbed equation which
vanishes on an open, nonempty set $\Omega$. First we localize and convexify
the situation and we introduce spherical coordinates (see Figure \ref{f-ucp}).
Without loss of generality we may assume that $\Omega$ is maximal, namely the
union of all open subsets on which $\psi$ vanishes; i.e., $\Omega
=M\setminus\operatorname*{supp}$ $\psi$.%

\begin{figure}
[h]
\begin{center}
\includegraphics[
height=2.2667in,
width=2.9334in
]%
{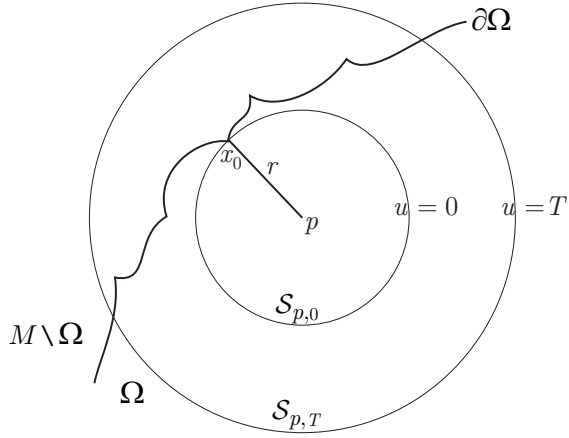}%
\caption{Local specification for the Carleman estimate}%
\label{f-ucp}%
\end{center}
\end{figure}
Since $M$ is connected, to prove that $\Omega=M$ it suffices to show that
$\Omega$ is closed (i.e., $\overline{\Omega}=\Omega$). If $\overline{\Omega
}\neq\Omega,$ then let $y_{0}\in\partial\Omega:=\overline{\Omega}%
\setminus\Omega,$ and let $B$ be an open, normal coordinate ball about $y_{0}%
$. Let $p\in\Omega\cap B$ and let $x_{0}\in\operatorname*{supp}\psi$ be a
point of the non-empty, compact set $\overline{B}\setminus\Omega=\overline
{B}\cap\operatorname*{supp}\psi$ which is closest to $p$. Let
\[
r:=d\left(  p,x_{0}\right)  =\min_{x\in\overline{B}\setminus\Omega}d\left(
p,x\right)  .
\]
For $z\in B$, let $u\left(  z\right)  :=d\left(  p,z\right)  -r$ be a ``radial
coordinate''. Note that $u=0$ defines a sphere, say $\mathcal{S}_{p,0}$, of
radius $r$ about $x_{0}$. We have larger hyperspheres $\mathcal{S}%
_{p,u}\subset M$ for $0\leq u\leq T$ with $T>0$ sufficiently \textit{small}.
In such a way we have parameterized an annular region ${N}_{T}:=\{\mathcal{S}%
_{p,u}\}_{u\in\lbrack0,T]}$ around $p$ of width $T$ and inner radius $r$,
ranging from the hypersphere $\mathcal{S}_{p,0}$ which is contained in
$\overline{\Omega}$, to the hypersphere $\mathcal{S}_{p,T}$. Note that
${N}_{T}$ contains some points where $\psi\neq0$, for otherwise $x_{0}%
\in\Omega$. Let $y$ denote a variable point in $\mathcal{S}_{p,0}$ and note
that points in ${N}_{T}$ may be identified with $(u,y)\in\left[  0,T\right]
\times\mathcal{S}_{p,0}$. Next, we replace the solution $\psi|{N}_{T}$ by a
cutoff
\[
v(u,y):=\varphi(u)\psi(u,y)
\]
with a smooth bump function $\varphi$ with $\varphi(u)=1$ for $u\leq0.8\,T $
and $\varphi(u)=0$ for $u\geq0.9\,T$. Then $\operatorname*{supp}v$ is
contained in ${N}_{T}$\thinspace. More precisely, it is contained in the
annular region ${N}_{0.9\,T}$\thinspace. Now our proof goes in two steps:
first we establish a Carleman inequality for any spinor section $v$ in the
domain of $\mathcal{D}$\ which satisfies $\operatorname*{supp}(v)\subset
N_{T}$. More precisely, we are going to show that for $T$ sufficiently small
there exists a constant $C$, such that
\begin{equation}
R\int_{u=0}^{T}\int_{\mathcal{S}_{p,u}}e^{R(T-u)^{2}}\,|v(u,y)|^{2}%
\,dy\,du\leq C\int_{u=0}^{T}\int_{\mathcal{S}_{p,u}}e^{R(T-u)^{2}%
}\,|\mathcal{D}v(u,y)|^{2}\,dy\,du \label{e:8.1}%
\end{equation}
holds for any real $R$ sufficiently large. In the second step we apply
\eqref{e:8.1} to our cutoff section $v$ and conclude that then $\psi$ is equal
$0$ on ${N}_{T/2}$\thinspace.

\medskip

\textbf{Step 1.} First consider a few technical points. The Dirac operator
$\mathcal{D}$\ has the form $G(u)(\partial_{u}+\mathcal{B}_{u})$ on the
annular region $[0,T]\times\mathcal{S}_{p,0}$, and it is obvious that we may
consider the operator $(\partial_{u}+\mathcal{B}_{u})$ instead of
$\mathcal{D}$. Moreover, we have by Lemma \ref{l:product} that $\mathcal{B}%
_{u}=B_{u}+C_{u}$ with a self-adjoint elliptic differential operator $B_{u}$
and an anti-symmetric operator $C_{u}$ of order zero, both on $\mathcal{S}%
_{p,u}$. Note that the metric structures depend on the normal variable $u$.

Now make the substitution
\[
v=:e^{-R(T-u)^{2}/2}v_{0}%
\]
which replaces (\ref{e:8.1}) by
\begin{equation}
R\int_{0}^{T}\int_{\mathcal{S}_{p,u}}\left|  v_{0}(u,y)\right|  ^{2}%
\,dy\,du\leq C\int_{0}^{T}\int_{\mathcal{S}_{p,u}}\left|  \dfrac{\partial
v_{0}}{\partial u}+\mathcal{B}_{u}v_{0}+R(T-u)v_{0}\right|  ^{2}\,dy\,du.
\label{e:8.5}%
\end{equation}
We denote the integral on the left side by $J_{0}$ and the integral on the
right side by $J_{1}$. Now we prove \eqref{e:8.5}. Decompose $\dfrac{\partial
}{\partial u}+\mathcal{B}_{u}+R(T-u)$ into its symmetric part $B_{u}+R(T-u)$
and anti-symmetric part $\partial_{u}+C_{u}$. This gives
\begin{align*}
J_{1}  &  =\int_{0}^{T}\int_{\mathcal{S}_{p,u}}\left|  \dfrac{\partial v_{0}%
}{\partial u}+\mathcal{B}_{u}v_{0}+R(T-u)v_{0}\right|  ^{2}\,dy\,du\\
&  =\int_{0}^{T}\int_{\mathcal{S}_{p,u}}\left|  \dfrac{\partial v_{0}%
}{\partial u}+C_{u}v_{0}\right|  ^{2}\,dy\,du+\int_{0}^{T}\int_{\mathcal{S}%
_{p,u}}\left|  \left(  B_{u}+R(T-u)\right)  v_{0}\right|  ^{2}\,dy\,du\\
&  \qquad+2\Re\int_{0}^{T}\int_{\mathcal{S}_{p,u}}\left\langle \dfrac{\partial
v_{0}}{\partial u}+C_{u}v_{0},B_{u}v_{0}+R(T-u)v_{0}\right\rangle \,dy\,du\,.
\end{align*}
Integrate by parts and use the identity for the real part
\[
\Re\left\langle f,Pf\right\rangle =\frac{1}{2}\left\langle f,(P+P^{\ast
})f\right\rangle
\]
in order to investigate the last and critical term which will be denoted by
$J_{2}$. This yields (where we drop domains of integration)
\begin{align*}
J_{2}  &  =2\Re\int\int\left\langle \dfrac{\partial v_{0}}{\partial u}%
+C_{u}v_{0},B_{u}v_{0}+R(T-u)v_{0}\right\rangle \,dy\,du\\
&  =2\Re\int\int\left\langle \dfrac{\partial v_{0}}{\partial u},B_{u}%
v_{0}+R(T-u)v_{0}\right\rangle \,dy\,du+2\Re\int\int\left\langle C_{u}%
v_{0},B_{u}v_{0}\right\rangle \,dy\,du\\
&  =-2\Re\int\int\left\langle v_{0},\left\{  \dfrac{\partial}{\partial
u}\left(  B_{u}+R(T-u)\right)  \right\}  v_{0}\right\rangle \,dy\,du-2\Re
\int\int\left\langle v_{0},C_{u}B_{u}v_{0}\right\rangle \,dy\,du\\
&  =2\int\int\left\langle v_{0},-\dfrac{\partial B_{u}}{\partial u}%
v_{0}+Rv_{0}\right\rangle \,dy\,du+\int\int\left(  v_{0};\left[  B_{u}%
,C_{u}\right]  v_{0}\right)  \,dy\,du\\
&  =2R\int_{0}^{T}\left\|  v_{0}\right\|  _{0}^{2}\,du+\int\int\left\langle
v_{0},-2\dfrac{\partial B_{u}}{\partial u}v_{0}+\left[  B_{u},C_{u}\right]
v_{0}\right\rangle \,dy\,du\\
&  =2RJ_{0}+J_{3},
\end{align*}
where $\Vert\cdot\Vert_{m}$ denotes the $m$-th Sobolev norm on
$E|_{\mathcal{S}_{p,u}}$ and $J_{3}$ requires a careful analysis. It follows
from the preceding decompositions of $J_{1}$ and $J_{2}$ that the proof of
\eqref{e:8.5} will be completed with $C=\tfrac{1}{2}$ when $J_{3}\geq0$. If
$J_{3}<0$ and $C=\tfrac{1}{2}$, it suffices to show that
\begin{equation}
\left|  J_{3}\right|  \leq\frac{1}{2}\left(  R\int_{0}^{T}\Vert v_{0}\Vert
_{0}^{2}\,du+\int_{0}^{T}\left\|  (B_{u}+R(T-u))v_{0}\right\|  ^{2}%
\,du\right)  . \label{e:8.8}%
\end{equation}
Since the operators $B_{u}$ are elliptic of order 1, G\aa rding's inequality
(\ref{e:gaarding}) yields
\[
\Vert f\Vert_{1}\leq c\left(  \Vert f\Vert_{0}+\Vert B_{u}f\Vert_{0}\right)
\]
for any section $f$ of $E$ on $\mathcal{S}_{p,u}$ (and $0\leq u\leq T$). Then,
also using the fact that $-2\tfrac{\partial B_{u}}{\partial u}+[B_{u},C_{u}]$
is a \textit{first}-order differential operator on $E|\mathcal{S}_{p,u}$, we
obtain
\begin{align*}
\left|  J_{3}\right|   &  \leq\int_{0}^{T}\Vert v_{0}\Vert_{0}\,\left\|
-2\tfrac{\partial B_{u}}{\partial u}v_{0}+[B_{u},C_{u}]v_{0}\right\|
_{0}\,du\leq c_{1}\int_{0}^{T}\Vert v_{0}\Vert_{0}\,\Vert v_{0}\Vert_{1}\,du\\
&  \leq c_{1}c\int_{0}^{T}\Vert v_{0}\Vert_{0}\left(  \Vert B_{u}v_{0}%
\Vert_{0}+\Vert v_{0}\Vert_{0}\right)  \,du\\
&  \leq c_{1}c\int_{0}^{T}\Vert v_{0}\Vert_{0}\left\{  \Vert(B_{u}%
+R(T-u))v_{0}\Vert_{0}+(R(T-u)+1)\Vert v_{0}\Vert_{0}\right\}  \,du\\
&  \leq c_{1}c(RT+1)\int_{0}^{T}\Vert v_{0}\Vert_{0}^{2}\,du+c_{1}c\int
_{0}^{T}\left\|  (B_{u}+R(T-u))v_{0}\right\|  _{0}\,\Vert v_{0}\Vert
_{0}\,du\,.
\end{align*}
The integrand of the second summand is equal to
\begin{equation}%
\begin{array}
[c]{c}%
\dfrac{\Vert(B_{u}+R(T-u))v_{0}\Vert_{0}}{\sqrt{c_{1}c}}\,\left(  \sqrt
{c_{1}c}\Vert v_{0}\Vert_{0}\right) \\
\leq\frac{1}{2}\left\{  \frac{1}{c_{1}c}\left\|  (B_{u}+R(T-u))v_{0}\right\|
_{0}^{2}\,+\,c_{1}c\Vert v_{0}\Vert_{0}^{2}\right\}
\end{array}
\label{e:8.10}%
\end{equation}
with the inequality due to the estimate $ab\leq\frac{1}{2}(a^{2}+b^{2})$. By
inserting \eqref{e:8.10} in the preceding inequality for $|J_{3}|$ we obtain
\begin{align*}
\left|  J_{3}\right|   &  \leq c_{1}c(RT+1)\int_{0}^{T}\Vert v_{0}\Vert
_{0}^{2}\,du\\
&  +c_{1}c\int_{0}^{T}\left(  \tfrac{1}{2}\left\{  \tfrac{1}{c_{1}c}\left\|
(B_{u}+R(T-u))v_{0}\right\|  _{0}^{2}\,+\,c_{1}c\Vert v_{0}\Vert_{0}%
^{2}\right\}  \right)  \,du\,\\
&  =c_{1}c(RT+1)\int_{0}^{T}\Vert v_{0}\Vert_{0}^{2}\,du+c_{1}c\int_{0}%
^{T}\tfrac{1}{2}\,c_{1}c\Vert v_{0}\Vert_{0}^{2}\,du\,\\
&  +\int_{0}^{T}\left(  \tfrac{1}{2}\left\|  (B_{u}+R(T-u))v_{0}\right\|
_{0}^{2}\,\right)  \,du\,\\
&  =\int_{0}^{T}\tfrac{1}{2}\left\|  (B_{u}+R(T-u))v_{0}\right\|  _{0}%
^{2}\,\,du+c_{1}c\left(  (RT+1)+\tfrac{1}{2}\,c_{1}c\right)  \int_{0}^{T}\Vert
v_{0}\Vert_{0}^{2}\,du\\
&  =\tfrac{1}{2}\int_{0}^{T}\left\|  (B_{u}+R(T-u))v_{0}\right\|  _{0}%
^{2}\,\,du+Rc_{1}c\left(  T+\frac{c_{1}c+2}{2R}\,\right)  \int_{0}^{T}\Vert
v_{0}\Vert_{0}^{2}\,du
\end{align*}
So (\ref{e:8.8}) holds for $T$ and $\frac{1}{R}$ sufficiently small, and we
then have the Carleman inequality (\ref{e:8.1}) for $C=\tfrac{1}{2}$.

\medskip

\textbf{Step 2}. To begin with, we have
\begin{equation}
e^{RT^{2}/4}\,\int_{0}^{\frac{T}{2}}\int_{\mathcal{S}_{p,u}}\left|
\psi(u,y)\right|  ^{2}\,dy\,du\leq\int_{0}^{T}\int_{\mathcal{S}_{p,u}%
}e^{R(T-u)^{2}}\,|\left(  \varphi\psi\right)  (u,y)|^{2}\,dy\,du=:I
\label{e-main1}%
\end{equation}
We apply our Carleman type inequality (\ref{e:8.1}):
\begin{equation}
I=\int_{0}^{T}\int_{\mathcal{S}_{p,u}}e^{R(T-u)^{2}}\,|\left(  \varphi
\psi\right)  (u,y)|^{2}\,dy\,du\leq\frac{C}{R}\int_{u=0}^{T}\int
_{\mathcal{S}_{p,u}}e^{R(T-u)^{2}}\,|\mathcal{D}(\varphi\psi)(u,y)|^{2}%
\,dy\,du. \label{e-main2}%
\end{equation}
Assuming that $\psi$ is a solution of the perturbed equation $\mathcal{D}%
\psi+\frak{P}_{A}(\psi)=0$, we get
\[
\mathcal{D}(\varphi\psi)=\varphi\mathcal{D}\psi+\mathbf{c}(du)\varphi^{\prime
}\psi=-\varphi\frak{P}_{A}(\psi)+\mathbf{c}(du)\varphi^{\prime}\psi.
\]
Using this in \eqref{e-main2} and noting that $\left(  a+b\right)  ^{2}%
\leq2\left(  a^{2}+b^{2}\right)  $, yields
\[
I\leq\frac{2C}{R}\int_{0}^{T}\int_{\mathcal{S}_{p,u}}e^{R(T-u)^{2}}\left(
|\varphi(u)\frak{P}_{A}(\psi)(u,y)|^{2}+|\mathbf{c}(du)\varphi^{\prime}%
(u)\psi(u,y)|^{2}\right)  \,dy\,du.
\]
Now we exploit our assumption
\begin{equation}
\left|  \frak{P}_{A}(\psi)(x)\right|  \leq P(\psi,x)|\psi(x)|\quad\text{for
$x\in M$} \label{e-main3}%
\end{equation}
about the perturbation with locally bounded $P(\psi,\cdot)$, say
\[
|P(\psi,(u,y))|\leq C_{0}:=\max_{x\in K}|P(\psi,x)|\quad\text{for all
$y\in\mathcal{S}_{p,u},\,u\in\lbrack0,T]$}%
\]
where $K$ is a suitable compact set. We obtain at once
\begin{align*}
\left(  1-\frac{2CC_{0}}{R}\right)  I  &  \leq\frac{2C}{R}\int_{0}^{T}%
\int_{\mathcal{S}_{p,u}}e^{R(T-u)^{2}}\,|\mathbf{c}(du)\varphi^{\prime}%
(u)\psi(u,y)|^{2}\,dy\,du\\
&  \leq\frac{2C}{R}e^{RT^{2}/25}\,\int_{0}^{T}\int_{\mathcal{S}_{p,u}%
}|\mathbf{c}(du)\varphi^{\prime}(u)\psi(u,y)|^{2}\,dy\,du.
\end{align*}
Here we use that $\varphi^{\prime}(u)=0$ for $0\leq u\leq0.8T$ so that we can
estimate the exponential and pull it in front of the integral. Using
(\ref{e-main1}),
\begin{align*}
&  \int_{0}^{\frac{T}{2}}\int_{\mathcal{S}_{p,u}}\left|  \psi(u,y)\right|
^{2}\,dy\,du\leq e^{-RT^{2}/4}I\\
&  \leq\frac{\frac{2C}{R}}{1-\frac{2CC_{0}}{R}}e^{RT^{2}\left(  \frac{1}%
{25}-\frac{1}{4}\right)  }\int_{0}^{T}\int_{\mathcal{S}_{p,u}}|\mathbf{c}%
(du)\varphi^{\prime}(u)\psi(u,y)|^{2}\,dy\,du.
\end{align*}
As $R\rightarrow\infty$, we get $\int_{0}^{\frac{T}{2}}\int_{\mathcal{S}%
_{p,u}}\left|  \psi(u,y)\right|  ^{2}\,dy\,du=0$ which contradicts $x_{0}%
\in\operatorname*{supp}\psi$.
\end{proof}

\bigskip

\section{The Index of Elliptic Operators on Partitioned Manifolds}

\subsection{Examples and the Hellwig--Vekua Index Theorem}

We begin with some elementary examples.

\medskip

\begin{example}
Consider the (trivially elliptic) ordinary differential operator on the unit
interval $I=[0,1]$ defined by
\[
P:C^{\infty}(I)\times C^{\infty}(I)\rightarrow C^{\infty}(I)\times C^{\infty
}(I),\text{where }\left(  f,g\right)  \mapsto\left(  f^{\prime},-g^{\prime
}\right)  ,
\]
with the boundary conditions $C^{\infty}(I)\times C^{\infty}(I)\rightarrow
C^{\infty}(\partial I)\cong\mathbb{C\times C}$ $\ $(where $\partial
I:=\left\{  0,1\right\}  $)
\begin{align*}
(i)\text{\ }R_{1}  &  :\left(  f,g\right)  \mapsto\left(  f-g\right)
|_{\partial I}\\
(ii)\text{\ }R_{2}  &  :\left(  f,g\right)  \mapsto f|_{\partial I}\\
(iii)\;R_{3}  &  :(f,g)\mapsto(f+g^{\prime})|_{\partial I}%
\end{align*}
We determine the index of the operators (for $i=1,2,3$)
\[
(P,R_{i}):C^{\infty}(I)\times C^{\infty}(I)\rightarrow C^{\infty}(I)\times
C^{\infty}(I)\times C^{\infty}(\partial I).
\]
Clearly, $\dim$ $\operatorname{Ker}(P,R_{i})=1$. To determine the cokernel,
one writes $P(f,g)=(F,G)$ and $R_{i}(f,g)=h$, with $F,G\in C^{\infty}(I)$ and
$h=\left(  h_{0},h_{1}\right)  \in\mathbb{C\times C}$ obtaining,
\[
f(t)=\int_{0}^{t}F(\tau)d\tau+c_{1},\;\;g(t)=-\int_{0}^{t}G(\tau)d\tau+c_{2}%
\]
and two more equations for the boundary condition. The dimension of
$\operatorname{Coker}(P,R_{i})$ is then the number of linearly independent
conditions on $F$, $G$, and $h$ which must be imposed in order to eliminate
the constants of integration. For each $i\in\left\{  1,2,3\right\}  $, there
is only one condition, namely $h_{0}=h_{1}-\int_{0}^{1}F(\tau)d\tau-\int
_{0}^{1}G(\tau)d\tau;$ $h_{0}=h_{1}-\int_{0}^{1}F(\tau)d\tau;$ resp.,
$h_{0}=h_{1}-\int_{0}^{1}F(\tau)d\tau-G\left(  0\right)  +G(1)$. So, the index
vanishes in all three cases.
\end{example}

\medskip

For a more comprehensive treatment of the existence and uniqueness of
boundary-value problems for ordinary differential equations (including
systems), we refer to \cite{CodLev} and \cite[p. 322-403]{Har64}.\smallskip

\medskip

We now consider the Laplace operator $\Delta:=\partial^{2}/\partial
x^{2}+\partial^{2}/\partial y^{2}$, as a linear elliptic differential operator
from $C^{\infty}(X)$ to $C^{\infty}(X)$, where $X$ is the unit disk $\left\{
z=x+iy\mid\left|  z\right|  \leq1\right\}  \subset\mathbb{C}$.

\medskip

\begin{example}
\label{DirichEx} For the Dirichlet boundary condition
\[
R:C^{\infty}(X)\rightarrow C^{\infty}(\partial X),\text{ with }R\left(
u\right)  =u|_{\partial X}\text{\ (}\partial X:=\left\{  z\in\mathbb{C}%
\mid\left|  z\right|  =1\right\}  \text{)},
\]
we show that
\[
\text{a) }\operatorname{Ker}(\Delta,R)=\left\{  0\right\}  \text{\ \ and\ \ b)
}\operatorname{Im}(\Delta,R)^{\bot}=\left\{  0\right\}  ,
\]
where $\bot$ is orthogonal complement in $L^{2}(X)\times L^{2}(\partial
X)$.\footnote{Here, consider that the intersection of the orthogonal
complement of the range $\operatorname{Im}(\Delta,R)$ relative to the usual
inner product in $L^{2}(X)\times L^{2}(\partial X)$ with the space $C^{\infty
}(X)\times C^{\infty}(\partial X)$ is isomorphic to $\operatorname{Coker}%
(\Delta,R)$. This is true, since it turns out that the image of the natural
Sobolev extension of $(\Delta,R)$ is closed in the $L^{2}$-norm, and its
$L^{2}$-orthogonal complement is contained in $C^{\infty}(X)\times C^{\infty
}(\partial X)$.} In particular, it follows that $\operatorname{index}%
(\Delta,R)=0$.\newline \qquad\textbf{For (a)}: $\operatorname{Ker}%
(\Delta,R)\ $consists of functions of the form $u+iv$, where $u$ and $v$ are
real-valued. Since the coefficients of the operators $\Delta$ and $R$ are
real, we may assume $v=0$ without loss of generality. Thus, consider a real
solution $u$ with $\Delta u=0$ in $X$ and $u=0$ on $\partial X$. Then (where
$\nabla u:=\left(  u_{x},u_{y}\right)  :=(\tfrac{\partial u}{\partial
x},\tfrac{\partial u}{\partial y})$)
\begin{equation}
0=-\int_{X}u\Delta u\,dxdy=\int_{X}\left|  \nabla u\right|  ^{2}dxdy,
\label{DirichStokes}%
\end{equation}
whence $\nabla u=0$. Thus, $u$ is constant, and indeed zero since $u=0$ on
$\partial X$. The trick lies in the equality (\ref{DirichStokes}), which
follows from Stokes' formula, namely $\int_{X}d\omega=\int_{\partial X}\omega
$, where $\omega$ is a 1-form. Indeed, setting $\omega:=u\wedge\ast du$, where
$\ast$ is the Hodge star operator ($\ast du=\ast\left(  u_{x}dx+u_{y}%
dy\right)  :=u_{x}dy-u_{y}dx$), we obtain
\[
d\omega=du\wedge\ast du+u\wedge d\ast du=\left|  \nabla u\right|  ^{2}dx\wedge
dy+\left(  u\Delta u\right)  dx\wedge dy.
\]
Using Stokes' formula and $u|_{\partial X}=0$, we have
\[
\int_{X}\left|  \nabla u\right|  ^{2}dxdy+\int_{X}\left(  u\Delta u\right)
dxdy=\int_{X}d\omega=\int_{\partial X}\omega=\int_{\partial X}u\wedge\ast
du=0.
\]
From this and $\Delta u=0$, we conclude that $\nabla u=0$ and $u$ is constant.
\newline \qquad\textbf{For (b)}: Choose $L\in C^{\infty}\left(  X\right)  $
and $l\in C^{\infty}\left(  \partial X\right)  $ with $\left(  L,l\right)  $
orthogonal to $\operatorname{Im}(\Delta,R)$, whence (relative to the usual
measures on $X$ and $\partial X$)
\begin{equation}
\int_{X}\left(  \Delta u\right)  L+\int_{\partial X}ul=0\text{ for all }u\in
C^{\infty}(X). \label{Llorth}%
\end{equation}
Using a 2-fold integration by parts (in the exterior calculus), we obtain
\begin{equation}
\int_{X}u\Delta L-\int_{X}\left(  \Delta u\right)  L=\int_{X}\left(  u(d\ast
dL)-d\left(  \ast du\right)  L\right)  =\int_{\partial X}\left(  u\ast
dL-L\ast du\right)  . \label{GrnForm}%
\end{equation}
First we consider $u$ with support $\operatorname*{supp}\left(  u\right)  :=$
the closure of $\left\{  z\in X\mid u(z)\neq0\right\}  $ contained in the
interior of $X$. Then
\[
\int_{X}u\Delta L=\int_{X}\left(  \Delta u\right)  L=-\int_{\partial X}ul=0,
\]
and so $\Delta L=0$. Now for $u\in C^{\infty}(X)$ we apply (\ref{Llorth}) and
(\ref{GrnForm}) to deduce that
\begin{align*}
\int_{\partial X}ul  &  =-\int_{X}\left(  \Delta u\right)  L=\int_{\partial
X}\left(  u\ast dL-L\ast du\right) \\
&  =\int_{\partial X}\left(  u\left(  xL_{x}+yL_{y}\right)  -L\left(
xu_{x}+yu_{y}\right)  \right)  .
\end{align*}
Thus, $l=xL_{x}+yL_{y}$ and $L|_{\partial X}=0$, and we finally apply (a).
Details are in \cite[p. 264]{Ho63}.
\end{example}

The preceding result $\operatorname{index}(\Delta,R)=0$ (for $Ru=u|_{\partial
X}$) can also be obtained by proving the symmetry of $\Delta$ and that the
$L^{2}$ extension on the domain defined by $Ru=0$ is a self-adjoint Fredholm extension.\bigskip

We now consider a $C^{\infty}$ vector field $\nu:\partial X\rightarrow
\mathbb{C}$ on the boundary $\partial X=\left\{  z\in\mathbb{C}\mid\left|
z\right|  =1\right\}  $. For $u\in C^{\infty}(X)$, $z\in\partial X$, and
$\nu(z)=\alpha(z)+i\beta(z)$, the ``directional derivative'' of the function
$u$ relative to $\nu$ at the point $z$ is
\[
\tfrac{\partial u}{\partial\nu}\left(  z\right)  :=\alpha(z)\tfrac{\partial
u}{\partial x}\left(  z\right)  +\beta(z)\tfrac{\partial u}{\partial y}\left(
z\right)  .
\]
From the standpoint of differential geometry it is better, either to denote
the vector field by $\tfrac{\partial}{\partial\nu}$ or to write the
directional derivative as simply as $\nu\left[  u\right]  \left(  z\right)  $,
since $\tfrac{\partial}{\partial x}$ and $\tfrac{\partial}{\partial y}$ can be
regarded as vector fields. The pair $(\Delta,\tfrac{\partial}{\partial\nu})$
defines a linear operator
\[
(\Delta,\tfrac{\partial}{\partial\nu}):C^{\infty}(X)\rightarrow C^{\infty
}(X)\oplus C^{\infty}(\partial X)\text{\ given by\ }u\mapsto(\Delta
u,\tfrac{\partial u}{\partial\nu}).
\]

\begin{theorem}
\label{VekuaThm}\textrm{(I. N. Vekua 1952)}. For $p\in\mathbb{Z}$ and
$\nu(z):=z^{p} $, we have that $(\Delta,\tfrac{\partial}{\partial\nu})$ is an
operator with finite-dimensional kernel and cokernel, and
\[
\operatorname{index}(\Delta,\tfrac{\partial}{\partial\nu})=2\left(
1-p\right)  .
\]%
\begin{figure}
[h]
\begin{center}
\includegraphics[
height=3.4861in,
width=4.3656in
]%
{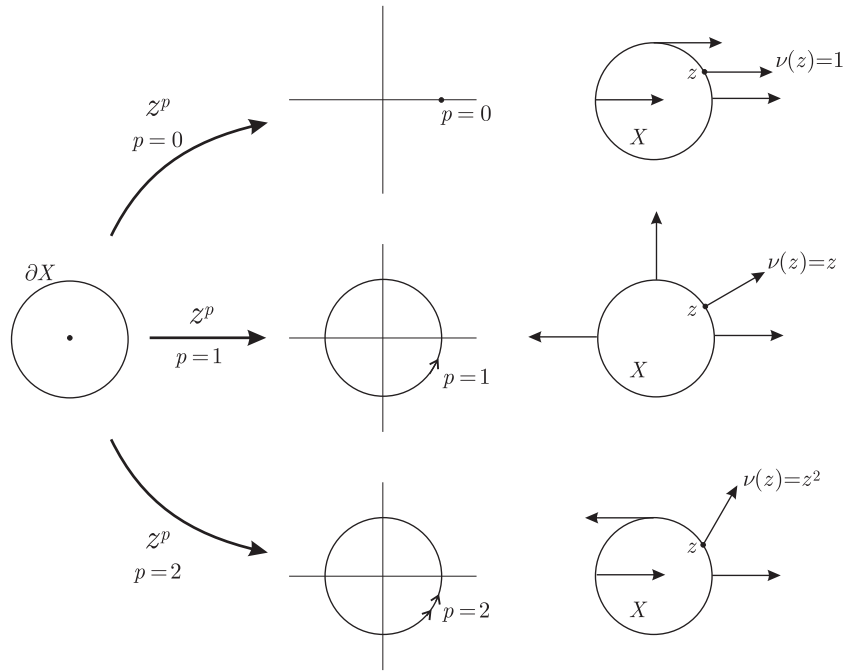}%
\caption{The vector field $\nu:\partial X\rightarrow\mathbb{C}$ with winding
number $p=0,1,2$}%
\label{f:vekua-wind}%
\end{center}
\end{figure}
\end{theorem}

\begin{remark}
The theorem of Vekua remains true, if we replace $z^{p}$ by any nonvanishing
``vector field'' $\nu:\partial X\rightarrow\mathbb{C}\setminus\left\{
0\right\}  $ with ``winding number'' $p$:%
\begin{figure}
[ptb]
\begin{center}
\includegraphics[
height=1.1684in,
width=4.7738in
]%
{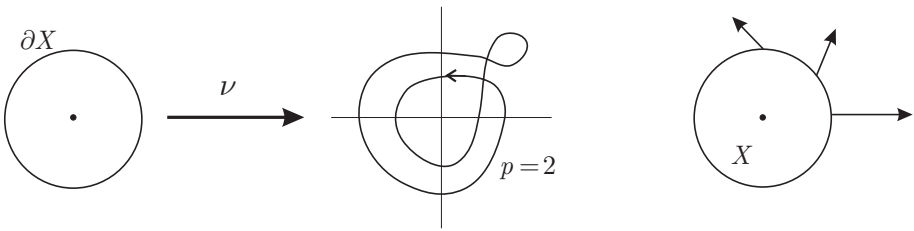}%
\caption{Another vector field $\nu$ with winding number 2 }%
\label{f:vekua-wind-mod}%
\end{center}
\end{figure}

Moreover, in place of the disk, we can take $X$ to be any simply-connected
domain in $\mathbb{C}$ with a ``smooth'' boundary $\partial X$. The reason is
the homotopy invariance of the index.
\end{remark}

\begin{remark}
In the theory of Riemann surfaces (e.g., in Riemann--Roch Theorem), one also
encounters the number $2(1-p)$ as the Euler characteristic of a closed surface
of genus $p$. This is no accident, but rather it is connected with the
relation between elliptic boundary-value problems and elliptic operators on
closed manifolds. Specifically, there is a relation between the index of
$(\Delta,\tfrac{\partial}{\partial\nu})$ and the index of Cauchy-Riemann
operator for complex line bundles over $S^{2}=\mathbb{P}^{1}\left(
\mathbb{C}\right)  $ with Chern number $1-p$ (e.g., see Example \ref{CREx} for
a start).
\end{remark}

\begin{remark}
Motivated by the method of replacing a differential equation by difference
equations, David Hilbert and Richard Courant expected ``linear problems of
mathematical physics which are correctly posed to behave like a system of $N$
linear algebraic equations in $N$ unknowns... If for a correctly posed problem
in linear differential equations the corresponding homogeneous problem
possesses only the trivial solution zero, then a uniquely determined solution
of the general inhomogeneous system exists. However, if the homogeneous
problem has a nontrivial solution, the solvability of the inhomogeneous system
requires the fulfillment of certain additional conditions.'' This is the
``heuristic principle'' which Hilbert and Courant saw in the Fredholm
Alternative \cite{CoHi}. G\"{u}nter Hellwig \cite{He52} (nicely explained in
\cite{Ha52}) in the real setting and Ilya Nestorovich Vekua \cite{Ve56} in
complex setting disproved it with their independently found example where the
principle fails for $p\neq1$. \newline We remark that in addition to these
``oblique-angle'' boundary-value problems, ``coupled'' oscillation equations,
as well as restrictions of boundary-value problems, even with vanishing index,
to suitable half-spaces, furnish further more or less elementary examples for
index $\neq0$. The simplest example of a system of first order differential
operators on the disc is provided in Example \ref{CRBVEx} below. A world of
more advanced, and for differential geometry much more meaningful examples, is
approached by the Atiyah--Patodi--Singer Index Theorem, see Section
\ref{ss:APSIndexThm} below.
\end{remark}

\begin{proof}
[Proof of Theorem \ref{VekuaThm}](After \cite[p. 266 f.]{Ho63}) Since the
coefficients of the differential operators $(\Delta,\tfrac{\partial}%
{\partial\nu})$ are real, we may restrict ourselves to real functions. Thus,
$u\in C^{\infty}(X)$ denotes a single real-valued function, rather than a
complex-valued function.

$\operatorname{Ker}(\Delta,\tfrac{\partial}{\partial\nu})$: It is well-known
that $\operatorname{Ker}(\Delta)$ consists of real (or imaginary) parts of
holomorphic functions on $X$ (e.g., see \cite[p. 175 ff.]{Ah53}). Hence,
$u\in\operatorname{Ker}(\Delta)$, exactly when $u=\frak{Re}(f)$ where $f=u+iv$
is holomorphic; i.e., the Cauchy-Riemann equation $\frac{\partial f}%
{\partial\bar{z}}=0$ holds, where $\frac{\partial}{\partial\bar{z}}:=\tfrac
{1}{2}(\frac{\partial}{\partial x}+i\frac{\partial}{\partial y})$.
Explicitly,
\[
0=\tfrac{\partial f}{\partial\bar{z}}=\tfrac{1}{2}\left(  \tfrac{\partial
}{\partial x}+i\tfrac{\partial}{\partial y}\right)  \left(  u+iv\right)
=\tfrac{1}{2}\left(  \tfrac{\partial u}{\partial x}-\tfrac{\partial
v}{\partial y}\right)  +\tfrac{i}{2}\left(  \tfrac{\partial u}{\partial
y}+\tfrac{\partial v}{\partial x}\right)  .
\]
Every holomorphic (= complex differentiable) function $f$ is twice complex
differentiable and its derivative is given by
\begin{align*}
&  \tfrac{\partial f}{\partial z}:=\tfrac{1}{2}\left(  \tfrac{\partial
}{\partial x}-i\tfrac{\partial}{\partial y}\right)  \left(  u+iv\right) \\
&  =\tfrac{1}{2}\left(  \tfrac{\partial u}{\partial x}+\tfrac{\partial
v}{\partial y}\right)  +\tfrac{i}{2}\left(  -\tfrac{\partial u}{\partial
y}+\tfrac{\partial v}{\partial x}\right)  =\tfrac{\partial u}{\partial
x}-i\tfrac{\partial u}{\partial y}.
\end{align*}
In this way we have a holomorphic function $\phi:=f^{\prime}$ for each
$u\in\operatorname{Ker}(\Delta)$. Since
\[
\tfrac{\partial u}{\partial\nu}=\frak{Re}\left(  z^{p}\right)  \tfrac{\partial
u}{\partial x}+\mathfrak{Im}\left(  z^{p}\right)  \tfrac{\partial u}{\partial
y}=\frak{Re}\left(  \left(  \tfrac{\partial u}{\partial x}-i\tfrac{\partial
u}{\partial y}\right)  z^{p}\right)  =\frak{Re}\left(  \phi\left(  z\right)
z^{p}\right)  ,
\]
the boundary condition $\tfrac{\partial u}{\partial\nu}=0$ ($\nu=z^{p}$) then
means that the real part $\frak{Re}(\phi(z)z^{p})$ vanishes for $\left|
z\right|  =1$. For $p\geq0$, $\phi(z)z^{p}$ is holomorphic as well as $\phi$,
and hence for $\phi(z):=\frac{\partial u}{\partial x}-i\frac{\partial
u}{\partial y},$ we have
\begin{align*}
u  &  \in\operatorname{Ker}\left(  \Delta,\tfrac{\partial}{\partial\nu
}\right)  \text{ with }\nu=z^{p},\text{ }p\geq0\\
&  \Rightarrow\frak{Re}\left(  \phi(z)z^{p}\right)  \in\operatorname{Ker}%
\left(  \Delta,R\right)  \text{ where }R\left(  \cdot\right)  =\left(
\cdot\right)  |_{\partial X}.
\end{align*}
Thus, we succeed in associating with the ``oblique-angle'' boundary-value
problem for $u$ a Dirichlet boundary-value problem for $\frak{Re}\left(
\phi(z)z^{p}\right)  $, which has only the trivial solution by Example
\ref{DirichEx}a. Since $\phi(z)z^{p}$ is holomorphic with $\frak{Re}\left(
\phi(z)z^{p}\right)  =0$, the partial derivatives of the imaginary part
vanish, and so there is a constant $C\in\mathbb{R}$ such that $\phi
(z)z^{p}=iC$ for all $z\in X$. If $p>0$, then we have $C=0$ (set $z=0$). Hence
$\phi=0$, and (by the definition of $\phi$) the function $u$ is constant
(i.e., $\dim\operatorname{Ker}\left(  \Delta,\tfrac{\partial}{\partial\nu
}\right)  =1$). If $p=0$, then $\tfrac{\partial u}{\partial x}-i\tfrac
{\partial u}{\partial y}=\phi(z)=iC$, and so $u(x,y)=-Cy+\widetilde{C}$,
whence $\dim\operatorname{Ker}\left(  \Delta,\tfrac{\partial}{\partial\nu
}\right)  =2$ in this case.

We now come to the case $p<0$, which curiously is not immediately reducible to
the case $q>0$ where $q:=-p$. One might try to look for a solution by simply
turning $\tfrac{\partial}{\partial\nu}$ around to $-\tfrac{\partial}%
{\partial\nu}$, but this is futile since the winding numbers of $\nu$ and
$-\nu$ about $0$ are the same. Besides, if $\nu_{p}\left(  z\right)  =z^{p}$,
we do \textit{not} have $\tfrac{\partial}{\partial\nu_{-p}}=-\tfrac{\partial
}{\partial\nu_{p}}$. In order to reduce the boundary-value problem with $p<0$
to the elementary Dirichlet problem, we must go through a more careful
argument. Note that $\phi(z)z^{p}$ can have a pole at $z=0$, whence
$\frak{Re}\left(  \phi(z)z^{p}\right)  $ is not necessarily harmonic. We write
the holomorphic function $\phi(z)$ as
\[
\phi(z)=\sum\nolimits_{j=0}^{q}a_{j}z^{j}+g\left(  z\right)  z^{q+1}%
\]
where $q:=-p$ and $g$ is holomorphic. We define a holomorphic function $\psi$
by
\[
\psi(z):=g\left(  z\right)  z+\sum\nolimits_{j=0}^{q-1}\overline{a}_{j}%
z^{q-j},
\]
with $\psi(0)=0$. Then one can write
\[
\phi(z)z^{p}=a_{q}+\sum\nolimits_{j=0}^{q-1}\left(  a_{j}z^{j-q}-\overline
{a}_{j}z^{q-j}\right)  +\psi(z)
\]
The boundary condition $\tfrac{\partial u}{\partial\nu}=0$ implies
$\frak{Re}\left(  \phi(z)z^{p}\right)  =0$ for $\left|  z\right|  =1$. By the
above equation, we have $0=\frak{Re}\left(  \phi(z)z^{p}\right)
=\frak{Re}(\psi(z)+a_{q})$ for $\left|  z\right|  =1$ since then $z^{-1}%
=\bar{z}.$ Since $\psi$ is holomorphic, we have again arrived at a Dirichlet
boundary value problem; this time for the function $\frak{Re}(\psi(z)+a_{q})$.
From Example \ref{DirichEx}a, it follows again that $\psi(z)+a_{q}$ is an
imaginary constant, whence $\psi(z)=\psi(0)=0$ and $a_{q}$ is pure imaginary.
We have
\[
\phi(z)=\phi(z)z^{p}z^{q}=a_{q}z^{q}+\sum\nolimits_{j=0}^{q-1}\left(
a_{j}z^{j}-\overline{a}_{j}z^{2q-j}\right)
\]
for arbitrary $a_{0},a_{1},...,a_{q-1}\in\mathbb{C}$ and $a_{q}\in
i\mathbb{R}$. As a vector space over $\mathbb{R}$, the set
\[
\left\{  \tfrac{\partial u}{\partial x}-i\tfrac{\partial u}{\partial y}\mid
u\in\operatorname{Ker}\left(  \Delta,\tfrac{\partial}{\partial\nu}\right)
\right\}
\]
has dimension $2q+1$; here we have restricted ourselves to real $u$, according
to our convention above. Since $u$ is uniquely determined by $\phi$ up to an
additive constant, it follows that for $\nu=z^{p}$ and $p<0$,
\[
\dim\operatorname{Ker}\left(  \Delta,\tfrac{\partial}{\partial\nu}\right)
=2q+2=2-2p
\]
$\operatorname{Coker}\left(  \Delta,\tfrac{\partial}{\partial\nu}\right)  $:
As Example \ref{DirichEx}b shows the equation $\Delta u=F$ has a solution for
each $F\in C^{\infty}(X)$. In view of this we can show
\[
\operatorname{Coker}\left(  \Delta,\tfrac{\partial}{\partial\nu}\right)
=\frac{C^{\infty}(X)\times C^{\infty}(\partial X)}{\operatorname{Im}\left(
\Delta,\tfrac{\partial}{\partial\nu}\right)  }\cong\frac{C^{\infty}(\partial
X)}{\tfrac{\partial}{\partial\nu}\left(  \operatorname{Ker}\Delta\right)  }%
\]
as follows. We assign to each representative pair $(F,h)\in C^{\infty
}(X)\times C^{\infty}(\partial X)$ the class of $h-\tfrac{\partial u}%
{\partial\nu}\in C^{\infty}(\partial X)$, where $u$ is chosen so that $\Delta
u=F$. This map is clearly well defined on the quotient space of pairs, and the
inverse map is given by $h\mapsto(0,h)$. Hence, we have found a representation
for $\operatorname{Coker}(\Delta,\tfrac{\partial}{\partial\nu})$ in terms of
the ``boundary functions'' $\left\{  \tfrac{\partial u}{\partial\nu}\mid
u\in\operatorname{Ker}\Delta\right\}  $, rather than the cumbersome pairs in
$\operatorname{Im}\left(  \Delta,\tfrac{\partial}{\partial\nu}\right)  $.
(This trick can always be applied for the boundary-value problems $(P,R)$,
when the operator $P$ is surjective.)

We therefore investigate the existence of solutions of the equation $\Delta
u=0$ with the ``inhomogeneous'' boundary condition $\tfrac{\partial
u}{\partial\nu}=h$, where $h$ is a given $C^{\infty}$ function on $\partial X
$. According to the trick introduced in the first part of our proof, it is
equivalent to ask for the existence of a holomorphic function $\phi$ with the
boundary condition $\frak{Re}\left(  \phi(z)z^{p}\right)  =h$, $\left|
z\right|  =1$, i.e., for a solution of a Dirichlet problem for $\frak{Re}%
\left(  \phi(z)z^{p}\right)  $. By Example \ref{DirichEx}b, there is a unique
harmonic function which restricts to $h$ on the boundary $\partial X$; hence,
we have a (unique up to an additive imaginary constant) holomorphic function
$\theta$ with $\frak{Re}\theta(z)=h$ for $\left|  z\right|  =1$.

In the case $p<0$, the boundary problem for $\phi$ is always solvable; namely,
set $\phi(z):=z^{-p}\theta(z)$. Hence, we have
\[
\dim\operatorname{Coker}\left(  \Delta,\tfrac{\partial}{\partial\nu}\right)
=0\text{ for }\nu(z)=z^{p}\text{ and }p\leq0.
\]
For $p>0$, we can construct a solution of the boundary-value problem in for
$\phi$ from $\theta$, if and only if there is a constant $C\in\mathbb{R}$,
such that $(\theta(z)-iC)/z^{p}$ is holomorphic (i.e., the holomorphic
function $\theta(z)-iC$ has a zero of order at least $p$ at $z=0$. Using the
Cauchy Integral Formula, these conditions on the derivatives of $\theta$ at
$z=0$ correspond to conditions on line integrals around $\partial X$. In this
way, we have $2p-1$ linear (real) equations that $h$ must satisfy in order
that the boundary-value problem have a solution. We summarize our results in
the following table ($\nu(z)=z^{p}$) and Figure \ref{f:vekua}:
\[%
\begin{tabular}
[c]{|c|c|c|c|}\hline\hline
$p%
\genfrac{}{}{0pt}{}{\mathstrut}{\mathstrut}%
$ & $\dim\operatorname{Ker}\left(  \Delta,\tfrac{\partial}{\partial\nu
}\right)  $ & $\dim\operatorname{Coker}\left(  \Delta,\tfrac{\partial
}{\partial\nu}\right)  $ & $\operatorname{index}\left(  \Delta,\tfrac
{\partial}{\partial\nu}\right)  $\\\hline\hline
$>0$ & $1$ & $2p-1$ & $2-2p$\\\hline
$\leq0$ & $2-2p$ & $0$ & $2-2p$\\\hline
\end{tabular}
\]%

\begin{figure}
[h]
\begin{center}
\includegraphics[
height=2.4958in,
width=1.7037in
]%
{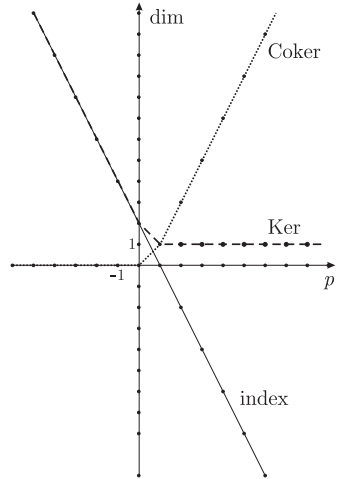}%
\caption{The dimensions of kernel, cokernel and the index of the Laplacian
with boundary condition given by $\nu\left(  z\right)  :=z^{p}$ for varying
$p$}%
\label{f:vekua}%
\end{center}
\end{figure}
\end{proof}

\textbf{Note 1}: We already noted in the proof the peculiarity that the case
$p<0$ cannot simply be played back to the case $p>0$. This is reflected here
in the asymmetry of the dimensions of kernel and cokernel and the index. It
simply reflects the fact that there are ``more'' rational functions with
prescribed poles than there are polynomials with ``corresponding'' zeros.

\textbf{Note 2}: In contrast to the Dirichlet Problem, which we could solve
via integration by parts (i.e., via Stokes' Theorem), the above proof is
function-theoretic in nature and cannot be used in higher dimensions. This is
no loss in our special case, since the index of the ``oblique-angle''
boundary-value problem must vanish anyhow in higher dimensions for topological
reasons; see \cite[p. 265 f.]{Ho63}. The actual mathematical challenge of the
function-theoretic proof arises less from the restriction $\dim X=2$ than from
a certain arbitrariness, namely the tricks and devices of the definitions of
the auxiliary functions $\phi,\psi,\theta$, by means of which the
oblique-angle problem is reduced to the Dirichlet problem.

\textbf{Note 3}: The theory of ordinary differential equations easily conveys
the impression that partial differential equations also possess a ``general
solution'' in the form of a functional relation between the unknown function
(``quantity'') $u$, the independent variables $x$ and some arbitrary constants
or functions, and that every ``particular solution'' is obtained by
substituting certain constants or functions $f,h,$ etc. for the arbitrary
constants and functions. (Corresponding to the higher degree of freedom in
partial differential equations, we deal not only with constants of integration
but with arbitrary functions.) The preceding calculations, regarding the
boundary value problem of the Laplace operator, clearly indicate how limited
this notion is which was conceived in the 18th century on the basis of
geometric intuition and physical considerations. The classical recipe of first
searching for general solutions and only at the end determining the arbitrary
constants and functions fails. For example, the specific form of boundary
conditions must enter the analysis to begin with.

\begin{example}
\label{CRBVEx}Let $X:=\left\{  z=x+iy\mid\left|  z\right|  <1\right\}  $ be
the unit disk and define an operator
\begin{align*}
&  T:C^{\infty}(X)\times C^{\infty}(X)\rightarrow C^{\infty}(X)\oplus
C^{\infty}(X)\oplus C^{\infty}(\partial X)\text{ by}\\
&  T\left(  u,v\right)  :=\left(  \tfrac{\partial u}{\partial\bar{z}}%
,\tfrac{\partial v}{\partial z},\left(  u-v\right)  |_{\partial X}\right)
\end{align*}
where $\frac{\partial}{\partial z}=\tfrac{1}{2}(\frac{\partial}{\partial
x}-i\frac{\partial}{\partial y})$ is complex differentiation and
$\frac{\partial}{\partial\bar{z}}=\tfrac{1}{2}(\frac{\partial}{\partial
x}+i\frac{\partial}{\partial y})$ is the Cauchy-Riemann differential operator
``formally adjoint'' to $\frac{\partial}{\partial z}$. We show that
$\dim(\operatorname{Ker}T)=1$ and $\operatorname{Coker}(T)=\left\{  0\right\}
$ and hence that $\operatorname{index}(T)=1$. Suppose that $\left(
u,v\right)  \in\operatorname{Ker}T$. Then $\frac{\partial u}{\partial\bar{z}%
}=0$ and $\frac{\partial v}{\partial z}=0$ in which case $u$ and $v$ are
harmonic. Then since $\left(  u-v\right)  |_{\partial X}=0,$ we have $u=v$ on
$X$. Now $\frac{\partial v}{\partial\bar{z}}=\frac{\partial u}{\partial\bar
{z}}=0\Rightarrow v$ is holomorphic, and $v^{\prime}\left(  z\right)
=\frac{\partial v}{\partial z}=0\Rightarrow v\;(=u)$ is constant. Thus,
$\dim(\operatorname{Ker}T)=1$. To show that $\operatorname{Coker}(T)=\left\{
0\right\}  $, or more precisely $(\operatorname{Im}T)^{\bot}=\left\{
0\right\}  $; see the footnote to Example \ref{DirichEx}, we choose arbitrary
$f,g\in C^{\infty}(X)$ and $h\in C^{\infty}(\partial X)$ and prove that $f$,
$g$ and $h$ must identically vanish, if
\begin{equation}
\int_{X}\left(  \tfrac{\partial u}{\partial\bar{z}}f+\tfrac{\partial
v}{\partial z}g\right)  +\int_{\partial X}(u-v)h=0\;\text{for all }u,v\in
C^{\infty}(X). \label{CokerCond}%
\end{equation}
Note that for $P,Q\in C^{\infty}(X)$
\begin{align*}
&  d\left(  Pdz+Qd\overline{z}\right)  =\tfrac{\partial P}{\partial\bar{z}%
}d\overline{z}\wedge dz+\tfrac{\partial Q}{\partial z}dz\wedge d\overline
{z}=\left(  \tfrac{\partial P}{\partial\bar{z}}-\tfrac{\partial Q}{\partial
z}\right)  d\overline{z}\wedge dz\\
&  =\left(  \tfrac{\partial P}{\partial\bar{z}}-\tfrac{\partial Q}{\partial
z}\right)  \left(  dx-idy\right)  \wedge\left(  dx+idy\right)  =2i\left(
\tfrac{\partial P}{\partial\bar{z}}-\tfrac{\partial Q}{\partial z}\right)
dx\wedge dy
\end{align*}
Thus we obtain the complex version of Stokes' Theorem,
\[
\int_{\partial X}Pdz+Qd\overline{z}=\int_{X}d\left(  Pdz+Qd\overline
{z}\right)  =2i\int_{X}\left(  \tfrac{\partial P}{\partial\bar{z}}%
-\tfrac{\partial Q}{\partial z}\right)  .
\]
From this, we get
\begin{align*}
\int_{X}\tfrac{\partial u}{\partial\bar{z}}f  &  =\int_{X}\tfrac{\partial
}{\partial\bar{z}}\left(  uf\right)  -\int_{X}u\tfrac{\partial f}{\partial
\bar{z}}=\frac{1}{2i}\int_{\partial X}uf\,dz-\int_{X}u\tfrac{\partial
f}{\partial\bar{z}}\text{\ and}\\
\int_{X}\tfrac{\partial v}{\partial z}g  &  =\int_{X}\tfrac{\partial}{\partial
z}\left(  vg\right)  -\int_{X}v\tfrac{\partial g}{\partial z}=\frac{-1}%
{2i}\int_{\partial X}vgd\bar{z}-\int_{X}v\tfrac{\partial g}{\partial z}.
\end{align*}
Hence,
\[
\int_{X}\left(  \tfrac{\partial u}{\partial\bar{z}}f+\tfrac{\partial
v}{\partial z}g\right)  =-\int_{X}\left(  u\tfrac{\partial f}{\partial\bar{z}%
}+v\tfrac{\partial g}{\partial z}\right)  +\frac{1}{2i}\int_{\partial
X}\left(  uf\,dz-vgd\bar{z}\right)  .
\]
Assuming (\ref{CokerCond}), we have
\begin{align*}
0  &  =\int_{X}\left(  \tfrac{\partial u}{\partial\bar{z}}f+\tfrac{\partial
v}{\partial z}g\right)  +\int_{\partial X}(u-v)h\\
&  =-\int_{X}\left(  u\tfrac{\partial f}{\partial\bar{z}}+v\tfrac{\partial
g}{\partial z}\right)  +\frac{1}{2i}\int_{\partial X}\left(  uf\,dz-vg\,d\bar
{z}\right)  +\int_{\partial X}(u-v)h
\end{align*}
By considering $u$ and $v$ with compact support inside the open disk, we
deduce that $\frac{\partial f}{\partial\bar{z}}=0$ and $\frac{\partial
g}{\partial z}=0$ (i.e., $f$ and $g$ are analytic and conjugate analytic
respectively). Thus, (\ref{CokerCond}) implies
\[
0=\frac{1}{2i}\int_{\partial X}\left(  uf\,dz-vg\,d\bar{z}\right)
+\int_{\partial X}(u-v)h,
\]
for all $u,v\in C^{\infty}(X)$. Choosing $v=u,$ we have
\begin{align*}
0  &  =\frac{1}{2i}\int_{\partial X}u\left(  f\,dz-g\,d\bar{z}\right)  \text{
for all }u\Rightarrow f\,dz=g\,d\bar{z}\text{ on }\partial X\\
&  \Rightarrow f\left(  e^{i\theta}\right)  \,ie^{i\theta}d\theta=-g\left(
e^{i\theta}\right)  ie^{-i\theta}d\theta\Rightarrow f\left(  e^{i\theta
}\right)  e^{i\theta}=-g\left(  e^{i\theta}\right)  e^{-i\theta}.
\end{align*}
However, since $f$ is analytic, the Fourier series of $f\left(  e^{i\theta
}\right)  e^{i\theta}$ has a nonzero coefficient for $e^{im\theta}$ only when
$m>0,$ and since $g$ is conjugate analytic, $g\left(  e^{i\theta}\right)
e^{-i\theta}$ only has a nonzero coefficient for $e^{im\theta}$ only when
$m<0$. Thus, $f=g=0$. Choosing $v=-u$, (\ref{CokerCond}) then yields
\[
0=\int_{\partial X}(u-v)h=2\int_{\partial X}uh\text{\ \ for all }u\in
C^{\infty}(X)\Rightarrow h=0.
\]
\end{example}

\begin{remark}
In engineering one calls a system of separate differential equations
\[%
\begin{array}
[c]{c}%
Pu=f\\
Qv=g,
\end{array}
\]
which are related by a ``transfer condition'' $R(u,v)=h$, a ``coupling
problem''; when the domains of $u$ and $v$ are different, but have a common
boundary (or boundary part) on which the transfer condition is defined, then
we have a ``transmission problem''; e.g., see \cite[p.7 ff]{Boo72}. Thus, we
may think of $T$ as an operator for a problem on the spherical surface
$X\cup_{\partial X}X$ with different behavior on the upper and lower
hemispheres, but with a fixed coupling along the equator.\bigskip
\end{remark}

\subsection{The Index of Twisted Dirac Operators on Closed
Manifolds\label{ss:IndexTDOonClosedMfds}}

Recall that the index of a Fredholm operator is a measure of its asymmetry: it
is defined by the difference between the dimension of the kernel (the null
space) of the operator and the dimension of the kernel of the adjoint operator
(= the codimension of the range). So, the index vanishes for self-adjoint
Fredholm operators. For an elliptic differential or pseudo-differential
operator on a closed manifold $M$, the index is finite and depends only on the
homotopy class of the principal symbol $\sigma$ of the operator over the
cotangent sphere bundle $S^{\ast}M$. It follows (see \cite[p. 257]{LaMi}) that
the index for elliptic \textit{differential} operators always vanishes on
\textit{closed} odd-dimensional manifolds. On even-dimensional manifolds one
has the Atiyah--Singer Index Theorem which expresses the index in explicit
topological terms, involving the Todd class defined by the Riemannian
structure of $M$ and the Chern class defined by gluing two copies of a bundle
over $S^{\ast}M$ by $\sigma$.

The original approach in proving the Index Theorem in the work of Atiyah and
Singer \cite{AtSi69} is based on the following clever strategy. The invariance
of the index under homotopy implies that the index (say, the \textit{analytic
index}) of an elliptic operator is stable under rather dramatic, but
continuous, changes in its principal symbol while maintaining ellipticity.
Moreover, the index is functorial with respect to certain operations, such as
addition and composition. Thus, the indices of elliptic operators transform
predictably under various global operations (or relations) such as direct
sums, embedding and cobordism. Using $K$-theory, a topological invariant (say,
the \textit{topological index}) with the same transformation properties under
these global operation is built from the symbol of the elliptic operator. It
turns out that the global operations are sufficient to construct enough vector
bundles and elliptic operators to deduce the Atiyah--Singer Index Theorem
(i.e., \textit{analytic index = topological index})\textit{.} With the aid of
Bott periodicity, it suffices to check that the two indices are the same in
the trivial case where the base manifold is just a single point. A
particularly nice exposition of this approach is found in E. Guentner's
article \cite{Gu93} following an argument of P. Baum.

Not long after this first proof (given in quite different variants), there
emerged a fundamentally different means of proving the Atiyah--Singer Index
Theorem, namely the \textit{heat kernel method}. This is outlined here in the
important case of the chiral half $\mathcal{D}^{+}$ of a twisted Dirac
operator $\mathcal{D}$. (In the index theory of closed manifolds, one usually
studies the index of a chiral half $\mathcal{D}^{+}$ instead of the total
Dirac operator $\mathcal{D}$, since $\mathcal{D}$ is symmetric for compatible
connections and then $\operatorname{index}$ $\mathcal{D}=0$.) The heat kernel
method had its origins in the late 1960s (e.g., in \cite{McK-Si}) and was
pioneered in the works \cite{Pa71}, \cite{Gi73}, \cite{ABP73}, etc.. In the
final analysis, it is debatable as to whether this method is really much
shorter or better. This depends on the background and tastes of the beholder.
Geometers and analysts (as opposed to topologists) are likely to find the heat
kernel method\textit{\ }appealing, because K-theory, Bott periodicity and
cobordism theory are avoided, not only for geometric operators which are
expressible in terms of twisted Dirac operators, but also largely for more
general elliptic pseudo-differential operators, as Melrose has done in
\cite{Me93}. Moreover, the heat method gives the index of a ``geometric''
elliptic differential operator naturally as the integral of a characteristic
form (a polynomial of curvature forms) which is expressed solely in terms of
the geometry of the operator itself (e.g., curvatures of metric tensors and
connections). One does not destroy the geometry of the operator by taking
advantage of the fact that it can be suitably deformed without altering the
index. Rather, in the heat kernel approach, the invariance of the index under
changes in the geometry of the operator is a consequence of the index formula
itself rather than a means of proof. However, considerable analysis and effort
are needed to obtain the heat kernel for $e^{-t\mathcal{D}^{2}}$ and to
establish its asymptotic expansion as $t\rightarrow0^{+}$. Also, it can be
argued that in some respects the K-theoretical embedding/cobordism methods are
more forceful and direct. Moreover, in \cite{LaMi}, we are cautioned that the
index theorem for families (in its strong form) generally involves torsion
elements in K-theory that are not detectable by cohomological means, and hence
are not computable in terms of local densities produced by heat asymptotics.
Nevertheless, when this difficulty does not arise, the K-theoretical
expression for the topological index may be less appealing than the integral
of a characteristic form, particularly for those who already understand and
appreciate the geometrical formulation of characteristic classes. All disputes
aside, the student who learns \textit{both} approaches and formulations will
be more accomplished (and probably older).

The classical geometric operators such as the Hirzebruch signature operator,
the de Rham operator, the Dolbeaut operator and even the Yang-Mills operator
can all be locally expressed in terms of chiral halves of twisted Dirac
operators. Thus, we will focus on index theory for such operators. The index
of any of these classical operators (and \textit{their} twists) can then be
obtained from the Local Index Theorem for twisted Dirac operators. This
theorem supplies a well-defined $n$-form on $M$, whose integral is the index
of the twisted Dirac operator. This $n$-form (or ``index density'') is
expressed in terms of forms for characteristic classes which are polynomials
in curvature forms. The Index Theorem thus obtained then becomes a formula
that relates a global invariant quantity, namely the index of an operator, to
the integral of a local quantity involving curvature. This is in the spirit of
the Gauss-Bonnet Theorem which can be considered a special case.\bigskip

\begin{definition}
\label{TwDirOpDefn}Let $M$ be an oriented Riemannian $n$-manifold ($n=2m$
even) with metric $h$, and oriented orthonormal frame bundle $FM$. Assume that
$M$ has a spin structure $P\rightarrow FM$, where $P$ is a principal
$\operatorname{Spin}\left(  n\right)  $-bundle and the projection
$P\rightarrow FM$ is a two-fold cover, equivariant with respect to
$\operatorname{Spin}\left(  n\right)  \rightarrow\operatorname{SO}\left(
n\right)  $. Furthermore, let $E\rightarrow M$ be a Hermitian vector bundle
with unitary connection $\varepsilon$. The \emph{twisted Dirac operator}
$\mathcal{D}$ associated with the above data is
\begin{equation}
\mathcal{D}:=\left(  1\otimes\mathbf{c}\right)  \circ\nabla:C^{\infty}\left(
E\otimes\Sigma\left(  M\right)  \right)  \rightarrow C^{\infty}\left(
E\otimes\Sigma\left(  M\right)  \right)  . \label{TwDirOp}%
\end{equation}
Here, $\Sigma\left(  M\right)  $ is the spin bundle over $M$ associated to
$P\rightarrow FM\rightarrow M$ via the spinor representation
$\operatorname{Spin}\left(  n\right)  \rightarrow\operatorname{End}\left(
\Sigma_{n}\right)  $,
\[
\mathbf{c}:C^{\infty}\left(  \Sigma\left(  M\right)  \otimes TM^{\ast}\right)
\rightarrow C^{\infty}\left(  \Sigma\left(  M\right)  \right)
\]
is Clifford multiplication, and
\[
\nabla:C^{\infty}\left(  E\otimes\Sigma\left(  M\right)  \right)  \rightarrow
C^{\infty}\left(  E\otimes\Sigma\left(  M\right)  \otimes TM^{\ast}\right)
\]
is the covariant derivative determined by the connection $\varepsilon$ and the
spinorial lift to $P$ of the Levi-Civita connection form, say $\theta$, on
$FM$.
\end{definition}

\bigskip

Note that $\mathcal{D}$ here is a special case of an operator of Dirac type
introduced in Subsection \ref{ss: Ops of Dirac Type} with $\frak{Cl}(M)$
acting on the second factor of $E\otimes\Sigma\left(  M\right)  $ which is
playing the role of $S$. Let $\Sigma^{\pm}\left(  M\right)  $ denote the
$\pm1$ eigenbundles of the complex Clifford volume element in $C^{\infty
}\left(  \frak{Cl}(M)\right)  $, given at a point $x\in M$ by $i^{m}%
e_{1}\cdots e_{n}$, where $e_{1},\ldots,e_{n}$ is an oriented, orthonormal
basis of $T_{x}M$. The $\Sigma^{\pm}\left(  M\right)  $ are the so-called
chiral halves of $\Sigma\left(  M\right)  =\Sigma^{+}\left(  M\right)
\oplus\Sigma^{+}\left(  M\right)  $. Since
\begin{align*}
&  \nabla\left(  C^{\infty}\left(  E\otimes\Sigma^{\pm}\left(  M\right)
\right)  \right)  \subseteq C^{\infty}\left(  E\otimes\Sigma^{\pm}\left(
M\right)  \otimes TM^{\ast}\right)  \text{ and}\\
&  \left(  1\otimes\mathbf{c}\right)  \left(  C^{\infty}\left(  E\otimes
\Sigma^{\pm}\left(  M\right)  \otimes TM^{\ast}\right)  \right)  \subseteq
C^{\infty}\left(  E\otimes\Sigma^{\mp}\left(  M\right)  \right)  ,
\end{align*}
we have
\[
\mathcal{D}=\mathcal{D}^{+}\oplus\mathcal{D}^{-}\text{, where }\mathcal{D}%
^{\pm}:C^{\infty}\left(  E\otimes\Sigma^{\pm}\left(  M\right)  \right)
\rightarrow C^{\infty}\left(  E\otimes\Sigma^{\mp}\left(  M\right)  \right)
.
\]
The symbol of the first-order differential operator $\mathcal{D}$ is computed
as follows. For $\phi\in C^{\infty}\left(  M\right)  $ with $\phi\left(
x\right)  =0$ and $\psi\in C^{\infty}\left(  E\otimes\Sigma\left(  M\right)
\right)  ,$ we have at $x$%
\begin{align*}
\left(  1\otimes\mathbf{c}\right)  \circ\nabla\left(  \phi\psi\right)   &
=\left(  1\otimes\mathbf{c}\right)  \circ\left(  \left(  d\phi\right)
\psi+\phi\nabla\psi\right)  =\left(  1\otimes\mathbf{c}\right)  \circ\left(
d\phi\right)  \psi\\
&  =\left(  1\otimes\mathbf{c}\left(  d\phi\right)  \right)  \psi.
\end{align*}
Thus, the symbol $\sigma\left(  \mathcal{D}\right)  :T_{x}M^{\ast}%
\rightarrow\operatorname{End}\left(  \Sigma\left(  M\right)  \right)  $ at the
covector $\xi_{x}\in T_{x}M^{\ast}$ is given by
\[
\sigma\left(  \mathcal{D}\right)  \left(  \xi_{x}\right)  =1\otimes
\mathbf{c}\left(  \xi_{x}\right)  \in\operatorname{End}\left(  E_{x}%
\otimes\Sigma_{x}\right)  .
\]
For $\xi_{x}\neq0$, $\sigma\left(  \mathcal{D}\right)  \left(  \xi_{x}\right)
$ is an isomorphism, since
\[
\sigma\left(  \mathcal{D}\right)  \left(  \xi_{x}\right)  \circ\sigma\left(
\mathcal{D}\right)  \left(  \xi_{x}\right)  =1\otimes\mathbf{c}\left(  \xi
_{x}\right)  ^{2}=-\left|  \xi_{x}\right|  ^{2}\operatorname{I}.
\]
Thus, $\mathcal{D}$ is an elliptic operator. Moreover, since $\sigma\left(
\mathcal{D}^{+}\right)  $ and $\sigma\left(  \mathcal{D}^{-}\right)  $ are
restrictions of $\sigma\left(  \mathcal{D}\right)  $, it follows that
$\mathcal{D}^{+}$ and $\mathcal{D}^{-}$ are elliptic. It can be shown that
$\mathcal{D}$ is formally self-adjoint, and $\mathcal{D}^{+}$ and
$\mathcal{D}^{-}$ are formal adjoints of each other (see \cite{LaMi}). We also
have a pair of self-adjoint elliptic operators
\begin{align*}
\mathcal{D}_{+}^{2}  &  :=\mathcal{D}^{2}|C^{\infty}\left(  E\otimes\Sigma
^{+}\left(  M\right)  \right)  =\mathcal{D}^{-}\circ\mathcal{D}^{+}\text{ }\\
\mathcal{D}_{-}^{2}  &  :=\mathcal{D}^{2}|C^{\infty}\left(  E\otimes\Sigma
^{-}\left(  M\right)  \right)  =\mathcal{D}^{+}\circ\mathcal{D}^{-}.
\end{align*}
For $\lambda\in\mathbb{C}$, let
\[
V_{\lambda}\left(  \mathcal{D}_{\pm}^{2}\right)  :=\left\{  \psi\in C^{\infty
}\left(  E\otimes\Sigma^{\pm}\left(  M\right)  \right)  \mid\mathcal{D}_{\pm
}^{2}\psi=\lambda\psi\right\}  .
\]
From the general theory of formally self-adjoint, elliptic operators on
compact manifolds, we know that
\[
\operatorname{Spec}\left(  \mathcal{D}_{\pm}^{2}\right)  =\left\{  \lambda
\in\mathbb{C}\mid V_{\lambda}\left(  \mathcal{D}_{\pm}^{2}\right)
\neq\left\{  0\right\}  \right\}  \text{.}%
\]
consists of the eigenvalues of $\mathcal{D}_{\pm}^{2}$ and is a discrete
subset of $\left[  0,\infty\right)  $, the eigenspaces $V_{\lambda}\left(
\mathcal{D}_{\pm}^{2}\right)  $ are finite-dimensional, and an $L^{2}\left(
E\otimes\Sigma^{\pm}\left(  M\right)  \right)  $-complete orthonormal set of
vectors can be selected from the $V_{\lambda}\left(  \mathcal{D}_{\pm}%
^{2}\right)  $. Note that $\mathcal{D}^{+}\left(  V_{\lambda}\left(
\mathcal{D}_{+}^{2}\right)  \right)  \subseteq V_{\lambda}\left(
\mathcal{D}_{-}^{2}\right)  $, since for $\psi\in V_{\lambda}\left(
\mathcal{D}_{+}^{2}\right)  $
\begin{align*}
\mathcal{D}_{-}^{2}\left(  \mathcal{D}^{+}\psi\right)   &  =\left(
\mathcal{D}^{+}\circ\mathcal{D}^{-}\right)  \left(  \mathcal{D}^{+}%
\psi\right)  =\mathcal{D}^{+}\left(  \left(  \mathcal{D}^{-}\circ
\mathcal{D}^{+}\right)  \left(  \psi\right)  \right) \\
&  =\mathcal{D}^{+}\left(  \mathcal{D}_{+}^{2}\left(  \psi\right)  \right)
=\mathcal{D}^{+}\left(  \lambda\psi\right)  =\lambda\mathcal{D}^{+}\left(
\psi\right)  ,
\end{align*}
and similarly $\mathcal{D}^{-}\left(  V_{\lambda}\left(  \mathcal{D}_{-}%
^{2}\right)  \right)  \subseteq V_{\lambda}\left(  \mathcal{D}_{+}^{2}\right)
$. For $\lambda\neq0$,
\[
\mathcal{D}^{\pm}|V_{\lambda}\left(  \mathcal{D}_{\pm}^{2}\right)
:V_{\lambda}\left(  \mathcal{D}_{\pm}^{2}\right)  \rightarrow V_{\lambda
}\left(  \mathcal{D}_{\mp}^{2}\right)
\]
is an isomorphism, since it has inverse $\frac{1}{\lambda}\mathcal{D}^{\mp} $.
Thus the set of nonzero eigenvalues (and their multiplicities) of
$\mathcal{D}_{+}^{2}$ coincides with that of $\mathcal{D}_{-}^{2}$. However,
in general
\[
\dim V_{0}\left(  \mathcal{D}_{+}^{2}\right)  -\dim V_{0}\left(
\mathcal{D}_{-}^{2}\right)  =\dim\operatorname{Ker}\left(  \mathcal{D}_{+}%
^{2}\right)  -\dim\operatorname{Ker}\left(  \mathcal{D}_{-}^{2}\right)
=\operatorname{index}\left(  \mathcal{D}^{+}\right)  \neq0.
\]
Since $\dim V_{\lambda}\left(  \mathcal{D}_{+}^{2}\right)  -\dim V_{\lambda
}\left(  \mathcal{D}_{-}^{2}\right)  =0$ for $\lambda\neq0$, obviously
\begin{align*}
\operatorname{index}\left(  \mathcal{D}^{+}\right)   &  =\dim V_{0}\left(
\mathcal{D}_{+}^{2}\right)  -\dim V_{0}\left(  \mathcal{D}_{-}^{2}\right) \\
&  =\sum\nolimits_{\lambda\in\operatorname{Spec}\left(  \mathcal{D}_{+}%
^{2}\right)  }e^{-t\lambda}\left(  \dim V_{\lambda}\left(  \mathcal{D}_{+}%
^{2}\right)  -\dim V_{\lambda}\left(  \mathcal{D}_{-}^{2}\right)  \right)  .
\end{align*}
This may seem like a very inefficient way to write $\operatorname{index}%
\left(  \mathcal{D}^{+}\right)  $, but the point is that the sum can be
expressed as the integral of the supertrace of the heat kernel for the
spinorial heat equation $\tfrac{\partial\psi}{\partial t}=-\mathcal{D}^{2}%
\psi$, from which the Local Index Theorem (Theorem \ref{LocIndThm} below) for
$\mathcal{D}^{+}$ will eventually follow. However, first the existence of the
heat kernel needs to be established.

Let the \textit{positive} eigenvalues of $\mathcal{D}_{\pm}^{2}$ be placed in
a sequence $0<\lambda_{1}\leq\lambda_{2}\leq\lambda_{3}\leq\ldots$ where each
eigenvalue is repeated according to its multiplicity. Let $u_{1}^{\pm},$
$u_{2}^{\pm},\ldots$ be an $L^{2}$-orthonormal sequence in $C^{\infty}\left(
E\otimes\Sigma^{+}\left(  M\right)  \right)  $ with $\mathcal{D}_{\pm}%
^{2}\left(  u_{j}^{\pm}\right)  =\lambda_{j}u_{j}^{\pm}$ (i.e., $u_{j}^{\pm
}\in V_{\lambda_{j}}\left(  \mathcal{D}_{\pm}^{2}\right)  $). We let
$u_{0_{1}}^{+},\ldots,u_{0_{n^{+}}}^{+}$ be an $L^{2}$-orthonormal basis of
$\operatorname{Ker}\mathcal{D}_{+}^{2}=\operatorname{Ker}\mathcal{D}_{+}$, and
$u_{0_{1}}^{-},\ldots,u_{0_{n^{-}}}^{-}$ be an $L^{2}$-orthonormal basis of
$\operatorname{Ker}\mathcal{D}_{-}^{2}=\operatorname{Ker}\mathcal{D}_{-}$. We
can pull back the bundle $E\otimes\Sigma^{\pm}\left(  M\right)  $ via either
of the projections $M\times M\times\left(  0,\infty\right)  \rightarrow M$
given by $\pi_{1}\left(  x,y,t\right)  :=x$ and $\pi_{2}\left(  x,y,t\right)
:=y$ and take the tensor product of the results to form a bundle
\[
\mathcal{K}^{\pm}:=\pi_{1}^{\ast}\left(  E\otimes\Sigma^{\pm}\left(  M\right)
\right)  \otimes\pi_{2}^{\ast}\left(  E\otimes\Sigma^{\pm}\left(  M\right)
\right)  \rightarrow M\times M\times\left(  0,\infty\right)  .
\]
Note that for $x\in M$, the Hermitian inner product $\left\langle
\;,\;\right\rangle _{x}$ on $\left(  E\otimes\Sigma^{\pm}\left(  M\right)
\right)  _{x}$ gives us a conjugate-linear map $\psi\mapsto\psi^{\ast}\left(
\cdot\right)  :=\left\langle \cdot,\psi\right\rangle _{x}$ from $\left(
E\otimes\Sigma^{\pm}\left(  M\right)  \right)  _{x}$ to its dual $\left(
E\otimes\Sigma^{\pm}\left(  M\right)  \right)  _{x}^{\ast}$. Thus, we can (and
do) make the identifications
\begin{align*}
&  \pi_{1}^{\ast}\left(  E\otimes\Sigma^{\pm}\left(  M\right)  \right)
\otimes\pi_{2}^{\ast}\left(  E\otimes\Sigma^{\pm}\left(  M\right)  \right)
\cong\left(  \pi_{1}^{\ast}\left(  E\otimes\Sigma^{\pm}\left(  M\right)
\right)  \right)  ^{\ast}\otimes\pi_{2}^{\ast}\left(  E\otimes\Sigma^{\pm
}\left(  M\right)  \right) \\
&  \cong\operatorname{Hom}\left(  \pi_{1}^{\ast}\left(  E\otimes\Sigma^{\pm
}\left(  M\right)  \right)  ,\pi_{2}^{\ast}\left(  E\otimes\Sigma^{\pm}\left(
M\right)  \right)  \right)  .
\end{align*}
The full proof of the following Proposition \ref{HeatEqnKerProp} will be found
in \cite{BlBo03}, but it is already contained in \cite{Gi95} for readers of
sufficient background.

\begin{proposition}
\label{HeatEqnKerProp}For $t>t_{0}>0$, the series $k^{\prime\pm},$ defined by
\[
k^{\prime\pm}\left(  x,y,t\right)  :=\sum_{j=1}^{\infty}e^{-\lambda_{j}t}%
u_{j}^{\pm}\left(  x\right)  \otimes u_{j}^{\pm}\left(  y\right)  ,
\]
converges uniformly in $C^{q}(\mathcal{K}^{\pm}|M\times M\times(t_{0}%
,\infty))$ for all $q\geq0$. Hence $k^{\prime\pm}\in C^{\infty}\left(
\mathcal{K}^{\pm}\right)  $, and $($for $t>0)$
\begin{equation}
\dfrac{\partial}{\partial t}k^{\prime\pm}\left(  x,y,t\right)  =-\sum
_{j=1}^{\infty}\lambda_{j}e^{-\lambda_{j}t}u_{j}^{\pm}\left(  x\right)
\otimes u_{j}^{\pm}\left(  y\right)  =-\mathcal{D}_{\pm}^{2}k^{\prime\pm
}\left(  x,y,t\right)  . \label{HeatEqnKer}%
\end{equation}
\end{proposition}

\begin{definition}
\label{SpinHeatKerDefn}The \emph{positive and negative twisted spinorial heat
kernels} $($or the \emph{heat kernels for}\textbf{\ }$\mathcal{D}_{\pm}^{2})$
$k^{\pm}\in C^{\infty}\left(  \mathcal{K}^{\pm}\right)  $ are given by
\begin{align*}
&  k^{\pm}\left(  x,y,t\right)  :=\sum\nolimits_{i=1}^{n^{\pm}}u_{0_{i}}^{\pm
}\left(  x\right)  \otimes u_{0_{i}}^{\pm}\left(  y\right)  +k^{\prime\pm
}\left(  x,y,t\right) \\
&  =\sum\nolimits_{i=1}^{n^{\pm}}u_{0_{i}}^{\pm}\left(  x\right)  \otimes
u_{0_{i}}^{\pm}\left(  y\right)  +\sum\nolimits_{j=1}^{\infty}e^{-\lambda
_{j}t}u_{j}^{\pm}\left(  x\right)  \otimes u_{j}^{\pm}\left(  y\right)
\text{\ for }t>0\text{.}%
\end{align*}
The \emph{total twisted spinorial heat kernel} (or the \emph{heat kernel
for}\textbf{\ }$\mathcal{D}^{2}$) is
\begin{align}
k  &  =\left(  k^{+},k^{-}\right)  \in C^{\infty}\left(  \mathcal{K}%
^{+}\right)  \oplus C^{\infty}\left(  \mathcal{K}^{-}\right)  \cong C^{\infty
}\left(  \mathcal{K}^{+}\oplus\mathcal{K}^{-}\right)  \subseteq C^{\infty
}\left(  \mathcal{K}\right)  ,\nonumber\\
&  \text{where }\mathcal{K}:=\mathcal{K}^{+}\oplus\mathcal{K}^{-}%
=\operatorname{Hom}\left(  \pi_{1}^{\ast}\left(  E\otimes\Sigma\left(
M\right)  \right)  ,\pi_{2}^{\ast}\left(  E\otimes\Sigma\left(  M\right)
\right)  \right)  . \label{HeatKerD2}%
\end{align}
\end{definition}

\bigskip

The terminology is justified in view of the following, whose proof is to be
found in \cite{BlBo03}.

\begin{proposition}
\label{HeatKerInitProp}Let $\psi_{0}^{\pm}\in C^{\infty}\left(  E\otimes
\Sigma^{\pm}\left(  M\right)  \right)  $ and let
\[
\psi^{\pm}\left(  x,t\right)  =\int_{M}\left\langle k^{\pm}\left(
x,y,t\right)  ,\psi_{0}^{\pm}\left(  y\right)  \right\rangle _{y}\text{ }%
\nu_{y}.
\]
Then for $t>0$, $\psi^{\pm}$ solves the heat equation with initial spinor
field $\psi_{0}^{\pm}:$
\begin{align*}
&  \dfrac{\partial\psi^{\pm}}{\partial t}=-\mathcal{D}_{\pm}^{2}%
\psi\;\text{and}\\
&  \lim_{t\rightarrow0^{+}}\psi^{\pm}\left(  \cdot,t\right)  =\psi_{0}^{\pm
}\text{ in }C^{q}\text{ for all }q\geq0\text{.}%
\end{align*}
Moreover, for $\psi_{0}\in C^{\infty}\left(  E\otimes\Sigma\left(  M\right)
\right)  $ and
\[
\psi\left(  x,t\right)  :=\int_{M}\left\langle k\left(  x,y,t\right)
,\psi_{0}\right\rangle \text{ }\nu_{y},
\]
we have $\dfrac{\partial\psi}{\partial t}=-\mathcal{D}^{2}\psi$ and
$\lim_{t\rightarrow0^{+}}\psi\left(  \cdot,t\right)  =\psi_{0}\left(
\cdot\right)  $ in $C^{q}$ for all $q\geq0$.
\end{proposition}

\bigskip

For any finite dimensional Hermitian vector space $\left(  V,\left\langle
\cdot,\cdot\right\rangle \right)  $ with orthonormal basis $e_{1},\ldots
,e_{N}$, we have (for $v\in V$)
\begin{align*}
&  \operatorname{Tr}\left(  v^{\ast}\otimes v\right)  =\sum\nolimits_{i=1}%
^{N}\left\langle \left(  v^{\ast}\otimes v\right)  \left(  e_{i}\right)
,e_{i}\right\rangle =\sum\nolimits_{i=1}^{N}\left\langle v^{\ast}\left(
e_{i}\right)  v,e_{i}\right\rangle \\
&  =\sum\nolimits_{i=1}^{N}\left\langle \left\langle e_{i},v\right\rangle
v,e_{i}\right\rangle =\sum\nolimits_{i=1}^{N}\left\langle e_{i},v\right\rangle
\left\langle v,e_{i}\right\rangle =\sum\nolimits_{i=1}^{N}\left|  \left\langle
e_{i},v\right\rangle \right|  ^{2}=\left|  v\right|  ^{2}.
\end{align*}
In particular, $k^{\pm}\left(  x,x,t\right)  \in\operatorname{End}\left(
\left(  E\otimes\Sigma^{\pm}\left(  M\right)  \right)  _{x}\right)  $ and
\[
\operatorname{Tr}\left(  k^{\pm}\left(  x,x,t\right)  \right)  =\sum
\nolimits_{i=1}^{n^{\pm}}\left|  u_{0_{i}}^{\pm}\left(  x\right)  \right|
^{2}+\sum\nolimits_{j=1}^{\infty}e^{-\lambda_{j}t}\left|  u_{j}^{\pm}\left(
x\right)  \right|  ^{2}.
\]
Since this series converges uniformly and $\left\|  u_{0_{i}}^{\pm}\right\|
_{2,0}=\left\|  u_{j}^{\pm}\right\|  _{2,0}=1$, we have
\[
\int_{M}\operatorname{Tr}\left(  k^{\pm}\left(  x,x,t\right)  \right)  \text{
}\nu_{x}=n^{\pm}+\sum\nolimits_{j=1}^{\infty}e^{-\lambda_{j}t}<\infty.
\]
For $t>0,$ we define the bounded operator $e^{-t\mathcal{D}_{\pm}^{2}}%
\in\operatorname{End}\left(  L^{2}\left(  E\otimes\Sigma^{\pm}\left(
M\right)  \right)  \right)  $ by
\[
e^{-t\mathcal{D}_{\pm}^{2}}\left(  \psi^{\pm}\right)  =\sum\nolimits_{i=1}%
^{n^{\pm}}\left(  u_{0_{i}}^{\pm},\psi_{0}^{\pm}\right)  u_{0_{i}}^{\pm}%
+\sum\nolimits_{j=1}^{\infty}e^{-\lambda_{j}t}\left(  u_{j}^{\pm},\psi
_{0}^{\pm}\right)  u_{j}^{\pm}.
\]
Note that $e^{-t\mathcal{D}_{\pm}^{2}}$ is of trace class, since
\[
\operatorname{Tr}\left(  e^{-t\mathcal{D}_{\pm}^{2}}\right)  =n^{\pm}%
+\sum\nolimits_{j=1}^{\infty}e^{-\lambda_{j}t}=\int_{M}\operatorname{Tr}%
\left(  k^{\pm}\left(  x,x,t\right)  \right)  \text{ }\nu_{x}<\infty.
\]
Now, we have
\begin{align}
\operatorname{index}\left(  \mathcal{D}^{+}\right)   &  =\dim V_{0}\left(
\mathcal{D}_{+}^{2}\right)  -\dim V_{0}\left(  \mathcal{D}_{-}^{2}\right)
\nonumber\\
&  =n^{+}-n^{-}+\sum\nolimits_{j=1}^{\infty}\left(  e^{-\lambda_{j}%
t}-e^{-\lambda_{j}t}\right) \nonumber\\
&  =n^{+}+\sum\nolimits_{j=1}^{\infty}e^{-\lambda_{j}t}-\left(  n^{-}%
+\sum\nolimits_{j=1}^{\infty}e^{-\lambda_{j}t}\right) \nonumber\\
&  =\int_{M}\left(  \operatorname{Tr}\left(  k^{+}\left(  x,x,t\right)
\right)  -\operatorname{Tr}\left(  k^{-}\left(  x,x,t\right)  \right)
\right)  \text{ }\nu_{x}. \label{IndexStr}%
\end{align}
Since $\mathcal{D}^{2}=\mathcal{D}_{+}^{2}\oplus\mathcal{D}_{-}^{2}$, we also
have the trace-class operator $e^{-t\mathcal{D}^{2}}\in\operatorname{End}%
\left(  L^{2}\left(  E\otimes\Sigma\left(  M\right)  \right)  \right)  $,
whose trace is given by
\[
\operatorname{Tr}\left(  e^{-t\mathcal{D}^{2}}\right)  =\int_{M}%
\operatorname{Tr}\left(  k\left(  x,x,t\right)  \right)  \text{ }\nu_{x}%
=\int_{M}\left(  \operatorname{Tr}\left(  k^{+}\left(  x,x,t\right)  \right)
+\operatorname{Tr}\left(  k^{-}\left(  x,x,t\right)  \right)  \right)  \text{
}\nu_{x}.
\]
The \textit{supertrace} of $k\left(  x,x,t\right)  $ is defined by
\[
\operatorname{Str}\left(  k\left(  x,x,t\right)  \right)  :=\operatorname{Tr}%
\left(  k^{+}\left(  x,x,t\right)  \right)  -\operatorname{Tr}\left(
k^{-}\left(  x,x,t\right)  \right)  ,
\]
and in view of (\ref{IndexStr}), we have
\begin{equation}
\operatorname{index}\left(  \mathcal{D}^{+}\right)  =\int_{M}%
\operatorname{Str}\left(  k\left(  x,x,t\right)  \right)  \text{ }\nu_{x}.
\label{IndexStr2}%
\end{equation}
The left side is independent of $t$ and so the right side is also independent
of $t$. The main task now is to determine the behavior of $\operatorname{Str}%
\left(  k\left(  x,x,t\right)  \right)  $ as $t\rightarrow0^{+}$. We suspect
that for each $x\in M$, as $t\rightarrow0^{+}$, $k\left(  x,x,t\right)  $ and
$\operatorname{Str}\left(  k\left(  x,x,t\right)  \right)  $ are influenced
primarily by the geometry (e.g., curvature form $\Omega^{\theta} $ of $M$ with
metric $h$ and Levi-Civita connection $\theta$, and the curvature
$\Omega^{\varepsilon}$ of the unitary connection for $E$) near $x $, since the
heat sources of points far from $x$ are not felt very strongly at $x$ for
small $t$. Indeed, we will give an outline a proof of the following Local
Index Formula, the full proof of which will appear in \cite{BlBo03}.

\begin{theorem}
[The Local Index Theorem]\label{LocIndThm} In the notation of Definitions
\ref{TwDirOpDefn} and \ref{SpinHeatKerDefn}, let $\mathcal{D}:C^{\infty
}\left(  E\otimes\Sigma\left(  M\right)  \right)  \rightarrow C^{\infty
}\left(  E\otimes\Sigma\left(  M\right)  \right)  $ be a twisted Dirac
operator and let $k\in C^{\infty}\left(  \mathcal{K}\right)  $ be the heat
kernel for\textbf{\ }$\mathcal{D}^{2}$. If $\Omega^{\varepsilon}$ is the
curvature form of the unitary connection $\varepsilon$ for $E$ and
$\Omega^{\theta}$ is the curvature form of the Levi-Civita connection $\theta$
for $\left(  M,h\right)  $ with volume element $\nu$, then
\begin{equation}
\lim_{t\rightarrow0^{+}}\operatorname{Str}\left(  k\left(  x,x,t\right)
\right)  =\left\langle \operatorname{Tr}\left(  e^{i\Omega^{\varepsilon}/2\pi
}\right)  \wedge\det\left(  \frac{i\Omega^{\theta}/4\pi}{\sinh\left(
i\Omega^{\theta}/4\pi\right)  }\right)  ^{\tfrac{1}{2}},\nu_{x}\right\rangle .
\label{StrLim2}%
\end{equation}
\end{theorem}

\begin{remark}
As will be explained below, the right side is really the inner product, with
volume form $\nu$ at $x$, of the canonical form $\mathbf{ch}\left(
E,\varepsilon\right)  \smallsmile\widehat{\mathbf{A}}\left(  M,\theta\right)
$ (depending on the connections $\varepsilon$ for $E$ and the Levi-Civita
connection $\theta$ for the metric $h$) which represents $\mathbf{ch}\left(
E\right)  \smallsmile\widehat{\mathbf{A}}\left(  M\right)  $. As a
consequence, we obtain the Index Theorem for twisted Dirac operators from the
Local Index Formula in Corollary \ref{IndFormTwDiracCor} below. Thus the Local
Index Formula is stronger than the Index Theorem for twisted Dirac operators.
Indeed, the Local Index Formula yields the Index Theorem for elliptic
operators which are locally expressible as twisted Dirac operators or direct
sums of such.
\end{remark}

\begin{corollary}
[Index formula for twisted Dirac operators]\label{IndFormTwDiracCor}For an
oriented Riemannian $n$-manifold $M$ ($n$ even) with spin structure, and a
Hermitian vector bundle $E\rightarrow M$ with unitary connection, let
$\mathcal{D}=\mathcal{D}^{+}\oplus\mathcal{D}^{-}$ be the twisted Dirac
operator, with $\mathcal{D}^{+}:C^{\infty}\left(  E\otimes\Sigma^{+}\left(
M\right)  \right)  \rightarrow C^{\infty}\left(  E\otimes\Sigma^{-}\left(
M\right)  \right)  $. We have
\begin{equation}
\operatorname{index}\left(  \mathcal{D}^{+}\right)  =\left(  \mathbf{ch}%
\left(  E\right)  \smallsmile\widehat{\mathbf{A}}\left(  M\right)  \right)
\left[  M\right]  , \label{IndFormTwDirac}%
\end{equation}
where $\mathbf{ch}\left(  E\right)  $ is the total Chern character class of
$E$ and $\widehat{\mathbf{A}}\left(  M\right)  $ is the total $\widehat
{\mathbf{A}} $ class of $M$, both defined below. In particular, we obtain:
\begin{align*}
&
\begin{array}
[c]{c}%
n=2\Rightarrow\text{ }\operatorname{index}\left(  \mathcal{D}^{+}\right)
=ch_{1}\left(  E\right)  \left[  M\right]  =c_{1}\left(  E\right)  \left[
M\right]  \text{ and}%
\end{array}
\\
&
\begin{array}
[c]{l}%
n=4\Rightarrow\left\{
\begin{array}
[c]{l}%
\operatorname{index}\left(  \mathcal{D}^{+}\right)  =\left(  \mathbf{ch}%
\left(  E\right)  \smallsmile\widehat{\mathbf{A}}\left(  M\right)  \right)
\left[  M\right] \\
=\left(  -\dim E\cdot\tfrac{1}{24}p_{1}\left(  TM\right)  +ch_{2}\left(
E\right)  \right)  \left[  M\right] \\
=\left(  -\frac{\dim E}{24}p_{1}\left(  TM\right)  +\tfrac{1}{2}c_{1}\left(
E\right)  ^{2}-c_{2}\left(  E\right)  \right)  \left[  M\right]  .
\end{array}
\right.  \text{ }%
\end{array}
\end{align*}
\newline 
\end{corollary}

\begin{proof}
By (\ref{IndexStr2}), (\ref{StrLim2}) and the above Remark, we have
\begin{align*}
\operatorname{index}\left(  \mathcal{D}^{+}\right)   &  =\int_{M}%
\operatorname{Str}\left(  k\left(  x,x,t\right)  \right)  \nu_{x}\\
&  =\lim_{t\rightarrow0^{+}}\int_{M}\operatorname{Str}\left(  k\left(
x,x,t\right)  \right)  \nu_{x}=\int_{M}\lim_{t\rightarrow0^{+}}%
\operatorname{Str}\left(  k\left(  x,x,t\right)  \right)  \nu_{x}\\
&  =\int_{M}\left\langle \mathbf{ch}\left(  E,\varepsilon\right)
_{x}\smallsmile\widehat{\mathbf{A}}\left(  M,\theta\right)  _{x},\nu
_{x}\right\rangle \nu_{x}=\left(  \mathbf{ch}\left(  E\right)  \smallsmile
\widehat{\mathbf{A}}\left(  M\right)  \right)  \left[  M\right]  .
\end{align*}
\end{proof}

\bigskip

We now explain the meaning of the form
\[
\operatorname{Tr}\left(  e^{i\Omega^{\varepsilon}/2\pi}\right)  \wedge
\det\left(  \frac{i\Omega^{\theta}/4\pi}{\sinh\left(  i\Omega^{\theta}%
/4\pi\right)  }\right)  ^{\tfrac{1}{2}}.
\]
The first part $\operatorname{Tr}\left(  e^{i\Omega^{\varepsilon}/2\pi
}\right)  $ is relatively easy. We have (recall $2m=\dim M$)
\begin{equation}
e^{i\Omega^{\varepsilon}/2\pi}:=\sum_{k=0}^{\infty}\frac{1}{k!}\left(
\frac{i}{2\pi}\right)  ^{k}\Omega^{\varepsilon}\wedge\overset{k}{\cdots}%
\wedge\Omega^{\varepsilon}=\sum_{k=0}^{m}\frac{1}{k!}\left(  \frac{i}{2\pi
}\right)  ^{k}\Omega^{\varepsilon}\wedge\overset{k}{\cdots}\wedge
\Omega^{\varepsilon}\text{,} \label{ChernForm}%
\end{equation}
where $\Omega^{\varepsilon}\wedge\overset{k}{\cdots}\wedge\Omega^{\varepsilon
}\in\Omega^{2k}\left(  \operatorname{End}\left(  E\right)  \right)  $. Also
$\operatorname{Tr}\left(  i^{k}\Omega^{\varepsilon}\wedge\overset{k}{\cdots
}\wedge\Omega^{\varepsilon}\right)  \in\Omega^{2k}\left(  M\right)  $ and
\[
\operatorname{Tr}\left(  e^{i\Omega^{\varepsilon}/2\pi}\right)  \in
\bigoplus_{k=1}^{m}\Omega^{2k}\left(  M\right)  .
\]
This (by one of many equivalent definitions) is a representative of the total
Chern character $\mathbf{ch}\left(  E\right)  \in\bigoplus_{k=1}^{m}%
H^{2k}\left(  M,\mathbb{Q}\right)  $. The curvature $\Omega^{\theta}$ of the
Levi-Civita connection $\theta$ for the metric $h$ has values in the
skew-symmetric endomorphisms of $TM$; i.e., $\Omega^{\theta}\in\Omega
^{2}\left(  \operatorname{End}\left(  TM\right)  \right)  $. A skew-symmetric
endomorphism of $\mathbb{R}^{2m}$, say $B\in\mathfrak{so}\left(  n\right)  $,
has pure imaginary eigenvalues $\pm ir_{k},$ where $r_{k}\in\mathbb{R}$
$\left(  1\leq k\leq m\right)  $. Thus, $iB$ has real eigenvalues $\pm r_{k}$.
Now $\frac{z/2}{\sinh\left(  z/2\right)  }$ is a power series in $z$ with
radius of convergence $2\pi$. Thus, $\frac{isB/2}{\sinh\left(  isB/2\right)
}$ is defined for $s$ sufficiently small and has eigenvalues $\frac{r_{k}%
s/2}{\sinh\left(  r_{k}s/2\right)  }$ each repeated twice. Hence
\begin{align*}
\det\left(  \frac{isB/2}{\sinh\left(  isB/2\right)  }\right)   &  =\prod
_{k=1}^{m}\left(  \frac{r_{k}s/2}{\sinh\left(  r_{k}s/2\right)  }\right)
^{2}\text{ and }\\
\det\left(  \frac{isB/2}{\sinh\left(  isB/2\right)  }\right)  ^{\tfrac{1}{2}}
&  =\prod_{k=1}^{m}\frac{r_{k}s/2}{\sinh\left(  r_{k}s/2\right)  }.
\end{align*}
The last product is a power series in $s$ of the form
\begin{equation}
\prod_{k=1}^{m}\frac{r_{k}s/2}{\sinh\left(  r_{k}s/2\right)  }=\sum
_{k=0}^{\infty}a_{k}\left(  r_{1}^{2},\ldots,r_{m}^{2}\right)  s^{2k},
\label{asubk}%
\end{equation}
where the coefficient $a_{k}\left(  r_{1}^{2},\ldots,r_{m}^{2}\right)  $ is a
homogeneous, symmetric polynomial in $r_{1}^{2},\ldots,r_{m}^{2}$ of degree
$k$. One can always express any such a symmetric polynomial as a polynomial in
the elementary symmetric polynomials $\sigma_{1},\ldots,\sigma_{m}$ in
$r_{1}^{2},\ldots,r_{m}^{2}$, where
\[
\sigma_{1}=\sum\nolimits_{i=1}^{m}r_{i}^{2},\text{ }\sigma_{2}=\sum
\nolimits_{i<j}^{m}r_{i}^{2}r_{j}^{2},\text{ }\sigma_{2}=\sum\nolimits_{i<j<k}%
^{m}r_{i}^{2}r_{j}^{2}r_{k}^{2},\ldots.
\]
These in turn may be expressed in terms of $\operatorname{SO}\left(  n\right)
$-invariant polynomials in the entries of $B\in\frak{so}\left(  n\right)  $
via
\begin{align*}
\det\left(  \lambda I-B\right)   &  =\prod_{j=1}^{m}\left(  \lambda
+ir_{j}\right)  \left(  \lambda-ir_{j}\right)  =\prod_{j=1}^{m}\left(
\lambda^{2}+r_{j}^{2}\right) \\
&  =\sum_{k=1}^{m}\sigma_{k}\left(  r_{1}^{2},\ldots,r_{m}^{2}\right)
\lambda^{2\left(  m-k\right)  }.
\end{align*}
On the other hand,
\begin{align*}
\det\left(  \lambda I-B\right)   &  =\sum_{k=1}^{m}\left(  \frac{1}{\left(
2k\right)  !}\sum_{(i),(j)}\delta_{i_{1}\cdots i_{2k}}^{j_{1}\cdots j_{2k}%
}B_{\;j_{1}}^{i_{1}}\cdots B_{\;j_{2k}}^{i_{2k}}\right)  \lambda^{2\left(
m-k\right)  }\text{, and so}\\
\sigma_{k}\left(  r_{1}^{2},\ldots,r_{m}^{2}\right)   &  =\frac{1}{\left(
2k\right)  !}\sum_{(i),(j)}\delta_{i_{1}\cdots i_{2k}}^{j_{1}\cdots j_{2k}%
}B_{\;j_{1}}^{i_{1}}\cdots B_{\;j_{2k}}^{i_{2k}},
\end{align*}
where $(i)=\left(  i_{1},\cdots,i_{2k}\right)  $ is an ordered $2k$-tuple of
distinct elements of $\left\{  1,\ldots,2m\right\}  $ and $(j)$ is a
permutation of $(i)$ with sign $\delta_{i_{1}\cdots i_{2k}}^{j_{1}\cdots
j_{2k}}$. If we replace $B_{\;j}^{i}$ with the 2-form $\frac{1}{2\pi}\left(
\Omega^{\theta}\right)  _{\;j}^{i}$ relative to an orthonormal frame field, we
obtain the Pontryagin forms
\[
p_{k}\left(  \Omega^{\theta}\right)  :=\frac{1}{\left(  2\pi\right)
^{2k}\left(  2k\right)  !}\sum_{(i),(j)}\delta_{i_{1}\cdots i_{2k}}%
^{j_{1}\cdots j_{2k}}\Omega_{i_{1}j_{1}}^{\theta}\wedge\cdots\wedge
\Omega_{i_{2k}j_{2k}}^{\theta},
\]
which represent the Pontryagin classes of the $\operatorname{SO}\left(
n\right)  $ bundle $FM.$ Note that $p_{k}\left(  \Omega^{\theta}\right)  $ is
independent of the choice of framing by the $ad$-invariance of the polynomials
$\sigma_{k}$. If we express the $a_{k}\left(  r_{1}^{2},\ldots,r_{m}%
^{2}\right)  $ as polynomials, say $\mathcal{A}_{k}\left(  \sigma_{1}%
,\ldots,\sigma_{k}\right)  $, in the $\sigma_{j}$ ($j\leq k$), we can
ultimately write
\[
\det\left(  \frac{isB/2}{\sinh\left(  isB/2\right)  }\right)  ^{\tfrac{1}{2}%
}=\sum_{k=0}^{\infty}\mathcal{A}_{k}\left(  \sigma_{1},\ldots,\sigma
_{k}\right)  s^{2k}.
\]
Formally replacing $B$ by $\frac{1}{2\pi}\Omega^{\theta}$, we finally have the
reasonable definition
\[
\det\left(  \frac{i\Omega^{\theta}/4\pi}{\sinh\left(  i\Omega^{\theta}%
/4\pi\right)  }\right)  ^{\tfrac{1}{2}}:=\sum_{k=0}^{\infty}\mathcal{A}%
_{k}\left(  p_{1}\left(  \Omega^{\theta}\right)  ,\ldots,p_{k}\left(
\Omega^{\theta}\right)  \right)  ,
\]
where the $p_{j}\left(  \Omega^{\theta}\right)  $ are multiplied via wedge
product when evaluating the terms in the sum; the order of multiplication does
not matter since $p_{j}\left(  \Omega^{\theta}\right)  $ is of even degree
$4j$. Also, since $\mathcal{A}_{k}\left(  p_{1}\left(  \Omega^{\theta}\right)
,\ldots,p_{k}\left(  \Omega^{\theta}\right)  \right)  $ is a $4k$-form, there
are only a finite number of nonzero terms in the infinite sum. Abbreviating
$p_{j}\left(  \Omega^{\theta}\right)  $ simply by $p_{j}$, one finds
\begin{align}
\det\left(  \frac{i\Omega^{\theta}/4\pi}{\sinh\left(  i\Omega^{\theta}%
/4\pi\right)  }\right)  ^{\tfrac{1}{2}}  &  =1-\frac{1}{24}p_{1}+\frac
{1}{5760}\left(  7p_{1}^{2}-4p_{2}\right) \nonumber\\
&  -\frac{1}{967\,680}\left(  31p_{1}^{3}-44p_{1}p_{2}+16p_{3}\right)
+\cdots. \label{AhatForm}%
\end{align}
This (by one definition) represents the total $\widehat{A}$-class of $M$,
denoted by
\begin{equation}
\widehat{\mathbf{A}}\left(  M\right)  \in\bigoplus_{k=1}^{m}H^{2k}\left(
M,\mathbb{Q}\right)  , \label{AhatTotalClass}%
\end{equation}
where actually $\widehat{\mathbf{A}}\left(  M\right)  $ has only nonzero
components in $H^{2k}\left(  M,\mathbb{Q}\right)  $ when $k$ is even (or
$2k\equiv0\operatorname{mod}4$). In (\ref{StrLim2}) the multi-degree forms
(\ref{ChernForm}) and (\ref{AhatForm}) have been wedged, and the top
($2m$-degree) component (relative to the volume form) has been harvested.

\bigskip

We now turn to our outline of the proof of Theorem \ref{LocIndThm}. The
well-known heat kernel (or fundamental solution) for the ordinary heat
equation $u_{t}=\Delta u$ in Euclidean space $\mathbb{R}^{n}$, is given by
\begin{equation}
e\left(  x,y,t\right)  =\left(  4\pi t\right)  ^{-n/2}\exp\left(  -\left|
x-y\right|  ^{2}/4t\right)  . \label{FundHeatEuc1}%
\end{equation}
Since $H\left(  x,y,t\right)  $ only depends on $r=\left|  x-y\right|  $ and
$t,$ it is convenient to write
\begin{equation}
e\left(  x,y,t\right)  =\mathcal{E}\left(  r,t\right)  :=\left(  4\pi
t\right)  ^{-n/2}\exp\left(  -r^{2}/4t\right)  . \label{FundHeatEuc2}%
\end{equation}
We do not expect such a simple expression for the heat kernel $k=\left(
k^{+},k^{-}\right)  $ of Definition \ref{SpinHeatKerDefn}. However, it can be
shown that for $x,y\in M$ (of even dimension $n=2m$) with $r=d(x,y):=$
Riemannian distance from $x$ to $y$ sufficiently small, we have an asymptotic
expansion as $t\rightarrow0^{+}$ for $k\left(  x,y,t\right)  $ of the form
\begin{equation}
k\left(  x,y,t\right)  \sim H_{Q}\left(  x,y,t\right)  :=\mathcal{E}\left(
d(x,y),t\right)  \sum\nolimits_{j=0}^{Q}h_{j}\left(  x,y\right)  t^{j},
\label{AsyExpxy}%
\end{equation}
for any fixed integer $Q>m+4$, where
\[
h_{j}\left(  x,y\right)  \in\operatorname{Hom}\left(  \left(  E\otimes
\Sigma\left(  M\right)  \right)  _{x},\left(  E\otimes\Sigma\left(  M\right)
\right)  _{y}\right)  ,\text{ }j\in\left\{  0,1,\ldots,Q\right\}  .
\]
The meaning of $k\left(  x,y,t\right)  \sim H_{Q}\left(  x,y,t\right)  $ is
that for $d(x,y)$ and $t$ sufficiently small,
\[
\left|  k\left(  x,y,t\right)  -H_{Q}\left(  x,y,t\right)  \right|  \leq
C_{Q}\mathcal{E}\left(  d(x,y),t\right)  t^{Q+1}\leq C_{Q}t^{Q-m+1},
\]
where $C_{Q}$ is a constant, independent of $(x,y,t)$. We then have
\begin{equation}
k\left(  x,x,t\right)  \sim\left(  4\pi t\right)  ^{-m}\sum_{j=0}^{Q}%
h_{j}\left(  x,x\right)  t^{j}=\left(  4\pi\right)  ^{-m}\sum_{j=0}^{Q}%
h_{j}\left(  x,x\right)  t^{j-m}. \label{AsyExp}%
\end{equation}
Using (\ref{IndexStr2}), i.e., $\int_{M}\operatorname{Str}\left(  k\left(
x,x,t\right)  \right)  \,\nu_{x}=\operatorname{index}\left(  \mathcal{D}%
^{+}\right)  $ and (\ref{AsyExp}), we deduce that
\begin{align*}
&  \int_{M}\operatorname{Str}\left(  h_{j}\left(  x,x\right)  \right)
\,\nu_{x}=0\;\;\text{for}\;\;j\in\left\{  0,1,\ldots,m-1\right\}  ,\text{
while}\\
&  \left(  4\pi\right)  ^{-m}\int_{M}\operatorname{Str}\left(  h_{m}\left(
x,x\right)  \right)  \,\nu_{x}=\int_{M}\operatorname{Str}\left(  k\left(
x,x,t\right)  \right)  \,\nu_{x}=\operatorname{index}\left(  \mathcal{D}%
^{+}\right)  .
\end{align*}
Thus, to prove the Local Index Formula, it suffices to show that
\[
\left(  4\pi\right)  ^{-m}\operatorname{Str}\left(  h_{m}\left(  x,x\right)
\right)  =\left\langle \operatorname{Tr}\left(  e^{i\Omega^{\varepsilon}/2\pi
}\right)  \wedge\det\left(  \frac{i\Omega^{\theta}/4\pi}{\sinh\left(
i\Omega^{\theta}/4\pi\right)  }\right)  ^{\tfrac{1}{2}},\nu_{x}\right\rangle
.
\]
While this may not be the intellectual equivalent of climbing Mount Everest,
it is not for the faint of heart.

We choose a normal coordinate system $\left(  y^{1},\ldots,y^{n}\right)  $ in
a coordinate ball $\mathcal{B}$ centered at the fixed point $x\in M$, so that
$\left(  y^{1},\ldots,y^{n}\right)  =0$ at $x$. The coordinate fields
$\partial_{1}:=\partial/\partial y^{1},\ldots,\partial_{n}:=\partial/\partial
y^{n}$ are orthonormal at $x$, and for any fixed $y_{0}\in\mathcal{B}$ with
coordinates $\left(  y_{0}^{1},\ldots,y_{0}^{n}\right)  $, the curve $t\mapsto
t\left(  y_{0}^{1},\ldots,y_{0}^{n}\right)  $ is a geodesic through $x$. By
parallel translating the frame $\left(  \partial_{1},\ldots,\partial
_{n}\right)  $ at $x$ along these radial geodesics, we obtain an orthonormal
frame field $\left(  E_{1},\ldots,E_{n}\right)  $ on $\mathcal{B}$ which
generally does not coincide with $\left(  \partial_{1},\ldots,\partial
_{n}\right)  $ at points $y\in\mathcal{B}$ other than at $x$. The framing
$\left(  E_{1},\ldots,E_{n}\right)  $ defines a particularly nice section
$\mathcal{B}\rightarrow FM|B$ and we may lift this to a section $\mathcal{B}%
\rightarrow P|B$ of the spin structure, which enables us to view the space
$C^{\infty}\left(  \Sigma\left(  M\right)  |\mathcal{B}\right)  $\ of spinor
fields on $\mathcal{B}$ as $C^{\infty}\left(  \mathcal{B},\Sigma_{n}\right)
$, i.e., functions on $\mathcal{B}$ with values in the fixed spinor
representation vector space $\Sigma_{n}=\Sigma_{n}^{+}\oplus\Sigma_{n}^{-}$.
By similar radial parallel translation (with respect to the connection
$\varepsilon$) of an orthonormal basis of the twisting bundle fiber $E_{x}$,
we can identify $C^{\infty}\left(  E|\mathcal{B}\right)  $ with $C^{\infty
}\left(  \mathcal{B},\mathbb{C}^{N}\right)  $, where $N=\dim_{\mathbb{C}}E$.
The coordinate expressions for the curvatures $\Omega^{\theta}$,
$\Omega^{\varepsilon}$ and $\mathcal{D}^{2}$ are as simple as possible in this
so-called radial gauge.

With the above identifications, we proceed as follows. For $0\leq
Q\in\mathbb{Z}$, let $\Psi_{Q}\in C^{\infty}\left(  \mathcal{B}\times\left(
0,\infty\right)  ,\mathbb{C}^{N}\otimes\Sigma_{2m}\right)  $ be of the form
\[
\Psi_{Q}\left(  y,t\right)  :=\mathcal{E}\left(  r,t\right)  \sum
\nolimits_{k=0}^{Q}U_{k}\left(  y\right)  t^{k},
\]
where $U_{k}\in C^{\infty}\left(  \mathcal{B},\mathbb{C}^{N}\otimes\Sigma
_{2m}\right)  $. If $U_{0}\left(  0\right)  \in\mathbb{C}^{N}\otimes
\Sigma_{2m} $ is arbitrarily specified, we seek a formula for $U_{k}\left(
y\right)  $, $k=0,\ldots,Q$, such that
\begin{equation}
\left(  \mathcal{D}^{2}+\partial_{t}\right)  \Psi_{Q}\left(  y,t\right)
=\mathcal{E}\left(  r,t\right)  t^{Q}\mathcal{D}^{2}\left(  U_{Q}\right)
\left(  y\right)  , \label{AsyCond}%
\end{equation}
where the square $\mathcal{D}^{2}$ of the Dirac operator $\mathcal{D}$ can be
written (where ``$\cdot$'' is Clifford multiplication) as
\[
\mathcal{D}^{2}\psi=-\Delta\psi+\tfrac{1}{2}\sum\nolimits_{j,k}\Omega
_{jk}^{\varepsilon}E_{j}\cdot E_{k}\cdot\psi+\tfrac{1}{4}S\psi,
\]
by virtue of the generalized Lichnerowicz formula (see \cite[p. 164]{LaMi}).
It is convenient to define the 0-th order operator $\mathcal{F}$ on
$C^{\infty}\left(  \mathcal{B},\mathbb{C}^{N}\otimes\Sigma_{2m}\right)  $ via
\begin{align*}
&  \mathcal{F}\left[  \psi\right]  :=\tfrac{1}{2}\sum\nolimits_{j,k}%
\Omega_{jk}^{\varepsilon}E_{j}\cdot E_{k}\cdot\psi,\text{ so that }\\
&  \mathcal{D}^{2}=-\Delta\psi+\left(  \mathcal{F}+\tfrac{1}{4}S\right)
\left[  \psi\right]  .
\end{align*}
The desired formula for the $U_{k}\left(  y\right)  $ involves the operator
$A$ on $C^{\infty}\left(  \mathcal{B},\mathbb{C}^{N}\otimes\Sigma_{2m}\right)
$ given by
\[
A\left[  \psi\right]  :=-h^{1/4}\mathcal{D}^{2}[h^{-1/4}\psi]=h^{1/4}%
\Delta\lbrack h^{-1/4}\psi]-\left(  \mathcal{F}+\tfrac{1}{4}S\right)  \left[
\psi\right]  ,
\]
where $h^{1/4}:=(\sqrt{\det h})^{1/2}$. For $s\in\left[  0,1\right]  $, let
\[
A_{s}\left[  \psi\right]  \left(  y\right)  :=A\left[  \psi\right]  \left(
sy\right)  .
\]
As is proved in \cite{Ble92} or in the forthcoming \cite{BlBo03}, we have

\begin{proposition}
\label{AsySolnProp}Let $U_{0}\left(  0\right)  \in\mathbb{C}^{N}\otimes
\Sigma_{2m}$, and let $V_{0}\in C^{\infty}\left(  \mathcal{B},\mathbb{C}%
^{N}\otimes\Sigma_{2m}\right)  $ be the constant function $V_{0}\left(
y\right)  \equiv U_{0}\left(  0\right)  $. Then the $U_{k}\left(  y\right)  $
which satisfy $($\ref{AsyCond}$)$ are given by
\begin{align}
&  \!\!\!\!\!\!\!\!U_{k}\left(  y\right)  =h\left(  y\right)  ^{-1/4}%
V_{k}\left(  y\right)  ,\text{ where}\nonumber\\
&  \!\!\!\!\!\!\!\!V_{k}\left(  y\right)  =\int_{I^{k}}\prod\nolimits_{i=0}%
^{k-1}\left(  s_{i}\right)  ^{i}\left(  A_{s_{k-1}}\circ\cdots\circ A_{s_{0}%
}\left[  V_{0}\right]  \right)  \left(  y\right)  ds_{0}\ldots ds_{k-1},
\label{UkForm}%
\end{align}
and where $I^{k}=\left\{  \left(  s_{0},\ldots s_{k}\right)  :s_{i}\in\left[
0,1\right]  ,i\in\left\{  0,\ldots,k-1\right\}  \text{ }\right\}  .$
\end{proposition}

\bigskip

Note that $U_{0}\left(  0\right)  \in\mathbb{C}^{N}\otimes\Sigma_{2m}$ may be
arbitrarily specified, and once $U_{0}\left(  0\right)  $ is chosen, the
$U_{m}\left(  y\right)  $ are uniquely determined via (\ref{UkForm}). Let
$h_{k}\left(  y\right)  \in\operatorname{End}\left(  \mathbb{C}^{N}%
\otimes\Sigma_{2m}\right)  $ be given by
\begin{equation}
h_{k}\left(  y\right)  \left(  U_{0}\left(  0\right)  \right)  :=U_{k}\left(
y\right)  \label{hkdefn}%
\end{equation}
(in particular, $h_{0}\left(  0\right)  =\operatorname{I}\in\operatorname{End}%
\left(  \mathbb{C}^{N}\otimes\Sigma_{2m}\right)  $), and
\[
H_{Q}\left(  0,y,t\right)  :=\mathcal{E}\left(  r,t\right)  \sum_{k=0}%
^{Q}h_{k}\left(  y\right)  t^{k}\in C^{\infty}\left(  \mathcal{B}%
,\operatorname{End}\left(  \mathbb{C}^{N}\otimes\Sigma_{2m}\right)  \right)
.
\]
We may regard $H_{Q}\left(  0,y,t\right)  $ as
\[
H_{Q}\left(  x,y,t\right)  \in\operatorname{Hom}\left(  E_{x}\otimes
\Sigma\left(  M\right)  _{x},E_{y}\otimes\Sigma\left(  M\right)  _{y}\right)
,
\]
where we recall that $x\in M$ is the point about which we have chosen normal
coordinates. For $y$ sufficiently close to $x$, we set
\begin{equation}
H_{Q}\left(  x,y,t\right)  :=\mathcal{E}\left(  d\left(  x,y\right)
,t\right)  \sum_{k=0}^{Q}h_{k}\left(  x,y\right)  t^{k}. \label{HQDefnGen}%
\end{equation}
Of course, one expects that $H_{Q}\left(  x,y,t\right)  $ provides the desired
asymptotic expansion (\ref{AsyExpxy}). Although this is very plausible, it is
not at all easy to prove honestly. When proofs are attempted in the
literature, often steps are skipped, hands are waved, and errors are made. An
extremely careful (and hence nearly unbearable) proof will be provided in
\cite{BlBo03}, but we must forgo this here. Thus, we will only state here
without proof that
\[
k\left(  x,y,t\right)  \sim H_{Q}\left(  x,y,t\right)  :=\mathcal{E}\left(
d\left(  x,y\right)  ,t\right)  \sum\nolimits_{j=0}^{Q}h_{j}\left(
x,y\right)  t^{j}%
\]
for $d\left(  x,y\right)  $ sufficiently small, where the $h_{j}$ are given in
(\ref{hkdefn}).

Using normal coordinates $\left(  y^{1},\ldots,y^{2m}\right)  \in B\left(
r_{0},0\right)  $ about $x\in M$ and the radial gauge, and selecting $V_{0}%
\in\mathbb{C}^{N}\otimes\Sigma_{2m}$, by Proposition \ref{AsySolnProp}, we
have
\begin{equation}
h_{m}\left(  x,x\right)  \left(  V_{0}\right)  =\int_{I^{m}}\prod
\nolimits_{i=0}^{m-1}\left(  s_{i}\right)  ^{i}\left(  \left(  A_{s_{m-1}%
}\circ\cdots\circ A_{s_{0}}\right)  \left[  \widetilde{V}_{0}\right]  \right)
\left(  0\right)  ds_{0}\ldots ds_{m-1}, \label{hmxx}%
\end{equation}
where $\widetilde{V}_{0}\in C^{\infty}\left(  B\left(  r_{0},0\right)
,\mathbb{C}^{N}\otimes\Sigma_{2m}\right)  $ is the constant extension of
$V_{0}$. Recall that for $\psi\in C^{\infty}\left(  B\left(  r_{0},0\right)
,\mathbb{C}^{N}\otimes\Sigma_{2m}\right)  $, we have
\[
A_{s}\left[  \psi\right]  \left(  y\right)  :=A\left[  \psi\right]  \left(
sy\right)  ,\text{ where }A\left[  \psi\right]  :=h^{1/4}\Delta\lbrack
h^{-1/4}\psi]-\left(  \mathcal{F}+\tfrac{1}{4}S\right)  \left[  \psi\right]
.
\]
While the right side of (\ref{hmxx}) may seem unwieldy, there is substantial
simplification due to facts that $\left(  A_{s_{m-1}}\circ\cdots\circ
A_{s_{0}}\right)  [\widetilde{V}_{0}]\left(  y\right)  $ is evaluated at $y=0$
in (\ref{hmxx}). Also, if $\gamma^{1},\ldots,\gamma^{n}$ denote the so-called
gamma matrices for Clifford multiplication by $\partial_{1},\ldots
,\partial_{n}$, only those terms of $A_{s_{m-1}}\circ\cdots\circ A_{s_{0}%
}[\widetilde{V}_{0}]\left(  0\right)  $ which involve the product
$\gamma_{n+1}:=\gamma^{1}\cdots\gamma^{n}$ will survive when the supertrace
$\operatorname{Str}\left(  h_{m}\left(  x,x\right)  \right)  $ is taken. As a
consequence, we have the following simplification (essentially contained in
\cite{Ble92}, or better yet, to appear in \cite{BlBo03})

\begin{proposition}
\label{Asuper0Prop}Let $R_{klji}\left(  0\right)  =h\left(  \Omega^{\theta
}\left(  \partial_{i},\partial_{j}\right)  \partial_{l},\partial_{k}\right)  $
denote the components of the Riemann curvature tensor of $h$ at $x,$ and let
$\Omega_{ij}^{\varepsilon}\left(  0\right)  :=\Omega^{\varepsilon}\left(
\partial_{i},\partial_{j}\right)  $ at $x$. Set
\begin{align}
&  \widetilde{\theta}^{1}\left(  \partial_{j}\right)  :=\tfrac{1}{8}%
\sum\nolimits_{k,l,i}R_{klji}\left(  0\right)  \gamma^{k}\gamma^{l}%
y^{i},\text{ }\nonumber\\
&  \mathcal{F}^{0}:=\tfrac{1}{2}\sum\nolimits_{i,j}F_{ij}\otimes\gamma
^{i}\gamma^{j}=\tfrac{1}{2}\sum\nolimits_{i,j}\Omega_{ij}^{\varepsilon}\left(
0\right)  \otimes\gamma^{i}\gamma^{j},\text{ and}\nonumber\\
&  A^{0}:=\sum\nolimits_{i}\left(  \partial_{i}^{2}+\widetilde{\theta}%
^{1}\left(  \partial_{i}\right)  ^{2}\right)  -\mathcal{F}_{0}.
\label{Asuper0}%
\end{align}
For $V_{0}\in\mathbb{C}^{N}\otimes\Sigma_{2m}$, define
\begin{equation}
h_{m}^{0}\left(  0,0\right)  \left(  V_{0}\right)  :=\int_{I^{k}}%
\prod\nolimits_{i=0}^{m-1}\left(  s_{i}\right)  ^{i}\left(  \left(
A_{s_{m-1}}^{0}\circ\cdots\circ A_{s_{0}}^{0}\right)  \left[  \widetilde
{V}_{0}\right]  \right)  \left(  0\right)  ds_{0}\ldots ds_{m-1}, \label{h0m}%
\end{equation}
where $\widetilde{V}_{0}\in C^{\infty}\left(  B\left(  r_{0},0\right)
,\mathbb{C}^{N}\otimes\Sigma_{2m}\right)  $ is the constant extension of
$V_{0}\in\mathbb{C}^{N}\otimes\Sigma_{2m}$. Then
\begin{equation}
\operatorname{Str}\left(  h_{m}\left(  0,0\right)  \right)
=\operatorname{Str}\left(  h_{m}^{0}\left(  0,0\right)  \right)  .
\label{StrStr0}%
\end{equation}
In other words, in the computation of $\operatorname{Str}\left(  h_{m}\left(
0,0\right)  \right)  $ given by (\ref{hmxx}), we may replace $A$ by $A^{0}$.
\end{proposition}

\bigskip

This is a substantial simplification, not only in that $A^{0}$ is a
second-order differential operator with coefficients which are at most
quadratic in $y$, but it also shows that $\operatorname{Str}\left(
h_{m}\left(  x,x\right)  \right)  $ only depends on the curvatures
$\Omega^{\theta}$ and $\Omega^{\varepsilon}$ at the point $x$. One might
regard the gist of the Index Formula for twisted Dirac operator as exhibiting
the global quantity $\operatorname{index}\left(  \mathcal{D}^{+}\right)  $ as
the integral of a form which may be locally computed. From this perspective,
Proposition \ref{Asuper0Prop} does the job. Also, knowing in advance that
$\operatorname{index}\left(  \mathcal{D}^{+}\right)  $ is insensitive to
perturbations in $h$ and $\varepsilon$, one suspects that $\operatorname{Str}%
\left(  h_{m}\left(  x,x\right)  \right)  \nu_{x}$ can be expressed in terms
of the standard forms which represent characteristic classes for $TM$ and $E$.
The Local Index Formula confirms this. Moreover, for low values of $m,$ say
$m=1$ or $2$ (i.e., for 2 and 4-manifolds), one can directly compute
$\operatorname{Str}\left(  h_{m}^{0}\left(  0,0\right)  \right)  $ using
(\ref{h0m}), and thereby verify Theorem \ref{LocIndThm} and hence obtain
Corollary \ref{IndFormTwDiracCor} rather easily. For readers who have no use
for the Local Index Theorem beyond dimension 4, this is sufficient. It
requires more effort to prove Theorem \ref{LocIndThm} for general $m$. For
lack of space, we cannot go into the details of this here, but they can be
found in \cite{Ble92} or \cite{BlBo03}. It is well worth mentioning that the
appearance of the $\sinh$ function in the Local Index Formula has its roots in
Mehler's formula for the heat kernel
\begin{equation}
e_{a}\left(  x,y,t\right)  =\frac{1}{\sqrt{4\pi\frac{\sinh\left(  2at\right)
}{2a}}}\exp\left(  -\frac{1}{4\frac{\sinh\left(  2at\right)  }{2a}}\left(
\cosh\left(  2at\right)  \left(  x^{2}+y^{2}\right)  -2xy\right)  \right)  .
\label{MehlerForm}%
\end{equation}
of the generalized 1-dimensional heat problem
\begin{align*}
&  u_{t}=u_{xx}-a^{2}x^{2}u,\text{\ \ \ }u\left(  x,t\right)  \in
\mathbb{R},\;\left(  y,t\right)  \in\mathbb{R\times}\left(  0,\infty\right)
,\\
&  u\left(  x,0\right)  =f\left(  x\right)  ,
\end{align*}
where $0\neq a\in\mathbb{R}$ is a given constant. A solution of this problem
is given by
\[
u\left(  x,t\right)  =\int_{-\infty}^{\infty}e_{a}\left(  x,y,t\right)
f\left(  y\right)  \,dy,
\]
and this reduces to the usual formula as $a\rightarrow0$. The nice idea of
using Mehler's formula in a rigorous derivation of the Local Index Theorem
appears to be due to Getzler in \cite{Get83} and \cite{Get86}, although it was
at least implicitly involved in earlier heuristic supersymmetric path integral
arguments for the Index Theorem. In the same vein, further simplifications and
details can be found in \cite{BeGeVe92}, \cite{Yu01}, and \cite{Ble92}, but
the treatment to be found in \cite{BlBo03} will be more self-contained and
less demanding.

While the Local Index Theorem (Theorem \ref{LocIndThm}) is stated for twisted
Dirac operators, the same proof may be applied to obtain the index formulas
for an elliptic operator, possibly on a \textit{nonspin} manifold, which is
only locally of the form of twisted Dirac operator $\mathcal{D}^{+}$. Indeed,
if $\mathcal{A}$ is such an operator and $k$ is the heat kernel for
$\mathcal{A}^{\ast}\mathcal{A}\oplus\mathcal{AA}^{\ast}$, then from the
spectral resolution of $\mathcal{A}$, we can still deduce from the asymptotic
expansion of $k$ that
\[
\operatorname{index}\left(  \mathcal{A}\right)  =\left(  4\pi\right)
^{-m}\int_{M}\operatorname{Str}\left(  h_{m}\left(  x,x\right)  \right)
\,\nu_{x},
\]
where the supertrace $\operatorname{Str}$ is defined in the natural way. The
crucial observation is that since $\mathcal{A}$ is locally in the form of a
twisted Dirac operator, we can compute $\operatorname{Str}\left(  h_{m}\left(
x,x\right)  \right)  $ in exactly the same way (i.e., locally) as we have
done. Since it is not easy to find first-order elliptic operators of
geometrical significance which are not expressible in terms of locally twisted
Dirac operators (or 0-th order perturbations thereof), the Local Index Formula
for twisted Dirac operators is much more comprehensive than it would appear at
first glance.

\bigskip

\subsection{Dirac Type Operators on Manifolds with Boundary \label{ss:Index of
Elliptic Boundary Value Problems}}

\bigskip

We fix the notation and recall basic properties of operators of Dirac type on
manifolds with boundary.

\subsubsection{The General Setting for Partitioned Manifolds}

\begin{figure}
[ptb]
\begin{center}
\includegraphics[
height=1.4736in,
width=2.9776in
]%
{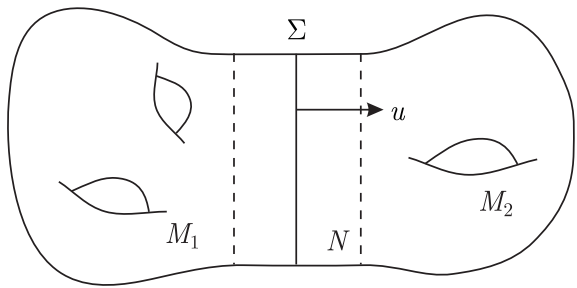}%
\caption{The partition of $M=M_{1}\cup_{\Sigma}M_{2}$}%
\label{f:split}%
\end{center}
\end{figure}

Let $M$ be a compact smooth Riemannian \textit{partitioned} manifold
\[
M=M_{1}\cup_{\Sigma}M_{2}:=M_{1}\cup M_{2}\ ,\quad\text{where $M_{1}\cap
M_{2}=\partial M_{1}=\partial M_{2}=\Sigma$}%
\]
and $\Sigma$ a hypersurface (see Figure \ref{f:split}).%
We assume that $M\setminus\Sigma$ does not have a closed connected component
(i.e. $\Sigma$ intersects any connected component of $M_{1}$ and $M_{2}$).
Let
\[
\mathcal{D}:C^{\infty}(M;S)\rightarrow C^{\infty}(M;S)
\]
be an \textit{operator of Dirac type} acting on sections of a Hermitian bundle
$S$ of Clifford modules over $M$, i.e. $\mathcal{D}=\mathbf{c}\circ\nabla$
where $\mathbf{c}$ denotes the Clifford multiplication and $\nabla$ is a
connection for $S$. Unlike the more general case introduced in Subsection
\ref{ss: Ops of Dirac Type}, we assume here that $\nabla$ is compatible with
$\mathbf{c}$ in the sense that $\nabla\mathbf{c}=0$. From this compatibility
assumption it follows that $\mathcal{D}$\ is symmetric and essentially
self-adjoint over $M$.

For even $n=\dim M$, the splitting $\mathfrak{Cl}(M)=\mathfrak{Cl}%
^{+}(M)\oplus\mathfrak{Cl}^{-}(M)$ of the Clifford bundles induces a
corresponding splitting of $S=S^{+}\oplus S^{-}$ and a \emph{chiral
decomposition}
\[
\mathcal{D}=\left(
\begin{array}
[c]{cc}%
0 & \mathcal{D}^{-}=(\mathcal{D}^{+})^{\ast}\\
\mathcal{D}^{+} & 0
\end{array}
\right)
\]
of the \emph{total} Dirac operator. The \emph{chiral Dirac operators}
$\mathcal{D}^{\pm}$ are elliptic but not symmetric, and for that reason they
may have nontrivial indices which provide us with important topological and
geometric invariants, as was found in the special case of twisted Dirac
operators in Subsection \ref{ss:IndexTDOonClosedMfds}.

Here we assume that all metric structures of $M$ and $S$ are product in a
collar neighborhood ${N}=[-1,1]\times\Sigma$ of $\Sigma$. If $u$ denotes the
normal coordinate (running from $M_{1}$ to $M_{2}$), then
\[
\mathcal{D}|_{{N}}=\sigma(\partial_{u}+\mathcal{B})\,,\text{ where }%
\sigma:=\mathbf{c}\left(  du\right)  \text{, }\partial_{u}:=\dfrac{\partial
}{\partial u},
\]
and $\mathcal{B}$ denotes the canonically associated Dirac operator over
$\Sigma$, called the \textit{tangential operator}. We have a similar product
formula for the chiral Dirac operator. Here the point of the product structure
is that then $\sigma$ and $\mathcal{B}$ do not depend on the normal variable.
Note that $\sigma$ (Clifford multiplication by $du$) is a unitary mapping
$L^{2}(\Sigma;S|_{\Sigma})\rightarrow L^{2}(\Sigma;S|_{\Sigma})$ with
$\sigma^{2}=-\operatorname{I}$ and $\sigma\mathcal{B}=-\mathcal{B}\sigma$. In
the non-product case, there are certain ambiguities in defining a `tangential
operator' which we shall not discuss here (but see also Formula
\ref{e-product}).

\medskip

\subsubsection{Analysis tools: Green's Formula}

For notational economy, we set $X:=M_{2}$. For greater generality, we consider
the chiral Dirac operator
\[
\mathcal{D}^{+}:C^{\infty}(X;S^{+})\rightarrow C^{\infty}(X;S^{-})
\]
and write $\mathcal{D}^{-}$ for its formally adjoint operator. The
corresponding results follow at once for the total Dirac operator.

\medskip

\begin{lemma}
\label{l:green} Let $\left\langle .,.\right\rangle _{\pm}$ denote the scalar
product in $L^{2}(X;S^{\pm})$. Then we have
\[
\left\langle \mathcal{D}^{+}f_{+},f_{-}\right\rangle _{-}-\left\langle
f_{+},\mathcal{D}^{-}f_{-}\right\rangle _{+}=-\int_{\Sigma}(\sigma
\gamma_{\infty}f_{+},\gamma_{\infty}f_{-})\,d\operatorname*{vol}%
\nolimits_{\Sigma}%
\]
for any $f_{\pm}\in C^{\infty}(X;S^{\pm})$.
\end{lemma}

Here
\begin{equation}
\gamma_{\infty}:C^{\infty}(X;S^{\pm})\rightarrow C^{\infty}(\Sigma;S^{\pm
}|_{\Sigma}) \label{e:gamma}%
\end{equation}
denotes the restriction of a section to the boundary $\Sigma$.

\bigskip

\subsubsection{Cauchy Data Spaces and the Calder{\'{o}}n Projection}

\label{s:calderon}

To explain the $L^{2}$ Cauchy data spaces we recall three additional, somewhat
delicate and not widely known properties of operators of Dirac type on compact
manifolds with boundary from \cite{BoWo93}:

\begin{enumerate}
\item  the invertible extension to the double;

\item  the Poisson type operator and the Calder{\'o}n projection; and

\item  the twisted orthogonality of the Cauchy data spaces for chiral and
total Dirac operators which gives the Lagrangian property in the symmetric
case (i.e., for the total Dirac operator).
\end{enumerate}

The idea and the properties of the Calder{\'{o}}n projection were announced in
Calder\'{o}n \cite{Ca63} and proved in Seeley \cite{Se66} in great generality.
In the following, we restrict ourselves to constructing the Calder\'{o}n
projection for operators of Dirac type (or, more generally, elliptic
differential operators of first order) which simplifies the presentation substantially.\bigskip

\subsubsection{Invertible Extension\label{sss:InvertibleExtension}}

First we construct the \textit{invertible double.} Clifford multiplication by
the inward normal vector gives a natural clutching of $S^{+}$ over one copy of
$X$ with $S^{-}$ over a second copy of $X$ to a smooth bundle $\widetilde
{S^{+}}$ over the closed double $\widetilde{X}$. The product forms of
$\mathcal{D}^{+}$ and $\mathcal{D}^{-}=\left(  \mathcal{D}^{+}\right)  ^{\ast
}$ fit together over the boundary and provide a new operator of Dirac type,
namely
\begin{equation}
\widetilde{\mathcal{D}^{+}}:=\mathcal{D}^{+}\cup\mathcal{D}^{-}:C^{\infty
}\left(  \widetilde{X},\widetilde{S^{+}}\right)  \rightarrow C^{\infty}\left(
\widetilde{X},\widetilde{S^{-}}\right)  . \label{DirDblOp}%
\end{equation}
Clearly $\left(  \mathcal{D}^{+}\cup\mathcal{D}^{-}\right)  ^{\ast
}=\mathcal{D}^{-}\cup\mathcal{D}^{+},$ and so $\operatorname{index}%
\widetilde{\mathcal{D}^{+}}=0$. It turns out that $\widetilde{\mathcal{D}^{+}%
}$ is invertible with a pseudo-differential elliptic inverse $(\widetilde
{\mathcal{D}^{+}})^{-1}$. Clearly, \textit{local} solutions of the homogenous
equation (here, solutions on \textit{one} copy of $X$) do not extend to
\textit{global} solutions on $\widetilde{X}$ in general. As a matter of fact,
the operator $\mathcal{D}^{+}$, even being the restriction of an invertible,
locally defined differential operator, is not invertible in general, and we
have $r^{+}(\widetilde{\mathcal{D}^{+}})^{-1}e^{+}\mathcal{D}^{+}%
\neq\operatorname{I}$, where $e^{+}:L^{2}\left(  X;S^{+}\right)  \rightarrow
L^{2}(\widetilde{X};\widetilde{S^{+}})$ denotes the extension-by-zero operator
and $r^{+}:H^{s}\left(  \widetilde{X};\widetilde{S^{+}}\right)  \rightarrow
H^{s}\left(  X;S^{+}\right)  $ the natural restriction operator for Sobolev
spaces for $s$ real. The precise decomposition $L^{2}(\widetilde{X}%
;\widetilde{S^{+}})=L^{2}\left(  X;S^{+}\right)  \times L^{2}\left(
X;S^{-}\right)  $ gives a different picture in case that a given operator $T$
on one component can be extended to an invertible operator $\widetilde{T}$ on
the whole space. This, indeed, would imply $T$ invertible with $T^{-1}%
=r^{+}\left(  \widetilde{T}\right)  ^{-1}e^{+}$. The $L^{2}$ extension of
$\widetilde{\mathcal{D}^{+}}$, however, is not a precise extension of the
$L^{2}$ extension of $\mathcal{D}^{+}$. Therefore, in our case, the $L^{2}%
$-argument breaks down.

\begin{example}
\label{CREx}In the simplest possible two-dimensional case we consider the
Cauchy-Riemann operator $\overline{\partial}:C^{\infty}(D^{2})\rightarrow
C^{\infty}(D^{2})$ over the disc $D^{2}$, where $\overline{\partial}=\frac
{1}{2}(\partial_{x}+i\partial_{y})$. In polar coordinates, this operator has
the form $\frac{1}{2}e^{i\varphi}(\partial_{r}+\frac{i}{r}\partial_{\varphi}%
)$. Therefore, after some small smooth perturbations (and modulo the factor
$\frac{1}{2}$), we assume that $\overline{\partial}$ has the following form in
a certain collar neighborhood of the boundary:
\[
\overline{\partial}=e^{i\varphi}(\partial_{r}+i\partial_{\varphi})
\]
Now we construct the invertible double of $\overline{\partial}$. By $E^{k}$,
$k\in\mathbb{Z}$, we denote the bundle, which is obtained from two copies of
$D^{2}\times\mathbb{C}$ by the identification $(z,w)=(z,z^{k}w)$ near the
equator. We obtain the bundle $E^{1}$ by gluing two halves of $D^{2}%
\times\mathbb{C}$ by $\sigma(\varphi)=e^{i\varphi}$ and $E^{-1}$ by gluing
with the adjoint symbol. In such a way we obtain the operator
\[
\widetilde{\overline{\partial}}:=\overline{\partial}\cup(\overline{\partial
})^{\ast}:C^{\infty}(S^{2};E^{1})\rightarrow C^{\infty}(S^{2};E^{-1})
\]
over the whole 2-sphere.\newline \qquad Let us analyze the situation more
carefully. We fix $N:=\left(  -\varepsilon,+\varepsilon\right)  \times S^{1}$,
a bicollar neighborhood of the equator. The operator formally adjoint to
$\overline{\partial}$ has the form
\[
(\overline{\partial})^{\ast}=e^{-i\varphi}\left(  -\partial_{u}+i\partial
_{\varphi}+1\right)
\]
($u=r-1$) in this cylinder. A section of $E^{1}$ is a couple $\left(
s_{1},s_{2}\right)  $ such that in $N$
\[
s_{2}\left(  u,\varphi\right)  =e^{i\varphi}s_{1}\left(  u,\varphi\right)  .
\]
The couple $\left(  \overline{\partial}s_{1},(\overline{\partial})^{\ast}%
s_{2}\right)  $ is a smooth section of $E^{-1}$. To show this, we check that
$(\overline{\partial})^{\ast}s_{2}=e^{-i\varphi}\overline{\partial}s_{1}$. In
the neighborhood $N$, we have
\begin{align*}
(\overline{\partial})^{\ast}s_{2}  &  =(\overline{\partial})^{\ast}\left(
e^{i\varphi}s_{1}\right)  =e^{-i\varphi}\left(  -\partial_{u}+i\partial
_{\varphi}+1\right)  \left(  e^{i\varphi}s_{1}\right) \\
&  =\partial_{u}s_{1}+ie^{-i\varphi}\partial_{\varphi}\left(  e^{i\varphi
}s_{1}\right)  +s_{1}=\left(  \partial_{u}+i\partial_{\varphi}\right)  s_{1}\\
&  =e^{-i\varphi}e^{i\varphi}\left(  \partial_{u}+i\partial_{\varphi}\right)
s_{1}=e^{-i\varphi}\left(  \overline{\partial}s_{1}\right)  .
\end{align*}
Then the operator $\overline{\partial}\cup\overline{\partial}^{\ast}$ becomes
injective and $\operatorname{index}\overline{\partial}\cup\overline{\partial
}^{\ast}=0$.
\end{example}

\subsubsection{The Poisson Operator and the Calder{\'{o}}n
Projection\label{sss:PoissonOperator}}

Next we investigate the solution spaces and their traces at the boundary. For
a total or chiral operator of Dirac type over a smooth compact manifold with
boundary $\Sigma$ and for any real $s$ we define the \emph{null space}
\[
\operatorname{Ker}(\mathcal{D}^{+},s):=\{f\in H^{s}(X;S^{+})\mid
\mathcal{D}^{+}f=0\text{ in $X\setminus\Sigma$}\}.
\]
The null spaces consist of sections which are distributional for negative $s$;
by elliptic regularity they are smooth in the interior; in particular they
possess a smooth restriction on the hypersurface $\Sigma_{\varepsilon
}=\{\varepsilon\}\times\Sigma$ parallel to the boundary $\Sigma$ of $X$ at a
distance $\varepsilon>0$. By a Riesz operator argument they can be shown to
also possess a trace over the boundary. Of course, that trace is no longer
smooth but belongs to $H^{s-\frac{1}{2}}(\Sigma;S^{+}|_{\text{$\Sigma$}})$.
More precisely, we have the following well-known \emph{General Restriction
Theorem} (for a proof see e.g. \cite{BoWo93}, Chapters 11 and 13):

\begin{theorem}
\label{t:res} \textrm{(a)} Let $s>\frac{1}{2}$. Then the restriction map
$\gamma_{\infty}$ of \eqref{e:gamma} extends to a bounded map
\begin{equation}
\gamma_{s}:H^{s}(X;S^{+})\rightarrow H^{s-\frac{1}{2}}(\text{$\Sigma$}%
;S^{+}|_{\text{$\Sigma$}}).
\end{equation}

\noindent\textrm{(b)} For $s\leq\frac{1}{2}$\thinspace, the preceding
reduction is no longer defined for arbitrary sections but only for solutions
of the operator $\mathcal{D}^{+}$: let $f\in\operatorname{Ker}(\mathcal{D}%
^{+},s)$ and let $\gamma_{(\varepsilon)}f$ denote the well-defined trace of
$f$ in $C^{\infty}(\Sigma_{\varepsilon};S^{+}|_{\text{$\Sigma$}})$. Then, as
$\varepsilon\rightarrow0_{+}$, the sections $\gamma_{(\varepsilon)}f$ converge
to an element $\gamma_{s}f\in H^{s-\frac{1}{2}}(\Sigma;S^{+}|_{\text{$\Sigma$%
}})$.

\noindent\textrm{(c)} Let $\widetilde{\mathcal{D}^{+}}$ denote the invertible
double of $\mathcal{D}^{+}$, $r_{+}$ denote the restriction operator
$r_{+}:H^{s}(\widetilde{X};\widetilde{S}^{+})\rightarrow H^{s}(X;S^{+})$ and
let $\gamma_{\infty}^{\ast}$ be the dual of $\gamma_{\infty}$ in the
distributional sense. For any $s\in\mathbb{R}$ the mapping \emph{(Poisson type
operator)}
\[
\mathcal{K}:=r_{+}\left(  \widetilde{\mathcal{D}^{+}}\right)  ^{-1}%
\gamma_{\infty}^{\ast}\sigma:C^{\infty}(\text{$\Sigma$};S^{+}|_{\Sigma
})\rightarrow C^{\infty}(X;S^{+})
\]
extends to a continuous map $\mathcal{K}^{(s)}:H^{s-1/2}(\Sigma;S^{+}%
|_{\text{$\Sigma$}})\rightarrow H^{s}(X;S^{+})$ with
\[
\operatorname{range}\mathcal{K}^{(s)}=\operatorname{Ker}(\mathcal{D}^{+},s).
\]
\end{theorem}

For $s=0$, Theorem \ref{t:res} can be reformulated in the following way:

\bigskip

\begin{corollary}
\label{c:restriction} For a constant $C$ independent of $f$, we have
\[
\left\|  \gamma\left(  f\right)  \right\|  _{-\frac{1}{2}}\leq C\left(
\left\|  \mathcal{D}^{+}f\right\|  _{0}+\left\|  f\right\|  _{0}\right)
\text{ for all }f\in D_{\max}\left(  \mathcal{D}^{+}\right)  ,\text{ }%
\]
where $D_{\max}\left(  \mathcal{D}^{+}\right)  :=\left\{  f\in L^{2}\left(
X;S\right)  \mid\mathcal{D}^{+}f\in L^{2}\left(  X;S\right)  \right\}  $. So,
the restriction
\[
\gamma:D_{\max}\left(  \mathcal{D}^{+}\right)  \longrightarrow H^{-1/2}%
(\Sigma;S^{+}|_{\text{$\Sigma$}})
\]
is well defined and bounded.
\end{corollary}

Proofs of Theorem and Corollary can be found, e.g., in Booss--Bavnbek and
Wojciechowski \cite{BoWo93}, Theorems 13.1 and 13.8 for our situation
($\mathcal{D}^{+}$ is of order 1); and in H\"{o}rmander \cite{Ho66} in greater
generality (Theorem 2.2.1 and the Estimate (2.2.8), p. 194).

The composition
\begin{equation}
\mathcal{P}(\mathcal{D}^{+}):=\gamma_{\infty}\circ\mathcal{K}:C^{\infty
}(\text{$\Sigma$};S^{+}|_{\text{$\Sigma$}})\rightarrow C^{\infty
}(\text{$\Sigma$};S^{+}|_{\text{$\Sigma$}}) \label{e:calderon}%
\end{equation}
is called the (\textit{Szeg\"{o}--})\textit{Calder\'{o}n projection}. It is a
pseudo-differential projection (idempotent, but in general not orthogonal). We
denote by $\mathcal{P}(\mathcal{D}^{+})^{(s)}$ its extension to the $s$-th
Sobolev space over $\Sigma$. It has the following geometric meaning.

We now have three options of defining the corresponding \emph{Cauchy data} (or
\emph{Hardy}) \emph{spaces}:

\begin{definition}
\label{d:cauchy-data} For all real $s$ we define
\begin{align*}
\Lambda(\mathcal{D}^{+},s)  &  :=\gamma_{s}(\operatorname{Ker}(\mathcal{D}%
^{+},s)),\\
\Lambda^{\mathrm{clos}}(\mathcal{D}^{+},s)  &  :=\overline{\gamma_{\infty
}\{f\in C^{\infty}(X;S^{+})\mid\mathcal{D}^{+}f=0\text{ in $X\setminus\Sigma$%
}\}}^{H^{s-\frac{1}{2}}(\text{$\Sigma$};S^{+}|_{\text{$\Sigma$}})},\text{
and}\\
\Lambda^{\mathrm{Cald}}(\mathcal{D}^{+},s)  &  :=\operatorname{range}%
\mathcal{P}(\mathcal{D}^{+})^{(s-\frac{1}{2})}\,.
\end{align*}
\end{definition}

The range of a projection is closed; the inclusions of the Sobolev spaces are
dense; and $\operatorname{range}\mathcal{P}(\mathcal{D}^{+})=\gamma_{\infty
}\{f\in C^{\infty}(X;S^{+})\mid\mathcal{D}^{+}f=0\text{ in $X\setminus$}%
\Sigma\}$, as shown in \cite{BoWo93}. So, the second and the third definition
of the Cauchy data space coincide. Moreover, for $s>\frac{1}{2}$ one has
$\Lambda(\mathcal{D}^{+},s)=\Lambda^{\mathrm{Cald}}(\mathcal{D}^{+},s)$. This
equality can be extended to the $L^{2}$ case ($s=\frac{1}{2}$, see also
Theorem \ref{t:criss-cross} below), and remains valid for any real $s$, as
proved in Seeley, \cite{Se66}, Theorem 6. For $s\leq\frac{1}{2} $, the result
is somewhat counter-intuitive (see also Example \ref{ex:cylinder}b in the
following Subsection).

We have:

\begin{proposition}
\label{p:all-the-same} For all $s\in\mathbb{R}$
\[
\Lambda(\mathcal{D}^{+},s)=\Lambda^{\mathrm{clos}}(\mathcal{D}^{+}%
,s)=\Lambda^{\mathrm{Cald}}(\mathcal{D}^{+},s).
\]
\end{proposition}

\subsubsection{Calder{\'{o}}n and Atiyah--Patodi--Singer Projection
\label{sss:Calderon APS proj}}

The Calder{\'{o}}n projection is closely related to another projection
determined by the `tangential' part of $\mathcal{D}^{+}$, described as
follows. Let $\mathcal{B}$ denote the \textit{tangential} symmetric elliptic
differential operator over $\Sigma$ in the product form
\[
\mathcal{D}\text{$^{+}\ =\sigma(\partial_{u}+\mathcal{B}):$}C^{\infty}\left(
N,S^{+}|N\right)  \rightarrow C^{\infty}\left(  N,S^{-}|N\right)  \text{ }%
\]
in a collar neighborhood $N$ of $\Sigma$ in $X$. It has discrete real
eigenvalues and a complete system of $L^{2}$ orthonormal eigensections. Let
$P_{\geq}(\mathcal{B})$ denote the spectral (Atiyah--Patodi--Singer)
projection onto the subspace $L_{+}(\mathcal{B})$ of $L^{2}(\Sigma
;S^{+}|_{\Sigma})$ spanned by the eigensections corresponding to the
nonnegative eigenvalues of $\mathcal{B}$. It is a pseudo-differential operator
and its principal symbol $p_{+}$ is the projection onto the eigenspaces of the
principal symbol $b(y,\zeta)$ of $\mathcal{B}$ corresponding to nonnegative
eigenvalues. It turns out that $p_{+}$ coincides with the principal symbol of
the Calder\'{o}n projection.

We call the space of pseudo-differential projections with the same principal
symbol $p_{+}$ the \emph{Grassmannian} $\mathcal{G}\mathrm{r}_{p_{+}}$ and
equip it with the operator norm corresponding to $L^{2}(\Sigma;S^{+}|_{\Sigma
})$. It has countably many connected components; two projections $P_{1}$,
$P_{2}$ belong to the same component, if and only if the \emph{virtual
codimension}
\begin{equation}
\mathbf{i}(P_{2},P_{1}):=\operatorname{index}\left\{  P_{2}P_{1}%
:\operatorname{range}P_{1}\rightarrow\operatorname{range}P_{2}\right\}
\label{d:iP1P2}%
\end{equation}
of $P_{2}$ in $P_{1}$ vanishes; the higher homotopy groups of each connected
component are given by Bott periodicity.

\begin{example}
\label{ex:cylinder} (\textit{a}) For the Cauchy-Riemann operator on the disc
$D^{2}=\{|z|\leq1\}$, the Cauchy data space is spanned by the eigenfunctions
$e^{ik\theta}$ of the tangential operator $\partial_{\theta}$ over
$S^{1}=[0,2\pi]/\{0,2\pi\}$ for nonnegative $k$. So, the Calder{\'{o}}n
projection and the Atiyah--Patodi--Singer projection coincide in this
case.\medskip\newline \medskip\noindent(\textit{b}) Next we consider the
cylinder $X^{R}=[0,R]\times\Sigma_{0}$ with {$\mathcal{D}$}${_{R}}%
=\sigma(\partial_{u}+B)$.%
\begin{figure}
[h]
\begin{center}
\includegraphics[
height=0.6555in,
width=1.6864in
]%
{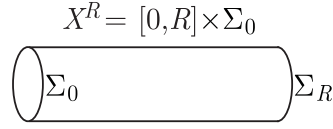}%
\caption{Cylinder of length $R$}%
\label{f:cylinder}%
\end{center}
\end{figure}
Here $B$ denotes a symmetric elliptic differential operator of first order
acting on sections of a bundle $E$ over $\Sigma_{0}$\thinspace, and $\sigma$ a
unitary bundle endomorphism with $\sigma^{2}=-\operatorname{I}$ and $\sigma
B=-B\sigma$. Let $B$ be invertible (for the ease of presentation). Let
$\{\varphi_{k},\lambda_{k}\}$ denote $B$'s spectral resolution of
$L^{2}(\Sigma_{0};E)$ with
\[
\dots\lambda_{-k}\leq...\lambda_{-1}<0<\lambda_{1}\leq...\lambda_{k}\leq\dots
\]
Then
\begin{equation}%
\begin{array}
[c]{ll}%
B\varphi_{k}=\lambda_{k}\varphi_{k} & \text{for all $k\in\mathbb{Z}%
\setminus\{0\}$}\,,\\
\lambda_{-k}=-\lambda_{k}\,,\ \sigma(\varphi_{k})=\varphi_{-k}\,,\text{
and}\ \sigma(\varphi_{-k})=-\varphi_{k} & \text{for $k>0$}.
\end{array}
\label{e:spectrum}%
\end{equation}
We consider
\begin{align*}
f\in\operatorname{Ker}({\mathcal{D}_{R}},0)  &  =\operatorname{span}%
\{e^{-\lambda_{k}u}\varphi_{k}\}_{k\in\mathbb{Z}\setminus\{0\}}\text{ in
$L^{2}(X^{R})$}\\
&  =\operatorname{Ker}{\mathcal{D}_{R}}_{\max}\text{ (kernel of maximal
extension)}.
\end{align*}
It can be written in the form
\begin{equation}
f(u,y)=f_{>}(u,y)+f_{<}(u,y),\text{ }u\in\left[  0,R\right]  , \label{e:(0)}%
\end{equation}
where
\[
f_{<}(u,y)=\sum_{k<0}a_{k}e^{-\lambda_{k}u}\varphi_{k}(y)\text{\ and\ }%
f_{>}(u,y)=\sum_{k>0}a_{k}e^{-\lambda_{k}u}\varphi_{k}(y).
\]
Because of
\[
\left\langle f,f\right\rangle _{L^{2}(X^{R})}<+\infty\iff\left\langle
f_{<},f_{<}\right\rangle <+\infty\;\text{and\ }\left\langle f_{>}%
,f_{>}\right\rangle <+\infty\,,
\]
the coefficients $a_{k}$ satisfy the conditions
\begin{equation}
\sum_{k<0}|a_{k}|^{2}\frac{e^{-2\lambda_{k}R}-1}{2|\lambda_{k}|}<+\infty
\quad\text{or, equivalently,}\quad\sum_{k<0}|a_{k}|^{2}\frac{e^{2|\lambda
_{k}|R}}{|\lambda_{k}|}<+\infty\label{e:(1)}%
\end{equation}
and
\begin{equation}
\sum_{k>0}|a_{k}|^{2}\frac{1-e^{-2\lambda_{k}R}}{2\lambda_{k}}<+\infty
\quad\text{or, equivalently,}\quad\sum_{k>0}|a_{k}|^{2}/{\lambda_{k}}%
<+\infty\,. \label{e:(2)}%
\end{equation}
We consider the Cauchy data space $\Lambda(${$\mathcal{D}$}${_{R}},0)$
consisting of all $\mathbb{\gamma}(f)$ with $f\in\operatorname{Ker}%
(${$\mathcal{D}$}${_{R}},0)$. Here $\mathbb{\gamma}(f)$ denotes the trace of
$f$ at the boundary
\[
\Sigma=\partial X^{R}=-{\Sigma}_{0}\sqcup{\Sigma}_{R}\,,
\]
where $\Sigma_{R}$ denotes a second copy of $\Sigma_{0}$\thinspace. According
to the spectral splitting $f=f_{>}+f_{<}$, we have
\[
\mathbb{\gamma}(f)=(s_{<}^{0},s_{<}^{R})+(s_{>}^{0},s_{>}^{R})
\]
where
\[
s_{>}^{0}=f_{>}(0),\ s_{<}^{0}=f_{<}(0),\ s_{>}^{R}=f_{>}(R),\ s_{<}^{R}%
=f_{<}(R).
\]
Because of \eqref{e:(1)} and \eqref{e:(2)}, we have
\[
(s_{<}^{0},s_{>}^{R})\in C^{\infty}({\Sigma}_{0}\cup{\Sigma}_{R}%
)\text{\ and\ }(s_{>}^{0},s_{<}^{R})\in H^{-1/2}({\Sigma}_{0}\cup{\Sigma}%
_{R}).
\]
Recall that
\[
\sum a_{k}\varphi_{k}\in H^{s}({\Sigma_{0}})\iff\sum|a_{k}|^{2}|k|^{2s/(m-1)}%
<+\infty
\]
and $|\lambda_{k}|\sim|k|^{\frac{1}{m-1}}$\thinspace\ for $k\rightarrow
\pm\infty$, where $m-1$ denotes the dimension of {$\Sigma$}${_{0}}$.

One notices that the estimate \eqref{e:(1)} for the coefficients of $s_{<}%
^{0}$ is stronger than the assertion that $\sum_{k<0}\lvert a_{k}\rvert
^{2}\lvert\lambda_{k}\rvert^{N}<+\infty$ for all natural $N$. Thus our
estimates confirm that not every smooth section can appear as initial value
over $\Sigma_{0}$ of a solution of {$\mathcal{D}$}${_{R}}f=0$ over the cylinder.

To sum up the example, the space $\Lambda(${$\mathcal{D}$}${_{R}},0)$ can be
written as the graph of an unbounded, densely defined, closed operator
$T:\operatorname{Dom}T\rightarrow H^{-\frac{1}{2}}(\Sigma_{R})$, mapping
$s_{<}^{0}+s_{>}^{0}=:s^{0}\mapsto s^{R}:=s_{<}^{R}+s_{>}^{R}$ with
$\operatorname{Dom}T\subset H^{-\frac{1}{2}}(\Sigma_{0})$\thinspace. To obtain
a closed subspace of $L^{2}(\Sigma)$ one takes the range $\Lambda
(${$\mathcal{D}$}${_{R}},\frac{1}{2})$ of the $L^{2}$ extension $\mathcal{P}%
(${$\mathcal{D}$}${_{R}})^{(0)}$ of the Calder{\'{o}}n projection. It
coincides with $\Lambda(${$\mathcal{D}$}${_{R}},0)\cap L^{2}(\Sigma)$ by
Proposition \ref{p:all-the-same}. In Theorem \ref{t:criss-cross} we show
without use of the pseudo-differential calculus why the intersection
$\Lambda(${$\mathcal{D}$}${_{R}},0)\cap L^{2}(\Sigma)$ must be closed in
$L^{2}(\Sigma)$. See \cite[p. 1214]{ScWo00} for another description of the
Cauchy data space $\Lambda(${$\mathcal{D}$}${_{R}},\frac{1}{2})$, namely as
the graph of a unitary elliptic pseudo-differential operator of order 0.

Since $\Sigma=-\Sigma_{0}\sqcup\Sigma_{R}$\thinspace, the tangential operator
takes the form $\mathcal{B}=B\oplus(-B)$ and we obtain from \eqref
{e:spectrum}
\[
\operatorname{range}P_{>}(\mathcal{B})^{(0)}=L_{+}(\mathcal{B}%
)=\operatorname{span}_{L^{2}(\Sigma)}\{(\varphi_{k},\sigma(\varphi
_{k}))\}_{k>0}\,.
\]
For comparison, we have in this example
\[
\operatorname{range}\mathcal{P}({\mathcal{D}_{R}})^{(0)}=\Lambda
({\mathcal{D}_{R}},{\tfrac{1}{2}})=\operatorname{span}_{L^{2}(\Sigma
)}\{(\varphi_{k},e^{-\lambda_{k}R}\varphi_{k})\}_{k>0}\,,
\]
hence $L_{+}(\mathcal{B})$ and $\Lambda(${$\mathcal{D}$}${_{R}},\frac{1}{2}) $
are transversal subspaces of $L^{2}(\Sigma)$. On the semi-infinite cylinder
$[0,\infty)\times\Sigma$, however, we have only one boundary component
$\Sigma_{0}$\thinspace. Hence
\[
\operatorname{range}P_{>}(B)^{(0)}=\operatorname{span}_{L^{2}(\Sigma_{0}%
)}\{\varphi_{k}\}_{k>0}=\lim_{R\rightarrow\infty}\operatorname{range}%
\mathcal{P}({\mathcal{D}_{R}})^{(0)}.
\]
\end{example}

\bigskip

One can generalize the preceding example: For any smooth compact manifold $X$
with boundary $\Sigma$ and any real $R\geq0$, let $X^{R}$ denote the stretched
manifold
\[
X^{R}\ :=\left(  [-R,0]\times\Sigma\right)  \cup_{\Sigma}X.
\]
Assuming product structures with $\mathcal{D}=\sigma(\partial_{u}%
+\mathcal{B})$ near $\Sigma$ gives a well-defined extension {$\mathcal{D}$%
}${_{R}}$ of $\mathcal{D}$. Nicolaescu, \cite{Ni95} proved that the
Calder{\'{o}n} projection and the Atiyah--Patodi--Singer projection coincide
up to a finite-dimensional component in the \emph{adiabatic limit}
($R\rightarrow+\infty$ in a suitable setting). Even for finite $R$ and, in
particular for $R=0$, one has the following interesting result. It was first
proved in Scott, \cite{Sco95} (see also Grubb, \cite{Gr99} and Wojciechowski,
\cite{DaKi99}, Appendix who both offered different proofs).

\begin{lemma}
For all $R\geq0$, the difference $\mathcal{P}(${$\mathcal{D}$}${_{R}}%
)-P_{\geq}(\mathcal{B})$ is an operator with a smooth kernel.
\end{lemma}

\medskip

\subsubsection{Twisted Orthogonality of Cauchy Data Spaces}

\label{TO} Green's formula (in particular the Clifford multiplication $\sigma$
in the case of Dirac type operators) provides a symplectic structure for
$L^{2}(\Sigma;S|_{\Sigma})$ for linear symmetric elliptic differential
operators of first order on a compact smooth manifold $X$ with boundary
$\Sigma$. For elliptic systems of second-order differential equations, various
interesting results have been obtained in the 1970s by exploiting the
symplectic structure of corresponding spaces (see e.g. \cite{LaSnTu75}).
Restricting oneself to first-order systems, the geometry becomes very clear
and it turns out that the Cauchy data space $\Lambda(\mathcal{D},\frac{1}{2})
$ is a Lagrangian subspace of $L^{2}(\Sigma;S|_{\Sigma})$.

More generally, in \cite{BoWo93} we described the orthogonal complement of the
Cauchy data space of the chiral Dirac operator $\mathcal{D}^{+}$ by
\begin{equation}
\sigma^{-1}(\Lambda(\mathcal{D}^{-},{\tfrac{1}{2}}))=(\Lambda(\mathcal{D}%
^{+},{\tfrac{1}{2}}))^{\perp}\,. \label{e:perp}%
\end{equation}
We obtained a short exact sequence
\[
0\rightarrow\sigma^{-1}(\Lambda(\mathcal{D}^{-},s))\hookrightarrow
H^{s-\frac{1}{2}}({\Sigma};S^{+}|_{\Sigma})\overset{\mathcal{K}^{\left(
s\right)  }}{\rightarrow}\operatorname{Ker}(\mathcal{D}^{+},s)\rightarrow0.
\]

For the total (symmetric) Dirac operator this means:

\begin{proposition}
The Cauchy data space $\Lambda(\mathcal{D},\frac{1}{2})$ of the \emph{total}
Dirac operator is a Lagrangian subspace of the Hilbert space $L^{2}%
(\Sigma;S|_{\Sigma})$ equipped with the symplectic form $\omega(\varphi
,\psi):=( ${$\sigma$}$\varphi,\psi)$.
\end{proposition}

\subsubsection{``Admissible'' Boundary Value Problems\label{sss:AdmissableBVP}}

We refer to \cite[Chapter 18]{BoWo93}, \cite{BrLe99}, and \cite{Sc01} for a
rigorous definition and treatment of large classes of admissible boundary
value problems defined by pseudo-differential projections. Prominent examples
belong to the Grassmannian $\mathcal{G}\mathrm{r}_{p_{+}}$ (introduced above
in Section \ref{sss:Calderon APS proj} before Example \ref{ex:cylinder}). On
even--dimensional manifolds, other prominent examples are local chiral
projections and unitary modifications as explained in \cite[p. 273]{BoWo93}.

For all admissible boundary conditions defined by a pseudo-differential
projection $R$ over $\Sigma$ the following features are common (For
simplicity, we suppress the distinction between total and chiral spinor bundle
in this paragraph and denote the bundle by $E$):

\begin{enumerate}
\item  We have an estimate
\begin{equation}
\left\|  f\right\|  _{1}\leq C\left(  \left\|  \mathcal{D}^{+}f\right\|
_{0}+\left\|  f\right\|  _{0}+\left\|  R\circ i\circ\gamma\left(  f\right)
\right\|  _{\frac{1}{2}}\right)  \;\text{for\ }f\in H^{1}\left(  X;E\right)  .
\label{e:a_priori_boundary}%
\end{equation}

\item  Defining a domain by
\begin{equation}
D=\operatorname{Dom}(\mathcal{D}_{R}^{+}):=\left\{  f\in H^{1}\left(
X;E\right)  \mid R\left(  f|_{\Sigma}\right)  =0\right\}  ,
\label{e:pseudodiff}%
\end{equation}
we obtain a closed Fredholm extension.

\item  The restriction of $\operatorname{Dom}(\mathcal{D}_{R}^{+})$ to the
boundary makes a Fredholm pair with the Cauchy data space of $\mathcal{D}^{+}$.

\item  The composition $R\mathcal{PD}^{+}$ defines a Fredholm operator from
the Cauchy data space to the range of $R$ (see also Proposition 3.34).

\item  The space $\operatorname{Ker}\left(  \mathcal{D}_{R}^{+}\right)  $
consists only of smooth sections.
\end{enumerate}

Warning: A new feature of operators of Dirac type on manifolds with boundary
is that the index of admissible boundary value problems can jump under
continuous or even smooth deformation of the coefficients. E.g., this is the
case for the Atiyah-Patodi-Singer boundary problem, as follows from the next subsection.

A total compatible (and so symmetric) Dirac operator $\mathcal{D}$ and an
orthogonal admissible projection $R$ define a self-adjoint extension
$\mathcal{D}_{D}^{+}=\mathcal{D}^{+\ast}|_{D}$, if the projection $R$ defines
the domain $D$ and satisfies the symmetry condition $I-R=\sigma^{-1}R\sigma$.
We denote by $\mathcal{G}\mathrm{r}_{p_{+}}^{\operatorname{sa}}$ or, shortly,
$\mathcal{G}\mathrm{r}^{\operatorname{sa}}$ the subspace of $\mathcal{G}%
\mathrm{r}_{p_{+}}$ of orthogonal projections which satisfy the preceding
condition and differ from the Atiyah--Patodi--Singer projection only by an
infinitely smoothing operator.

Posing a suitable \emph{well-posed} boundary value problem provides for a
nicely spaced discrete spectrum near 0. Then, varying the coefficients of the
differential operator and the imposed boundary condition suggests the use of
the powerful topological concept of spectral flow. From
(\ref{e:a_priori_boundary}) and a careful analysis of the corresponding
parametrices we see in \cite[Section 3]{BoLePh01} under which conditions the
curves of the induced self-adjoint $L^{2}$ extensions become continuous curves
in $\mathcal{CF}^{\operatorname{sa}}\left(  L^{2}(X;E)\right)  $ in the gap
topology so that their spectral flow is well-defined and truly homotopy invariant.

We summarize the main results. They depend strongly on the weak unique
continuation property (either in the form of Section \ref{WeakUCP} or in the
weaker form of (\ref{e:1.9-new}) which is sufficient here) and the invertible
extension (Section \ref{sss:InvertibleExtension}).

\begin{lemma}
\label{l:grass} For fixed $\mathcal{D}$ the mapping
\[
\mathcal{G}\mathrm{r}^{\operatorname{sa}}\left(  \mathcal{D}\right)  \ni
P\longmapsto\mathcal{D}_{P}\in\mathcal{CF}^{\operatorname{sa}}\left(
L^{2}(X;E)\right)
\]
is continuous from the operator norm to the gap metric.
\end{lemma}

\begin{theorem}
\label{t:dependence} Let $X$ be a compact Riemannian manifold with boundary.
Let $\left\{  \mathcal{D}_{s}\right\}  _{s\in M}$, $M$ a metric space, be a
family of compatible operators of Dirac type. We assume that in each local
chart, the coefficients of $\mathcal{D}_{s}$ depend continuously on $s$. Then
we have\newline (a) The Poisson operator $K_{s}:L^{2}\left(  \Sigma
;E|_{\Sigma}\right)  \rightarrow H^{1/2}\left(  X;E\right)  $ of
$\mathcal{D}_{s}$ depends continuously on $s$ in the operator norm.\newline
(b) The Calder\'{o}n projector $P_{+}(s):L^{2}\left(  \Sigma;E_{|\Sigma
}\right)  \rightarrow L^{2}\left(  \Sigma;E_{|\Sigma}\right)  $ of
$\mathcal{D}_{s}$ depends continuously on $s$ in the operator norm.\newline
(c) The family
\[
M\ni s\longmapsto\left(  \mathcal{D}_{s}\right)  _{P_{+}\left(  s\right)  }%
\in\mathcal{CF}^{\operatorname{sa}}\left(  L^{2}(X;E)\right)
\]
is continuous.\newline (d) Let $P_{t\{t\in Y\}}$ be a norm-continuous path of
orthogonal projections in $L^{2}(X;E)$. If
\[
P_{t}\in\bigcap_{s\in M}\mathcal{G}\mathrm{r}^{\operatorname{sa}}\left(
\mathcal{D}_{s}\right)  ,\;t\in Y,
\]
then
\[
M\times Y\ni\left(  s,t\right)  \longmapsto\left(  \mathcal{D}_{s}\right)
_{P_{t}}\in\mathcal{CF}^{\operatorname{sa}}\left(  L^{2}(X;E)\right)
\]
is continuous.
\end{theorem}

Note that (b) is a pseudo-differential reformulation of the continuity of
Cauchy data spaces which is valid in much greater generality (see our Section
\ref{sss:Symmetric Operators}).

\medskip

We close this Section with a recent result of \cite{Gr02}:

\begin{lemma}
\label{l:grubb}\textrm{(G. Grubb 2002)} Let $\mathcal{D}$ be an operator of
Dirac type over a compact manifold with boundary. Let $P$ be an orthogonal
projection which defines an `admissible' self--adjoint boundary condition for
$\mathcal{D}$. Then there exists an invertible operator of Dirac type
$\mathcal{B}^{\prime}$ over the boundary such that $P=P_{>}(\mathcal{B}%
^{\prime})$.
\end{lemma}

Grubb's Lemma shows that the Atiyah--Patodi--Singer boundary projection is the
most general admissible self--adjoint boundary condition, in the specified sense.

\subsection{The Atiyah--Patodi--Singer Index Theorem\label{ss:APSIndexThm}}

\bigskip

Let $X$ be a compact, oriented Riemannian manifold with boundary $Y=\partial
X$ with $\dim X=n=2m$ even. Let $\mathcal{D}:C^{\infty}\left(  X,S\right)
\rightarrow C^{\infty}\left(  X,S\right)  $ be a compatible operator of Dirac
type where $S\rightarrow X$ is a bundle of Clifford modules. Relative to the
splitting $S=S^{+}\oplus S^{-}$ into chiral halves, we have the operators
$\mathcal{D}^{+}:C^{\infty}\left(  X,S^{+}\right)  \rightarrow C^{\infty
}\left(  X,S^{-}\right)  $ and $\mathcal{D}^{-}:C^{\infty}\left(
X,S^{-}\right)  \rightarrow C^{\infty}\left(  X,S^{+}\right)  $ which are
formal adjoints on sections with support in $X\setminus Y$. We assume that all
structures (e.g., Riemannian metric, Clifford module, connection) are products
on some collared neighborhood $N\cong\lbrack-1,1]\times Y$ of $Y$. Then
$\mathcal{D}^{+}|_{{N}}:=\mathcal{D}^{+}:C^{\infty}\left(  N,S^{+}|N\right)
\rightarrow C^{\infty}\left(  N,S^{-}|N\right)  $ has the form
\[
\mathcal{D}^{+}|_{{N}}=\sigma(\partial_{u}+\mathcal{B}).
\]
Here $u\in\lbrack-1,1]$ is the normal coordinate (i.e., $N=\left\{  \left(
u,y\right)  \mid y\in Y,\text{ }u\in\lbrack-1,1]\right\}  $) with
$\partial_{u}=\frac{\partial}{\partial u}$ the inward normal), $\sigma
=\mathbf{c}\left(  du\right)  $ is the (unitary) Clifford multiplication by
$du$ with $\sigma\left(  S^{+}|N\right)  =S^{-}|N$, and
\[
\mathcal{B}:C^{\infty}\left(  Y,S^{+}|_{Y}\right)  \rightarrow C^{\infty
}\left(  Y,S^{+}|_{Y}\right)
\]
denotes the canonically associated (elliptic, self-adjoint) Dirac operator
over $Y$, called the \textit{tangential operator}. Note that due to the
product structure, $\sigma$ and $\mathcal{B}$ do not depend on $u$. Let
$P_{\geq}(\mathcal{B})$ denote the spectral (Atiyah--Patodi--Singer)
projection onto the subspace $L_{+}(\mathcal{B})$ of $L^{2}(Y,S^{+}|_{\partial
X})$ spanned by the eigensections corresponding to the nonnegative eigenvalues
of $\mathcal{B}$. Let
\[
C^{\infty}\left(  X,S^{+};P_{\geq}\right)  :=\left\{  \psi\in C^{\infty
}\left(  X,S^{+}\right)  \mid P_{\geq}(\mathcal{B})\left(  \psi|_{Y}\right)
=0\right\}  ,\text{ and}%
\]%
\[
\mathcal{D}_{P_{\geq}}^{+}:=\mathcal{D}^{+}|_{C^{\infty}\left(  X,S^{+}%
;P_{\geq}\right)  }:C^{\infty}\left(  X,S^{+};P_{\geq}\right)  \rightarrow
C^{\infty}\left(  X,S^{-}\right)  .
\]
The eta function for $\mathcal{B}$ is defined by
\[
\eta_{\mathcal{B}}\left(  s\right)  :=\sum\nolimits_{\lambda\in
\operatorname{spec}\mathcal{B-}\left\{  0\right\}  }\left(
\operatorname*{sign}\lambda\right)  m_{\lambda}\left|  \lambda\right|
^{-s}\text{, }%
\]
for $\mathfrak{R}\left(  s\right)  $ sufficiently large, where $m_{\lambda}$
is the multiplicity of $\lambda$. Implicit in the following result
(originating in \cite{AtPaSi75}) is that $\eta_{\mathcal{B}}$ extends to a
meromorphic function on all $\mathbb{C}$, which is holomorphic at $s=0$ so
that $\eta_{\mathcal{B}}\left(  0\right)  $ is finite.

\begin{theorem}
[Atiyah-Patodi-Singer Index Formula]\label{tAPS}The above operator
$\mathcal{D}_{P_{\geq}}^{+}$ has finite index given by
\[
\operatorname{index}\mathcal{D}_{P_{\geq}}^{+}=\int_{X}\left(  \mathbf{ch}%
\left(  S,\varepsilon\right)  \wedge\widetilde{\mathbf{A}}\left(
X,\theta\right)  \right)  -\frac{m_{0}+\eta_{\mathcal{B}}\left(  0\right)
}{2}.
\]
where $m_{0}=\dim\left(  \operatorname{Ker}\mathcal{B}\right)  ,$
$\mathbf{ch}\left(  S,\varepsilon\right)  \in\Omega^{\ast}\left(
X,\mathbb{R}\right)  $ is the total Chern character form of the complex vector
bundle $S$ with compatible, unitary connection $\varepsilon$, and
$\widetilde{\mathbf{A}}\left(  X,\theta\right)  \in\Omega^{\ast}\left(
X,\mathbb{R}\right)  $ is closely related to the total $\widehat{\mathbf{A}%
}\left(  X,\theta\right)  $\textbf{\ }form\textbf{\ }relative to the
Levi-Civita connection $\theta$, namely $\widetilde{\mathbf{A}}\left(
X,\theta\right)  _{4k}=2^{2k-m}\widehat{\mathbf{A}}\left(  X,\theta\right)
_{4k}$.
\end{theorem}

\begin{proof}
(outline) The proof found in \cite{AtPaSi75} or \cite{BoWo93} is based on the
heat kernel method for computing the index, but the process is less
straightforward than in the closed case because of the boundary condition. The
appropriate heat kernel is constructed by means of Duhamel's method. Namely,
an exact kernel is obtained from an approximate one by an iterative process
initiated by writing the error as the integral of a derivative of the
convolution of the true and approximate kernel; see (\ref{Duhamel}) below. The
initial approximate heat kernel is obtained by patching together two heat
kernels, denoted by $\mathcal{E}_{c}$ and $\mathcal{E}_{d}$. Here,
$\mathcal{E}_{c}$ is a heat kernel for a Dirac operator over an infinite
extension $\left[  0,\infty\right)  \times Y$ of the collared neighborhood
$N=\left[  0,1\right]  \times Y$ of $Y$ in $X$, for which the boundary
condition $P_{\geq}(\mathcal{B})\left(  \psi|_{Y}\right)  =0$ is imposed. The
heat kernel $\mathcal{E}_{d}$ is the usual one (without boundary conditions)
for $e^{-t\widetilde{\mathcal{D}^{-}}\widetilde{\mathcal{D}^{+}}}$, where
$\widetilde{\mathcal{D}^{\pm}}$ are chiral halves of the invertible Dirac
operator (see (\ref{DirDblOp})), namely
\[
\widetilde{\mathcal{D}^{\pm}}:=\mathcal{D}^{\pm}\cup\mathcal{D}^{\mp
}:C^{\infty}\left(  \widetilde{X},\widetilde{S^{\pm}}\right)  \rightarrow
C^{\infty}\left(  \widetilde{X},\widetilde{S^{\pm}}\right)
\]
over the double $\widetilde{X}$, a closed manifold without boundary.

We begin with the construction of $\mathcal{E}_{c}$. We let
\[
D:=\partial_{u}+\mathcal{B}:C^{\infty}\left(  N,S^{+}|N\right)  \rightarrow
C^{\infty}\left(  N,S^{+}|N\right)  \text{,}%
\]
which has formal adjoint $D^{\ast}:=-\partial_{u}+\mathcal{B}$. Define
\begin{align*}
&  \mathcal{D}:=D|\operatorname{Dom}\mathcal{D}\text{, where}\\
&  \operatorname{Dom}\mathcal{D}:=\left\{  f\in H^{1}\left(  \mathbb{R}%
_{+}\times Y,\pi^{\ast}\left(  S^{+}|Y\right)  \right)  \mid P_{\geq}\left(
f|_{\left\{  0\right\}  \times Y}\right)  =0\right\}  \text{, and}\\
&  \mathcal{D}^{\ast}=D^{\ast}|\operatorname{Dom}\mathcal{D}^{\ast}\text{,
where}\\
&  \operatorname{Dom}\mathcal{D}^{\ast}:=\left\{  f\in H^{1}\left(
\mathbb{R}_{+}\times Y,\pi^{\ast}\left(  S^{-}|Y\right)  \right)  \mid
P_{<}\left(  f|_{\left\{  0\right\}  \times Y}\right)  =0\right\}  .
\end{align*}
We also have the Laplacians given by
\begin{align*}
&  \Delta_{c}:=\mathcal{D}^{\ast}\mathcal{D}|\operatorname{Dom}\mathcal{D}%
^{\ast}\mathcal{D}\text{, where}\\
&  \operatorname{Dom}\mathcal{D}^{\ast}\mathcal{D}:=\left\{
\begin{array}
[c]{l}%
f\in H^{1}\left(  \mathbb{R}_{+}\times Y,\pi^{\ast}\left(  S^{+}|Y\right)
\right)  \mid\\
P_{\geq}\left(  f|_{\left\{  0\right\}  \times Y}\right)  =0,\text{ }%
P_{<}\left(  Df|_{\left\{  0\right\}  \times Y}\right)  =0
\end{array}
\right\}  ,\text{ and}\\
&  \Delta_{c\ast}:=\mathcal{DD}^{\ast}|\operatorname{Dom}\mathcal{DD}^{\ast
}\text{, where}\\
&  \operatorname{Dom}\mathcal{DD}^{\ast}:=\left\{
\begin{array}
[c]{l}%
f\in H^{1}\left(  \mathbb{R}_{+}\times Y,\pi^{\ast}\left(  S^{+}|Y\right)
\right)  \mid\\
P_{<}\left(  f|_{\left\{  0\right\}  \times Y}\right)  =0,\text{ }P_{\geq
}\left(  D^{\ast}f|_{\left\{  0\right\}  \times Y}\right)  =0
\end{array}
\right\}  .
\end{align*}
Let $\varphi_{\lambda}\left(  y\right)  \in C^{\infty}\left(  Y,S^{+}%
|Y\right)  $ be an eigensection of $\mathcal{B}$ with eigenvalue $\lambda
\in\mathbb{R}$. Note that $g_{\lambda}\left(  t;u,y\right)  =f_{\lambda
}\left(  t,u\right)  \varphi_{\lambda}\left(  y\right)  $ is a solution of the
heat equation
\begin{align*}
0  &  =\left(  \partial_{t}+\Delta_{c}\right)  g_{\lambda}=\partial
_{t}g_{\lambda}+\left(  \partial_{u}+\mathcal{B}\right)  \left(  -\partial
_{u}+\mathcal{B}\right)  g_{\lambda}=\left(  \partial_{t}-\partial_{u}%
^{2}+\mathcal{B}^{2}\right)  g_{\lambda}\\
&  =\left(  \partial_{t}f_{\lambda}-\partial_{u}^{2}f_{\lambda}+\lambda
^{2}f_{\lambda}\right)  \varphi_{\lambda}\left(  y\right)  ,
\end{align*}
with $g_{\lambda}\left(  t;\cdot,\cdot\right)  \in\operatorname{Dom}%
\mathcal{D}^{\ast}\mathcal{D}$, when $f_{\lambda}:\left(  0,\infty\right)
\times\lbrack0,\infty)\rightarrow\mathbb{R}$ solves the heat problem
\begin{align}
&  \partial_{t}f_{\lambda}=\partial_{u}^{2}f_{\lambda}-\lambda^{2}f_{\lambda
}\text{ with boundary condition}\nonumber\\
&  f_{\lambda}\left(  t,0\right)  =0\text{ if }\lambda\geq0\nonumber\\
&  \partial_{u}f_{\lambda}\left(  t,0\right)  +\lambda f_{\lambda}\left(
t,0\right)  =Df_{\lambda}=0\text{ if }\lambda<0. \label{flam}%
\end{align}
Recall that the\textit{\ complementary error function} is
\[
\operatorname{erfc}\left(  x\right)  :=\frac{2}{\sqrt{\pi}}\int_{x}^{\infty
}e^{-\xi^{2}}d\xi\text{.}%
\]
Of use to us are the facts
\begin{align}
&  \operatorname{erfc}^{\prime}\left(  x\right)  =\frac{-2}{\sqrt{\pi}%
}e^{-x^{2}}\text{ and }\nonumber\\
&  \operatorname{erfc}\left(  x\right)  \leq\frac{2}{\sqrt{\pi}}%
\frac{e^{-x^{2}}}{x+\sqrt{x^{2}+\frac{4}{\pi}}}\leq e^{-x^{2}}\text{ for
}x\geq0. \label{ErfcEst}%
\end{align}
For $\lambda\geq0,$ the heat kernel for the problem (\ref{flam}) is (via the
method of images)
\[
e_{\lambda}\left(  t;u,v\right)  :=\frac{e^{-\lambda^{2}t}}{2\sqrt{\pi t}%
}\left(  e^{-\frac{\left(  u-v\right)  ^{2}}{4t}}-e^{-\frac{\left(
u+v\right)  ^{2}}{4t}}\right)  \text{ for }u,v\geq0\text{ and }t>0.
\]
For $\lambda<0,$ and $u,v\geq0$ and $t>0,$ the heat kernel is (using Laplace
transforms)
\[
e_{\lambda}\left(  t;u,v\right)  :=\frac{e^{-\lambda^{2}t}}{2\sqrt{\pi t}%
}\left(  e^{-\frac{\left(  u-v\right)  ^{2}}{4t}}+e^{-\frac{\left(
u+v\right)  ^{2}}{4t}}+\lambda e^{\lambda\left(  u+v\right)  }%
\operatorname{erfc}\left(  \tfrac{u+v}{2\sqrt{t}}-\lambda\sqrt{t}\right)
\right)  .
\]
The heat kernel for $\Delta_{c}$ (i.e., the kernel for $e^{-t\Delta_{c}}$) is
then
\begin{equation}
\mathcal{E}_{c}\left(  t;u,y;v,z\right)  =\sum\nolimits_{\lambda}e_{\lambda
}\left(  t;u,v\right)  \varphi_{\lambda}\left(  y\right)  \otimes
\varphi_{\lambda}\left(  z\right)  ^{\ast} \label{Kc}%
\end{equation}
Here $\lambda$ ranges over $\operatorname{spec}\mathcal{B}$ and (by a
convenient abuse of notation) $\varphi_{\lambda}\left(  y\right)
\otimes\varphi_{\lambda}\left(  z\right)  ^{\ast}$ is really a sum
$\sum\nolimits_{k}\varphi_{\lambda,k}\left(  y\right)  \otimes\varphi
_{\lambda,k}\left(  z\right)  ^{\ast}$ where $\left\{  \varphi_{\lambda
,1},\ldots,\varphi_{\lambda,k}\right\}  $ is an orthonormal eigenbasis of the
eigenspace for $\lambda$, and $\left\{  \varphi_{\lambda,1}^{\ast}%
,\ldots,\varphi_{\lambda,k}^{\ast}\right\}  $ is the dual basis. Similarly,
the heat kernel for $\Delta_{c\ast}$ is
\begin{equation}
\mathcal{E}_{c\ast}\left(  t;u,y;v,z\right)  =\sum\nolimits_{\lambda
}e_{\lambda\ast}\left(  t;u,v\right)  \varphi_{\lambda}\left(  y\right)
\otimes\varphi_{\lambda}\left(  z\right)  ^{\ast}, \label{Kcstar}%
\end{equation}
where $e_{\lambda\ast}$ is the heat kernel for the problem
\begin{align}
&  \partial_{t}f_{\lambda}=\partial_{u}^{2}f_{\lambda}-\lambda^{2}f_{\lambda
}\text{ with boundary condition}\nonumber\\
&  f_{\lambda}\left(  t,0\right)  =0\text{ if }\lambda<0\nonumber\\
&  \partial_{u}f_{\lambda}\left(  t,0\right)  -\lambda f_{\lambda}\left(
t,0\right)  =D^{\ast}f_{\lambda}=0\text{ if }\lambda\geq0, \label{flamstar}%
\end{align}
where $f_{\lambda}:\left(  0,\infty\right)  \times\lbrack0,\infty
)\rightarrow\mathbb{R}$. Thus, for $u,v\geq0$ and $t>0$, we have
\[
e_{\lambda\ast}\left(  t;u,v\right)  :=\frac{e^{-\lambda^{2}t}}{2\sqrt{\pi t}%
}\left(  e^{-\frac{\left(  u-v\right)  ^{2}}{4t}}-e^{-\frac{\left(
u+v\right)  ^{2}}{4t}}\right)  \text{ for }\lambda<0,
\]
while for $\lambda\geq0$ (note the sign change in passing from (\ref{flam}) to
(\ref{flamstar}))
\[
e_{\lambda\ast}\left(  t;u,v\right)  :=\frac{e^{-\lambda^{2}t}}{2\sqrt{\pi t}%
}\left(  e^{-\frac{\left(  u-v\right)  ^{2}}{4t}}+e^{-\frac{\left(
u+v\right)  ^{2}}{4t}}-\lambda e^{\lambda\left(  u+v\right)  }%
\operatorname{erfc}\left(  \tfrac{u+v}{2\sqrt{t}}+\lambda\sqrt{t}\right)
\right)  .
\]
Combining (\ref{Kc}) and (\ref{Kcstar}), we obtain the trace of kernel
$\mathcal{E}_{c}-\mathcal{E}_{c\ast}$ for $e^{-t\Delta_{c}}-e^{-t\Delta
_{c\ast}}$ evaluated at the point $\left(  u,y;u,y\right)  $ of the diagonal
\begin{align*}
&  \mathcal{K}\left(  t;u,y\right)  :=\operatorname{Tr}\sum\nolimits_{\lambda
}\left(  e_{\lambda}\left(  t;u,u\right)  -e_{\ast\lambda}\left(
t;u,u\right)  \right)  \varphi_{\lambda}\left(  y\right)  \otimes
\varphi_{\lambda}^{\ast}\left(  y\right) \\
&  =\sum_{\lambda\geq0}\left(  \frac{-e^{-\lambda^{2}t}}{2\sqrt{\pi t}%
}e^{-u^{2}/t}-\frac{e^{-\lambda^{2}t}}{2\sqrt{\pi t}}e^{-u^{2}/t}+\lambda
e^{2\lambda u}\operatorname{erfc}\left(  \tfrac{u}{\sqrt{t}}+\lambda\sqrt
{t}\right)  \right)  \left|  \varphi_{\lambda}\left(  y\right)  \right|
^{2}\\
&  +\sum_{\lambda<0}\left(  \frac{e^{-\lambda^{2}t}}{2\sqrt{\pi t}}%
e^{-u^{2}/t}+\lambda e^{-2\lambda u}\operatorname{erfc}\left(  \tfrac{u}%
{\sqrt{t}}-\lambda\sqrt{t}\right)  +\frac{e^{-\lambda^{2}t}}{2\sqrt{\pi t}%
}e^{-u^{2}/t}\right)  \left|  \varphi_{\lambda}\left(  y\right)  \right|
^{2}\\
&  =\sum_{\lambda}\operatorname*{sign}\left(  \lambda\right)  \left(
\frac{-e^{-\lambda^{2}t}e^{-u^{2}/t}}{\sqrt{\pi t}}+\left|  \lambda\right|
e^{2\left|  \lambda\right|  u}\operatorname{erfc}\left(  \tfrac{u}{\sqrt{t}%
}+\left|  \lambda\right|  \sqrt{t}\right)  \right)  \left|  \varphi_{\lambda
}\left(  y\right)  \right|  ^{2}\\
&  =\sum\nolimits_{\lambda}\operatorname*{sign}\left(  \lambda\right)
\frac{\partial}{\partial u}\left(  \tfrac{1}{2}e^{2\left|  \lambda\right|
u}\operatorname{erfc}\left(  \tfrac{u}{\sqrt{t}}+\left|  \lambda\right|
\sqrt{t}\right)  \right)  \left|  \varphi_{\lambda}\left(  y\right)  \right|
^{2}.
\end{align*}
We used $\operatorname{erfc}^{\prime}\left(  x\right)  =\frac{-2}{\sqrt{\pi}%
}e^{-x^{2}}$ to obtain the last equality, and we have set $\left.
\operatorname*{sign}\left(  0\right)  :=1\right.  $ for convenience. We
compute (where the interchange of the sum and integral can be justified)
\begin{align*}
&  \mathcal{K}\left(  t\right)  :=\int_{0}^{\infty}\int_{Y}\mathcal{K}\left(
t;u,y\right)  dydu\\
&  =\int_{0}^{\infty}\int_{Y}\sum\nolimits_{\lambda}\operatorname*{sign}%
\left(  \lambda\right)  \frac{\partial}{\partial u}\left(  \tfrac{1}%
{2}e^{2\left|  \lambda\right|  u}\operatorname{erfc}\left(  \tfrac{u}{\sqrt
{t}}+\left|  \lambda\right|  \sqrt{t}\right)  \right)  \left|  \varphi
_{\lambda}\left(  y\right)  \right|  ^{2}dydu\\
&  =\sum_{\lambda}m_{\lambda}\operatorname*{sign}\left(  \lambda\right)
\int_{0}^{\infty}\frac{\partial}{\partial u}\left(  \tfrac{1}{2}e^{2\left|
\lambda\right|  u}\operatorname{erfc}\left(  \tfrac{u}{\sqrt{t}}+\left|
\lambda\right|  \sqrt{t}\right)  \right)  \text{ }du\\
&  =\sum_{\lambda}m_{\lambda}\operatorname*{sign}\left(  \lambda\right)
\left.  \tfrac{1}{2}e^{2\left|  \lambda\right|  u}\operatorname{erfc}\left(
\tfrac{u}{\sqrt{t}}+\left|  \lambda\right|  \sqrt{t}\right)  \right|
_{u=0}^{\infty}\\
&  =-\sum_{\lambda}\frac{m_{\lambda}\operatorname*{sign}\lambda}%
{2}\operatorname{erfc}\left(  \left|  \lambda\right|  \sqrt{t}\right)
=-\tfrac{1}{2}m_{0}-\sum_{\lambda\neq0}m_{\lambda}\frac{\operatorname*{sign}%
\lambda}{2}\operatorname{erfc}\left(  \left|  \lambda\right|  \sqrt{t}\right)
.
\end{align*}
Recall that $m_{\lambda}$ is the multiplicity of $\lambda$. Thus,
\[
\mathcal{K}\left(  t\right)  +\tfrac{1}{2}m_{0}=-\sum_{\lambda\neq0}%
\frac{m_{\lambda}\operatorname*{sign}\lambda}{2}\operatorname{erfc}\left(
\left|  \lambda\right|  \sqrt{t}\right)
\]
For $\left|  \lambda\right|  >0$, one verifies using integration by parts and
substitution that
\[
\int_{0}^{\infty}\operatorname{erfc}\left(  \left|  \lambda\right|  \sqrt
{t}\right)  t^{s-1}dt=\frac{\left|  \lambda\right|  ^{-2s}}{s\sqrt{\pi}}%
\Gamma\left(  s+\tfrac{1}{2}\right)  .
\]
Using this and the fact that $\mathcal{K}\left(  t\right)  +\tfrac{1}{2}h\leq
Ce^{-\alpha t}$ for constants $C$ and $\alpha>0$,
\begin{align*}
&  \int_{0}^{\infty}\left(  \mathcal{K}\left(  t\right)  +\tfrac{1}{2}%
m_{0}\right)  t^{s-1}dt\\
&  =-\int_{0}^{\infty}\left(  \sum\nolimits_{\lambda\neq0}\frac{m_{\lambda
}\operatorname*{sign}\lambda}{2}\operatorname{erfc}\left(  \left|
\lambda\right|  \sqrt{t}\right)  \right)  t^{s-1}dt\\
&  =-\sum\nolimits_{\lambda\neq0}\frac{m_{\lambda}\operatorname*{sign}\lambda
}{2}\int_{0}^{\infty}\operatorname{erfc}\left(  \left|  \lambda\right|
\sqrt{t}\right)  t^{s-1}dt\\
&  =-\sum\nolimits_{\lambda\neq0}\frac{m_{\lambda}\operatorname*{sign}\lambda
}{2}\frac{\left|  \lambda\right|  ^{-2s}}{s\sqrt{\pi}}\Gamma\left(
s+\tfrac{1}{2}\right) \\
&  =-\frac{\Gamma\left(  s+\frac{1}{2}\right)  }{2s\sqrt{\pi}}\sum
\nolimits_{\lambda\neq0}m_{\lambda}\operatorname*{sign}\left(  \lambda\right)
\left|  \lambda\right|  ^{-2s}\\
&  =-\frac{\Gamma\left(  s+\frac{1}{2}\right)  }{2s\sqrt{\pi}}\eta_{B}\left(
2s\right)  .
\end{align*}
\textit{Suppose} that we have an asymptotic expansion
\begin{align*}
&  \mathcal{K}\left(  t\right)  \sim\sum\nolimits_{k=-n+1}^{N}a_{k}%
t^{k/2}\text{ as }t\rightarrow0^{+};\text{ i.e.,}\\
&  \mathcal{K}\left(  t\right)  -\sum\nolimits_{k=-n+1}^{N}a_{k}%
t^{k/2}=\operatorname{O}(t^{\frac{1}{2}\left(  N+1\right)  })\text{ as
}t\rightarrow0^{+}.
\end{align*}
In (\ref{KasyExp}) below, such an asymptotic expansion will eventually be
produced (as was done in \cite[p. 239f]{BoWo93}) from the known asymptotic
expansion (see \cite[p. 68]{Gi95}) of the heat kernel on a closed manifold,
namely the double of $M$. Since $\mathfrak{R}\left(  s\right)  >-\frac{N+1}%
{2}\Rightarrow\mathfrak
{R}\left(  \tfrac{1}{2}\left(  N+1\right)  +s-1\right)  >-1$, we then have
that
\[
f_{1}\left(  s\right)  :=\int_{0}^{1}\left(  \mathcal{K}\left(  t\right)
-\sum\nolimits_{k=-n+1}^{N}a_{k}t^{k/2}\right)  t^{s-1}dt=\int_{0}%
^{1}\operatorname{O}\left(  t^{\frac{1}{2}\left(  N+1\right)  +s-1}\right)
ds
\]
is holomorphic for $\mathfrak{R}\left(  s\right)  >-\frac{N+1}{2}$. We also
have the entire function
\[
f_{\infty}\left(  s\right)  :=\int_{1}^{\infty}\left(  \mathcal{K}\left(
t\right)  +\tfrac{1}{2}h\right)  t^{s-1}dt
\]
We claim that
\[
-\frac{\Gamma\left(  s+\frac{1}{2}\right)  }{2s\sqrt{\pi}}\eta_{\mathcal{B}%
}\left(  2s\right)  =\frac{m_{0}}{2s}+\sum\nolimits_{k=-n+1}^{N}\frac{a_{k}%
}{s+\frac{1}{2}k}+f_{1}\left(  s\right)  +f_{\infty}\left(  s\right)  .
\]
Indeed, we have
\begin{align*}
&  -\frac{\Gamma\left(  s+\frac{1}{2}\right)  }{2s\sqrt{\pi}}\eta
_{\mathcal{B}}\left(  2s\right) \\
&  =\int_{0}^{\infty}\left(  \mathcal{K}\left(  t\right)  +\tfrac{1}%
{2}h\right)  t^{s-1}dt=\int_{0}^{1}\left(  \mathcal{K}\left(  t\right)
+\tfrac{1}{2}h\right)  t^{s-1}dt+f_{\infty}\left(  s\right) \\
&  =\int_{0}^{1}\left(  \tfrac{1}{2}m_{0}+\sum\nolimits_{k=-n+1}^{N}%
a_{k}t^{k/2}\right)  t^{s-1}dt+f_{1}\left(  s\right)  +f_{\infty}\left(
s\right) \\
&  =\int_{0}^{1}\tfrac{1}{2}m_{0}t^{s-1}dt+\sum\nolimits_{k=-n+1}^{N}a_{k}%
\int_{0}^{1}t^{\left(  s+\frac{1}{2}k\right)  -1}dt+f_{1}\left(  s\right)
+f_{\infty}\left(  s\right) \\
&  =\tfrac{1}{2}m_{0}\left(  \left.  \frac{t^{s}}{s}\right|  _{t=0}%
^{t=1}\right)  +\sum\nolimits_{k=-n+1}^{N}a_{k}\left.  \frac{t^{s+\frac{1}%
{2}k}}{s+\frac{1}{2}k}\right|  _{t=0}^{t=1}+f_{1}\left(  s\right)  +f_{\infty
}\left(  s\right) \\
&  =\frac{m_{0}}{2s}+\sum\nolimits_{k=-n+1}^{N}\frac{a_{k}}{s+\frac{1}{2}%
k}+f_{1}\left(  s\right)  +f_{\infty}\left(  s\right)  .
\end{align*}
Thus, where $\theta_{N}\left(  s\right)  =f_{1}\left(  s\right)  +f_{\infty
}\left(  s\right)  $ is holomorphic for $\mathfrak{R}\left(  s\right)
>-\frac{N+1}{2}$,
\begin{equation}
\eta_{\mathcal{B}}\left(  2s\right)  =-\frac{2s\sqrt{\pi}}{\Gamma\left(
s+\tfrac{1}{2}\right)  }\left(  \frac{m_{0}}{2s}+\sum\nolimits_{k=-n+1}%
^{N}\frac{a_{k}}{\tfrac{1}{2}k+s}+\theta_{N}\left(  s\right)  \right)  .
\label{EtaBAsy}%
\end{equation}
The heat kernels, say $\mathcal{E}_{d}$ for $\widetilde{\mathcal{D}^{-}%
}\widetilde{\mathcal{D}^{+}}$ and $\mathcal{E}_{d\ast}$ for $\widetilde
{\mathcal{D}^{+}}\widetilde{\mathcal{D}^{-}}$, over the double $\widetilde{X}$
are more familiar. For $t>0$,
\[
\operatorname{Tr}e^{-t\widetilde{\mathcal{D}^{-}}\widetilde{\mathcal{D}^{+}}%
}-\operatorname{Tr}e^{-t\widetilde{\mathcal{D}^{+}}\widetilde{\mathcal{D}^{-}%
}}=\int_{\widetilde{X}}\operatorname{Tr}\left(  \mathcal{E}_{d}\left(
t;x,x\right)  -\mathcal{E}_{d\ast}\left(  t;x,x\right)  \right)  \text{ }dx,
\]
and there is the asymptotic expansion
\[
F\left(  t;x\right)  :=\operatorname{Tr}\left(  \mathcal{E}_{d}\left(
t;x,x\right)  -\mathcal{E}_{d\ast}\left(  t;x,x\right)  \right)  \sim
\sum\nolimits_{k\geq-n}\alpha_{k}\left(  x\right)  t^{k/2}.
\]
For $0<a<b<1$, let $\rho_{\left(  a,b\right)  }:\left[  0,1\right]
\rightarrow\left[  0,1\right]  $ be $C^{\infty}$ and increasing with
\[
\rho_{\left(  a,b\right)  }\left(  u\right)  =\left\{
\begin{array}
[c]{cl}%
0 & \text{for }u\leq a\\
1 & \text{for }u\geq b.
\end{array}
\right.
\]
Thinking of $u$ as the normal coordinate function on $N=\left[  0,1\right]
\times Y\subset X$, and extending by the constant values $0$ and $1$, we can
(and do) regard $\rho_{\left(  a,b\right)  }$ as a function on $X$. Let
\begin{align*}
&  Q\left(  t;x,x^{\prime}\right)  :=\varphi_{c}\left(  x\right)
\mathcal{E}_{c}\left(  t;x,x^{\prime}\right)  \psi_{c}\left(  x^{\prime
}\right)  +\varphi_{d}\left(  x\right)  \mathcal{E}_{d}\left(  t;x,x^{\prime
}\right)  \psi_{d}\left(  x^{\prime}\right) \\
&  :=\left(  1-\rho_{\left(  5/7,6/7\right)  }\left(  x\right)  \right)
\mathcal{E}_{c}\left(  t;x,x^{\prime}\right)  \left(  1-\rho_{\left(
3/7,4/7\right)  }\left(  x^{\prime}\right)  \right) \\
&  +\rho_{\left(  1/7,2/7\right)  }\left(  x\right)  \mathcal{E}_{d}\left(
t;x,x^{\prime}\right)  \rho_{\left(  3/7,4/7\right)  }\left(  x^{\prime
}\right)  .
\end{align*}
We have that $\psi_{c}+\psi_{d}=1$ and $\left\{  \psi_{c},\psi_{d}\right\}  $
is a partition of unity for the cover $\left\{  u^{-1}\left(  \tfrac{3}%
{7},\infty\right)  ,u^{-1}\left[  0,\tfrac{4}{7}\right)  \right\}  $ of $X$.
Moreover, for $j=c,d$
\begin{equation}
\varphi_{j}|\operatorname*{supp}\psi_{j}=1\text{\ and\ }\operatorname*{dist}%
\left(  \operatorname*{supp}\partial_{u}^{k}\varphi_{j},\operatorname*{supp}%
\psi_{j}\right)  \geq\tfrac{1}{7}\text{ (for }k\geq1\text{).}
\label{suppphipsi}%
\end{equation}
Of course, we do not expect $Q\left(  t;x,x^{\prime}\right)  $ to equal the
exact kernel for $\partial_{t}+\mathcal{D}^{\ast}\mathcal{D}$ throughout $X$.
However, note that for $x,x^{\prime}\in\left[  0,1/7\right)  \times Y$,
$Q\left(  t;x,x^{\prime}\right)  =\mathcal{E}_{c}\left(  t;x,x^{\prime
}\right)  $, so that $Q\left(  t;x,x^{\prime}\right)  $ meets the APS boundary
condition. Let $x^{\prime}\in X$ be fixed. For $x\in X\setminus\left[
0,6/7\right)  \times Y,$ we have $Q\left(  t;x,x^{\prime}\right)
=\mathcal{E}_{d}\left(  t;x,x^{\prime}\right)  \psi_{d}\left(  x^{\prime
}\right)  ,$ while for $x\in\left[  0,1/7\right)  \times Y,$ $Q\left(
t;x,x^{\prime}\right)  =\mathcal{E}_{c}\left(  t;x,x^{\prime}\right)  \psi
_{c}\left(  x^{\prime}\right)  $. Thus, for $x\in X\setminus\left(  \left[
\tfrac{1}{7},\tfrac{6}{7}\right]  \times Y\right)  $, $-\left(  \partial
_{t}+\mathcal{D}^{\ast}\mathcal{D}\right)  \left.  Q\left(  \cdot
;\cdot,x^{\prime}\right)  \right|  _{\left(  t;x,x^{\prime}\right)  }=0$.
Moreover, for $d\left(  x,x^{\prime}\right)  <1/7$,
\begin{align*}
&  \varphi_{c}\left(  x\right)  \psi_{c}\left(  x^{\prime}\right)
+\varphi_{d}\left(  x\right)  \psi_{d}\left(  x^{\prime}\right) \\
&  =\left(  1-\rho_{\left(  5/7,6/7\right)  }\left(  x\right)  \right)
\left(  1-\rho_{\left(  3/7,4/7\right)  }\left(  x^{\prime}\right)  \right)
+\rho_{\left(  1/7,2/7\right)  }\left(  x\right)  \rho_{\left(
3/7,4/7\right)  }\left(  x^{\prime}\right) \\
&  =1-\rho_{\left(  5/7,6/7\right)  }\left(  x\right)  -\rho_{\left(
3/7,4/7\right)  }\left(  x^{\prime}\right) \\
&  +\rho_{\left(  5/7,6/7\right)  }\left(  x\right)  \rho_{\left(
3/7,4/7\right)  }\left(  x^{\prime}\right)  +\rho_{\left(  1/7,2/7\right)
}\left(  x\right)  \rho_{\left(  3/7,4/7\right)  }\left(  x^{\prime}\right) \\
&  =1-\rho_{\left(  5/7,6/7\right)  }\left(  x\right)  -\rho_{\left(
3/7,4/7\right)  }\left(  x^{\prime}\right)  +\rho_{\left(  5/7,6/7\right)
}\left(  x\right)  +\rho_{\left(  3/7,4/7\right)  }\left(  x^{\prime}\right)
=1
\end{align*}
Note that $d\left(  x,x^{\prime}\right)  <1/7\Rightarrow\rho_{\left(
5/7,6/7\right)  }\left(  x\right)  \rho_{\left(  3/7,4/7\right)  }\left(
x^{\prime}\right)  =\rho_{\left(  5/7,6/7\right)  }\left(  x\right)  ,$ etc..
Because of this, we expect that
\begin{align*}
&  \lim_{t\rightarrow0^{+}}Q\left(  t;x,x^{\prime}\right) \\
&  =\lim_{t\rightarrow0^{+}}\left(  \varphi_{c}\left(  x\right)
\mathcal{E}_{c}\left(  t;x,x^{\prime}\right)  \psi_{c}\left(  x^{\prime
}\right)  +\varphi_{d}\left(  x\right)  \mathcal{E}_{d}\left(  t;x,x^{\prime
}\right)  \psi_{d}\left(  x^{\prime}\right)  \right) \\
&  =\varphi_{c}\left(  x\right)  \mathcal{\delta}\left(  x,x^{\prime}\right)
\psi_{c}\left(  x^{\prime}\right)  +\varphi_{d}\left(  x\right)
\mathcal{\delta}\left(  x,x^{\prime}\right)  \psi_{d}\left(  x^{\prime}\right)
\\
&  =\left(  \varphi_{c}\left(  x\right)  \psi_{c}\left(  x^{\prime}\right)
+\varphi_{d}\left(  x\right)  \psi_{d}\left(  x^{\prime}\right)  \right)
\mathcal{\delta}\left(  x,x^{\prime}\right)  =\mathcal{\delta}\left(
x,x^{\prime}\right)  .
\end{align*}
Thus, on the operator level, we expect that $\lim_{t\rightarrow0^{+}}%
Q_{t}=\operatorname{I}$; i.e.,
\[
\lim_{t\rightarrow0^{+}}Q_{t}\left(  f\right)  \left(  x\right)
=\lim_{t\rightarrow0^{+}}\int_{X}Q\left(  t;x,x^{\prime}\right)  f\left(
x^{\prime}\right)  \,dx^{\prime}=\int_{X}\mathcal{\delta}\left(  x,x^{\prime
}\right)  f\left(  x^{\prime}\right)  \,dx^{\prime}=f\left(  x\right)  .
\]
The extent to which $Q\left(  \cdot;\cdot,x^{\prime}\right)  $ fails to
satisfy the heat equation in $\left[  \tfrac{1}{7},\tfrac{6}{7}\right]  \times
Y$ is given by
\[
C\left(  t;x,x^{\prime}\right)  :=-\left(  \partial_{t}+\mathcal{D}^{\ast
}\mathcal{D}\right)  \left.  Q\left(  \cdot;\cdot,x^{\prime}\right)  \right|
_{\left(  t;x,x^{\prime}\right)  }.
\]
For $x^{\prime}\in X$ fixed and $x=\left(  u,y\right)  \in N=\left[
0,1\right]  \times Y$, we have
\begin{align}
&  -C\left(  t;x,x^{\prime}\right) \label{CEst}\\
&  =\left(  \partial_{t}+\mathcal{D}^{\ast}\mathcal{D}\right)  Q\left(
t;x,x^{\prime}\right) \nonumber\\
&  =\left(  \partial_{t}+\mathcal{D}^{\ast}\mathcal{D}\right)  \left(
\sum\nolimits_{j\in\left\{  c,d\right\}  }\varphi_{j}\left(  x\right)
\mathcal{E}_{j}\left(  t;x,x^{\prime}\right)  \psi_{j}\left(  x^{\prime
}\right)  \right) \nonumber\\
&  =\sum_{j\in\left\{  c,d\right\}  }\partial_{u}^{2}\left(  \varphi
_{j}\left(  x\right)  \right)  \mathcal{E}_{j}\left(  t;x,x^{\prime}\right)
\psi_{j}\left(  x^{\prime}\right)  +2\partial_{u}\left(  \varphi_{j}\left(
x\right)  \right)  \partial_{u}\mathcal{E}_{j}\left(  t;x,x^{\prime}\right)
\psi_{j}\left(  x^{\prime}\right) \nonumber
\end{align}
For $d\left(  x,x^{\prime}\right)  <\tfrac{1}{7}$, we have $C\left(
t;x,x^{\prime}\right)  =0$, since (\ref{suppphipsi}) implies that $\psi
_{c}\left(  x^{\prime}\right)  =0$ when $\partial_{u}^{k}\left(  \varphi
_{j}\left(  x\right)  \right)  \neq0$ and $d\left(  x,x^{\prime}\right)
<\tfrac{1}{7}$. Using this and the estimates
\begin{align*}
\mathcal{E}_{j}\left(  t;x,x^{\prime}\right)   &  \leq At^{-n/2}e^{-Bd\left(
x,x^{\prime}\right)  ^{2}/t}\text{ and }\\
\partial_{u}\mathcal{E}_{j}\left(  t;x,x^{\prime}\right)   &  \leq
At^{-\left(  n+1\right)  /2}e^{-Bd\left(  x,x^{\prime}\right)  ^{2}/t}%
\end{align*}
for positive constants $A$ and $B$, it follows from (\ref{CEst}) that
\begin{align}
\left|  C\left(  t;x,x^{\prime}\right)  \right|   &  \leq c_{1}t^{-\left(
n+1\right)  /2}e^{-Bd\left(  x,x^{\prime}\right)  ^{2}/t}\nonumber\\
&  \leq c_{1}t^{-\left(  n+1\right)  /2}e^{-B7^{-2}/t}=c_{1}e^{-c_{2}%
/t},\text{ for }0<t<T_{0}<\infty, \label{CEst2}%
\end{align}
for positive constants $c_{1}$ and $c_{2}$.

Let $\mathcal{E}\left(  t;x,x^{\prime}\right)  $ be the exact heat kernel for
$\mathcal{D}^{\ast}\mathcal{D}$ and let the corresponding operator be
$\mathcal{E}\left(  t\right)  =\exp\left(  -t\mathcal{D}^{\ast}\mathcal{D}%
\right)  $. On the operator level, we have
\begin{align}
\mathcal{E}\left(  t\right)  -Q\left(  t\right)   &  =\mathcal{E}\left(
t\right)  Q\left(  0\right)  -\mathcal{E}\left(  0\right)  Q\left(  t\right)
\nonumber\\
&  =\int_{0}^{t}\tfrac{d}{ds}\left(  \mathcal{E}\left(  s\right)  Q\left(
t-s\right)  \right)  ds\nonumber\\
&  =\int_{0}^{t}\tfrac{d\mathcal{E}}{ds}Q\left(  t-s\right)  +\mathcal{E}%
\left(  s\right)  \tfrac{d}{ds}Q\left(  t-s\right)  \text{\thinspace
}ds\nonumber\\
&  =\int_{0}^{t}-\mathcal{D}^{\ast}\mathcal{DE}\left(  s\right)  Q\left(
t-s\right)  +\mathcal{E}\left(  s\right)  \tfrac{d}{ds}Q\left(  t-s\right)
\text{\thinspace}ds\nonumber\\
&  =\int_{0}^{t}\mathcal{E}\left(  s\right)  \left(  -\mathcal{D}^{\ast
}\mathcal{D}-\tfrac{d}{d\left(  t-s\right)  }\right)  Q\left(  t-s\right)
\text{\thinspace}ds\nonumber\\
&  =\int_{0}^{t}\mathcal{E}\left(  s\right)  C\left(  t-s\right)
\text{\thinspace}ds. \label{Duhamel}%
\end{align}
Hence on the kernel level,
\begin{align}
\mathcal{E}\left(  t;x,x^{\prime}\right)   &  =Q\left(  t;x,x^{\prime}\right)
+\int_{0}^{t}\int_{X}\mathcal{E}\left(  s;x,z\right)  C\left(  t-s;z,x^{\prime
}\right)  \text{\thinspace}dzds\nonumber\\
&  =Q\left(  t;x,x^{\prime}\right)  +\left(  \mathcal{E}\ast C\right)  \left(
t;x,x^{\prime}\right)  , \label{EQconvEqn}%
\end{align}
where the convolution operation is defined by
\[
\left(  \alpha\ast\beta\right)  \left(  t;x,x^{\prime}\right)  :=\int_{0}%
^{t}\int_{X}\alpha\left(  s;x,z\right)  \beta\left(  t-s;z,x^{\prime}\right)
dzds.
\]
Note that (\ref{EQconvEqn}) can be written as
\[
\mathcal{E}\left(  \operatorname{I}-\ast C\right)  =\mathcal{E}-\mathcal{E}%
\ast C=Q.
\]
Thus, at least formally,
\begin{align}
\mathcal{E}  &  =Q\left(  \operatorname{I}-\ast C\left(  t\right)  \right)
^{-1}=Q\left(  \operatorname{I}+\sum\nolimits_{k=1}^{\infty}\left(  \ast
C\right)  ^{k}\right) \nonumber\\
&  =Q+Q\ast\sum\nolimits_{k=1}^{\infty}C_{k}=Q+Q\ast\mathcal{C},\text{ where}
\label{LeviSum}%
\end{align}%
\[
C_{1}:=C=-\left(  \partial_{t}+\mathcal{D}^{\ast}\mathcal{D}\right)
Q,\;C_{k+1}:=C_{1}\ast C_{k}\text{, and\ }\mathcal{C}:=\sum\nolimits_{k=1}%
^{\infty}C_{k}\text{.}%
\]
The proof of the validity of the Levi sum (\ref{LeviSum}) for $\mathcal{E}$
rests on the estimate (\ref{CEst2}) and we refer to \cite[Ch. 22]{BoWo93} for
the details. The above constructions can be also be applied to obtain a kernel
$\mathcal{E}_{\ast}\left(  t;x,x^{\prime}\right)  $ for $\exp\left(
-t\mathcal{DD}^{\ast}\right)  $ starting from
\[
Q_{\ast}\left(  t;x,x^{\prime}\right)  :=\varphi_{c}\left(  x\right)
\mathcal{E}_{c\ast}\left(  t;x,x^{\prime}\right)  \psi_{c}\left(  x^{\prime
}\right)  +\varphi_{d}\left(  x\right)  \mathcal{E}_{d\ast}\left(
t;x,x^{\prime}\right)  \psi_{d}\left(  x^{\prime}\right)  .
\]
Recall that for $\left(  u,y\right)  \in\left[  0,\infty\right)  \times Y$,
\begin{align}
&  \mathcal{K}\left(  t;u,y\right)  :=\operatorname{Tr}\left(  \mathcal{E}%
_{c}\left(  t;\left(  u,y\right)  ,\left(  u,y\right)  \right)  -\mathcal{E}%
_{c\ast}\left(  t;\left(  u,y\right)  ,\left(  u,y\right)  \right)  \right)
,\text{ and}\nonumber\\
&  \mathcal{K}\left(  t\right)  :=\int_{0}^{\infty}\int_{Y}\mathcal{K}\left(
t;u,y\right)  \,dydu. \label{KtDefn}%
\end{align}
For $x\in\widetilde{X}$, let
\[
F\left(  t;x\right)  :=\operatorname{Tr}\left(  \mathcal{E}_{d}\left(
t;x,x\right)  -\mathcal{E}_{d\ast}\left(  t;x,x\right)  \right)  .
\]
Now,
\[
\operatorname{index}\mathcal{D}_{P_{\geq}}^{+}=\operatorname{Tr}\left(
e^{-t\mathcal{D}^{\ast}\mathcal{D}}-e^{-t\mathcal{DD}^{\ast}}\right)
=\operatorname{Tr}\left(  \mathcal{E}\left(  t\right)  -\mathcal{E}_{\ast
}\left(  t\right)  \right)  .
\]
We use the following notation for functions which are ``exponentially close''
as $t\rightarrow0^{+}:$
\begin{equation}
f\left(  t\right)  \sim_{\exp}g\left(  t\right)  \text{ }\Leftrightarrow
\left|  f\left(  t\right)  -g\left(  t\right)  \right|  \leq ae^{-b/t}\text{
for }t\in\left(  0,\varepsilon\right)  \label{ExpAsyDefn}%
\end{equation}
for some positive constants $a,b,\varepsilon$. We claim (but must omit
important details here -- see \cite[Ch. 22]{BoWo93}) that
\[
\operatorname{Tr}\left(  \mathcal{E}\left(  t\right)  -\mathcal{E}_{\ast
}\left(  t\right)  \right)  \sim_{\exp}\operatorname{Tr}\left(  Q\left(
t\right)  -Q_{\ast}\left(  t\right)  \right)  .
\]
Since $\mathcal{E}=Q+Q\ast\mathcal{C}$ (resp., $\mathcal{E}_{\ast}=Q_{\ast
}+Q_{\ast}\ast\mathcal{C}_{\ast}$), this result ultimately rests on
(\ref{CEst2}) (resp., the corresponding result for $C_{\ast}$). Since
$\varphi_{c}\psi_{c}=\psi_{c}$ and $\varphi_{d}\psi_{d}=\psi_{d}$,
\begin{align*}
&  \operatorname{Tr}\left(  Q\left(  t\right)  -Q_{\ast}\left(  t\right)
\right)  =\int_{X}\operatorname{Tr}Q\left(  t;x,x\right)  \,dx\\
&  =\int_{X}\operatorname{Tr}\left(
\begin{array}
[c]{c}%
\varphi_{c}\left(  x\right)  \left(  \mathcal{E}_{c}-\mathcal{E}_{c\ast
}\right)  \left(  t;x,x\right)  \psi_{c}\left(  x\right) \\
+\varphi_{d}\left(  x\right)  \left(  \mathcal{E}_{d}-\mathcal{E}_{c\ast
}\right)  \left(  t;x,x\right)  \psi_{d}\left(  x\right)
\end{array}
\right)  \,dx\\
&  =\int_{0}^{1}\int_{Y}\mathcal{K}\left(  t;u,y\right)  \psi_{c}\left(
u\right)  \,dydu+\int_{X}F\left(  t;x\right)  \psi_{d}\left(  x\right)  \,dx.
\end{align*}
Let
\[
\mathcal{K}_{a}\left(  t\right)  :=\int_{0}^{a}\int_{Y}\mathcal{K}\left(
t;u,y\right)  dydu.
\]
Note that
\begin{align*}
\mathcal{K}\left(  t\right)  -\mathcal{K}_{a}\left(  t\right)   &  =\int
_{a}^{\infty}\int_{Y}\mathcal{K}\left(  t;u,y\right)  dydu\\
&  =\sum\nolimits_{\lambda}m_{\lambda}\operatorname*{sign}\left(
\lambda\right)  \left.  \tfrac{1}{2}e^{2\left|  \lambda\right|  u}%
\operatorname{erfc}\left(  \tfrac{u}{\sqrt{t}}+\left|  \lambda\right|
\sqrt{t}\right)  \right|  _{u=a}^{\infty}\\
&  =-\tfrac{1}{2}\sum\nolimits_{\lambda}m_{\lambda}\operatorname*{sign}\left(
\lambda\right)  e^{2\left|  \lambda\right|  a}\operatorname{erfc}\left(
\tfrac{a}{\sqrt{t}}+\left|  \lambda\right|  \sqrt{t}\right)  .
\end{align*}
Thus (using (\ref{ErfcEst})) for some constant $c>0$%
\begin{align*}
\left|  \mathcal{K}\left(  t\right)  -\mathcal{K}_{a}\left(  t\right)
\right|   &  \leq\sum\nolimits_{\lambda}m_{\lambda}e^{2\left|  \lambda\right|
a}\exp\left(  -\left(  \tfrac{a}{\sqrt{t}}+\left|  \lambda\right|  \sqrt
{t}\right)  ^{2}\right) \\
&  \leq\tfrac{1}{2}e^{-a^{2}/t}\sum\nolimits_{\lambda}m_{\lambda}e^{-\left|
\lambda\right|  ^{2}t}\sim_{t\rightarrow0^{+}}ce^{-a^{2}/t}t^{-n/2},
\end{align*}
which says that $\mathcal{K}\left(  t\right)  \sim_{\exp}\mathcal{K}%
_{a}\left(  t\right)  $. Hence,
\begin{align}
&  \operatorname{index}\mathcal{D}_{P_{\geq}}^{+}=\operatorname{Tr}\left(
\mathcal{E}\left(  t\right)  -\mathcal{E}_{\ast}\left(  t\right)  \right)
\sim_{\exp}\operatorname{Tr}\left(  Q\left(  t\right)  -Q_{\ast}\left(
t\right)  \right) \nonumber\\
&  =\int_{0}^{1}\int_{Y}\mathcal{K}\left(  t;u,y\right)  \psi_{c}\left(
u\right)  \,dydu+\int_{X}F\left(  t;x\right)  \psi_{d}\left(  x\right)
\,dx\nonumber\\
&  \sim_{\exp}\left(  \int_{0}^{\infty}\int_{Y}\mathcal{K}\left(
t;u,y\right)  \,dydu+\int_{X}F\left(  t;x\right)  \psi_{d}\left(  x\right)
\,dx\right) \nonumber\\
&  =\mathcal{K}\left(  t\right)  +\int_{X}F\left(  t;x\right)  \,dx.
\label{IndKF}%
\end{align}
The last equality follows from (\ref{KtDefn}) and the fact that
\begin{equation}
F\left(  t;x\right)  =\operatorname{Tr}\left(  \mathcal{E}_{d}\left(
t;x,x\right)  -\mathcal{E}_{d\ast}\left(  t;x,x\right)  \right)  =0\text{ for
}x\in N=\left[  0,1\right]  \times Y, \label{FNeq0}%
\end{equation}
which is seen as follows. Since $\sigma^{\ast}=-\sigma$,
\begin{align*}
\widetilde{\mathcal{D}^{+}}|_{{N}}  &  =\sigma(\partial_{u}+\mathcal{B}%
):C^{\infty}\left(  N,\widetilde{S^{+}}\right)  \rightarrow C^{\infty}\left(
N,\widetilde{S^{-}}\right)  \text{ implies}\\
\widetilde{\mathcal{D}^{-}}|_{{N}}  &  =\left(  \sigma(\partial_{u}%
+\mathcal{B})\right)  ^{\ast}=(\partial_{u}-\mathcal{B})\sigma:C^{\infty
}\left(  N,\widetilde{S^{-}}\right)  \rightarrow C^{\infty}\left(
N,\widetilde{S^{+}}\right)  ,
\end{align*}
Thus, over $N$,
\begin{align*}
\widetilde{\mathcal{D}^{-}}\widetilde{\mathcal{D}^{+}}  &  =(\partial
_{u}-\mathcal{B})\sigma\sigma(\partial_{u}+\mathcal{B})=\left(  -\partial
_{u}^{2}+\mathcal{B}^{2}\right)  \left|  C^{\infty}\left(  N,\widetilde{S^{+}%
}\right)  \right.  \text{ and}\\
\widetilde{\mathcal{D}^{+}}\widetilde{\mathcal{D}^{-}}  &  =\sigma
(\partial_{u}+\mathcal{B})(\partial_{u}-\mathcal{B})\sigma=\sigma^{2}%
(\partial_{u}-\mathcal{B})(\partial_{u}+\mathcal{B})\\
&  =\left(  -\partial_{u}^{2}+\mathcal{B}^{2}\right)  \left|  C^{\infty
}\left(  N,\widetilde{S^{-}}\right)  .\right.
\end{align*}
Since $\left[  \left(  -\partial_{u}^{2}+\mathcal{B}^{2}\right)
,\sigma\right]  =0$, $\sigma$ maps the eigenspaces of $\widetilde
{\mathcal{D}^{-}}\widetilde{\mathcal{D}^{+}}$ pointwise isometrically (over
$N$) onto those of $\widetilde{\mathcal{D}^{+}}\widetilde{\mathcal{D}^{-}}$.
Hence, we obtain (\ref{FNeq0}) upon taking the trace of the eigensection
expansion of $\mathcal{E}_{d}\left(  t;x,x\right)  -\mathcal{E}_{d\ast}\left(
t;x,x\right)  $ for $x\in N$. From (\ref{IndKF}), we obtain
\begin{align}
\mathcal{K}\left(  t\right)   &  =\operatorname{index}\mathcal{D}_{P_{\geq}%
}^{+}-\int_{X}F\left(  t;x\right)  \,dx\nonumber\\
&  \sim\operatorname{index}\mathcal{D}_{P_{\geq}}^{+}-\int_{X}\sum
\nolimits_{k\geq-n}\alpha_{k}\left(  x\right)  t^{k/2}\,dx\nonumber\\
&  =\operatorname{index}\mathcal{D}_{P_{\geq}}^{+}-\sum\nolimits_{k\geq
-n}A_{k}t^{k/2}\text{ }(\text{for }A_{k}:=\int_{X}\alpha_{k}\left(  x\right)
\,dx), \label{KasyExp}%
\end{align}
where $F\left(  t;x\right)  \sim\sum\nolimits_{k\geq-n}\alpha_{k}\left(
x\right)  t^{k/2}$ is the asymptotic expansion of the trace $F\left(
t;x\right)  $ of the kernel for $e^{-t\widetilde{\mathcal{D}^{-}}%
\widetilde{\mathcal{D}^{+}}}-e^{-t\widetilde{\mathcal{D}^{+}}\widetilde
{\mathcal{D}^{-}}}$, which is known for elliptic differential operators (e.g.,
$\widetilde{\mathcal{D}^{-}}\widetilde{\mathcal{D}^{+}}$ and $\widetilde
{\mathcal{D}^{+}}\widetilde{\mathcal{D}^{-}}$) over closed manifolds such as
$\widetilde{X}$ (see \cite[p. 68]{Gi95}). Thus, according to (\ref{EtaBAsy})
with $a_{k}=-A_{k}$ for $k\neq0$ and $a_{0}=\operatorname{index}%
\mathcal{D}_{P_{\geq}}^{+}-A_{0}$,
\[
\eta_{\mathcal{B}}\left(  2s\right)  =-\frac{2s\sqrt{\pi}}{\Gamma\left(
s+\tfrac{1}{2}\right)  }\left(  \frac{m_{0}}{2s}+\sum\nolimits_{k=-n+1}%
^{N}\frac{a_{k}}{\tfrac{1}{2}k+s}+\theta_{N}\left(  s\right)  \right)
\]%
\[
\eta_{\mathcal{B}}\left(  2s\right)  =-\frac{2s\sqrt{\pi}}{\Gamma\left(
s+\tfrac{1}{2}\right)  }\left(  \frac{m_{0}+2\operatorname{index}%
\mathcal{D}_{P_{\geq}}^{+}}{2s}-\sum\nolimits_{k=-n+1}^{N}\frac{A_{k}}%
{\tfrac{1}{2}k+s}+\theta_{N}\left(  s\right)  \right)  .
\]
Setting $s=0$, yields $\eta_{\mathcal{B}}\left(  0\right)  =-\left(
m_{0}+2\operatorname{index}\mathcal{D}_{P_{\geq}}^{+}-2A_{0}\right)  $ or
\[
\operatorname{index}\mathcal{D}_{P_{\geq}}^{+}=A_{0}-\tfrac{1}{2}\left(
m_{0}+\eta_{\mathcal{B}}\left(  0\right)  \right)  .
\]
\end{proof}

\subsection{Symplectic Geometry of Cauchy Data
Spaces\label{ss:SymplecticGeometryOfCauchyDataSpaces}}

As we have seen in our Section \ref{sss:Symmetric Operators} there exists a
concept of Cauchy data spaces which solely is based on the concepts of minimal
and maximal domain and which is more elementary and more general than the
definitions provided in Section \ref{sss:PoissonOperator} which are based on
pseudo-differential calculus.

In this Section we stay in the real category and do not assume product
structure near $\Sigma=\partial X$ unless otherwise stated. The operator
$\mathcal{D}$ need not be of Dirac type. We only assume that it is a linear,
elliptic, symmetric, differential operator of first order.

\medskip

\subsubsection{The Natural Cauchy Data Space}

Let $\mathcal{D}_{0}$ denote the restriction of $\mathcal{D}$ to the space
$C_{0}^{\infty}(X;S)$ of smooth sections with support in the interior of $X$.
As mentioned above, there is no natural choice of the order of the Sobolev
spaces for the boundary reduction. Therefore, a systematic treatment of the
boundary reduction may begin with the minimal closed extension {$\mathcal{D}$%
}$_{\min}:=\overline{{\mathcal{D}}_{0}}$ and the adjoint {$\mathcal{D}$%
}$_{\max}:=(${$\mathcal{D}$}$_{0})^{\ast}$ of {$\mathcal{D}$}$_{0}$. Clearly,
{$\mathcal{D}$}$_{\max}$ is the maximal closed extension. This gives
\[
D_{\min}:=\operatorname{Dom}({\mathcal{D}}_{\min})=\overline{C_{0}^{\infty
}(X;S)}^{\mathcal{G}}=\overline{C_{0}^{\infty}(X;S)}^{H^{1}(X;S)}%
\]
and
\begin{multline*}
D_{\max}:=\operatorname{Dom}({\mathcal{D}}_{\max})\\
=\{u\in L^{2}(X;S)\mid\mathcal{D}u\in L^{2}(X;S)\text{ in the sense of
distributions}\}.
\end{multline*}
Here, the superscript $\mathcal{G}$ means the closure in the graph norm which
coincides with the first Sobolev norm on $C_{0}^{\infty}(X;S)$. We form the
space $\mathbf{\beta}$ of \emph{natural boundary values} with the
\emph{natural trace map} $\gamma$ as described in Section \ref{ss:Symmetric
Operators}.

There we defined also the \emph{natural Cauchy data space} $\Lambda
(\mathcal{D}):=\gamma(\operatorname{Ker}\mathcal{D}_{\max})$ under the
assumption that $\mathcal{D}$ has a self-adjoint $L^{2}$ extension with a
compact resolvent. Such an extension always exists. Take for instance
{$\mathcal{D}$}$_{\mathcal{P}({\mathcal{D}})}$, the operator {$\mathcal{D}$}
with domain
\[
\operatorname{Dom}_{\mathcal{P}({\mathcal{D}})}:=\{f\in H^{1}(X;S)\mid
\mathcal{P}(\mathcal{D})^{(\frac{1}{2})}(f|_{\Sigma})=0\},
\]
where $\mathcal{P}(\mathcal{D})$ denotes the Calder{\'{o}}n projection defined
in (\ref{e:calderon}).

Clearly, $D_{\max}$ and $D_{\min}$ are $C^{\infty}(X)$ modules, and so the
space $\mathbf{\beta}$ is a $C^{\infty}(\Sigma)$ module. This shows that
$\mathbf{\beta}$ is \textit{local} in the following sense: If $\Sigma$
decomposes into $r$ connected components $\Sigma=\Sigma_{1}\sqcup\dots
\sqcup\Sigma_{r}\,,$ then $\mathbf{\beta}$ decomposes into
\[
\mathbf{\beta}=\bigoplus_{j=1}^{r}\mathbf{\beta}_{j},
\]
where
\[
\mathbf{\beta}_{j}:=\gamma\left(  \{f\in D_{\max}\mid\operatorname*{supp}%
f\subset{N}_{j}\}\right)
\]
with a suitable collared neighborhood ${N}_{j}$ of $\Sigma_{j}$\thinspace.
Note that each $\mathbf{\beta}_{j}$ is a closed symplectic subspace of
$\mathbf{\beta}$ and therefore a symplectic Hilbert space.

\noindent By Theorem \ref{t:res}a and, alternatively and in greater
generality, by H\"{o}rmander \cite{Ho66} (Theorem 2.2.1 and the Estimate
(2.2.8), p. 194), the space $\mathbf{\beta}$ is naturally embedded in the
distribution space $H^{-\frac{1}{2}}(\Sigma;S|_{\Sigma})$. Under this
embedding we have $\Lambda(\mathcal{D})=\Lambda(\mathcal{D},0)$, where the
last space was defined in Definition \ref{d:cauchy-data}.

If the metrics are product close to $\Sigma$, we can give a more precise
description of the embedding of $\mathbf{\beta}$, namely as a \textit{graded}
space of distributions. Let $\{\varphi_{k},\Lambda_{k}\}$ be a spectral
resolution of $L^{2}(\Sigma)$ by eigensections of $\mathcal{B}$. (Here and in
the following we do not mention the bundle $S$). Once again, for simplicity,
we assume $\operatorname{Ker}\mathcal{B}=\{0\}$ and have $\mathcal{B}%
\varphi_{k}=\lambda_{k}\varphi_{k}$ for all $k\in\mathbb{Z}\setminus\{0\}$,
and $\lambda_{-k}=-\lambda_{k}\,$, $\sigma(\varphi_{k})=\varphi_{-k}\,$, and
$\sigma(\varphi_{-k})=-\varphi_{k}$ for $k>0$. In \cite{BoFu99}, Proposition
7.15 (see also \cite{BrLe99} for a more general setting) it was shown that
\begin{align}
\mathbf{\beta}  &  =\mathbf{\beta}_{-}\oplus\mathbf{\beta}_{+}\quad\text{
with}\nonumber\\
\mathbf{\beta}_{-}  &  :=\overline{[\{\varphi_{k}\}_{k<0}]}^{H^{\frac{1}{2}%
}(\Sigma)}\text{\ and\ }\mathbf{\beta}_{+}:=\overline{[\{\varphi_{k}\}_{k>0}%
]}^{H^{-\frac{1}{2}}(\Sigma)}\,. \label{e:beta-split}%
\end{align}
Then $\mathbf{\beta}_{-}$ and $\mathbf{\beta}_{+}$ are Lagrangian and
transversal subspaces of $\mathbf{\beta}$.

\subsubsection{Criss--cross Reduction}

Let us define two Lagrangian and transversal subspaces $L_{\pm}$ of
$L^{2}(\Sigma)$ in a similar way, namely by the closure in $L^{2}(\Sigma)$ of
the linear span of the eigensections with negative, resp. with positive
eigenvalue. Then $L_{+}$ is dense in $\mathbf{\beta}_{+}$\thinspace, and
$\mathbf{\beta}_{-}$ is dense in $L_{-}$. This anti--symmetric relation may
explain some of the well--observed delicacies of dealing with spectral
invariants of continuous families of Dirac operators.

Moreover, $\gamma(D_{\operatorname{aps}})=\mathbf{\beta}_{-}$\thinspace,
where
\begin{equation}
D_{\operatorname{aps}}:=\{f\in H^{1}(X)\mid P_{>}(f|_{\Sigma})=0\}
\label{e:def_aps}%
\end{equation}
denotes the domain corresponding to the Atiyah--Patodi--Singer boundary
condition. Note that a series $\sum_{k<0}c_{k}\varphi_{k}$ may converge to an
element $\varphi\in L^{2}(\Sigma)$ without converging in $H^{\frac{1}{2}%
}(\Sigma)$. Therefore such $\varphi\in L_{-}$ can not appear as trace at the
boundary of any $f\in D_{\max}$.

Recall Proposition \ref{p:lagr} and note that $(\mathbf{\beta}_{\mathbf{-}%
},\Lambda(\mathcal{D}))$ is a Fredholm pair.

This can all be achieved without the symbolic calculus of pseudo-differential
operators. Therefore one may ask how the preceding approach to Cauchy data
spaces and boundary value problems via the maximal domain and our symplectic
space $\mathbf{\beta}$ is related to the approach via the Calder{\'{o}}n
projection, which we reviewed in the preceding section. How can results from
the distributional theory be translated into $L^{2}$ results?%

\begin{figure}
[h]
\begin{center}
\includegraphics[
height=2.5408in,
width=3.1055in
]%
{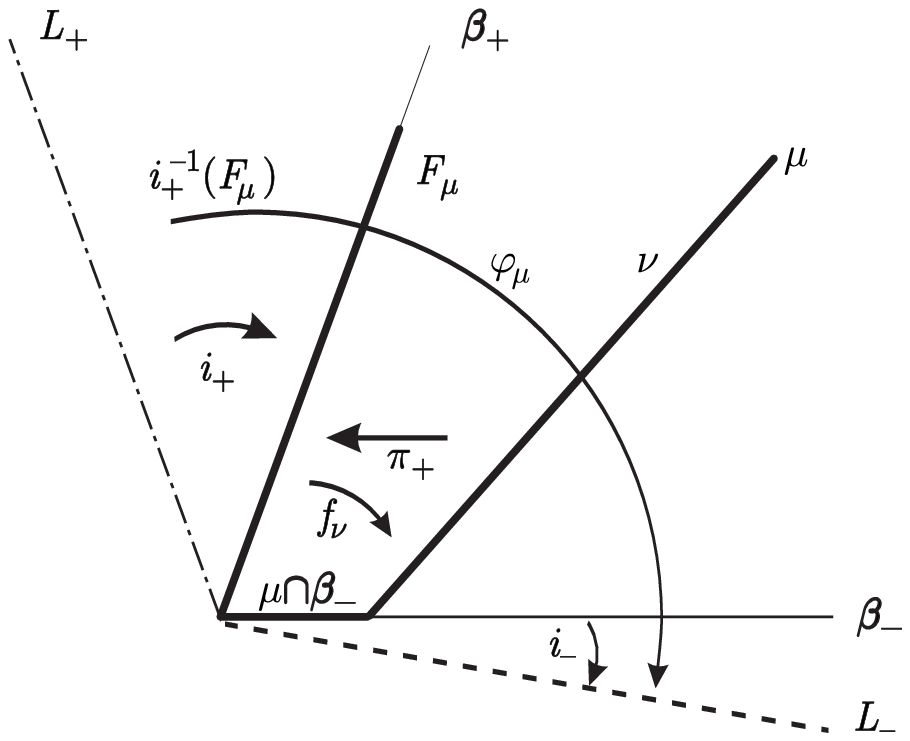}%
\caption{The mapping $\tau:\mathcal{FL}_{\mathbf{\beta}_{-}}\left(
\mathbf{\beta}\right)  \rightarrow\mathcal{FL}_{L_{-}}\left(  L\right)  $}%
\label{f:criss-cross}%
\end{center}
\end{figure}
To relate the two approaches we recall a fairly general symplectic
`Criss--Cross' Reduction Theorem from \cite{BoFuOt01} (Theorem 1.2). Let
$\mathbf{\beta}$ and $L$ be symplectic Hilbert spaces with symplectic forms
$\omega_{\beta}$ and $\omega_{L}$\thinspace, respectively. Let
\[
\mathbf{\beta}=\mathbf{\beta}_{-}\dotplus\mathbf{\beta}_{+}\qquad
\text{and}\qquad L=L_{-}\dotplus L_{+}%
\]
be direct sum decompositions by transversal (not necessarily orthogonal) pairs
of Lagrangian subspaces. We assume that there exist continuous, injective
mappings
\[
i_{-}:\mathbf{\beta}_{-}\longrightarrow L_{-}\qquad\text{and}\qquad
i_{+}:L_{+}\longrightarrow\mathbf{\beta}_{+}%
\]
with dense images and which are compatible with the symplectic structures,
i.e.
\[
\omega_{L}(i_{-}(x),a)=\omega_{\mathbf{\beta}}(x,i_{+}(a))\quad\text{ for all
$a\in L_{+}$ and $x\in\mathbf{\beta}_{-}$}\,.
\]
Let $\mu\in\mathcal{FL}_{\mathbf{\beta}_{-}}(\mathbf{\beta})$, e.g. $\mu
=(\mu\cap\mathbf{\beta}_{-})\dotplus\nu$ with a suitable closed $\nu$. Let us
define (see also Figure \ref{f:criss-cross})
\begin{equation}
\tau(\mu):=i_{-}(\mu\cap\mathbf{\beta}_{-})+\operatorname{graph}(\varphi_{\mu
})\,, \label{e:tau-operator}%
\end{equation}
where
\[%
\begin{matrix}
\varphi_{\mu}: & i_{+}^{-1}(F_{\mu}) & \longrightarrow &  L_{-}\\
\  & x & \mapsto &  i_{-}\circ f_{\nu}\circ i_{+}(x)
\end{matrix}
\qquad.
\]
Here $F_{\mu}$ denotes the image of $\mu$ under the projection $\pi_{+}$ from
$\mu$ to $\mathbf{\beta}_{+}$ along $\mathbf{\beta}_{-}$ and $f_{\nu}:F_{\mu
}\longrightarrow\mathbf{\beta}_{-}$ denotes the uniquely determined bounded
operator which yields $\nu$ as its graph. Then:

\begin{theorem}
\label{t:criss-cross}The mapping (\ref{e:tau-operator}) defines a continuous
mapping
\[
\tau:\mathcal{FL}_{\mathbf{\beta}_{-}}(\mathbf{\beta})\longrightarrow
\mathcal{FL}_{\mathbf{L}_{-}}(L)
\]
which maps the Maslov cycle $\mathcal{M}_{\mathbf{\beta}_{-}}(\mathbf{\beta})
$ of $\mathbf{\beta}_{-}$ into the Maslov cycle $\mathcal{M}_{L_{-}}(L)$ of
$L_{-}$ and preserves the Maslov index
\[
\mathbf{mas}(\{\mu_{s}\}_{s\in\lbrack0,1]},\mathbf{\beta}_{-})=\mathbf{mas}%
(\{\tau(\mu_{s})\}_{s\in\lbrack0,1]},L_{-})
\]
for any continuous curve $[0,1]\ni s\mapsto\mu_{s}\in\mathcal{FL}%
_{\mathbf{\beta}_{-}}(\mathbf{\beta})$.
\end{theorem}

In the product case, the `Criss--Cross' Reduction Theorem implies for our two
types of Cauchy data that all results proved in the theory of natural boundary
values ($\mathbf{\beta}$ theory) remain valid in the $L^{2}$ theory. In
particular we have:

\begin{corollary}
The $L^{2}(\Sigma)$ part $\Lambda(\mathcal{D})\cap L^{2}(\Sigma)$ of the
natural Cauchy data space $\Lambda(\mathcal{D})$ is closed in $L^{2}(\Sigma)$.
Actually, it is a Lagrangian subspace of $L^{2}(\Sigma)$ and it forms a
Fredholm pair with the component $L_{-}$\thinspace, defined at the beginning
of this subsection.
\end{corollary}

\subsection{Non-Additivity of the Index\label{ss:non_add_index}}

\medskip

\subsubsection{The Bojarski Conjecture}

The \emph{Bojarski Conjecture} gives quite a different description of the
index of an elliptic operator over a closed partitioned manifold $M=M_{1}%
\cup_{\Sigma}M_{2}$\thinspace. It relates the `quantum' quantity index with a
`classical' quantity, the Fredholm intersection index of the Cauchy data
spaces from both sides of the hypersurface $\Sigma$. It was suggested in
\cite{Bo79} and proved in \cite{BoWo93} for operators of Dirac type.

\begin{proposition}
\label{p:boj} Let $M$ be a partitioned manifold as before and let
$\Lambda(\mathcal{D}_{j}^{+},\tfrac{1}{2})$ denote the $L^{2}$ Cauchy data
spaces, $j=1,2$ (see Definition \ref{d:cauchy-data}). Then
\[
\operatorname{index}\mathcal{D}^{+}=\operatorname{index}\left(  \Lambda\left(
\mathcal{D}_{1}^{+},\tfrac{1}{2}\right)  ,\Lambda\left(  \mathcal{D}_{2}%
^{+},\tfrac{1}{2}\right)  \right)  .
\]
\end{proposition}

Recall that
\begin{align*}
\operatorname{index}\left(  \Lambda\left(  \mathcal{D}_{1}^{+},\tfrac{1}%
{2}\right)  ,\Lambda\left(  \mathcal{D}_{2}^{+},\tfrac{1}{2}\right)  \right)
&  :=\dim\left(  \Lambda\left(  \mathcal{D}_{1}^{+},\tfrac{1}{2}\right)
\cap\Lambda\left(  \mathcal{D}_{2}^{+},\tfrac{1}{2}\right)  \right) \\
&  -\dim\left(  \frac{L^{2}\left(  \Sigma;S|_{\Sigma}\right)  }{\Lambda\left(
\mathcal{D}_{1}^{+},\tfrac{1}{2}\right)  +\Lambda\left(  \mathcal{D}_{2}%
^{+},\tfrac{1}{2}\right)  }\right)  .
\end{align*}
It is equal to $\mathbf{i}(\operatorname{I}-\mathcal{P}(\mathcal{D}_{2}%
^{+}),\mathcal{P}(\mathcal{D}_{1}^{+}))$, where $\mathcal{P}(\mathcal{D}%
_{j}^{+})$ denotes the corresponding Calder\'{o}n projections, and
$\mathbf{i}\left(  \cdot,\cdot\right)  $ was defined in (\ref{d:iP1P2}).

The proof of the Proposition depends on the unique continuation property for
Dirac operators and the Lagrangian property of the Cauchy data spaces, more
precisely the chiral twisting property \eqref{e:perp} .

\medskip

\subsubsection{Generalizations for Global Boundary Conditions}

On a smooth compact manifold $X$ with boundary $\Sigma$, the solution spaces
$\operatorname{Ker}(\mathcal{D},s)$ depend on the order $s$ of
differentiability and they are infinite-dimensional. To obtain a finite index
one must apply suitable boundary conditions (see \cite{BoWo93} for local and
global boundary conditions for operators of Dirac type). In this report, we
restrict ourselves to boundary conditions of Atiyah--Patodi--Singer type
(i.e., $P\in\mathcal{G}\mathrm{r}(\mathcal{D})$), and consider the extension
\begin{equation}
\mathcal{D}_{P}:\operatorname{Dom}(\mathcal{D}_{P})\longrightarrow L^{2}(X;S)
\end{equation}
of $\mathcal{D}$\ defined by the domain
\begin{equation}
\operatorname{Dom}(\mathcal{D}_{P}):=\{f\in H^{1}(X;S)\mid P^{(0)}(f|_{\Sigma
})=0\}.
\end{equation}
It is a closed operator in $L^{2}(X;S)$ with finite-dimensional kernel and
cokernel. We have an explicit formula for the adjoint operator
\begin{equation}
(\mathcal{D}_{P})^{\ast}=\mathcal{D}_{\sigma(\operatorname{I}-P)\sigma^{\ast}%
}\,. \label{e:adjoint-domain}%
\end{equation}
In agreement with Proposition \ref{p:lagr}b, the preceding equation shows that
an extension $\mathcal{D}_{P}$ is self-adjoint, if and only if
$\operatorname{Ker}P^{(0)}$ is a Lagrangian subspace of the symplectic Hilbert
space $L^{2}(\Sigma;S|_{\Sigma})$.

Let us recall the \emph{Boundary Reduction Formula for the Index} of (global)
elliptic boundary value problems discussed in \cite{BoWo85} (inspired by
\cite{Se69}, see \cite{BoWo93} for a detailed proof for Dirac operators). Like
the Bojarski Conjecture, the point of the formula is that it gives an
expression for the index in terms of the geometry of the Cauchy data in the
symplectic space of all (here $L^{2}$) boundary data.

\begin{proposition}
\label{p:seeley-brf}
\[
\operatorname{index}\mathcal{D}_{P}=\operatorname{index}\left\{
P\mathcal{P}\left(  \mathcal{D}\right)  :\Lambda(\mathcal{D},\tfrac{1}%
{2})\rightarrow\operatorname{range}(P^{\left(  0\right)  })\right\}
\]
\end{proposition}

\medskip

\subsubsection{Pasting Formulas}

We shall close this section by mentioning a slight modification of the
Bojarski Conjecture/Theorem, namely a non-additivity formula for the splitting
of the index over partitioned manifolds.

From the Atiyah--Singer Index Theorem (here the expression of the index on the
closed manifold $M$ by an integral of the \textit{index density}) and the
Atiyah--Patodi--Singer Index Theorem (Theorem \ref{tAPS}), we obtain at once
\[
\operatorname{index}\mathcal{D}=\operatorname{index}\left(  \mathcal{D}%
_{1}\right)  _{P_{\leq}}+\operatorname{index}\left(  \mathcal{D}_{2}\right)
_{P_{\geq}}-\dim\operatorname{Ker}(\mathcal{B}).
\]
Then an Agranovic--Dynin type correction formula (based on Proposition
\ref{p:seeley-brf}) yields:

\begin{theorem}
\label{t:non-add} Let $P_{j}$ be projections belonging to $\mathcal{G}%
\mathrm{r}(\mathcal{D}_{j})$, $j=1,2$. Then
\[
\operatorname{index}\mathcal{D}=\operatorname{index}\left(  \mathcal{D}%
_{1}\right)  _{P_{1}}+\operatorname{index}\left(  \mathcal{D}_{2}\right)
_{P_{2}}-\mathbf{i}\left(  P_{2},\operatorname{I}-P_{1}\right)
\]
\end{theorem}

It turns out that the boundary correction term $\mathbf{i}(P_{2}%
,\operatorname{I}-P_{1})$ equals the index of the operator $\sigma
(\partial_{t}+\mathcal{B}) $ on the cylinder $[0,1]\times\Sigma$ with the
boundary conditions $P_{1} $ at $t=0$ and $P_{2}$ at $t=1$. A direct proof of
Theorem \ref{t:non-add} can be derived from Proposition \ref{p:seeley-brf} by
elementary operations with the virtual indices $\mathbf{i}\left(
P_{1},\mathcal{P}\left(  \mathcal{D}_{1}\right)  \right)  $ and $\mathbf{i}%
\left(  P_{2},\mathcal{P}\left(  \mathcal{D}_{2}\right)  \right)  $; see,
e.g., \cite{DaZh96} in a more general setting.

\begin{remark}
\label{r:index} (\emph{a}) In this section we have not always distinguished
between the \textit{total} and the \textit{chiral} Dirac operator because all
the discussed index formulas are valid in both cases.\newline \noindent
(\emph{b}) Important index formulas for (global) elliptic boundary value
problems for operators of Dirac type can also be obtained without analyzing
the concept and the geometry of the Cauchy data spaces (see e.g. the
celebrated Atiyah--Patodi--Singer Index Theorem \ref{tAPS} or \cite{Sc01} for
a recent survey of index formulas where there is no mention of the
Calder{\'{o}}n projection). The basic reason is that the index is an invariant
represented by a local density inside the manifold plus a correction term
which lives on the boundary and may be local or non-local as well. However,
these formulas do not explain the simple origin of the index or the spectral
flow, namely that \emph{all} index information is naturally coded by the
geometry of the Cauchy data spaces. To us it seems necessary to use the
Calder{\'{o}}n projection in order to understand (not calculate) the index of
an elliptic boundary problem and the reason for the locality or non-locality.
\end{remark}

\subsection{Pasting of Spectral Flow}

\smallskip

\subsubsection{Spectral Flow and the Maslov Index}

Let $\{\mathcal{D}_{t}\}_{t\in\lbrack0,1]}$ be a continuous family of (from
now on always \textit{total}) Dirac operators with the same principal symbol
and the same domain $D$. To begin with, we do not distinguish between the case
of a closed manifold (when $D$ is just the first Sobolev space and all
operators are essentially self-adjoint) and the case of a manifold with
boundary (when $D$ is specified by the choice of a suitable boundary value condition).

We consider the \emph{spectral flow} $\operatorname{SF}\left\{  \mathcal{D}%
_{t,D}\right\}  $ (see Section \ref{sss:self-adjoint fred ops}). We want a
pasting formula for the spectral flow. To achieve that, one replaces the
spectral flow of a continuous one-parameter family of self-adjoint Fredholm
operators, which is a `quantum' quantity, by the \textit{Maslov index} of a
corresponding path of Lagrangian Fredholm pairs, which is a `quasi-classical'
quantity. The idea is due to Floer and was worked out subsequently by Yoshida
in dimension 3, by Nicolaescu in all odd dimensions, and pushed further by
Cappell, Lee and Miller, Daniel and Kirk, and many other authors. For a
survey, see \cite{BoFu98}, \cite{BoFu99}, \cite{DaKi99}, \cite{KiLe00}.

In this review we give two spectral flow formulas of that type. To begin with,
we consider the case of a manifold with boundary. Then weak UCP for Dirac type
operators (established in Section \ref{WeakUCP}) implies weak inner UCP in the
sense of (\ref{e:1.9-new}), if the manifold is connected. Actually, it would
suffice that there is no connected component without boundary. Let us fix the
space $\mathbf{\beta}$ for the family. By Proposition \ref{p:lagr}c the
corresponding family $\{\Lambda(\mathcal{D}_{t})\}$ of natural Cauchy data
spaces is continuous. Applying the \textit{General Boundary Reduction Formula}
for the spectral flow (Theorem \ref{t:sff-old}) gives a family version of the
Bojarski conjecture (our Proposition \ref{p:boj}):

\begin{theorem}
\label{t:sff-mwb} The spectral flow of the family $\{\mathcal{D}_{t,D}\}$ and
the Maslov index\quad$\mathbf{mas}\left(  \{\Lambda(\mathcal{D}_{t}%
)\},\gamma(D)\right)  $ are well-defined and we have
\begin{equation}
\operatorname{SF}\left\{  \mathcal{D}_{t,D}\right\}  =\mathbf{mas}\left(
\{\Lambda(\mathcal{D}_{t})\},\gamma(D)\right)  {.} \label{e:sff-mwb}%
\end{equation}
\end{theorem}

We have various corollaries for the spectral flow on closed manifolds with
fixed hypersurface (see \cite{BoFu99}). Note that a partitioned manifold
$M=M_{1}\cup_{\Sigma}M_{2}$ can be considered as a manifold $M^{\#}%
=M_{1}\sqcup M_{2}$ with boundary $\partial M^{\#}=-\Sigma\sqcup\Sigma$. Here
$\sqcup$ denotes taking the disjoint union. Then any elliptic operator
$\mathcal{D}$ over $M$ defines an operator $\mathcal{D}^{\#}$ over $M^{\#}$
with natural Fredholm extension $\mathcal{D}_{D}^{\#}$ by fixing the domain
\[
D:=\{(f_{1},f_{2})\in H^{1}(M_{1})\times H^{1}(M_{2})\mid(f_{1}|_{\Sigma
},f_{2}|_{\Sigma})\in\Delta\},
\]
where $\Delta$ denotes the diagonal of $L^{2}\left(  -\Sigma;S|_{\Sigma
}\right)  \times L^{2}\left(  \Sigma;S|_{\Sigma}\right)  $. By elliptic
regularity we have
\[
\operatorname{Ker}\mathcal{D}_{D}^{\#}=\operatorname{Ker}\mathcal{D\;}%
\text{and\ } \operatorname{SF} \left\{  \mathcal{D}_{t,D}^{\#}\right\}  =
\operatorname{SF} \left\{  \mathcal{D}_{t}\right\}  .
\]

For product structures near $\Sigma$, one can apply Theorem
\ref{t:criss-cross} and obtain an $L^{2}$ version of the preceding Theorem
which gives a new proof and a slight generalization of the Yoshida--Nicolaescu
Formula (for details see \cite{BoFuOt01}, Section 3). For safety reasons, we
assume that $\mathcal{D}_{t}$ is a zero-order perturbation of $\mathcal{D}%
_{0}$, induced by a continuous change of the defining connection, or
alternatively, $\mathcal{D}_{t}:=\mathcal{D}_{0}+C_{t}$ where $C_{t}$ is a
self-adjoint bundle morphism.

\begin{theorem}%
\begin{multline*}
\operatorname{SF}\{\mathcal{D}_{t}\}=\mathbf{mas}\left(  \left\{  \Lambda
_{t}^{1}\cap L^{2}\left(  -\Sigma\right)  \overset{\cdot}{+}\Lambda_{t}%
^{2}\cap L^{2}\left(  \Sigma\right)  \right\}  ,\Delta\right) \\
=:\mathbf{mas}\left(  \{\Lambda_{t}^{1}\cap L^{2}(-\Sigma)\},\{\Lambda_{t}%
^{2}\cap L^{2}(\Sigma)\}\right)  ,
\end{multline*}
where the Cauchy data spaces $\Lambda_{t}^{j}$ are taken on each side $j=1,2 $
of the hypersurface $\Sigma$.
\end{theorem}

\begin{rem}
Here it is not compelling to use the symplectic geometry of the Cauchy data
spaces (see Remark \ref{r:index}b). Actually, deep gluing formulas can and
have been obtained for the spectral flow by coding relevant information not in
the full infinite-dimensional Cauchy data spaces but in families of Lagrangian
subspaces of suitable finite-dimensional symplectic spaces, like the kernel of
the tangential operator (see \cite{CaLeMi96} and \cite{CaLeMi00}).
\end{rem}

\medskip

\subsubsection{Correction Formula for the Spectral Flow}

Let $D,\,D^{\prime}$ with $D_{\min}<D,D^{\prime}<D_{\max}$ be two domains such
that both $\{\mathcal{D}_{t,D}\}$ and $\{\mathcal{D}_{t,D^{\prime}}\}$ become
families of self-adjoint Fredholm operators. We assume that $D$ and
$D^{\prime}$ differ only by a finite dimension, more precisely that
\begin{equation}
\dim\frac{\gamma(D)}{\gamma(D)\cap\gamma(D^{\prime})}=\dim\frac{\gamma
(D^{\prime})}{\gamma(D)\cap\gamma(D^{\prime})}<+\infty\,.
\label{e:corr-finite}%
\end{equation}
Then we find from Theorem \ref{t:sff-mwb} (for details see \cite{BoFu99},
Theorem 6.5):
\begin{align}
&  \mathbf{sf}\left\{  \mathcal{D}_{t,D}\right\}  -\mathbf{sf}\left\{
\mathcal{D}_{t,D^{\prime}}\right\} \label{e:hormander}\\
&  =\mathbf{mas}\left(  \left\{  \Lambda\left(  \mathcal{D}_{t}\right)
\right\}  ,\gamma(D^{\prime})\right)  -\mathbf{mas}\left(  \left\{
\Lambda\left(  \mathcal{D}_{t}\right)  \right\}  ,\gamma(D)\right) \nonumber\\
&  =\sigma_{\text{H\"{o}r}}\left(  \Lambda(\mathcal{D}_{0}),\Lambda
(\mathcal{D}_{1});\gamma(D^{\prime}),\gamma(D)\right) \nonumber
\end{align}
(see Remark \ref{r:horm}b). The assumption \eqref{e:corr-finite} is rather
restrictive. The pair of domains, for instance, defined by the
Atiyah--Patodi--Singer projection and the Calder\'{o}n projection, may not
always satisfy this condition. For the present proof, however, it seems indispensable.

\bigskip

\section{The Eta Invariant\label{s:eta}}

As mentioned in the Introduction, a systematic functional analytical frame is
missing for the eta invariant in contrast to the index and the spectral flow
(see, however, \cite{CoMo95} for an ambitious approach to establish analogues
of Sobolev spaces, pseudo-differential operators, and zeta and eta functions
in the context of noncommutative spectral geometry). Basically, however, the
concept of the eta invariant of a (total and compatible, hence self-adjoint)
Dirac operator is rather an immediate generalization of the index. Instead of
measuring the chiral asymmetry of the zero eigenvalues we now measure the
asymmetry of the \textit{whole} spectrum.

Let $\mathcal{D}$ be an operator of Dirac type; i.e., roughly speaking, an
operator with a real discrete spectrum which is nicely spaced without finite
accumulation points and with an infinite number of eigenvalues on both sides
of the real line. In close analogy with the definition of the zeta-function
for essentially positive elliptic operators like the Laplacian, we set
\[
\eta_{\mathcal{D}}(s):=\sum_{\lambda\in\operatorname{spec}(\mathcal{D}%
)\setminus\left\{  0\right\}  }\operatorname*{sign}(\lambda)\,\lambda^{-s}\,.
\]
Clearly, the formal sum $\eta_{\mathcal{D}}(s)$ is well defined for complex
$s$ with $\Re(s)$ sufficiently large, and it vanishes for a symmetric spectrum
(i.e., if for each $\lambda\in\operatorname{spec}(\mathcal{D})$ also
$-\lambda\in\operatorname{spec}(\mathcal{D})$).

In Subsections \ref{ss:IndexTDOonClosedMfds} and \ref{ss:APSIndexThm}, we
expressed the index by the difference of the traces of two related heat
operators. Similarly, we also have a heat kernel expression for the eta
function. Let $\mathcal{D}$ be any self-adjoint operator with compact
resolvent and let $\{\lambda_{k}\}$ denote its eigenvalues ordered so that
$\cdots\leq\lambda_{k-1}\leq\lambda_{k}\leq\lambda_{k+1}\leq\cdots$, each
repeated according its multiplicity. Formally, we have
\begin{align}
\eta_{\mathcal{D}}(s)\cdot\Gamma(\frac{s+1}{2})  &  =\sum\nolimits_{\lambda
_{k}\neq0}\operatorname*{sign}(\lambda_{k})\cdot|\lambda_{k}|^{-s}\cdot
\int_{0}^{\infty}r^{\frac{s-1}{2}}e^{-r}\,dr\label{e:eta-def}\\
&  =\sum\nolimits_{\lambda_{k}\neq0}\lambda_{k}(\lambda_{k}^{2})^{-\frac
{s+1}{2}}\int_{0}^{\infty}(t\lambda_{k}^{2})^{\frac{s-1}{2}}e^{-t\lambda
_{k}^{2}}\,d(t\lambda_{k}^{2})\nonumber\\
&  =\sum\nolimits_{\lambda_{k}\neq0}\int_{0}^{\infty}t^{\frac{s-1}{2}}%
\lambda_{k}e^{-t\lambda_{k}^{2}}\,dt=\int_{0}^{\infty}t^{\frac{s-1}{2}%
}\operatorname{Tr}e^{-t\mathcal{D}^{2}}\,dt.\nonumber
\end{align}
For comparison we give the corresponding formula for the zeta function of the
(positive) Dirac Laplacian $\mathcal{D}^{2}$:
\begin{equation}
\zeta_{\mathcal{D}^{2}}(s):=\operatorname{Tr}(\mathcal{D}^{2})^{-s}=\frac
{1}{\Gamma(s)}\,\int_{0}^{\infty}t^{s-1}\operatorname{Tr}e^{-t{\mathcal{D}%
}^{2}}\,dt. \label{e:zeta-def}%
\end{equation}
For the zeta function, we must assume that $\mathcal{D}$ has no vanishing
eigenvalues (i.e. $\mathcal{D}^{2}$ is positive). Otherwise the integral on
the right side is divergent. (The situation, however, can be cured by
subtracting the orthogonal projection onto the kernel of $\mathcal{D}^{2}$
from the heat operator before taking the trace.) For the eta function, on the
contrary, it clearly does not matter whether there are 0--eigenvalues and
whether the summation is over all or only over the nonvanishing eigenvalues.

The derivation of \eqref{e:eta-def} and \eqref{e:zeta-def} is completely
elementary for $\Re(s)>\frac{1+\dim X}{2}$\thinspace, resp. $\Re(s)>\frac{\dim
X}{2}$\thinspace, where $X$ denotes the underlying manifold. It follows at
once that $\eta(s)$ (and $\zeta(s)$) admit a meromorphic extension to the
whole complex plane. However, it is not clear at all how to characterize the
operators for which the eta function (and the zeta function) have a finite
value at $s=0$, the \textit{eta invariant} (resp. the \textit{zeta invariant}).

Historically, the eta invariant appeared for the first time in the 1970s as an
error term showing up in the index formula for the APS spectral boundary value
problem of a Dirac operator $\mathcal{D}$ on a compact manifold $X$ with
smooth boundary $\Sigma$ (see our Subsection \ref{ss:APSIndexThm}). More
precisely, what arose was the eta invariant of the tangential operator (i.e.,
the induced Dirac operator over the closed manifold $\Sigma$). Even in that
case it was hard to establish the existence and finiteness of the eta invariant.

Strictly speaking, one can define the eta invariant as the constant term in
the Laurent expansion of the eta function around the point $s=0$. For various
applications this suffices. Many practical calculations, however, are much
facilitated when we know the regularity \textit{a priori}.

Basically, there are three different approaches to establish it: the original
proof by {Atiyah}, {Patodi}, and {Singer}, worked out in \cite[Corollary
22.9]{BoWo93} and summarized in our Subsection \ref{ss:APSIndexThm}; it is
based on an assumption about the existence of a suitable asymptotic expansion
for the corresponding heat kernel on the infinite cylinder $\mathbb{R}%
_{+}\times\Sigma$. An intrinsic proof can be found in {Gilkey} \cite[Section
3.8]{Gi95}. It does not exploit that $\Sigma$ bounds $X$, but requires strong
topological means.

For a compatible (!) Dirac operator over a closed manifold {Bismut} and
{Freed} \cite{BiFr86} have shown that the eta function is actually a
holomorphic function of $s$ for $\Re(s)>-2$. They used the heat kernel
representation \eqref{e:eta-def} which implies that the eta invariant, when it
exists, can be expressed as
\[
\eta_{\mathcal{D}}(0)=\frac{1}{\sqrt{\pi}}\int_{0}^{\infty}\frac{1}{\sqrt{t}%
}\cdot\operatorname{Tr}\left(  \mathcal{D}e^{-t\mathcal{D}^{2}}\right)  \,dt.
\]
It follows from \eqref{e:eta-def} that the estimate
\[
|\operatorname{Tr}\mathcal{D}e^{-t\mathcal{D}^{2}}|<c\sqrt{t}%
\]
implies the regularity of the eta function at $s=0$. In fact, {Bismut} and
{Freed} proved a sharper result, using nontrivial results from stochastic
analysis. Inspired by calculations presented in {Bismut} and {Cheeger}
\cite[Section 3]{BiCh89}, {Wojciechowski} gave a purely analytic reformulation
of the details of their proof. This provides a third and completely elementary
way of proving the regularity of the eta function at $s=0$. The essential
steps are:

\begin{theorem}
\label{t:bismut_freed} Let $\mathcal{D}:C^{\infty}(\Sigma;E)\rightarrow
C^{\infty}(\Sigma;E)$ denote a compatible Dirac operator over a closed
manifold $\Sigma$ of odd dimension $m$. Let $\operatorname{e}(t;x,x^{\prime})$
denote the integral kernel of the heat operator $e^{-t\mathcal{D}^{2}}$. Then
there exists a positive constant $C$ such that
\[
\lvert\left.  \operatorname{Tr}\left(  \mathcal{D}_{x}\/e(t;x,x^{\prime
})\right)  \right|  _{x=x^{\prime}}\rvert<C\sqrt{t}%
\]
for all $x\in\Sigma$ and $0<t<1$.
\end{theorem}

We recall the definition of the `local' $\eta$ function.

\begin{definition}
Let $\{f_{k};\lambda_{k}\}_{k\in\mathbb{Z}}$ be a discrete spectral resolution
of $\mathcal{D}$. Then we define
\begin{align*}
\eta_{\mathcal{D}}(s;x)  &  :=\sum_{\lambda_{k}\neq0}\operatorname*{sign}%
(\lambda_{k})\/\lvert\lambda_{k}\rvert^{-s}\/\left\langle f_{k}(x),f_{k}%
(x)\right\rangle _{E_{x}}\\
&  =\frac{1}{\Gamma(\frac{s+1}{2})}\int_{0}^{\infty}t^{\frac{s-1}{2}}\left(
\sum_{\lambda_{k}\neq0}\lambda_{k}e^{-t\lambda_{k}^{2}}\/\left\langle
f_{k}(x),f_{k}(x)\right\rangle \right)  dt\\
&  =\frac{1}{\Gamma(\frac{s+1}{2})}\int_{0}^{\infty}t^{\frac{s-1}{2}%
}\operatorname{Tr}\mathcal{D}\/\operatorname{e}(t;x,x)\/dt.
\end{align*}
\end{definition}

\begin{corollary}
Under the assumptions of the preceding theorem the `local' $\eta$ function
$\eta_{\mathcal{D}}(s;x)$ is holomorphic in the half plane $\Re(s)\geq-2$ for
any $x\in\Sigma$.
\end{corollary}

\bigskip

In the decade or so following 1975, it was generally believed that the
existence of a finite eta invariant was a very special feature of operators of
Dirac type on closed manifolds which are boundaries, and then, more generally,
of Dirac type operators on all closed manifolds. Only after the seminal paper
by {Douglas} and {Wojciechowski} \cite{DoWo91} was it gradually realized that
globally elliptic self-adjoint boundary value problems for operators of Dirac
type also have a finite eta invariant. Once again, there are quite different
approaches to obtain that result.

The work by {Wojciechowski} and collaborators is based on the Duhamel
Principle and provides complete asymptotic expansions of the heat kernels for
the self-adjoint Fredholm extension $\mathcal{D}_{P}$ where $P$ belongs to the
smooth self-adjoint Grassmannian. We recall from our review of the
Atiyah--Patodi--Singer Index Theorem that the Duhamel Principle allows one to
study the interior contribution and the boundary contribution separately and
identify the singularities caused by the boundary contribution. It seems that
Wojciechowski's method is only applicable if the metric structures close to
the boundary are product.

Based on joint work with {Seeley}, {Grubb} \cite{Gr99} also obtained nice
asymptotic expansions of the trace of the heat kernels in \eqref{e:eta-def}
and \eqref{e:zeta-def}. Contrary to the qualitative arguments of the
Duhamel--Wojciechowski approach, Grubb's approach requires the explicit
computation of certain coefficients in the expansion. Some of them have to
vanish to guarantee the desired regularity.

For a slightly larger class of self-adjoint Fredholm extensions, {Br{\"{u}%
}ning} and {Lesch} \cite{BrLe99a} also studied the eta invariant. However, the
authors had to deal with the residue of the eta function at $s=0$ which is not
present in the Duhamel--Wojciechowski approach.

In spite of the differences between the various approaches and types of
results it seems that one consequence can be drawn immediately, namely there
must be a deeper meaning of the eta invariant beyond its role as error term in
the APS index formula. Happily, such a meaning was found by {Singer} in
\cite{Si85} for the eta invariant on closed manifolds and successively
generalized by {Wojciechowski} and collaborators for manifolds with boundary,
namely the identification of the eta invariant as the phase of the zeta
function regularized determinant. We shall explain this now, and postpone our
main topic, the pasting formulas for the eta invariant on partitioned manifolds.

\bigskip

\subsection{Functional Integrals and Spectral Asymmetry\label{ss:physics}}

Several important quantities in quantum mechanics and quantum field theory are
expressed in terms of \emph{quadratic functionals} and \emph{functional
integrals}. The concept of the determinant for Dirac operators arises
naturally when one wants to evaluate the corresponding path integrals. As
{Itzykson} and {Zuber} report in the chapter on \emph{Functional Methods} of
their monograph \cite{ItZu80}: ``The path integral formalism of {Feynman} and
{Kac} provides a unified view of quantum mechanics, field theory, and
statistical models. The original suggestion of an alternative presentation of
quantum mechanical amplitudes in terms of path integrals stems from the work
of {Dirac} (1933) and was brilliantly elaborated by {Feynman} in the 1940s.
This work was first regarded with some suspicion due to the difficult
mathematics required to give it a decent status. In the 1970s it has, however,
proved to be the most flexible tool in suggesting new developments in field
theory and therefore deserves a thorough presentation.''

We shall restrict our discussion to the easiest variant of that complex matter
by focusing on the partition function of a quadratic functional given by the
Euclidean action of a \emph{Dirac operator} which is assumed to be elliptic
with imaginary time due to Wick rotation and coupled to continuously varying
vector potentials (sources, fields, connections), for the ease of presentation
in vacuum. We refer to {Bertlmann}, \cite{Be96} and {Schwarz}, \cite{Sch93}
for an introduction to the quantum theoretic language for mathematicians and
for a more extensive treatment of general aspects of quadratic functionals and
functional integrals involving the relations to the Lagrangian and Hamiltonian formalism.

There are various alternative notions of ``path integral'' around, some more
sophisticated than others. A mathematically rigorous formulation of the
concept of ``path integrals'', as physicists typically ``understand'' it in
quantum field theory is flimsy at best. For fields on Minkowski space, they
are not mathematical integrals at all, because no measure is defined. For
fields on Euclidean space (with imaginary time), one can construct genuine
measures in limited settings which are not entirely realistic. Even then there
is the issue of continuing the integrals back to real time which is generally
ill-defined on a curved space-time. Many physicists do not care about such
matters. Indeed, such physicists use path integrals primarily as compact
generators of recipes to produce Feynman integrals which when regularized and
renormalized yield coefficients in a formal power series in coupling constants
for physical quantities of interest. However, no one has ever proved that
these series converge, even for quantum electrodynamics (QED). Indeed the
consensus of those who care is that these are only asymptotic series. Adding
first few terms of these renormalized perturbation series yields remarkable
11-decimal point agreement with experiment, and QED is thereby hailed as a
huge success. For many mathematicians and a few physicists, the great tragedy
is that these recipes work so well without a genuine mathematical foundation.
Although path integrals may not make precise sense per se, not only do they
generate successful recipes in physics, but they can motivate fruitful ideas
and precise concepts. In mathematics, path integrals have motivated the
$\zeta$-regularized determinant for the (Euclidean) Dirac operator, as a
mathematically genuine canonical object, independent of particular choices
made for regularization, which can be precisely calculated in principle.

A special feature of Dirac operators is that their determinants involve a
\emph{phase}, the imaginary part of the determinant's logarithm. As we will
see now, this is a consequence of the fact that, unlike second-order
semi-bounded Laplacians, first-order Dirac operators have an infinite number
of both positive and negative eigenvalues. Then the phase of the determinant
reflects the spectral asymmetry of the corresponding Dirac operator.

The simplest path integral we meet in quantum field theory takes the form of
the \emph{partition function} and can be written formally as the integral%

\begin{equation}
Z(\beta):=\int_{\Gamma}\,e^{-\beta S(\omega)}\,d\omega\,, \label{e:partitionf}%
\end{equation}
where $d\omega$ denotes functional integration over the space $\Gamma
:=\Gamma(M;E)$ of sections of a Euclidean vector bundle $E$ over a Riemannian
manifold $M$.

In quantum theoretic language, $M$ is space or space-time; a $\omega\in\Gamma$
is a position function of a particle or a spinor field. The scaling parameter
$\beta$ is a real or complex parameter, most often $\beta=1$. The functional
$S$ is a quadratic real-valued functional on $\Gamma$ defined by
$S(\omega):=\left\langle \omega,T\omega\right\rangle $ with a fixed linear
symmetric operator $T:\Gamma\rightarrow\Gamma$. Typically $T:=\mathcal{D}$ is
a Dirac operator and $S(\omega)$ is the action $S(\omega)=\int_{M}\left\langle
\omega,\mathcal{D}\omega\right\rangle $.

Mathematically speaking, the integral (\ref{e:partitionf}) is an oscillating
integral like the \emph{Gaussian integral}. It is ill-defined in general because

\begin{itemize}
\item [(I)]as it stands, it is meaningless when $\dim\Gamma(M;E)=+\infty$
(i.e., when $\dim M\geq1$); and,

\item[(II)] even when $\dim\Gamma(M;E)<\infty$ (i.e. when $\dim M=0$ and $M$
consists of a finite set of points), the integral $Z(\beta)$ diverges unless
$\beta S(\omega)$ is positive and nondegenerate.
\end{itemize}

Nevertheless, these expressions have been used and construed in quantum field
theory. As a matter of fact, reconsidering the physicists' use and
interpretation of these mathematically ill-defined quantities, one can
describe certain formal manipulations which lead to normalizing and evaluating
$Z(\beta)$ in a mathematically precise way.

We begin with a few calculations in \emph{Case II}, inspired by {Adams} and
{Sen}, \cite{AdSe96}, to show how \emph{spectral asymmetry} is naturally
entering into the calculations even in the finite-dimensional case and how
this suggests a definition of the determinant in the infinite-dimensional case
for the Dirac operator.

Then, let $\dim\Gamma=d<\infty$ and, for a symmetric endomorphism $T$, let
\[
S(\omega):=\left\langle \omega,T\omega\right\rangle \text{ for all }\omega
\in\Gamma.
\]
\textit{Case 1}. We assume that $S$ positive and nondegenerate, i.e. $T$ is
strictly positive, say $\operatorname{spec}T=\{\lambda_{1},\dots\lambda_{d}\}$
with all $\lambda_{j}>0$. This is the classical case. We choose an orthonormal
system of eigenvectors $(e_{1},\dots,e_{d})$ of $T$ as basis for $\Gamma$. We
have $S(\omega)=\sum\lambda_{j}x_{j}^{2}$ for $\omega=\sum x_{j}e_{j}$ and,
for real $\beta>0,$ we get
\begin{align*}
Z(\beta)  &  =\int_{G}e^{-\beta S(\omega)}\,d\omega=\int_{\mathbb{R}^{d}%
}\/dx_{1}\/\dots dx_{d}\,e^{-\beta\sum\lambda_{j}x_{j}^{2}}\\
&  =\int_{-\infty}^{\infty}e^{-\beta\lambda_{1}x_{1}^{2}}\/dx_{1}\int
_{-\infty}^{\infty}e^{-\beta\lambda_{2}x_{2}^{2}}\/dx_{2}\,\dots\int_{-\infty
}^{\infty}e^{-\beta\lambda_{d}x_{d}^{2}}\/dx_{d}\\
&  =\sqrt{\frac{\pi}{\beta\lambda_{1}}}\,\sqrt{\frac{\pi}{\beta\lambda_{2}}%
}\dots\sqrt{\frac{\pi}{\beta\lambda_{d}}}=\pi^{d/2}\cdot\beta^{-d/2}\cdot(\det
T)^{-1/2}\,.
\end{align*}
In that way the determinant appears when evaluating the simplest quadratic integral.

\textit{Case 2}. If the functional $S$ is positive and degenerate, $T\geq0$,
the partition function is given by
\[
Z(\beta)=\pi^{\zeta/2}\cdot\beta^{-\mathbb{\zeta}/2}\cdot(\det\omega
T)^{-1/2}\cdot\operatorname*{Vol}(\operatorname{Ker}T),
\]
where $\mathbb{\zeta}:=\dim\Gamma-\dim\operatorname{Ker}T$ and $\widetilde
{T}:=T|_{(\operatorname{Ker}T)^{\perp}}$, but, of course $\operatorname*{Vol}%
(\operatorname{Ker}T)=\infty$. For approaches to ``renormalizing'' this
quantity in quantum chromodynamics, we refer to \cite{AdSe96}, \cite{BMSW97},
\cite{Sch93}. One approach customary in physics is to take $\pi^{\zeta/2}%
\beta^{-\mathbb{\zeta}/2}(\det\widetilde{T})^{-1/2}$ as the definition of the
integral by setting the factor $\operatorname*{Vol}(\operatorname{Ker}T)$
equal to 1.

\textit{Case 3}. Now we assume that the functional $S$ is nondegenerate, i.e.
$T$ invertible, but $S$ is neither positive nor negative. We decompose
$\Gamma=\Gamma_{+}\times\Gamma_{-}$ and $T=T_{+}\oplus T_{-}$ with
$T_{+},-T_{-}$ strictly positive on $\Gamma_{\pm}$. Formally, we obtain
\begin{align*}
Z(\beta)  &  =\left(  \int_{\Gamma_{+}}d\omega_{+}e^{-\beta\left\langle
\omega_{+},T_{+}\omega\right\rangle }\right)  \left(  \int_{\Gamma_{-}}%
d\omega_{-}e^{-(-\beta)\left\langle \omega_{-},-T_{-}\omega\right\rangle
}\right) \\
&  =\pi^{d_{+}/2}\beta^{-d_{+}/2}(\det T_{+})^{-1/2}\,\pi^{d_{-}/2}%
(-\beta)^{-d_{-}/2}(\det-T_{-})^{-1/2}\\
&  =\pi^{\mathbb{\zeta}/2}\beta^{-d_{+}/2}(-\beta)^{-d_{-}/2}(\det
|T|)^{-1/2}\,
\end{align*}
where $d_{\pm}:=\dim\Gamma_{\pm}$, hence $\mathbb{\zeta}=d_{+}+d_{-}$ and
$|T|:=\sqrt{\widetilde{T}^{2}}=T_{+}\oplus-T_{-}$.

\textit{Case 4}. In the preceding formula, the term $(\beta)^{-d_{+}/2}%
(-\beta)^{-d_{-}/2}$ is undefined for $\beta\in\mathbb{R}_{\pm}$. We shall
replace it by a more intelligible term for $\beta=1$ by first expanding
$Z(\beta)$ in the upper complex half plane and then formally setting $\beta=1
$. More precisely, let $\beta\in\mathbb{C}_{+}=\{z\in\mathbb{C}\mid\Im z>0\}$
and write $\beta=|\beta|e^{i\theta}$ with $\theta\in\lbrack0,\pi]$, hence
$-\beta=|\beta|e^{i(\theta-\pi)} $ with $\theta-\pi\in\lbrack-\pi,0]$. We set
$\beta^{a}:=|\beta|^{a}e^{i\theta a}$ and get
\begin{align*}
\beta^{-d_{+}/2}(-\beta)^{-d_{-}/2}  &  =(|\beta|e^{i\theta})^{-d_{+}%
/2}(|\beta|e^{i(\theta-\pi)})^{-d_{-}/2}\\
&  =|\beta|^{-\mathbb{\zeta}/2}\,e^{-i\frac{d_{+}}{2}\theta}\,e^{-i\frac
{d_{-}}{2}\theta}\,e^{i\pi\frac{d_{-}}{2}}\,.
\end{align*}
Moreover,
\begin{align*}
-\frac{d_{+}}{2}\theta-\frac{d_{-}}{2}\theta+\pi\frac{d_{-}}{2}  &
=-\frac{\theta}{2}\left(  d_{+}+d_{-}\right)  +\frac{\pi}{2}\left(
\frac{d_{-}}{2}+\frac{d_{+}}{2}+\frac{d_{-}}{2}-\frac{d_{+}}{2}\right) \\
&  =-\frac{\theta}{2}\mathbb{\zeta}+\frac{\pi}{4}\left(  \mathbb{\zeta}%
-\eta\right)  =-\frac{\pi}{4}\left(  \frac{2\theta\mathbb{\zeta}}{\pi}%
+(\eta-\mathbb{\zeta})\right)  ,
\end{align*}
where $\zeta:=d_{+}+d_{-}$ is the finite-dimensional equivalent of the $\zeta$
invariant, counting the eigenvalues, and $\eta:=d_{+}-d_{-}$ the
finite-dimensional equivalent of the $\eta$ invariant, measuring the spectral
asymmetry of $T$. We obtain
\[
Z\left(  \beta\right)  =\pi^{\mathbb{\zeta}/2}\left|  \beta\right|
^{-\mathbb{\zeta}/2}e^{-i\frac{\pi}{4}(\frac{2\mathbb{\zeta}\theta}{\pi}%
+(\eta-\mathbb{\zeta}))}\,(\det|T|)^{-1/2}%
\]
and, formally, for $\beta=1$, i.e. $\theta=0$,
\begin{equation}
Z(1)=\pi^{\mathbb{\zeta}/2}\,\underbrace{e^{-i\frac{\pi}{4}(\eta
-\mathbb{\zeta})}\,(\det|T|)^{-1/2}}_{=:\det T}\,. \label{e:finite}%
\end{equation}
Equation (\ref{e:finite}) suggests a nonstandard definition of the determinant
for the infinite-dimensional case.\bigskip

\begin{remark}
(a) The methods and results of this section also apply to real-valued
quadratic functionals on complex vector spaces. Since the integration in
(\ref{e:partitionf}) in this case is over the real vector space underlying
$\Gamma$, the expressions for the partition functions in this case become the
square of those above.\newline (b) In the preceding calculations we worked
with ordinary commuting numbers and functions. The resulting Gaussian
integrals are also called \emph{bosonic} integrals. If we consider
\emph{fermionic} integrals, we work with Grassmannian variables and obtain the
determinant not in the denominator but in the nominator (see e.g.
\cite{BeGeVe92} or \cite{Be96}).
\end{remark}

\subsection{The $\mathbb{\zeta}$--Determinant for Operators of Infinite
Rank\label{s:zetadet}}

Once again, our point of departure is finite-dimensional linear algebra. Let
$T:\mathbb{C}^{d}\rightarrow\mathbb{C}^{d}$ be an invertible, positive
operator with eigenvalues $0<\lambda_{1}\leq\lambda_{2}\leq...\leq\lambda_{d}$
. We have the equality
\begin{align*}
\det T  &  =\prod\lambda_{j}=\exp\{\sum\ln\lambda_{j}e^{-s\ln\lambda_{j}%
}|_{s=0}\}\\
&  =\exp(-\frac{d}{ds}(\sum\lambda_{j}^{-s})|_{s=0})=e^{-\frac{d}{ds}\zeta
_{T}(s)|_{s=0}}\,,
\end{align*}
where $\mathbb{\zeta}_{T}(s):=\sum_{j=1}^{d}\lambda_{j}^{-s}$.

We show that the preceding formula generalizes naturally, when $T$ is replaced
by a \emph{positive definite self-adjoint elliptic} operator $L$ (for the ease
of presentation, of second order, like the Laplacian) acting on sections of a
Hermitian vector bundle over a closed manifold $M$ of dimension $m$. Then $L$
has a discrete spectrum $\operatorname{spec}L=\{\lambda_{j}\}_{j\in\mathbb{N}%
}$ with $0<\lambda_{1}\leq\lambda_{2}\leq\dots$, satisfying the asymptotic
formula $\lambda_{n}\sim Cn^{m/2}$ for a constant $C>0$ depending on $L$ (see
e.g. \cite{Gi95}, Lemma 1.6.3). We extend $\zeta_{L}(s):=\sum_{j=1}^{\infty
}\lambda_{j}^{-s}$ in the complex plane by
\[
\zeta_{L}(s):=\frac{1}{\Gamma(s)}\int_{0}^{\infty}t^{s-1}\,\operatorname{Tr}%
e^{-tL}\,dt
\]
with $\Gamma(s):=\int_{0}^{\infty}t^{s-1}e^{-t}dt$. Note that $e^{-tL}$ is the
heat operator transforming any initial section $f_{0}$ into a section $f_{t}$
satisfying the heat equation $\tfrac{\partial}{\partial t}f+Lf=0$. Clearly
$\operatorname{Tr}e^{-tL}=\sum e^{-t\lambda_{j}}$\thinspace.

One shows that the original definition of $\mathbb{\zeta}_{L}(s)$ yields a
holomorphic function for $\Re(s)$ large and that its preceding extension is
meromorphic in the entire complex plane with simple poles only. The point
$s=0$ is a regular point and $\mathbb{\zeta}_{L}(s)$ is a holomorphic function
at $s=0$. From the asymptotic expansion of $\Gamma(s)\sim\frac{1}{s}%
+\gamma+sh(s)$ close to $s=0$ with the Euler number $\gamma$ and a suitable
holomorphic function $h$ we obtain an explicit formula
\[
\mathbb{\zeta}_{L}^{\prime}(0)\sim\int_{0}^{\infty}\frac{1}{t}%
\operatorname{Tr}e^{-tL}\,dt-\gamma\mathbb{\zeta}_{L}(0)\,.
\]
This is explained in great detail, e.g., in \cite{Wo99}. Therefore, {Ray} and
{Singer} in \cite{RaSi71} could introduce $\det_{z}(L)$ by defining:
\[
{\det}_{z}L:=e^{-\frac{d}{ds}\zeta_{L}(s)|_{s=0}}=e^{-\zeta_{L}^{\prime}%
(0)}\,.
\]

\medskip The preceding definition does not apply immediately to the main hero
here, the Dirac operator $\mathcal{D}$ which has infinitely many positive
${\lambda_{j}}$ and negative eigenvalues ${-\mu_{j}}$. Clearly by the
preceding argument
\[
{\det}_{z}\mathcal{D}^{2}=e^{-\mathbb{\zeta}_{\mathcal{D}^{2}}^{\prime}}%
\qquad\text{and}\qquad{\det}_{z}|\mathcal{D}|=e^{-\mathbb{\zeta}%
_{|\mathcal{D}|}^{\prime}}=e^{-\frac{1}{2}\mathbb{\zeta}_{\mathcal{D}^{2}%
}^{\prime}}\,.
\]
For the Dirac operator we set
\[
\ln\det\mathcal{D}:=-\tfrac{d}{ds}\zeta_{\mathcal{D}}(s)|_{s=0}%
\]
with, choosing \footnote{Choosing the alternative representation, namely
$(-1)^{-s}=e^{-i\pi s}$ yields the opposite sign of the phase of the
determinant which may appear to be more natural for some quantum field models
and also in view of \eqref{e:finite}. However, we follow the more common
convention introduced by Singer in [Si85, p. 331] when defining the
determinant of operators of Dirac type.} the branch $(-1)^{-s}=e^{i{\pi}s}$,%

\begin{align*}
\mathbb{\zeta}_{\mathcal{D}}(s)  &  =\sum\lambda_{j}^{-s}+\sum(-1)^{-s}\mu
_{j}^{-s}=\sum\lambda_{j}^{-s}+e^{i\pi s}\sum\mu_{j}^{-s}\\
&  =\frac{\sum\lambda_{j}^{-s}+\sum\mu_{j}^{-s}}{2}+\frac{\sum\lambda_{j}%
^{-s}-\sum\mu_{j}^{-s}}{2}\\
&  \qquad\qquad+e^{i\pi s}\left\{  \frac{\sum\lambda_{j}^{-s}+\sum\mu_{j}%
^{-s}}{2}-\frac{\sum\lambda_{j}^{-s}-\sum\mu_{j}^{-s}}{2}\right\} \\
&  =\tfrac{1}{2}\left\{  \zeta_{\mathcal{D}^{2}}(\frac{s}{2})+\eta
_{\mathcal{D}}(s)\right\}  +\tfrac{1}{2}e^{i{\pi}s}\left\{  \zeta
_{\mathcal{D}^{2}}(\frac{s}{2})-\eta_{\mathcal{D}}(s)\right\}  ,
\end{align*}
where $\eta_{\mathcal{D}}(s):=\sum\lambda_{j}^{-s}-\sum\mu_{j}^{-s}$. Later we
will show that $\eta_{\mathcal{D}}(s)$ is a holomorphic function of $s$ for
$\Re(s)$ large with a meromorphic extension to the whole complex plane which
is holomorphic in the neighborhood of $s=0$ . We obtain
\begin{multline*}
\mathbb{\zeta}_{\mathcal{D}}^{\prime}(s)=\tfrac{1}{4}\mathbb{\zeta
}_{\mathcal{D}^{2}}^{\prime}(\tfrac{s}{2})+\tfrac{1}{2}\eta_{\mathcal{D}%
}^{\prime}(s)+\tfrac{1}{2}i\pi e^{i\pi s}\{\mathbb{\zeta}_{\mathcal{D}^{2}%
}(\tfrac{s}{2})-\eta_{\mathcal{D}}(s)\}\\
+\tfrac{1}{2}e^{i\pi s}\{\tfrac{1}{2}\mathbb{\zeta}_{\mathcal{D}^{2}}^{\prime
}(\tfrac{s}{2})-\eta_{\mathcal{D}}^{\prime}(s)\}.
\end{multline*}
It follows that
\[
\mathbb{\zeta}_{\mathcal{D}}^{\prime}(0)=\tfrac{1}{2}\mathbb{\zeta
}_{\mathcal{D}^{2}}^{\prime}(0)+\tfrac{i\pi}{2}\left\{  \mathbb{\zeta
}_{\mathcal{D}^{2}}(0)-\eta_{\mathcal{D}}(0)\right\}
\]
and
\begin{align*}
{\det}_{z}\mathcal{D}  &  =e^{-\frac{1}{2}\mathbb{\zeta}_{\mathcal{D}^{2}%
}^{\prime}(0)}\,e^{-\frac{i\pi}{2}\left\{  \mathbb{\zeta}_{\mathcal{D}^{2}%
}(0)-\eta_{\mathcal{D}}(0)\right\}  }=e^{-\frac{i\pi}{2}\left\{
\mathbb{\zeta}_{|\mathcal{D}|}(0)-\eta_{\mathcal{D}}(0)\right\}
}\,e^{-\mathbb{\zeta}_{|\mathcal{D}|}^{\prime}(0)}\\
&  =e^{-\frac{i\pi}{2}\left\{  \mathbb{\zeta}_{|\mathcal{D}|}(0)-\eta
_{\mathcal{D}}(0)\right\}  }\,{\det}_{z}|\mathcal{D}|\,.
\end{align*}
So, the Dirac operator's `partition function' in the sense of
(\ref{e:partitionf}) becomes
\[
Z(1)=\pi^{\mathbb{\zeta}_{|\mathcal{D}|}(0)}({\det}_{\mathbb{\zeta}%
}\mathcal{D})^{-\frac{1}{2}}\,.
\]

\bigskip

\subsection{Spectral Invariants of Different
`Sensitivity'\label{ss:SpectralInvariants}}

In the preceding formulas four spectral invariants of the Dirac operator
$\mathcal{D}$ enter:

\subsubsection{The Index\label{sss:The Index}}

First recall that the index of arbitrary elliptic operators on closed
manifolds and the spectral flow of 1-parameter families of self-adjoint
elliptic operators are topological invariants and so stable under small
variation of the coefficients and, by definition, solely depending on the
multiplicity of the eigenvalue 0. In the theory of bounded or closed (not
necessarily bounded) Fredholm operators and bounded or not necessarily bounded
self-adjoint Fredholm operators and the related $K$ and $K^{-1}$ theory, we
have a powerful functional analytical and topological frame for discussing
these invariants. Moreover, index and spectral flow are local invariants,
i.e., can be expressed by an integral where the integrand is locally expressed
by the coefficients of the operator(s). Consequently, we have simple, precise
pasting formulas for index and spectral flow on partitioned manifolds where
the error term is localized along the separating hypersurface.\newline
\ \ \ \ On manifolds with boundary the Calder\'{o}n projection and its range,
the Cauchy data space, change continuously when we vary the Dirac operator as
shown in Section \ref{sss:AdmissableBVP}, exploiting the unique continuation
property of operators of Dirac type. The same is not true for the
Atiyah--Patodi--Singer projection: it can jump from one connected component of
the Grassmannian to another component under small changes of the Dirac
operator. Correspondingly, the index of the APS boundary problem can jump
under small variation of the coefficients (i.e., of the defining connection or
the underlying Riemannian or Clifford structure).\newline \ \ \ \ Regarding
parity of the manifold, the index density (of all elliptic differential
operators) vanishes on odd-dimensional manifolds for symmetry reason. Then, in
the closed case the index vanishes, and on manifolds with boundary the APS
Index Theorem takes the simple form $\operatorname*{index}\mathcal{D}%
_{P_{\geq}}=-\dim\operatorname{Ker} B^{+}$ where we have the total Dirac
operator on the left and a chiral component of the induced tangential operator
on the right. Once again, the formula shows the instability of the index under
small changes of the Dirac operator. This is no contradiction to the stability
of the index on the spaces $\mathcal{F}$, respectively $\mathcal{CF}$,
discussed in Section \ref{sss:metrics} because the graph norm distance between
two APS realizations $\mathcal{D}_{P_{\geq}}$ and $\mathcal{D}_{P_{\geq}%
}^{\prime}$ can remain bounded away from zero when $\mathcal{D}$ runs to
$\mathcal{D}^{\prime}$. This is the case if and only if the dimension of the
kernel of the tangential operator changes under the deformation.

\subsubsection{The $\zeta$--invariant\label{sss:The Zeta}}

Similarly, $\mathbb{\zeta}_{\mathcal{D}^{2}}(0)$ is also local; i.e., it is
given by the integral $\int_{M}\alpha(x)\/dx$, where $\alpha(x)$ denotes a
certain coefficient in the heat kernel expansion and is locally expressed by
the coefficients of $\mathcal{D}$. In particular, $\mathbb{\zeta}%
_{\mathcal{D}^{2}}(0)$ remains unchanged for small changes of the spectrum.
Actually, $\mathbb{\zeta}_{L}(0)$ vanishes when $L$ is the square of a
self-adjoint elliptic operator on a closed manifold of odd dimension. It can
be defined (and it vanishes, see \cite[Appendix]{PaWo02a}) for a large class
of squares of operators of Dirac type with globally elliptic boundary
conditions on compact, smooth manifolds (of odd dimension) with boundary. So,
there are no nontrivial pasting formulas at all in such cases.

\subsubsection{The $\eta$--invariant\label{sss:The Eta}}

Unlike the index and spectral flow on closed manifolds, we have neither an
established functional analytical nor a topological frame for discussing
$\eta_{\mathcal{D}}(0)$, nor is it given by an integral of a locally defined
expression. On the circle, e.g., consider the operator
\begin{equation}
\mathcal{D}_{a}:=-i\frac{d}{dx}+a=e^{-ixa}\mathcal{D}_{0}e^{ixa},
\label{e:nonlocal_eta}%
\end{equation}
so that $\mathcal{D}_{a}$ and $\mathcal{D}_{0}$ have the same total symbol
(i.e., coincide locally), but $\eta_{\mathcal{D}_{a}}(0)=-2a$ depends on $a$.

The $\eta$--invariant depends, however, only on finitely many terms of the
symbol of the resolvent $(\mathcal{D}-\lambda)^{-1}$ and the real part (in
$\mathbb{R}/\mathbb{Z}$) will not change when one changes or removes a finite
number of eigenvalues. The integer part changes according to the net sign
change occurring under removing or modifying eigenvalues. Moreover, the first
derivative of the eta invariant of a smooth 1-parameter family of Dirac type
operators is local, namely the spectral flow, as noted in our Introduction.
This leads again to precise (though not so simple) pasting formulas for the
eta invariant on partitioned manifolds.\newline \ \ \ \ In even dimensions,
the eta invariant vanishes on any closed manifold $\Sigma$ for any Dirac type
operator which is the tangential operator of a Dirac type operator on a
suitable manifold which has $\Sigma$ as its boundary because of the induced
precise symmetry of the spectrum due to the anti-commutativity of the
tangential operator with Clifford multiplication. For the study of eta of
boundary value problems on even-dimensional manifolds see \cite{KlWo96}.

\subsubsection{The Modulus of the Determinant\label{sss:The Modulus of Det}}

The number $\mathbb{\zeta}_{\mathcal{D}^{2}}^{\prime}(0)$ is the most delicate
of the invariants involved: It is neither a local invariant, nor does it
depend only on the total symbol of the Dirac operator. Even small changes of
the eigenvalues will change the $\mathbb{\zeta}^{\prime}$ invariant and hence
the determinant. Moreover, no precise pasting formulas are obtained but only
adiabatic ones (i.e., by inserting a long cylinder around the separating
hypersurface (see \cite{PaWo02a}, \cite{PaWo02b}, \cite{PaWo02c}).

Without proof we present the main results by Wojciechowski and collaborators,
based on \cite{Wo99} where the $\zeta$--function regularized determinant was
established for pseudo-differential boundary value conditions belonging to the
smooth, self-adjoint Grassmannian. The first is a boundary correction formula,
proved in \cite{ScWo00} (see also the recent \cite{Sco02}):

\begin{theorem}
\label{t:scott-woj} \textrm{(Scott, Wojciechowski)}. Let $\mathcal{D}$ be a
Dirac operator over an odd--dimen\-sional compact manifold $M$ with boundary
$\Sigma$ and let $P\in\mathcal{G}\mathrm{r}^{\operatorname{sa}}(\mathcal{D}%
)$\,. Then the range of the Calder{\'o}n projection $\mathcal{P}(\mathcal{D})$
(the Cauchy data space $\Lambda(\mathcal{D},\frac12)$) and the range of $P$
can be written as the graphs of unitary, elliptic operators of order 0, $K$,
resp. $T$ which differ from the operator $(B^{+}B^{-})^{-1/2} B^{+}:
C^{\infty}(\Sigma;S^{+}|_{\Sigma}) \to C^{\infty}(\Sigma;S^{-}|_{\Sigma})$ by
a smoothing operator. Moreover,
\begin{equation}
{\det}_{\zeta}\mathcal{D}_{P} = {\det}_{\zeta}\mathcal{D}_{\mathcal{P}%
(\mathcal{D})} \,\cdot\ {\det}_{\operatorname{Fr}} \frac12 (\operatorname{I} +
KT^{-1})\,.
\end{equation}
\end{theorem}

The second result, in most simple form, is found in \cite{PaWo00}:

\begin{theorem}
\textrm{(Park, Wojciechowski)} Let $R\in\mathbb{R}$ be positive, let $M^{R}$
denote the stretched partitioned manifold $M^{R}=M_{1}\,\cup_{\Sigma
}\,[-R,0]\times\Sigma\,\cup_{\Sigma}\,[0,R]\times\Sigma\,\cup_{\Sigma}\,
M_{2}$\,, and let $\mathcal{D}_{R}$, $\mathcal{D}_{1,R}$, $\mathcal{D}_{2,R}$
denote the corresponding Dirac operators. We assume that the tangential
operator $\mathcal{B}$ is invertible. Then
\[
\lim_{R\to\infty} \frac{{\det}_{\zeta}\mathcal{D}_{R}^{\ 2}} {\Bigl({\det
}_{\zeta}(\mathcal{D}_{1,R})_{\operatorname{I}- P_{>}}^{\,2}\Bigr) \cdot
\Bigl({\det}_{\zeta}(\mathcal{D}_{2,R})_{P_{>}}^{\,2}\Bigr)} \ =\ 2^{-\zeta
_{\mathcal{B}^{2}}(0)}\,.
\]
\end{theorem}

\medskip

Although {Felix Klein} in \cite{Kl27} rated the determinant as \emph{simplest
example of an invariant}, today we must give an inverse rating. For the
present authors, it is not the invariants that are stable under the largest
transformation groups which deserve the highest interest, but rather
(according to {Dirac}'s approach to elementary particle physics) the
\emph{finest} invariants which exhibit \emph{anomalies} under small
perturbations. Correspondingly, the determinant and its amplitude are the most
subtle and the most fascinating objects of our study. They are much more
difficult to comprehend than the $\eta$-invariant \ref{sss:The Zeta}; and the
$\eta$-invariant is much more difficult to comprehend than the index.

\subsection{Pasting Formulas for the Eta--Invariant - Outlines}

\bigskip

In the rest of this review, i.e. over the next 33 pages, we shall prove a
strikingly simple (to state) additivity property of the $\eta$-invariant. We
fix the assumptions and the notation.

\subsubsection{Assumptions and Notation}

\noindent(a) Let $M$ be an odd--dimensional closed \emph{partitioned}
Riemannian manifold $M=M_{1}\cup_{\Sigma}M_{2}$ with $M_{1},M_{2}$ compact
manifolds with common boundary {$\Sigma$}. Let $S$ be a bundle of Clifford
modules over $M$.

\noindent(b) To begin with we assume that $\mathcal{D}$ is a \emph{compatible}
(= true) Dirac operator over $M$. Thus, in particular, $\mathcal{D}$ is
symmetric and has a unique self-adjoint extension in $L^{2}(M;S)$.

\noindent(c) We assume that there exists a bicollared cylindrical neighborhood
{(\emph{a neck}}) $N\simeq(-1,1)\times\Sigma$ of the separating hypersurface
{$\Sigma$}${,}$ such that the Riemannian structure on $M$ and the Hermitian
structure on $S$ are product in $N$; i.e., they do not depend on the normal
coordinate $u$, when restricted to {$\Sigma$}$_{u}=\{u\}\times\Sigma$. Our
convention for the orientation of the coordinate $u$ is that it runs from
$M_{1}$ to $M_{2}$; i.e., $M_{1}\cap N=(-1,0]\times\Sigma$ and $N\cap
M_{2}=[0,1)\times\Sigma$. Then the operator $\mathcal{D}$ takes the following
form on $N$:
\begin{equation}
\mathcal{D}|_{N}=\sigma(\partial_{u}+B), \label{e:product_form}%
\end{equation}
where the principal symbol in $u$--direction $\sigma:S|_{\Sigma}\rightarrow
S|_{\Sigma}$ is a unitary bundle isomorphism (Clifford multiplication by the
normal vector $du$) and the tangential operator $B:C^{\infty}(${$\Sigma$%
}$;S|_{\Sigma})\rightarrow C^{\infty}(${$\Sigma$}$;S|_{\Sigma})$ is the
corresponding Dirac operator on {$\Sigma$}. Note that $\sigma$ and $B$ do not
depend on the normal coordinate $u$ in $N$ and they satisfy the following
identities
\begin{equation}
\sigma^{2}=-\operatorname{I}\ ,\ \sigma^{\ast}=-\sigma\ ,\ \sigma\cdot
B=-B\cdot\sigma\ ,\ B^{\ast}=B. \label{e:tangential_identities}%
\end{equation}
Hence, $\sigma$ is a skew-adjoint involution and $S$, the bundle of spinors,
decomposes in $N$ into $\pm i$--eigenspaces of $\sigma$, $S|_{N}=S^{+}\oplus
S^{-}$. It follows that \eqref{e:product_form} leads to the following
representation of the operator $\mathcal{D}$ in $N$
\[
\mathcal{D}|_{N}=
\begin{pmatrix}
i & 0\\
0 & -i
\end{pmatrix}
\cdot\left(  \partial_{u}+
\begin{pmatrix}
0 & B_{-}=B_{+}^{\ast}\\
B_{+} & 0
\end{pmatrix}
\right)  ,
\]
where $B_{+}:C^{\infty}(${$\Sigma$}$;S^{+})\rightarrow C^{\infty}(${$\Sigma$%
}$;S^{-})$ maps the spinors of positive chirality into the spinors of negative chirality.

\noindent(d) To begin with we consider only the case of $\operatorname{Ker}
B=\{0\}$. That implies that $B$ is an \emph{invertible} operator. More
precisely, there exists a pseudo-differential elliptic operator $L$ of order
$-1$ such that $BL=\operatorname{I}_{S}=LB$ (see, for instance, \cite{BoWo93},
Proposition 9.5).

\noindent(e) For real $R>0$ we study the closed \emph{stretched} manifold
$M^{R}$ which we obtain from $M$ by inserting a cylinder of length $2R$, i.e.
replacing the collar $N$ by the cylinder $(-2R-1,+1)\times\Sigma$
\[
M^{R}=M_{1}\cup([-2R,0]\times\Sigma)\cup M_{2}\,.
\]
We extend the bundle $S$ to the stretched manifold $M^{R}$ in a natural way.
The extended bundle will be also denoted by $S$. The Riemannian structure on
$M$ and the Hermitian structure on $S$ are product on $N$. Hence we can extend
them to smooth metrics on $M^{R}$ in a natural way and, at the end, we can
extend the operator $\mathcal{D}$ to an operator $\mathcal{D}^{R}$ on $M^{R}$
by using formula \eqref{e:product_form}. Then $M^{R}$ splits into two
manifolds with boundary: $M^{R}=M_{1}^{R}\cup M_{2}^{R}$ with $M_{1}^{R}%
=M_{1}\cup\left(  (-2R,-R]\times\Sigma\right)  $, $M_{2}^{R}=\left(
[-R,0)\times\Sigma\right)  \cup M_{2}$, and $\partial M_{1}=\partial M_{2}%
^{R}=\{-R\}\times\Sigma$. Consequently, the operator $\mathcal{D}^{R}$ splits
into $\mathcal{D}^{R}=\mathcal{D}_{1}^{R}\cup\mathcal{D}_{2}^{R}$. We shall
impose spectral boundary conditions to obtain self-adjoint operators
$\mathcal{D}_{1,P_{<}}$, $\mathcal{D}_{1,P_{<}}^{R}$, $\mathcal{D}_{2,P_{>}}$,
and $\mathcal{D}_{2,P_{>}}^{R}$ in the corresponding $L^{2}$ spaces on the
parts (see \eqref{e:domain_m_2}).

\noindent(f) We also introduce the complete, noncompact Riemannian manifold
\emph{with cylindrical end}
\[
M_{2}^{\infty}:=\left(  (-\infty,0]\times\Sigma\right)  \cup M_{2}%
\]
by gluing the half--cylinder $(-\infty,0]\times\Sigma$ to the boundary
{$\Sigma$} of $M_{2}$. Clearly, the Dirac operator $\mathcal{D}$ extends also
to $C^{\infty}(M_{2}^{\infty},S)$.

\begin{remark}
Our presentation is somewhat simplified by our assumption (b) that
$\mathcal{D}$ is compatible and assumption (d) that the tangential operator
$B$ is invertible. Both assumptions can be lifted. This is done in the
literature; see \cite{Wo95} and \cite{Wo99}.
\end{remark}

We recall the following ideas in the big scheme from Section
\ref{ss:SymplecticGeometryOfCauchyDataSpaces} of this review. Let $P_{>}$
(respectively $P_{<}$) denote the spectral projection of $B$ onto the subspace
of $L^{2}(${$\Sigma$}$;S|_{\Sigma})$ spanned by the eigensections
corresponding to the positive (respectively negative) eigenvalues. Then
$P_{>}$ is a self-adjoint elliptic boundary condition for the operator
$\mathcal{D}_{2}=\mathcal{D}|_{M_{2}}$ (see \cite{BoWo93}, Proposition 20.3).
This means that the operator $\mathcal{D}_{2,P_{>}}$ defined by
\begin{equation}%
\begin{array}
[c]{lll}%
\mathcal{D}_{2,P_{>}} & =\mathcal{D}|_{M_{2}} & \\
\operatorname{Dom}(\mathcal{D}_{2,P_{>}}) & =\{s\in H^{1}(M_{2};S|_{M_{2}%
})\mid P_{>}(s|_{\Sigma})=0\} &
\end{array}
\label{e:domain_m_2}%
\end{equation}
is an unbounded self-adjoint operator in $L^{2}(M_{2};S|_{M_{2}})$ with
compact resolvent. In particular,
\[
\mathcal{D}_{2,P_{>}}:\operatorname{Dom}(\mathcal{D}_{2,P_{>}})\rightarrow
L^{2}(M_{2};S|_{M_{2}})
\]
is a Fredholm operator with discrete real spectrum and the kernel of
$\mathcal{D}_{2,P_{>}}$ consists of smooth sections of $S|_{M_{2}}$.

As mentioned before, the eta function of $\mathcal{D}_{2,P_{>}}$ is well
defined and enjoys all properties of the eta function of the Dirac operator
defined on a closed manifold. In particular, $\eta_{\mathcal{D}_{2,P_{>}}}%
(0)$, the eta invariant of $\mathcal{D}_{2,P_{>}}$, is well defined.
Similarly, $P_{<}$ is a self-adjoint boundary condition for the operator
$\mathcal{D}|_{M_{1}}$, and we define the operator $\mathcal{D}_{1,P_{<}}$
using a formula corresponding to \eqref{e:domain_m_2}. \emph{To keep track of
the various manifolds, operators, and integral kernels we refer to the
following table where we have collected the major notations}.
\[%
\begin{tabular}
[c]{|c|c|c|}\hline
$\text{manifolds}$ & $\text{operators}$ & $\text{integral kernels}%
$\\\hline\hline
$M=M_{1}\cup_{\Sigma}M_{2}%
\genfrac{}{}{0pt}{}{\mathstrut}{\mathstrut}%
$ & $e^{-t\mathcal{D}^{2}},\ \mathcal{D}e^{-t\mathcal{D}^{2}}$ &
$\mathcal{E}(t;x,x^{\prime})$\\\hline
$M^{R}=M_{1}^{R}\cup_{\Sigma}M_{2}^{R}%
\genfrac{}{}{0pt}{}{\mathstrut}{\mathstrut}%
$ & $e^{-t(\mathcal{D}^{R})^{2}},\ \mathcal{D}^{R}e^{-t(\mathcal{D}^{R})^{2}}$%
& $\mathcal{E}^{R}(t;x,x^{\prime})$\\\hline
$M_{2}%
\genfrac{}{}{0pt}{}{\mathstrut}{\mathstrut}%
$ & $e^{-t\mathcal{D}_{2,P_{>}}^{\ 2}},\ \mathcal{D}_{2}e^{-t\mathcal{D}%
_{2,P_{>}}^{\ 2}}$ & $\mathcal{E}_{2}(t;x,x^{\prime})$\\\hline
$M_{2}^{R}=\left(  [-R,0]\times\Sigma\right)  \cup M_{2}%
\genfrac{}{}{0pt}{}{\mathstrut}{\mathstrut}%
$ & $e^{-t(\mathcal{D}_{2,P_{>}}^{R})^{2}},\ \mathcal{D}_{2}^{R}%
e^{-t(\mathcal{D}_{2,P_{>}}^{R})^{2}}$ & $\mathcal{E}_{2}^{R}(t;x,x^{\prime}%
)$\\\hline
$M_{2}^{\infty}=\left(  (-\infty,0]\times\Sigma\right)  \cup M_{2}%
\genfrac{}{}{0pt}{}{\mathstrut}{\mathstrut}%
$ & $e^{-t(\mathcal{D}_{2}^{\infty})^{2}},\ \mathcal{D}_{2}^{\infty
}e^{-t(\mathcal{D}_{2}^{\infty})^{2}}$ & $\mathcal{E}_{2}^{\infty
}(t;x,x^{\prime})$\\\hline
$\Sigma_{\operatorname{cyl}}^{\infty}=(-\infty,+\infty)\times\Sigma%
\genfrac{}{}{0pt}{}{\mathstrut}{\mathstrut}%
$ & $e^{-tD_{\operatorname{cyl}}^{2}},\,D_{\operatorname{cyl}}%
e^{-tD_{\operatorname{cyl}}^{2}}$ & $\mathcal{E}_{\operatorname{cyl}%
}(t;x,x^{\prime})$\\\hline
$\Sigma_{\operatorname{cyl}/2}^{\infty}=[0,+\infty)\times\Sigma%
\genfrac{}{}{0pt}{}{\mathstrut}{\mathstrut}%
$ & $e^{-tD_{\operatorname{aps}}^{2}},\ D_{\operatorname{aps}}%
e^{-tD_{\operatorname{aps}}^{2}}$ & $\mathcal{E}_{\operatorname{aps}%
}(t;x,x^{\prime})$\\\hline
\end{tabular}
\]
In addition, on $M_{2}^{R}$ we have the operator $Q_{2}^{R}(t)$ with kernel
$Q_{2}^{R}(t;x,x^{\prime})$ and
\[
C^{R}(t)=\left(  (\mathcal{D}_{2,P_{>}}^{R})^{2}+\tfrac{d}{dt}\right)
Q_{2}^{R}(t)\text{ with kernel }C^{R}(t;x,x^{\prime})\text{.}%
\]

\subsubsection{The Gluing Formulas}

The most basic results for pasting $\eta$ are the following theorem on the
adiabatic limits of the $\eta$ invariants and its additivity corollary:

\begin{theorem}
\label{t:adiabatic_eta} Attaching a cylinder of length $R>0$ at the boundary
of the manifold $M_{2}$, we can approximate the eta invariant of the spectral
boundary condition on the prolonged manifold $M_{2}^{R}$ by the corresponding
integral of the `local' eta function of the closed stretched manifold $M^{R}%
$:
\[
\lim_{R\rightarrow\infty}\left\{  \eta_{\mathcal{D}_{2,P_{>}}^{R}}%
(0)-\int_{M_{2}^{R}}\eta_{\mathcal{D}^{R}}(0;x)\,dx\right\}  \equiv
0\operatorname{mod}\mathbb{Z}.
\]
\end{theorem}

\bigskip

\begin{corollary}
\label{c:0.1} $\eta_{\mathcal{D}}(0)\equiv\eta_{\mathcal{D}_{1,P_{<}}}%
(0)+\eta_{\mathcal{D}_{2,P_{>}}}(0)\operatorname{mod}\mathbb{Z}.$
\end{corollary}

\begin{remark}
(a) With hindsight, it is not surprising that modulo the integers the
preceding additivity formula for the $\eta$-invariant on a partitioned
manifold is precise. An intuitive argument runs as follows. ``Almost all''
eigensections and eigenvalues of the operator $\mathcal{D}$ on the closed
partitioned manifold $M=M_{1}\cup M_{2}$ can be traced back, either to
eigensections $\psi_{1,k}$ and eigenvalues $\mu_{1,k}$ of the spectral
boundary problem $\mathcal{D}_{1,P_{<}}$ on the part $M_{1},$ or to
eigensections $\psi_{2,\ell}$ and eigenvalues $\mu_{2,\ell}$ of the spectral
boundary problem $\mathcal{D}_{2,P_{>}}$ on the part $M_{2}$. While we have no
explicit exact correspondence, due to the product form of the Dirac operator
in a neighborhood of the separating hypersurface, eigensections on one part
$M_{1}$ or $M_{2}$ of the manifold $M$ can be extended to smooth sections on
the whole of $M$. These are not true eigensections of $\mathcal{D}$, but they
have a relative error which is rapidly decreasing as $R\rightarrow\infty$ when
we attach cylinders of length $R$ to the part manifolds or, equivalently,
insert a cylinder of length $2R$ in $M$. There is also a residual set
$\{\mu_{0,j}\}$ of eigenvalues of $\mathcal{D}$ which can neither be traced
back to eigenvalues of $\mathcal{D}_{1,P_{<}}$ nor to those of $\mathcal{D}%
_{2,P_{>}}$. These eigenvalues can, however, be traced back to the kernel of
the Dirac operators $\mathcal{D}_{1}^{\infty}$ and $\mathcal{D}_{2}^{\infty}$
on the part manifolds $M_{1}^{\infty}$ and $M_{2}^{\infty}$ with cylindrical
ends. Because of Fredholm properties the residual set is finite and, hence (as
noticed in Section \ref{ss:SpectralInvariants}) can be discarded for
calculating the eta invariant modulo $\mathbb{Z}$.\newline \ \ \ Therefore, no
$R$ (i.e., no prolongation of the bicollar neighborhood $N$) enters the
\emph{formula}. Nevertheless, our \emph{arguments} rely on an adiabatic
argument to separate the spectrum of $\mathcal{D}$ into its three parts
\begin{equation}
\operatorname{spec}\mathcal{D}\sim\{\mu_{0,j}\}\cup\{\mu_{1,k}\}\cup
\{\mu_{2,\ell}\}. \label{e:spec_strain}%
\end{equation}
For the most part, however, we need not make all arguments explicit on the
level of the single eigenvalue. It suffices to work on the level of the eta
invariant for the following reason. Unlike the index, the eta invariant cannot
be described by a local formula, as explained in Section
\ref{ss:SpectralInvariants}. Nevertheless, it can be described by an integral
over the manifold. The integrand, however, is not defined in local terms
solely. In particular, when writing the eta function in integral form and
decomposing the $\eta$ integral
\[
\eta_{\mathcal{D}}(s)=\int_{M}\eta_{\mathcal{D}}(s;x,x)\,dx=\int_{M_{1}}%
\eta_{\mathcal{D}}(s;x_{1},x_{1})\,dx_{1}+\int_{M_{2}}\eta_{\mathcal{D}%
}(s;x_{2},x_{2})\,dx_{2}%
\]
there is no geometrical interpretation of the integrals on the right over the
two parts of the manifold. This is very unfortunate. But for sufficiently
large $R$, the integrals become intelligible and can be read as the $\eta$
invariants of $\mathcal{D}_{1,P_{<}}^{R}$ and $\mathcal{D}_{2,P_{>}}^{R}$.
That is the meaning of the adiabatic limit.\newline (b) \label{r:finis}
Theorem \ref{t:0.5} can be generalized to larger classes of boundary
conditions by variational argument (see \cite{LeWo96}) yielding the general
gluing formula (in $\mathbb{R}/\mathbb{Z}$)
\[
\eta_{\mathcal{D}}(0)=\eta_{\mathcal{D}_{1,P_{1}}}(0)+\eta_{\mathcal{D}%
_{2,P_{2}}}(0)+\eta_{P_{1},\operatorname{I}-P_{2}}^{N}(0),
\]
where $\eta_{P_{1},\operatorname{I}-P_{2}}^{N}(0)$ denotes the $\eta
$--invariant on the cylinder $N$, see also \cite{BrLe99}, \cite{DaFr94},
\cite{MaMe95}, and the recent review \cite{PaWo02e}. G. Grubb's new result of
\cite{Gr02}, mentioned in our Introduction, may open an alternative route for
proving the general gluing formula. The integer jump was calculated in
\cite{KiLe00} yielding (among other formulas)
\[
\eta_{\mathcal{D}}(0)=\eta_{\mathcal{D}_{1,\operatorname{I}-P}}(0)+\eta
_{\mathcal{D}_{2,P}}(0)+2\operatorname{SF}\{\mathcal{D}_{1,\operatorname{I}%
-P_{t}}\}+2\operatorname{SF}\{\mathcal{D}_{2,P_{t}}\},
\]
where $\{P_{t}\}$ is a smooth curve in the Grassmannian from $P$ to the
Calder\'{o}n projection $\mathcal{P}(\mathcal{D}_{2})$. \newline (c) An
interesting feature of \cite{KiLe00} is that the pasting formula is derived
from the Scott--Wojciechowski Comparison Formula (Theorem \ref{t:scott-woj}).
This is quite analogous to the existence of two completely different proofs of
the pasting formula for the index (Theorem \ref{t:non-add}), where also the
derivation from the boundary reduction formula is much, much shorter than
arguing via the Atiyah-Singer Index Theorem and the Atiyah-Patodi-Singer Index
Theorem plus the Agranovic-Dynin Formula. However, for this review we prefer
the long way because of the many interesting insights about the gluing on the
eigenvalue level which can be gained.
\end{remark}

\bigskip

\subsubsection{Plan of the Proof}

Let $\mathcal{E}_{2}^{R}(t)$ denote the integral kernel of the operator
$\mathcal{D}_{2}^{R}e^{-t(\mathcal{D}_{2,P_{>}}^{R})^{2}}$ defined on the
manifold $M_{2}^{R}=([-R,0]\times\Sigma)\cup M_{2}$. Without proof, we
mentioned before that the eta invariant of the self-adjoint operator
$\mathcal{D}_{2,P_{>}}^{R}$ is well defined and we have
\begin{align}
\eta_{\mathcal{D}_{2,P_{>}}^{R}}(0)  &  =\frac{1}{\sqrt{\pi}}\int_{0}^{\infty
}\frac{dt}{\sqrt{t}}\int_{M_{2}^{R}}\operatorname{Tr}\mathcal{E}_{2}%
^{R}(t;x,x)\,dx\nonumber\\
&  =\frac{1}{\sqrt{\pi}}\int_{0}^{\sqrt{R}}\frac{dt}{\sqrt{t}}\int_{M_{2}^{R}%
}\operatorname{Tr}\mathcal{E}_{2}^{R}(t;x,x)\,dx\label{e:0_sqrt_r}\\
&  \qquad+\frac{1}{\sqrt{\pi}}\int_{\sqrt{R}}^{\infty}\frac{dt}{\sqrt{t}}%
\int_{M_{2}^{R}}\operatorname{Tr}\mathcal{E}_{2}^{R}(t;x,x)\,dx.
\label{e:sqrt_r_infty}%
\end{align}
We first deal with the integral of \eqref{e:0_sqrt_r} and show that it splits
into an interior contribution and a cylinder contribution as $R\rightarrow
\infty$. This will be done in Subsection \ref{ss:small-t-chopped} by first in
paragraph \ref{sss:duhamel} applying the Duhamel method which we introduced
before in the proof of Theorem \ref{tAPS} (pp. \pageref{tAPS}ff). By Lemma
\ref{l:crucial_estimate}, Lemma \ref{l:error_r_estimate}, Proposition
\ref{p:error est}, and Corollary \ref{c:7.5}, we can replace the heat kernel
$\mathcal{E}_{2}^{R}$ of the Atiyah-Patodi-Singer boundary problem on the
prolonged manifold $M_{2}^{R}$ with boundary $\Sigma$ by an artificially glued
integral kernel $Q_{2}^{R}$ which, near the boundary, is equal to the heat
kernel of the APS problem on the half-infinite cylinder and in the interior
equal to the heat kernel on the stretched closed manifold. We can do it in
such a way that the original small-$t$ integral (i.e., the integral from $0$
to $\sqrt{R}$ in (\ref{e:0_sqrt_r})) can be approximated by the new integral
up to an error of order $O(e^{-cR})$. Then we shall show in Lemma
\ref{l:trace_even_odd} and Lemma \ref{lem:3.1} (paragraph \ref{sss:discarding}%
) that the cylinder contribution is traceless. More precisely, we obtain that
the trace $\operatorname{Tr}Q_{2}^{R}(T;x,x)$ can be replaced pointwise (for
$x\in M_{2}^{R}$) by the trace $\operatorname{Tr}\mathcal{E}^{R}(t;x,x)$ of
the integral kernel of the operator $\mathcal{D}^{R}e^{-t(\mathcal{D}^{R}%
)^{2}}$ which is defined on the stretched manifold $M^{R}$. Consequently, the
small-$t$ \textit{chopped} $\eta$-invariant of the closed stretched manifold
$M^{R}$ coincides with the sum of the small-$t$ chopped $\eta$-invariants of
the APS problems on the two prolonged manifolds $M_{1}^{R}$ and $M_{2}^{R}$ up
to $O\left(  e^{-cR}\right)  $.

Then we will show that the integral of \eqref{e:sqrt_r_infty} vanishes as
$R\rightarrow\infty$. This is in Lemma \ref{l:lemma.7.1} (p.
\pageref{l:lemma.7.1}f) a direct consequence of Theorem \ref{t:6.1} which
states that the eigenvalues of $\mathcal{D}_{2,P_{>}}^{R}$ are uniformly
bounded away from 0. Theorem \ref{t:6.1} is of independent interest. The proof
is a long story stretching over Subsections \ref{ss:large-t-aps} and
\ref{ss:LowestEigenv}. First, in paragraph \ref{sss:cylinder} we shall
consider the operator $D_{\operatorname{cyl}}$ on the infinite cylinder
$\Sigma_{\operatorname{cyl}}^{\infty}$ in Definition \ref{d:sobolev} and Lemma
\ref{l:sobolev_d}. We obtain that $D_{\operatorname{cyl}}$ has no eigenvalues
in the interval $\left(  -\lambda_{1},\lambda_{1}\right)  $ where $\lambda
_{1}$ denotes the smallest positive eigenvalue of the tangential operator $B$
on the manifold $\Sigma$. Next, in paragraph \ref{sss:cylindrical-ends} we
investigate the operator $\mathcal{D}_{2}^{\infty}$ on the manifold
$M_{2}^{\infty}$ with infinite cylindrical end in Lemma
\ref{l:infty_eigen_estimate}, Lemma \ref{l:general_ess_self_adjoint}, Lemma
\ref{l:ess_self_adjoint}, and Proposition \ref{p:lemma_6.2}. Proposition
\ref{p:lemma_6.2} states that $\mathcal{D}_{2}^{\infty}$ has only finitely
many eigenvalues in the aforementioned interval $\left(  -\lambda_{1}%
,\lambda_{1}\right)  $, each of finite multiplicity. Its proof is somewhat
delicate and involves Lemma \ref{l:trace_compact} and Corollary
\ref{c:compact}. Then, in Subsection \ref{ss:LowestEigenv} the proof of
Theorem \ref{t:6.1} follows with Lemma \ref{l:eigenext}.

By then we will have that the sum of the $\eta$-invariants of the APS boundary
problems on the prolonged manifolds $M_{1}^{R}$ and $M_{2}^{R}$ can be
replaced by the small-$t$ chopped $\eta$-invariant on the stretched closed
manifold $M^{R}$ with an error exponentially vanishing as $R\rightarrow\infty
$. We then would like to repeat the preceding chain of arguments around
Theorem \ref{t:6.1} to show that the large-$t$ chopped $\eta$-invariant of the
operator $\mathcal{D}^{R}$ on $M^{R}$ also vanishes. However, this can be done
only in $\mathbb{R}/\mathbb{Z}$ because in general the eigenvalues of
$\mathcal{D}^{R}$ are not bounded away from 0. So we shall present a different
chain of arguments in Subsection \ref{ss:Spectrum on Closed Stretched Mfd}.
The main technical result is Theorem \ref{t:0.3}, once again of independent
interest. It describes a partition of the eigenvalues into two subsets, the
exponentially small ones and the eigenvalues bounded away from 0. The key for
our arguments is a gluing construction (Definition \ref{d:VR}). Then we first
show by Lemma \ref{l:R0}, Lemma \ref{l:1.3} and Proposition \ref{p:specproj}
that $\mathcal{D}^{R}$ has at least $q:=\dim\operatorname{Ker}\mathcal{D}%
_{1}^{\infty}+\dim\operatorname{Ker}\mathcal{D}_{2}^{\infty}$ exponentially
small eigenvalues belonging to eigensections which we can approximate by
pasting together $L^{2}$ solutions. Then we shall show in Lemma \ref{l:2.1},
Theorem \ref{t:2.2}, Lemma \ref{l:2.3} and Proposition \ref{p:2.4} that this
makes the list of eigenvalues approaching 0 as $R\rightarrow\infty$ complete.

It follows in Lemma \ref{l:etaR} (at the beginning of the closing Subsection
\ref{ss:LowestEigenv}) that the large-$t$ chopped $\eta$-invariant on $M^{R}$
vanishes asymptotically up to an integer error. This establishes Theorem
\ref{t:adiabatic_eta}. To arrive at Corollary \ref{c:0.1}, we must get rid of
the adiabatic limit. This is a simple consequence of the locality of the
derivative in $R$-direction of the $\eta$-invariants on the stretched part
manifolds with APS boundary condition and on the stretched closed manifold
(Proposition \ref{p:muller}).

\medskip

\subsection{The Adiabatic Additivity of the Small-$t$ Chopped $\eta
$-invariant\label{ss:small-t-chopped}}

\bigskip

\subsubsection{Applying Duhamel's Method to the Small-$t$ Chopped $\eta
$-invariant\label{sss:duhamel}}

The simplest construction of a parametrix for $\mathcal{E}_{2}^{R}(t)$ (i.e.,
of an approximate heat kernel) is the following: we glue the kernel
$\mathcal{E}$ of the operator $\mathcal{D}e^{-t\mathcal{D}^{2}}$ (given on the
whole, closed manifold $M$) and the kernel $\mathcal{E}_{\operatorname{aps}%
}^{\infty}$ of the $L^{2}$ extension of the operator $\sigma(\partial
_{u}+B)e^{-t(\sigma(\partial_{u}+B))^{2}}$, given on the semi-infinite
cylinder $[-R,\infty)\times\Sigma$ and subject to the Atiyah--Patodi--Singer
boundary condition at the end $u=-R$. In that construction the gluing happens
on the neck $N=[0,1)\times\Sigma$ with suitable cutoff functions.

Locally, the heat kernel is always of the form $(4\pi t)^{-m/2}e^{c_{1}%
t}e^{-\lvert x-x^{\prime}\rvert^{2}/4t}$\thinspace. By Duhamel's Principle we
get after gluing a similar global result for the kernel $\operatorname{e}%
_{2}^{R}(t;x,x^{\prime})$ of the operator $e^{-t(\mathcal{D}_{2,P_{>}}%
^{R})^{2}}$ and, putting a factor $t^{-1/2}$ in front, for the kernel of the
combined operator $\mathcal{D}e^{-t(\mathcal{D}_{2,P_{>}}^{R})^{2}}$ (e.g.,
see Gilkey \cite{Gi95}, Lemma 1.9.1). That yields two crucial estimates:

\begin{lemma}
\label{l:crucial_estimate} There exist positive reals $c_{1}$, $c_{2}$, and
$c_{3}$ which do not depend on $R$, such that for all $x,x^{\prime}\in
M_{2}^{R}$ and any $t>0$ and $R>0$,
\begin{gather}
\lvert\operatorname{e}_{2}^{R}(t;x,x^{\prime})\rvert\leq c_{1}\cdot
t^{-\frac{m}{2}}\cdot e^{c_{2}t}\cdot e^{-c_{3}\frac{d^{2}(x,x^{\prime})}{t}%
},\label{e:4.3.a}\\
\lvert\mathcal{E}_{2}^{R}(t;x,x^{\prime})\rvert\leq c_{1}t^{-\frac{1+m}{2}%
}\cdot e^{c_{2}t}\cdot e^{-c_{3}\frac{d^{2}(x,x^{\prime})}{t}}.
\label{e:4.3.b}%
\end{gather}
Here $d(x,x^{\prime})$ denotes the geodesic distance.
\end{lemma}

Notice that exactly the same type of estimate is also valid for the kernel
$\mathcal{E}^{R}(t;x,x^{\prime})$ on the stretched closed manifold $M^{R}$ and
for the kernel $\mathcal{E}_{\operatorname{aps}}^{\infty}(t;x,x^{\prime})$ on
the infinite cylinder. For details see also \cite{BoWo93}, Theorem 22.14.
There, however, the term $e^{c_{2}t}$ was suppressed in the final formula
because the emphasis was on small time asymptotics.

As mentioned before, as $R\rightarrow\infty$, we want to separate the
contribution to the kernel $\mathcal{E}_{2}^{R}$ which comes from the cylinder
and the contribution from the interior by a gluing process. Unfortunately, the
inequality \ref{e:4.3.b} does not suffice to show that the contribution to the
eta invariant, more precisely to the integral \eqref{e:0_sqrt_r}, which comes
from the `error' term vanishes with $R\rightarrow\infty$. Therefore, we
introduce a different parametrix for the kernel $\mathcal{E}_{2}^{R}$.

Instead of gluing over the fixed neck $N=[0,1)\times\Sigma,$ we glue over a
segment $N^{R}$ of growing length of the attached cylinder, say $N^{R}%
:=(-\frac{4}{7}R,-\frac{3}{7}R)\times\Sigma$ (the reason for choosing these
ratios will be clear soon). Thus, we choose a smooth partition of unity
$\{\chi_{\operatorname{aps}},\chi_{\operatorname{int}}\}$ on $M_{2}^{R}$
suitable for the covering $\{U_{\operatorname{aps}},U_{\operatorname{int}}\}$
with $U_{\operatorname{aps}}:=[-R,-\frac{3}{7}R)\times\Sigma$ and
$U_{\operatorname{int}}:=\left(  (-\frac{4}{7}R,0]\times\Sigma\right)  \cup
M_{2}$, hence $U_{\operatorname{aps}}\cap U_{\operatorname{int}}=N^{R}$.
Moreover, we choose nonnegative smooth cutoff functions $\{\psi
_{\operatorname{aps}},\psi_{\operatorname{int}}\}$ such that
\begin{align*}
\  &  \psi_{j}\equiv1\;\text{on $\{x\in M_{2}^{R}\mid\operatorname*{dist}(x,$%
}\operatorname*{supp}\text{$\chi_{j})<\frac{1}{7}R\}$ and}\\
&  \psi_{j}\equiv0\;\text{on $\{x\in M_{2}^{R}\mid\operatorname*{dist}%
(x,\operatorname*{supp}\chi_{j})\geq\frac{2}{7}R\}$}%
\end{align*}
for $j\in\{${$\operatorname{aps}$}$,${$\operatorname{int}$}$\}$. We notice
\begin{equation}
\operatorname*{dist}(\operatorname*{supp}\psi_{j}^{\prime}%
,\operatorname*{supp}\chi_{j})=\operatorname*{dist}(\operatorname*{supp}%
\psi_{j}^{\prime\prime},\operatorname*{supp}\chi_{j})\geq\frac{1}{7}R\,.
\label{e:distance_support}%
\end{equation}
Moreover, we may assume that
\[
\left|  {\frac{\partial^{k}\psi_{j}}{\partial u^{k}}}\right|  \leq c_{0}/R
\]
for all $k$, where $c_{0}$ is a certain positive constant.

For any parameter $t>0$, we define an operator $Q_{2}^{R}(t)$ on $C^{\infty
}(M_{2}^{R};S)$ with a smooth kernel, given by
\begin{equation}
Q_{2}^{R}(t;x,x^{\prime}):=\psi_{\operatorname{aps}}(x)\mathcal{E}%
_{{\operatorname{aps}}}^{\infty}(t;x,x^{\prime})\chi_{\operatorname{aps}%
}(x^{\prime})+\psi_{\operatorname{int}}(x)\mathcal{E}^{R}(t;x,x^{\prime}%
)\chi_{\operatorname{int}}(x^{\prime}). \label{e:parametrix_r}%
\end{equation}
Recall that $\mathcal{E}^{R}$ denotes the kernel of the operator
$\mathcal{D}^{R}e^{-t(\mathcal{D}^{R})^{2}}$, given on the stretched closed
manifold $M^{R}$. Notice that, by construction, $Q_{2}^{R}(t)$ maps
$L^{2}(M_{2}^{R};S) $ into the domain of the operator $\mathcal{D}_{2,P_{>}%
}^{R}$.

Then, for $x^{\prime}\in U_{\operatorname{aps}}$ with $\chi
_{\operatorname{aps}}(x^{\prime})=1$, we have by definition:
\begin{equation}
Q_{2}^{R}(t;x,x^{\prime})=\left\{
\begin{array}
[c]{ll}%
\mathcal{E}_{\operatorname{aps}}^{\infty}(t;x,x^{\prime}) & \text{if
$d(x,\operatorname*{supp}\chi_{\operatorname{aps}})<\frac{1}{7}R$, and}\\
0 & \text{if $d(x,\operatorname*{supp}\chi_{\operatorname{aps}})\geq\frac
{2}{7}R$}.
\end{array}
\right.  \label{e:error_vanishing_cyl}%
\end{equation}
Correspondingly, we have for $x^{\prime}\in U_{\operatorname{int}}$ with
$\chi_{\operatorname{int}}(x^{\prime})=1$,
\begin{equation}
Q_{2}^{R}(t;x,x^{\prime})=\left\{
\begin{array}
[c]{ll}%
\mathcal{E}^{R}(t;x,x^{\prime}) & \text{if $d(x,\operatorname*{supp}%
\chi_{\operatorname{int}})<\frac{1}{7}R$, and}\\
0 & \text{if $d(x,\operatorname*{supp}\chi_{\operatorname{int}})\geq\frac
{2}{7}R$}.
\end{array}
\right.  \label{e:error_vanishing_int}%
\end{equation}
For fixed $t>0$, we determine the difference between the precise kernel
$\mathcal{E}_{2}^{R}(t;x,x^{\prime})$ and the approximate one $Q_{2}%
^{R}(t;x,x^{\prime})$. Let $C^{R}(t)$ denote the operator $\left(
(\mathcal{D}_{2,P_{>}}^{R})^{2}+\tfrac{d}{dt}\right)  \circ Q_{2}^{R}(t)$ and
$C^{R}(t;x,x^{\prime})$ its kernel. By definition, we have $\left(
(\mathcal{D}_{2,P_{>}}^{R})^{2}+\tfrac{d}{dt}\right)  \circ\mathcal{E}_{2}%
^{R}(t)=0$. Thus, $C^{R}(t)$ `measures' the error we make when replacing the
precise kernel $\mathcal{E}_{2}^{R}(t;x,x^{\prime})$ by the glued, approximate one.

More precisely, we have by Duhamel's Formula
\[
\mathcal{E}_{2}^{R}(t;x,x^{\prime})-Q_{2}^{R}(t;x,x^{\prime})=-\int_{0}%
^{t}ds\int_{M_{2}^{R}}\/dz\,\mathcal{E}_{2}^{R}(s;x,z)C^{R}(t-s;z,x^{\prime})
\]
with
\begin{align*}
C^{R}  &  (t-s;z,x^{\prime})=\left(  (\mathcal{D}_{2,(z)}^{R})^{2}%
+\mathcal{D}\right)  Q_{2}^{R}(t-s;z,x^{\prime})\\
&  =\left(  (\mathcal{D}_{2\,(z)}^{R})^{2}-\tfrac{d}{ds}\right)  Q_{2}%
^{R}(t-s;z,x^{\prime})\\
&  =\psi_{\operatorname{aps}}^{\prime\prime}(z)\mathcal{E}_{\operatorname{aps}%
}^{R}(t-s;z,x^{\prime})\chi_{\operatorname{aps}}(x^{\prime})+2\psi
_{\operatorname{aps}}^{\prime}(z)\frac{\partial}{\partial u}(\mathcal{E}%
_{\operatorname{aps}}^{R}(t-s;z,x^{\prime}))\chi_{\operatorname{aps}%
}(x^{\prime})\\
&  \qquad\qquad\qquad+\psi_{\operatorname{aps}}(z)\underbrace{\left(
\mathcal{D}_{(z)}^{2}-\tfrac{d}{ds}\right)  \mathcal{E}_{\operatorname{aps}%
}^{R}(t-s;z,x^{\prime})}_{=0}\chi_{\operatorname{aps}}(x^{\prime})\\
&  \quad+\psi_{\operatorname{int}}^{\prime\prime}(z)\mathcal{E}^{R}%
(t-s;z,x^{\prime})\chi_{\operatorname{int}}(x^{\prime})+2\psi
_{\operatorname{int}}^{\prime}(z)\frac{\partial}{\partial u}(\mathcal{E}%
^{R}(t-s;z,x^{\prime}))\chi_{\operatorname{int}}(x^{\prime})\\
&  \qquad\qquad\qquad+\psi_{\operatorname{int}}(z)\underbrace{\left(
(\mathcal{D}_{(z)}^{R})^{2}-\tfrac{d}{ds}\right)  \mathcal{E}^{R}%
(t-s;z,x^{\prime})}_{=0}\chi_{\operatorname{int}}(x^{\prime}).
\end{align*}
Here, $\mathcal{D}_{(z)}$ denotes the operator $\mathcal{D}$ acting on the $z
$ variable; and in the partial derivative $\frac{\partial}{\partial u}$ the
letter $u$ denotes the normal coordinate of the variable $z$.

As stated in \eqref{e:distance_support}, the supports of $\chi_{j}$ and
$\psi_{j}^{\prime}$ (and, equally, $\psi_{j}^{\prime\prime}$) are disjoint and
separated from each other by a distance $R/7$ in the normal variable for
$j\in\{${$\operatorname{aps}$}$,${$\operatorname{int}$}$\}$. Then the error
term $C^{R}(t-s;z,x^{\prime})$ vanishes both for the distance in the normal
variable $d(z,x^{\prime})<R/7$ and, actually, whenever $z$ or $x^{\prime}$ are
outside the segment $[-\frac{6}{7}R,\frac{1}{7}R]\times\Sigma$.

Let $z$ and $x^{\prime}$ be on the cylinder and $\lvert u-v\rvert>R/7$ where
$u$ and $v$ denote their normal coordinates. We investigate the error term
$C^{R}(t-s;z,x^{\prime})$ which consists of six summands. Two of them vanish
as we have pointed out above. The remaining four summands involve the kernels
$\mathcal{E}_{\operatorname{aps}}^{\infty}(t-s;z,x^{\prime})$ on the infinite
cylinder $[-R,\infty)\times\Sigma$ and $\mathcal{E}^{R}(t-s;z,x^{\prime})$ on
the stretched closed manifold $M^{R}$. We shall use that both kernels can be
estimated according to inequality \eqref{e:4.3.b}.

We estimate the first summand
\begin{align*}
\lvert\psi_{\operatorname{aps}}^{\prime\prime}(z)\mathcal{E}%
_{\operatorname{aps}}^{\infty}(t-s;z,x^{\prime})\chi_{\operatorname{aps}%
}(x^{\prime})\rvert &  \leq\frac{c_{0}}{R}c_{1}(t-s)^{-\frac{1+m}{2}}%
e^{c_{2}t}e^{-c_{3}\frac{d^{2}(z,x^{\prime})}{t-s}}\\
&  \leq c_{1}^{\prime}e^{c_{2}^{\prime}t}e^{-c_{3}^{\prime}R^{2}/t}\,.
\end{align*}
Here we have used $t\geq s\geq0$ and
\[
(t-s)^{-(1+m)/2}e^{-c_{2}\frac{d^{2}(z,x^{\prime})}{(t-s)}}\leq ct^{-(1+m)/2}%
e^{-c_{2}\frac{d^{2}(z,x^{\prime})}{t}}\leq\tilde{c}e^{-c_{2}\frac
{d^{2}(z,x^{\prime})}{2t}}.
\]
Similarly we estimate the second summand
\begin{multline*}
2\left|  {\psi_{\operatorname{aps}}^{\prime}(z)\frac{\partial}{\partial
u}\mathcal{E}_{\operatorname{aps}}^{\infty}(t-s;z,x^{\prime})\chi
_{\operatorname{aps}}(x^{\prime})}\right| \\
\leq\frac{c_{0}}{R}c_{1}\frac{(t-s)^{-\frac{1+m}{2}}}{\sqrt{t}}e^{c_{2}%
t}e^{-c_{3}\frac{d^{2}(z,x^{\prime})}{t-s}}\leq c_{1}^{\prime}e^{c_{2}%
^{\prime}t}e^{-c_{3}^{\prime}R^{2}/t}\,,
\end{multline*}
where the factor $1/\sqrt{t}$ comes from the differentiation of the kernel as
explained before. The third and fourth summands, involving the kernel
$\mathcal{E}^{R}$ of the closed stretched manifold $M^{R}$, are treated in
exactly the same way. Altogether we have proved

\begin{lemma}
\label{l:error_r_estimate} The error kernel $C^{R}(t;u,v)$ vanishes for
$u\notin\lbrack-\frac{6}{7}R,-\frac{1}{7}R]$. Moreover, $C^{R}(t;u,v)$
vanishes whenever $\lvert u-v\rvert\leq R/7$. For arbitrary $x,x^{\prime}\in
M_{2}^{R}$ we have the estimate
\[
\lvert C^{R}(t;x,x^{\prime})\rvert\leq c_{1}e^{c_{2}t}e^{-c_{3}R^{2}/t}%
\]
with constants $c_{1},c_{2},c_{3}$ independent of $x,x^{\prime},t,R$.
\end{lemma}

We consider the pointwise error
\[
\mathcal{E}_{2}^{R}(t;x,x)-Q_{2}^{R}(t;x,x)=\int_{0}^{t}ds\int_{M_{2}^{R}%
}dz\text{ }\mathcal{E}_{2}^{R}(s;x,z)C^{R}(t-s;z,x).
\]
We obtain the following proposition as a consequence of the preceding lemma.

\begin{proposition}
\label{p:error est}For all $x\in M_{2}^{R}$ and all $t>0$ we have
\[
\operatorname{Tr}\mathcal{E}_{2}^{R}(t;x,x)-\operatorname{Tr}Q_{2}%
^{R}(t;x,x)=\operatorname{Tr}\left(  \mathcal{E}_{2}^{R}(t;x,x)-Q_{2}%
^{R}(t;x,x)\right)  .
\]
Moreover, there exist positive constants $c_{1},c_{2},c_{3}$, independent of
$R$, such that the `error' term satisfies the inequality
\[
\left|  {\mathcal{E}_{2}^{R}(t;x,x)-Q_{2}^{R}(t;x,x)}\right|  \leq c_{1}\cdot
e^{c_{2}t}\cdot e^{-c_{3}(R^{2}/t)}\,.
\]
\end{proposition}

\begin{proof}
We estimate the error term
\begin{align*}
\lvert{\mathcal{E}}_{2}^{R}(t;x,x)  &  -Q_{2}^{R}(t;x,x)\rvert\\
&  \leq\int_{0}^{t}ds\int_{M_{2}^{R}}dz\,\lvert{\mathcal{E}}_{2}%
^{R}(s;x,z)C^{R}(t-s;z,x)\rvert\\
&  \leq c_{1}^{2}e^{c_{2}t}\cdot\int_{0}^{t}ds\,\int_{M_{2}^{R}}dz\,\left\{
s^{-\frac{d+1}{2}}\cdot e^{-c_{3}\frac{d^{2}(x,z)}{s}}\right\}  \cdot
e^{-c_{3}\frac{d^{2}(x,z)}{t-s}}\\
&  \leq c_{1}^{2}e^{c_{2}t}\cdot\int_{0}^{t}ds\,\int_{\operatorname*{supp}%
_{z}C^{R}(t-s;z,x)}dz\,e^{-c_{4}\frac{d^{2}(x,z)}{s}}\cdot e^{-c_{3}%
\frac{d^{2}(x,z)}{t-s}}\\
&  \leq c_{1}^{2}e^{c_{2}t}\cdot\int_{0}^{t}ds\int_{\operatorname*{supp}%
_{z}C^{R}(t-s;z,x)}dz\,e^{-c_{5}\frac{t\cdot R^{2}}{s(t-s)}}\\
&  \leq c_{1}^{2}e^{c_{2}t}\cdot cR\cdot\int_{0}^{t}ds\,e^{-c_{5}\frac{t\cdot
R^{2}}{s(t-s)}}\leq c_{1}^{2}e^{c_{2}t}\cdot2cR\cdot\int_{0}^{t/2}%
ds\,e^{-c_{5}\frac{t\cdot R^{2}}{s(t/2)}}\\
&  =c_{1}^{2}e^{c_{2}t}\cdot2cR\cdot\int_{0}^{t/2}ds\,e^{-2c_{5}\frac{R^{2}%
}{s}}\,.
\end{align*}
Here we have used that $\operatorname*{Vol}(\operatorname*{supp}_{z}%
C^{R}(t-s;z,x))\sim\operatorname*{Vol}(\Sigma)\cdot R$ according to Lemma
\ref{l:error_r_estimate}. We investigate the last integral.
\begin{multline*}
\int_{0}^{t}e^{-\frac{c}{s}}\,ds=-\int_{0}^{t}\frac{s^{2}}{c}\cdot
e^{-\frac{c}{s}}\cdot(-\frac{c}{s^{2}})\,ds\\
<-\int_{0}^{t}\frac{t^{2}}{c}\cdot e^{-\frac{c}{s}}\cdot(-\frac{c}{s^{2}%
})\,ds=-\frac{t^{2}}{c}\int_{\infty}^{\frac{c}{t}}e^{-r}\,dr=\frac{t^{2}}%
{c}e^{-\frac{c}{t}}\,.
\end{multline*}
Thus we have
\[
\lvert\mathcal{E}_{2}^{R}(t;x,x)-Q_{2}^{R}(t;x,x)\rvert\leq c_{1}^{2}%
e^{c_{2}t}\cdot2cR\cdot\frac{t^{2}}{c_{6}R^{2}}e^{-\frac{c_{6}R^{2}}{t}}\leq
c_{7}e^{c_{2}t}\cdot e^{-c_{8}(R^{2}/t)}\,.
\]
\end{proof}

The preceding proposition shows that, for $t$ smaller than $\sqrt{R}$, the
trace $\operatorname{Tr}\mathcal{E}_{2}^{R}(t;x,x)$ of the kernel of the
operator $\mathcal{D}_{2}^{R}e^{-t(\mathcal{D}_{2,P_{>}}^{R})^{2}}$ approaches
the trace $\operatorname{Tr}Q_{2}^{R}(t;x,x)$ of the approximative kernel
pointwise as $R\rightarrow\infty$. In particular, we have:

\begin{corollary}
\label{c:7.5} The following equality holds, as $R\rightarrow\infty$,
\begin{multline*}
\frac{1}{\sqrt{\pi}}\int_{0}^{\sqrt{R}}\frac{dt}{\sqrt{t}}\int_{M_{2}^{R}%
}\operatorname{Tr}\mathcal{E}_{2}^{R}(t;x,x)\,dx\\
=\frac{1}{\sqrt{\pi}}\int_{0}^{\sqrt{R}}\frac{dt}{\sqrt{t}}\int_{M_{2}^{R}%
}\operatorname{Tr}Q_{2}^{R}(t;x,x)\,dx+\operatorname{O}(e^{-cR}).
\end{multline*}
\end{corollary}

\begin{proof}
We have
\begin{align*}
\frac{1}{\sqrt{\pi}}\int_{0}^{\sqrt{R}}  &  \frac{dt}{\sqrt{t}}\int_{M_{2}%
^{R}}\operatorname{Tr}\mathcal{E}_{2}^{R}(t;x,x)\,dx\\
&  =\frac{1}{\sqrt{\pi}}\int_{0}^{\sqrt{R}}\frac{dt}{\sqrt{t}}\int_{M_{2}^{R}%
}\operatorname{Tr}Q_{2}^{R}(t;x,x)\,dx\\
&  \quad+\frac{1}{\sqrt{\pi}}\int_{0}^{\sqrt{R}}\frac{dt}{\sqrt{t}}\int
_{M_{2}^{R}}\operatorname{Tr}\left(  \mathcal{E}_{2}^{R}(t;x,x)-Q_{2}%
^{R}(t;x,x)\right)  \,dx,
\end{align*}
and we have to show that the second summand on the right side is
$\operatorname{O}(e^{-cR})$ as $R\rightarrow\infty$. We estimate
\begin{align*}
&  \left|  {\frac{1}{\sqrt{\pi}}\int_{0}^{\sqrt{R}}\frac{dt}{\sqrt{t}}%
\int_{M_{2}^{R}}\operatorname{Tr}}\left(  {\mathcal{E}_{2}^{R}(t;x,x)-Q_{2}%
^{R}(t;x,x)}\right)  {\,dx}\right| \\
&  \leq\frac{1}{\sqrt{\pi}}\int_{0}^{\sqrt{R}}\frac{dt}{\sqrt{t}}\int
_{M_{2}^{R}}\lvert\mathcal{E}_{2}^{R}(t;x,x)-Q_{2}^{R}(t;x,x)\rvert\,dx\\
&  \leq\frac{1}{\sqrt{\pi}}\int_{0}^{\sqrt{R}}\frac{dt}{\sqrt{t}}\int
_{M_{2}^{R}}c_{1}\cdot e^{c_{2}t}\cdot e^{-c_{3}(R^{2}/t)}\,dx\\
&  \leq\frac{c_{1}\operatorname*{Vol}(M_{2}^{R})}{\sqrt{\pi}}\int_{0}%
^{\sqrt{R}}\frac{e^{c_{2}t}\cdot e^{-c_{3}(R^{2}/t)}}{\sqrt{t}}\,dt\\
&  \leq c_{4}R\int_{0}^{\sqrt{R}}e^{c_{2}\sqrt{R}}\cdot e^{-c_{5}R^{3/2}%
}\,dt\leq c_{4}R^{3/2}\cdot e^{-c_{6}R}\leq c_{7}\cdot e^{-c_{8}R}\,.
\end{align*}
\bigskip
\end{proof}

\subsubsection{Discarding the Cylinder Contributions\label{sss:discarding}%
\bigskip}

Corollary \ref{c:7.5} shows that the essential part of the local eta function
of the spectral boundary condition on the half manifold with attached cylinder
of length $R$, i.e., the `small--time' integral from 0 to $\sqrt{R}$ can be
replaced, as $R\rightarrow\infty$, by the corresponding integral over the
trace $\operatorname{Tr}Q_{2}^{R}(t;x,x)$ of the approximate kernel,
constructed in \eqref{e:parametrix_r}. Now we show that $\operatorname{Tr}%
Q_{2}^{R}(t;x,x)$ can be replaced pointwise (for $x\in M_{2}^{R}$) by the
trace $\operatorname{Tr}\mathcal{E}^{R}(t;x,x)$ of the kernel of the operator
$\mathcal{D}^{R}e^{-t(\mathcal{D}^{R})^{2}}$ which is defined on the stretched
closed manifold $M^{R}$.

Consider the Dirac operator
\[
\sigma(\partial_{u}+B):C^{\infty}([0,\infty)\times\Sigma;S)\rightarrow
C^{\infty}([0,\infty)\times\Sigma;S),
\]
on the semi-infinite cylinder with the domain
\[
\{s\in C_{0}^{\infty}([0,\infty)\times\Sigma;S)\mid P_{>}(s|_{\{0\}\times
\Sigma})=0\}.
\]
It has a unique self-adjoint extension which we denote by
$D_{\operatorname{aps}}$. Recall that the integral kernel $\mathcal{E}%
_{\operatorname{aps}}^{\infty}$ of the operator $D_{\operatorname{aps}%
}e^{-t(D_{\operatorname{aps}})^{2}}$ enters in the definition of the
approximative kernel $Q_{2}^{R}$ as given in \eqref{e:parametrix_r}. We show
that $\mathcal{E}_{\operatorname{aps}}^{\infty}(t;x,x)$ is traceless for all
$x\in\lbrack0,\infty)\times\Sigma$. Then
\begin{equation}
\operatorname{Tr}Q_{2}^{R}(t;x,x)=\operatorname{Tr}\mathcal{E}^{R}%
(t;x,x)\ \text{ for all $x\in M_{2}^{R}$}, \label{e:7.6}%
\end{equation}
follows.

To prove that a product $TV$ is traceless, the following easy result can be used.

\begin{lemma}
\label{l:trace_even_odd} Let $\sigma$ be unitary with $\sigma^{2}%
=-\operatorname{I}$. We consider an operator $V$ of trace class which is
`even', i.e. it commutes with $\sigma$. Moreover, $T$ is odd, i.e. it
anticommutes with $\sigma$. Then
\[
\operatorname{Tr}(TV)=0.
\]
\end{lemma}

\begin{proof}
We have, by unitary equivalence,
\[
\operatorname{Tr}(TV)=\operatorname{Tr}(-\sigma(TV)\sigma)=\operatorname{Tr}%
(-\sigma T\sigma V)=\operatorname{Tr}(\sigma^{2}TV)=\operatorname{Tr}(-TV).
\]
\end{proof}

\begin{lemma}
\label{lem:3.1} Let $\chi:[0,\infty)\rightarrow\mathbb{R}$ be a smooth
function with compact support and $t>0$. Then the trace of the operator
$\chi\cdot D_{\operatorname{aps}}e^{-t(D_{\operatorname{aps}})^{2}}$ vanishes.
In particular,
\[
\int_{\Sigma}\operatorname{Tr}\mathcal{E}_{\operatorname{aps}}^{\infty
}(t;u,y;u,y)\,dy=0
\]
for all $u\in\lbrack0,\infty)$.\bigskip
\end{lemma}

\begin{proof}
Clearly, $D_{\operatorname{aps}}^{2}=\left(  \sigma(\partial_{u}+B)\right)
^{2}=-\partial_{u}^{2}+B^{2}$ is even, hence also the power series
$e^{-t(D_{\operatorname{aps}})^{2}}$ is even. On the other hand, as with $B$,
$\sigma B$ is also odd. So,
\[
\operatorname{Tr}\left(  \chi\cdot\sigma Be^{-t(D_{\operatorname{aps}})^{2}%
}\right)  =0.
\]
To show that
\[
\operatorname{Tr}\left(  (\chi\cdot\sigma\partial_{u}%
e^{-t(D_{\operatorname{aps}})^{2}}\right)  =0,
\]
we need a slightly more specific argument: Let $\operatorname{e}%
_{\operatorname{aps}}^{\infty}$ denote the heat kernel of the operator
$D_{\operatorname{aps}}$. For $u,v\in\lbrack0,\infty)$ and $y,z\in\Sigma$ it
has the following form (see e.g. \cite{BoWo93}, Formulae 22.33 and 22.35):
\[
\operatorname{e}_{\operatorname{aps}}(t;u,y;v,z)=\sum_{k\in\mathbb{Z}}%
e_{k}(t;u,v)\varphi_{k}(y)\otimes\varphi_{k}^{\ast}(z)
\]
for an orthonormal system $\{\varphi_{k}\}$ of eigensections of $B$. Hence,
\[
\sigma\partial_{u}\operatorname{e}_{\operatorname{aps}}(t;u,y;v,z)=\sum
_{k\in\mathbb{Z}}e_{k}^{\prime}(t;u,v)\sigma\varphi_{k}(y)\otimes\varphi
_{k}^{\ast}(z).
\]
But $\left\langle \sigma\varphi_{k};\varphi_{k}\right\rangle =0$ on $\Sigma$
since $\sigma$ is skew-adjoint.
\end{proof}

\bigskip

\subsection{Asymptotic Vanishing of Large-$t$ Chopped $\eta$-invariant on
Stretched Part Manifold\label{ss:large-t-aps}}

So far we have found
\begin{multline*}
\eta_{\mathcal{D}_{2,P_{>}}^{R}}(0)=\frac{1}{\sqrt{\pi}}\int_{0}^{\sqrt{R}%
}\frac{dt}{\sqrt{t}}\int_{M_{2}^{R}}\operatorname{Tr}\mathcal{E}%
^{R}(t;x,x)\,dx+\operatorname{O}(e^{-cR})\\
+\frac{1}{\sqrt{\pi}}\int_{\sqrt{R}}^{\infty}\frac{dt}{\sqrt{t}}\int
_{M_{2}^{R}}\operatorname{Tr}\mathcal{E}_{2}^{R}(t;x,x)\,dx
\end{multline*}
as $R\rightarrow\infty$. To prove Theorem \ref{t:adiabatic_eta}, we still have
to show
\begin{gather}
\frac{1}{\sqrt{\pi}}\int_{\sqrt{R}}^{\infty}\frac{dt}{\sqrt{t}}\int_{M_{2}%
^{R}}\operatorname{Tr}\mathcal{E}_{2}^{R}(t;x,x)\,dx=\operatorname{O}%
(e^{-cR})\text{ and}\label{e:lemma.7.1}\\
\frac{1}{\sqrt{\pi}}\int_{\sqrt{R}}^{\infty}\frac{dt}{\sqrt{t}}\int_{M_{2}%
^{R}}\operatorname{Tr}\mathcal{E}^{R}(t;x,x)\,dx=\operatorname{O}(e^{-cR}),
\label{e:7.18}%
\end{gather}
as $R\rightarrow\infty$. Recall that $\mathcal{E}_{2}^{R}(t;x,x^{\prime})$
denotes the kernel of the operator $\mathcal{D}_{2}^{R}e^{-t(\mathcal{D}%
_{2,P_{>}}^{R})^{2}}$ on the compact manifold $M_{2}^{R}$ with boundary
$\{-R\}\times\Sigma$, and $\mathcal{E}^{R}(t;x,x^{\prime})$ the kernel of the
operator $\mathcal{D}^{R}e^{-t(\mathcal{D}^{R})^{2}}$ on the closed stretched
manifold $M^{R}$.

\medskip

In the following we show \eqref{e:lemma.7.1}, i.e. that we can neglect the
contribution to the eta invariant of $\mathcal{D}_{2,P_{>}}^{R}$ which comes
from the large $t$ asymptotic of $\mathcal{E}_{2}^{R}(t;x,x^{\prime})$. The
key to that is that the eigenvalue of $\mathcal{D}_{2,P_{>}}^{R}$ with the
smallest absolute value is uniformly bounded away from zero.

\begin{theorem}
\label{t:6.1} Let $\mu_{0}(R)$ denote the smallest (in absolute value)
nonvanishing eigenvalue of the operator $\mathcal{D}_{2,P_{>}}^{R}$ on the
manifold $M_{2}^{R}$. Let us assume, as always in this section, that
$\operatorname{Ker}B=\{0\}$. Then there exists a positive constant $c_{0}$,
which does not depend on $R$ such that $\mu_{0}(R)>c_{0}$ for $R$ sufficiently large.
\end{theorem}

\begin{remark}
\label{r:adiabatic} As we will discover, this result indicates that the
behavior of the small eigenvalues on $M_{2}^{R}$ differs from that on the
stretched, closed manifold $M^{R}$. On the manifold with boundary $M_{2}^{R}$
with the attached cylinder of length $R$, the eigenvalues are bounded away
from 0 when $R\rightarrow\infty$ due to the spectral boundary condition. That
is the statement of Theorem \ref{t:6.1} which we are going to prove in the
next two sections. However, on $M^{R}$, the set of eigenvalues splits into one
set of eigenvalues becoming exponentially small and another one of eigenvalues
being uniformly bounded away from 0 as $R\rightarrow\infty$. This we are going
to show further below. Roughly speaking, the reason for the different behavior
is that on $M_{2}^{R}$ the \emph{eigensections} must satisfy the spectral
boundary condition. Therefore they are exponentially decreasing on the
cylinder, and the \emph{eigenvalues} are bounded away from 0. But on $M^{R}$
we have to cope with \emph{eigensections} on a closed manifold which need not
decrease, but require part of the \emph{eigenvalues} to decrease exponentially
(for details see Theorem \ref{t:0.3} below).
\end{remark}

\bigskip

\subsubsection{The Cylindrical Dirac Operator\label{sss:cylinder}\bigskip}

To prove Theorem \ref{t:6.1} we first recall a few properties of the
cylindrical Dirac operator $D_{\operatorname{cyl}}:=\sigma(\partial_{u}+B)$ on
the infinite cylinder $\Sigma_{\operatorname{cyl}}^{\infty}:=(-\infty
,+\infty)\times\Sigma$. A special feature of the cylindrical manifold
$\Sigma_{\operatorname{cyl}}^{\infty}$ is that we may apply the theory of
\textit{Sobolev spaces} exactly as in the case of $\mathbb{R}^{m}$. The point
is that we can choose a covering of the open manifold $\Sigma
_{\operatorname{cyl}}^{\infty}$ by a \emph{finite} number of coordinate
charts. We can also choose a finite trivialization of the bundle
$S|_{\Sigma_{\operatorname{cyl}}^{\infty}}$. Let $\{U_{\iota},\kappa_{\iota
}\}_{\iota=1}^{K}$ be such a trivialization, where $\kappa_{\iota
}:S|_{U_{\iota}}\rightarrow V_{\iota}\times\mathbb{C}^{N}$ is a bundle
isomorphism and $V_{\iota}$ an open (possibly non--compact) subset of
$\mathbb{R}^{m}$. Let $\{f_{\iota}\}$ be a corresponding partition of unity.
We assume that for any $\iota$ the derivatives of the function $f_{\iota}$ are bounded.

\begin{definition}
\label{d:sobolev} We say that a section (or distribution) $s$ of the bundle
$S$ over $\Sigma_{\operatorname{cyl}}^{\infty}$ belongs to the $p$-th
\emph{Sobolev space} $\mathcal{H}^{p}(\Sigma_{\operatorname{cyl}}^{\infty}%
;S)$, $p\in\mathbb{R}$, if and only if $f_{\iota}\cdot s$ belongs to the
Sobolev space $\mathcal{H}^{p}(\mathbb{R}^{m};\mathbb{C}^{N})$ for any $\iota
$. We define the $p$-th \emph{Sobolev norm}
\[
\left\|  s\right\|  _{p}:=\sum\nolimits_{\iota=1}^{K}\left\|  \left(
\operatorname{I}+\Delta_{\iota}\right)  ^{p/2}\left(  f_{\iota}\cdot s\right)
\right\|  _{L^{2}(\mathbb{R}^{m})}\,,
\]
where $\Delta_{\iota}$ denotes the Laplacian on the trivial bundle $V_{\iota
}\times\mathbb{C}^{N}\subset\mathbb{R}^{m}\times\mathbb{C}^{N}$.
\end{definition}

\begin{lemma}
\label{l:sobolev_d} \textrm{(a)} For the unique self-adjoint $L^{2}$ extension
of $D_{\operatorname{cyl}}$ (denoted by the same symbol) we have
\[
\operatorname{Dom}(D_{\operatorname{cyl}})=\mathcal{H}^{1}(\Sigma
_{\operatorname{cyl}}^{\infty};S).
\]
\noindent\textrm{(b)} Let $\lambda_{1}$ denote the smallest positive
eigenvalue of the operator $B$ on the manifold $\Sigma$. Then we have
\begin{equation}
\left\langle (D_{\operatorname{cyl}})^{2}s;s\right\rangle \geq\lambda_{1}%
^{2}\left\|  s\right\|  ^{2} \label{e:bounded_away}%
\end{equation}
for all $s\in\operatorname{Dom}(D_{\operatorname{cyl}})$, and for any $\mu
\in(-\lambda_{1},+\lambda_{1})$ the operator
\[
D_{\operatorname{cyl}}-\mu:\mathcal{H}^{1}(\Sigma_{\operatorname{cyl}}%
^{\infty};S)\rightarrow L^{2}(\Sigma_{\operatorname{cyl}}^{\infty};S)
\]
is an isomorphism of Hilbert spaces.\newline \noindent\textrm{(c)} Let
$\mathcal{R}_{\operatorname{cyl}}(\mu)$ denote the inverse of the operator
$D_{\operatorname{cyl}}-\mu$. Then the family $\{\mathcal{R}%
_{\operatorname{cyl}}(\mu)\}_{\mu\in(-\lambda_{1},\lambda_{1})}$ is a smooth
family of elliptic pseudo-differential operators of order $-1$.
\end{lemma}

\begin{proof}
(a) follows immediately from the corresponding result on the model manifold
$\mathbb{R}^{m}$. To prove (b) we consider a spectral resolution
$\{\varphi_{k},\lambda_{k}\}_{{k}\in\mathbb{Z}\setminus0}$ of $L^{2}%
(\Sigma;S)$ generated by the tangential operator $B$. Because of
(\ref{e:tangential_identities}), we have $\lambda_{-k}=-\lambda_{k}$. We can
assume $\varphi_{-k}=\sigma\varphi_{k}$ for $k\in\mathbb{N}$. We consider a
section $s$ belonging to the dense subspace $C_{0}^{\infty}(\Sigma
_{\operatorname{cyl}}^{\infty};S)$ of $\operatorname{Dom}%
(D_{\operatorname{cyl}})$, and expand it in terms of the preceding spectral
resolution
\[
s(u,y)=\sum\nolimits_{k\in\mathbb{Z}\setminus\{0\}}f_{k}(u)\varphi_{k}(y).
\]
Since $(D_{\operatorname{cyl}})^{2}=-\partial_{u}^{2}+B^{2}$, we obtain
\[
(D_{\operatorname{cyl}})^{2}s=\sum\nolimits_{k}(\lambda_{k}^{2}f_{k}%
-f_{k}^{\prime\prime})\varphi_{k},
\]
hence
\begin{align*}
\left\langle (D_{\operatorname{cyl}})^{2}s;s\right\rangle  &  =\sum
\nolimits_{k}\int_{-\infty}^{\infty}(\lambda_{k}^{2}f_{k}(u)-f_{k}%
^{\prime\prime}(u))\bar{f}_{k}(u)du\\
&  \geq\lambda_{1}^{2}\left\|  s\right\|  ^{2}-\sum\nolimits_{k}\int_{-\infty
}^{\infty}f_{k}^{\prime\prime}(u)\bar{f}_{k}(u)du\\
&  =\lambda_{1}^{2}\left\|  s\right\|  ^{2}+\sum\nolimits_{k}\int_{-\infty
}^{\infty}f_{k}^{\prime}(u)\bar{f}_{k}^{\prime}(u)du\geq\lambda_{1}%
^{2}\left\|  s\right\|  ^{2}.
\end{align*}
It follows that $(D_{\operatorname{cyl}})^{2}$ (and therefore
$D_{\operatorname{cyl}}$) has bounded inverse in $L^{2}(\Sigma
_{\operatorname{cyl}}^{\infty};S)$ and, more generally, that
$(D_{\operatorname{cyl}})^{2}-\mu$ is invertible for $\mu\in(-\lambda
_{1},\lambda_{1})$. To prove (c) we apply the symbolic calculus and construct
a parametrix $S$ for the operator $D_{\operatorname{cyl}}$; i.e., $S$ is an
elliptic pseudo-differential operator of order $-1$ such that
$SD_{\operatorname{cyl}}=\operatorname{I}+T$, where $T$ is a smoothing
operator. Thus $D_{\operatorname{cyl}}^{-1}=S-TD_{\operatorname{cyl}}^{-1}$.
The operator $TD_{\operatorname{cyl}}^{-1}$ is a smoothing operator, hence
$D_{\operatorname{cyl}}^{-1}$ is an elliptic pseudo-differential operator of
order $-1$. The same argument can be applied to the resolvent $\mathcal{R}%
_{\operatorname{cyl}}(\mu)=(D_{\operatorname{cyl}}-\mu)^{-1}$ for arbitrary
$\mu\in(-\lambda_{1},\lambda_{1})$. The smoothness of the family follows by
standard calculation.
\end{proof}

\bigskip

\subsubsection{The Part Manifolds with Half-infinite Cylindrical
Ends\label{sss:cylindrical-ends}}

To prove Theorem \ref{t:6.1} we need to refine the preceding results on the
infinite cylinder $\Sigma_{\operatorname{cyl}}^{\infty}$ to the Dirac
operator, naturally extended to the manifold $M_{2}^{\infty}=\left(
(-\infty,0]\times\Sigma\right)  \cup M_{2}$ with cylindrical end. Let
$C_{0}^{\infty}(M_{2}^{\infty};S)$ denote the space of compactly supported
smooth sections of $S$ over $M_{2}^{\infty}$. Then
\begin{equation}
\mathcal{D}_{2}^{\infty}|_{C_{0}^{\infty}(M_{2}^{\infty};S)}:C_{0}^{\infty
}(M_{2}^{\infty};S)\rightarrow L^{2}(M_{2}^{\infty};S) \label{e:d_j_infty}%
\end{equation}
is symmetric. Moreover, we have

\begin{lemma}
\label{l:infty_eigen_estimate} Let $s\in C^{\infty}(M_{2}^{\infty};S)$ be an
eigensection of $\mathcal{D}_{2}^{\infty}$. Then there exist $C,c>0$ such
that, on $(-\infty,0]\times\Sigma$, one has $\lvert s(u,y)\rvert\leq Ce^{cu}$.
\end{lemma}

\begin{proof}
Let $\{\varphi_{k},\lambda_{k}\}_{{k}\in\mathbb{Z}\setminus0}$ be a spectral
resolution of the tangential operator $B$. Because of
(\ref{e:tangential_identities}) we have $\lambda_{-k}=-\lambda_{k}$ and we can
assume that $\varphi_{-k}=\sigma\varphi_{k}$ for $k\in\mathbb{N}$. Then
\begin{equation}
\left\{  \varphi_{k}^{\pm}=\tfrac{1}{\sqrt{2}}\left(  \varphi_{k}\pm
\sigma\varphi_{k}\right)  ,\pm\lambda_{k}\right\}  _{k\in\mathbb{N}}
\label{e:eigen_gb}%
\end{equation}
is a spectral resolution of the composed operator $\sigma B$ on $\Sigma$.
Notice that we have
\begin{equation}
\sigma\varphi_{k}^{+}=-\varphi_{k}^{-}\quad\text{and}\quad\sigma\varphi
_{k}^{-}=\varphi_{k}^{+}\,. \label{e:toggle}%
\end{equation}
Let $s\in C^{\infty}(M_{2}^{\infty};S)$ and
\begin{equation}
\mathcal{D}_{L^{2}}^{\infty}\psi=\mu\psi\label{e:psi_mu}%
\end{equation}
with $\mu\in\mathbb{R}$. We expand $s|_{(-\infty,0]\times\Sigma}$ in terms of
the spectral resolution of $\sigma B$ just constructed:
\[
s(u,y)=\sum_{k=1}^{\infty}f_{k}(u)\varphi_{k}^{+}(y)+g_{k}(u)\varphi_{k}%
^{-}(y)\,.
\]
Because of \eqref{e:eigen_gb}, \eqref{e:toggle}, and \eqref{e:psi_mu} the
coefficients $f_{k},g_{k}$ must satisfy the system of ordinary differential
equations
\[%
\begin{pmatrix}
\lambda_{k} & \partial/\partial u\\
-\partial/\partial u & -\lambda_{k}%
\end{pmatrix}%
\begin{pmatrix}
f_{k}\\
g_{k}%
\end{pmatrix}
=\mu%
\begin{pmatrix}
f_{k}\\
g_{k}%
\end{pmatrix}
\]
or, equivalently,
\[%
\begin{pmatrix}
f_{k}^{\prime}\\
g_{k}^{\prime}%
\end{pmatrix}
=\mathbf{A}
\begin{pmatrix}
f_{k}\\
g_{k}%
\end{pmatrix}
\qquad\text{with $\mathbf{A}:=
\begin{pmatrix}
0 & -(\mu+\lambda_{k})\\
\mu-\lambda_{k} & 0
\end{pmatrix}
$}\,.
\]
Since $s\in L^{2}$, of the eigenvalues $\pm\sqrt{\lambda_{k}^{2}-\mu^{2}}$ of
$\mathbf{A}$\textbf{,} only those which are on the positive real line enter in
the construction of $s$ by solving the preceding differential equation. In
particular, all coefficients $f_{k},g_{k}$ must vanish identically for
$\lambda_{k}\leq\mu$. Thus,
\begin{equation}
s(u,y)=\sum_{\lambda_{k}>\mu}a_{k}\left(  \exp\left(  \sqrt{\lambda_{k}%
^{2}-\mu^{2}}\,u\right)  \,\varphi_{k}^{+}-\dfrac{\lambda_{k}-\mu}%
{\sqrt{\lambda_{k}^{2}-\mu^{2}}}\exp\left(  \sqrt{\lambda_{k}^{2}-\mu^{2}%
}\,u\right)  \,\varphi_{k}^{-}\right)  \label{e:gb_expansion}%
\end{equation}
and, in particular,
\[
\lvert s(u,y)\rvert\leq C\exp\left(  \sqrt{\lambda_{k_{0}}^{2}-\mu^{2}}%
\,\frac{u}{2}\right)  ,\qquad u<0,
\]
for some constant $C$, where $\lambda_{k_{0}}$ denotes the smallest positive
eigenvalue of $B$ with $\lambda_{k_{0}}>\lvert\mu\rvert$.
\end{proof}

In spectral theory we are looking for self-adjoint $L^{2}$ extensions of a
symmetric operator. We recall: on a \emph{closed} manifold, the Dirac operator
is \emph{essentially self-adjoint}; i.e. its minimal closed extension is
self-adjoint (and therefore there do not exist other self-adjoint extensions)
and it is a Fredholm operator. On a compact manifold with boundary, the
situation is much more complicated. There is a huge variety of dense domains
to which the Dirac operator can be extended such that it becomes self-adjoint;
and there is a smaller, but still large variety where the extension of the
Dirac operator becomes self-adjoint \textit{and} Fredholm (see Section
\ref{sss:Symmetric Operators} above and Boo\ss-Bavnbek and Furutani
\cite{BoFu98}); a special type of self-adjoint and Fredholm domains are the
domains specified by the boundary conditions belonging to the Grassmannian of
all self-adjoint generalized Atiyah--Patodi--Singer projections.

Now we shall show that the situation on manifolds with (infinite) cylindrical
ends resembles the situation on closed manifolds.

We recall the following simple lemma (see also Reed and Simon \cite{ReSi72},
Theorem VIII.3, Corollary, p. 257).

\begin{lemma}
\label{l:general_ess_self_adjoint} Let $A$ be a densely defined symmetric
operator in a separable complex Hilbert space $\mathcal{H}$. We assume that
$\operatorname{range}(A+i${$\operatorname{I}$}$)$ is dense in $\mathcal{H}$.
Then $A$ is essentially self-adjoint.
\end{lemma}

\begin{proof}
Since $A$ is symmetric, the operator $A+i\operatorname{I}$ is injective and
the operator $(A+i\operatorname{I})^{-1}$ is well defined and bounded on the
dense subspace $\operatorname{range}(A+i\operatorname{I})$ of $\mathcal{H}$.
Then the closure $R_{i}$ of $(A+i\operatorname{I})^{-1}$ has the whole space
$\mathcal{H}$ as domain and $R_{i}$ is bounded and injective. Now a standard
argument of functional analysis (see e.g. Pedersen \cite{Pe89}, Proposition
5.1.7) says that the inverse $R_{i}^{-1}$ of a densely defined, closed, and
injective operator $R_{i}$ has the same properties. Thus our $R_{i}^{-1}$ is
closed; and by construction it is the minimal closed extension of
$A+i\operatorname{I}$. Therefore, $R_{i}^{-1}-i\operatorname{I}$ is symmetric
and the minimal closed extension of $A$, hence self-adjoint and equal
$A^{\ast}$.
\end{proof}

We apply the lemma for $\mathcal{H}=L^{2}(M_{2}^{\infty};S)$ and take for $A $
the operator of \eqref{e:d_j_infty}. To prove that the range $(\mathcal{D}%
_{2}^{\infty}+i${$\operatorname{I}$}$)(C_{0}^{\infty}(M_{2}^{\infty};S))$ is
dense in $L^{2}(M_{2}^{\infty};S)$ we consider a section $s\in L^{2}%
(M_{2}^{\infty};S)$ which is orthogonal to $(\mathcal{D}_{2}^{\infty}%
+i${$\operatorname{I}$}$)(C_{0}^{\infty}(M_{2}^{\infty};S))$; i.e. the
distribution $(\mathcal{D}_{2}^{\infty}-i${$\operatorname{I}$}$)s$ vanishes
when applied to any test function, hence
\begin{equation}
(\mathcal{D}_{2}^{\infty}-i{\operatorname{I}})s=0.
\label{e:psi_i_eigensection}%
\end{equation}
Since $\mathcal{D}_{2}^{\infty}-i$ is elliptic, by elliptic regularity $s$ is
smooth at all interior points, that is for our complete manifold in all
points. On the cylinder $(-\infty,0]\times\Sigma$ we expand $s$ in terms of
the eigensections of the composed operator {$\sigma$}$B$ on $\Sigma$. It
follows that $s$ satisfies an estimate of the form
\begin{equation}
\lvert s(u,y)\rvert\leq Ce^{cu},\qquad(u,y)\in(-\infty,0]\times\Sigma,
\label{e:psi_estimate}%
\end{equation}
for some constants $C,c>0$ (according to Lemma \ref{l:infty_eigen_estimate}).
On the manifold $M_{2}^{R}$ with cylindrical end of finite length $R$ we apply
Green's formula and get
\begin{equation}
\left\langle \mathcal{D}_{2}^{R}s^{R};s^{R}\right\rangle -\left\langle
s^{R};\mathcal{D}_{2}^{R}s^{R}\right\rangle =-\int_{\{-R\}\times\Sigma
}({\sigma}s|_{\{-R\}\times\Sigma}\,dy,s|_{\{-R\}\times\Sigma}),
\label{e:green_psi}%
\end{equation}
where $s^{R}$ denotes the restriction of $s$ to the manifold $M_{2}^{R}$ with
boundary $\{-R\}\times\Sigma$. For $R\rightarrow\infty$, the right side of
\eqref{e:green_psi} vanishes; and the left side becomes $2i\left|  s\right|
^{2}$ by \eqref{e:psi_i_eigensection}. Hence $s=0$. Thus we have proved

\begin{lemma}
\label{l:ess_self_adjoint} The operator \eqref{e:d_j_infty} is essentially self-adjoint.
\end{lemma}

We denote the (unique) self-adjoint $L^{2}$ extension by the same symbol
$\mathcal{D}_{2}^{\infty}$, and we define the Sobolev spaces on the manifold
$M_{2}^{\infty}$ as in Definition \ref{d:sobolev}. Once again, the point is
that manifolds with cylindrical ends, even if they are not compact but only
complete, are like the infinite cylinder sufficiently simple to be covered by
a finite system of local charts. Clearly
\[
\operatorname{Dom}(\mathcal{D}_{2}^{\infty})=\mathcal{H}^{1}(M_{2}^{\infty
};S)\text{\ and\ }\mathcal{D}_{2}^{\infty}:\mathcal{H}^{1}(M_{2}^{\infty
};S)\rightarrow L^{2}(M_{2}^{\infty};S)
\]
is bounded. There are, however, substantial differences between the properties
of the simple Dirac operator $D_{\operatorname{cyl}}$ on the infinite cylinder
and the Dirac operator $\mathcal{D}_{2}^{\infty}$ on the manifold with
cylindrical end. For instance, from $B$ the discreteness of the spectrum and
the regularity at 0 (i.e., 0 is not an eigenvalue) are passed on to
$D_{\operatorname{cyl}}$, but not to $\mathcal{D}_{2}^{\infty}$. Yet we can
prove the following result:

\begin{proposition}
\label{p:lemma_6.2} The operator
\[
\mathcal{D}_{2}^{\infty}:\operatorname{Dom}(\mathcal{D}_{2}^{\infty
})=\mathcal{H}^{1}(M_{2}^{\infty};S)\rightarrow L^{2}(M_{2}^{\infty};S)
\]
is a Fredholm operator and its spectrum in the interval $(-\lambda_{1}%
,\lambda_{1})$ consists of finitely many eigenvalues of finite multiplicity.
Here $\lambda_{1}$ denotes the smallest positive eigenvalue of $B$.
\end{proposition}

\begin{note}
Actually, using more advanced methods one can show that the essential spectrum
of $\mathcal{D}_{2}^{\infty}$ is equal to $(-\infty,-\lambda_{1}]\cup
\lbrack\lambda_{1},\infty)$ (see for instance M\"{u}ller \cite{Mu94}, Section 4).
\end{note}

\medskip

Before proving the proposition we shall collect various criteria for the
compactness of a bounded operator between Sobolev spaces on an open manifold.
Let $X$ be a complete (not necessarily compact) Riemannian manifold with a
fixed Hermitian bundle. Recall the three cornerstones of the Sobolev analysis
of Dirac operators for $X$ \emph{closed}.

\begin{description}
\item [Rellich Lemma]The inclusion $\mathcal{H}^{1}(X)\subset L^{2}(X)$ is compact.

\item[Compact Resolvent] To each Dirac operator $\mathcal{D}$ we have a
parametrix $\mathcal{R}$ which is an elliptic pseudo-differential operator of
order $-1$ with principal symbol equal to the inverse of the principal symbol
of $\mathcal{D}$. So $\mathcal{R}$ is a bounded operator from $L^{2}(X)$ to
$\mathcal{H}^{1}(X)$, and hence compact in $L^{2}(X)$. In particular, for
$\mu$ in the resolvent set the resolvent $(\mathcal{D}-\mu${$\operatorname{I}%
$}$)^{-1}$ is compact as operator in $L^{2}(X)$.

\item[Smoothing Operator] Any integral operator over $X$ with smooth kernel is
a smoothing operator, i.e. it maps distributional sections of arbitrary low
order into smooth sections. Moreover, it is of trace class and thus compact.
\end{description}

In the \emph{general} case, i.e. for not necessarily compact $X$, the Rellich
Lemma remains valid for sections with compact support. A compact resolvent is
not attainable, hence the essential spectrum appears. Operators with smooth
kernel remain smoothing operators, but in general they are no longer of trace
class nor compact. We recall:

\begin{lemma}
\label{l:trace_compact} Let $X$ be a complete (not necessarily compact)
Riemannian manifold with fixed Hermitian bundle. Let $K$ be a compact subset
of $X$.\newline \noindent\ (a) The injection $\mathcal{H}^{1}(X)\subset
L^{2}(X)$ defines a compact operator when restricted to sections with support
in $K$. In particular, for any cutoff function $\chi$ with support in $K$ and
any bounded operator $\mathcal{R}:L^{2}(X)\rightarrow\mathcal{H}^{1}(X)$ the
operator $\chi\mathcal{R}$ is compact in $L^{2}(X)$.\newline \noindent\ (b)
Let $T:L^{2}(X)\rightarrow L^{2}(X)$ be an integral operator with a kernel
$k(x,y)\in L^{2}(X^{2})$. Then the operator $T$ is a bounded, compact operator
(in fact it is of Hilbert--Schmidt class).\newline \noindent\ (c) Let
$T:L^{2}(X)\rightarrow L^{2}(X)$ be a bounded compact operator and
$\mathcal{H}^{\prime}$ a closed subspace of $L^{2}(X)$, e.g. $\mathcal{H}%
^{\prime}:=L^{2}(X^{\prime})$ where $X^{\prime}$ is a submanifold of $X$ of
codimension 0. Assume that $T(\mathcal{H}^{\prime})\subset\mathcal{H}^{\prime
}$. Then $T|_{\mathcal{H}^{\prime}}$ is compact as operator from
$\mathcal{H}^{\prime}$ to $\mathcal{H}^{\prime}$.
\end{lemma}

\begin{proof}
(a) follows immediately from the local Rellich Lemma. (b) is the famous
Hilbert--Schmidt Lemma. Also (c) is well known, see e.g. H\"{o}rmander
\cite[Proposition 19.1.13]{Ho85} where (c) is proved within the category of
trace class operators.
\end{proof}

In general an integral operator $T$ with smooth kernel is not compact even if
either $\operatorname*{supp}_{x}k(x,x^{\prime})$ or $\operatorname*{supp}%
_{x^{\prime}}k(x,x^{\prime})$ are contained in a compact subset $K\subset X$.
Consider for instance on $\Sigma_{\operatorname{cyl}}^{\infty}=(-\infty
,+\infty)\times\Sigma$ an integral operator $T$ with a smooth kernel of the
form
\[
k(x,x^{\prime})=\chi(x)d(x,x^{\prime}),
\]
where $d(x,x^{\prime})$ denotes the distance and $\chi$ is a function with
support in a ball of radius 1 (and equal 1 in a smaller ball). Then $T$ is not
a compact operator on $L^{2}(\Sigma_{\operatorname{cyl}}^{\infty})$: choose a
sequence $\{s_{n}\}$ of $L^{2}$ functions of norm 1 and with
$\operatorname*{supp}s_{n}$ contained in a ball of radius 1 such that
$d(\operatorname*{supp}\chi,\operatorname*{supp}s_{n})=n$. Then for any $n$ we
have $\left|  Ts_{n}\right|  >Cn$. Thus $T$ is not compact, in fact not even bounded.

For the bounded resolvent (see Lemma \ref{l:sobolev_d})
\[
\mathcal{R}_{\operatorname{cyl}}:L^{2}(\Sigma_{\operatorname{cyl}}^{\infty
};S)\rightarrow\mathcal{H}^{1}(\Sigma_{\operatorname{cyl}}^{\infty};S)
\]
we have, however, the following corollary to the preceding lemma. It provides
an example of a compact integral operator on an open manifold with a smooth
kernel which is compactly supported only in one variable.

\begin{corollary}
\label{c:compact} Let $\chi$ and $\psi$ be smooth cutoff functions on
$\Sigma_{\operatorname{cyl}}^{\infty}$ with support contained in the
half-cylinder $(-\infty,0)\times\Sigma$. Let $\operatorname*{supp}\chi$ be
compact. Then the operators $\chi\mathcal{R}_{\operatorname{cyl}}\psi$ and
$\psi\mathcal{R}_{\operatorname{cyl}}\chi$ are compact in $L^{2}%
(\Sigma_{\operatorname{cyl}}^{\infty};S)$.
\end{corollary}

\begin{proof}
The operator $\chi\mathcal{R}_{\operatorname{cyl}}\psi$ is compact according
to the preceding lemma, claim (a). Its adjoint operator is $\psi
\mathcal{R}_{\operatorname{cyl}}\chi$, since $\mathcal{R}_{\operatorname{cyl}%
}$ is self-adjoint. Thus it is also compact (even if its range is \emph{not}
compactly supported).
\end{proof}

\bigskip

\subsection{The Estimate of the Lowest Nontrivial
Eigenvalue\label{ss:LowestEigenv}\label{s:lowest_eigenvalue}}

In this section we prove Theorem \ref{t:6.1}. Recall that the tangential
operator $B$ is assumed to be nonsingular and that $\lambda_{1}$ denotes the
smallest positive eigenvalue of $B$. So far, we have established that

\begin{description}
\item [I]the operator $D_{\operatorname{cyl}}$ on the infinite cylinder
$\Sigma_{\operatorname{cyl}}^{\infty}$ has no eigenvalues in the interval
$(-\lambda_{1},+\lambda_{1})$, and

\item[II] the operator $\mathcal{D}_{2}^{\infty}$ on the manifold
$M_{2}^{\infty}$ with infinite cylindrical end has only finitely many
eigenvalues in the interval $(-\lambda_{1},+\lambda_{1})$, each of finite multiplicity.
\end{description}

We have to show that

\begin{description}
\item [III]the nonvanishing eigenvalues of $(\mathcal{D}_{2}^{R})_{P_{>}}$ are
bounded away from 0 by a bound independent of $R$.
\end{description}

\bigskip

\textbf{Proof of Theorem} \textbf{\ref{t:6.1}} The idea of the proof is the
following. We define a positive constant $\mu_{1}$ independent of $R$. Then
let $R$ be a positive real (more precisely $R>R_{0}$ for a suitable positive
$R_{0}$), and $s\in L^{2}(M_{2}^{R};S)$ any eigensection with eigenvalue
$\mu\in(-\lambda_{1}/\sqrt{2},+\lambda_{1}/\sqrt{2})$, i.e.
\[
s\in\operatorname{Dom}(\mathcal{D}_{2}^{R})_{P_{>}}\text{ i.e. $P_{>}%
(s|_{\{-R\}\times\Sigma})=0$\ and\ }\mathcal{D}_{2}^{R}s=\mu s.
\]
Then we show that $\mu^{2}>\mu_{1}/2$ for a certain real $\mu_{1}>0$ which is
independent of $R$ and $s$. A natural choice of $\mu_{1}$ is
\begin{equation}
\mu_{1}=\min\left\{  \left.  \frac{\left|  \mathcal{D}_{2}^{\infty}%
\Psi\right|  ^{2}}{\left|  \Psi\right|  ^{2}}\right|  \;\Psi\in\mathcal{H}%
^{1}(M_{2}^{\infty};S)\text{\ and\ }\Psi\perp\operatorname{Ker}\mathcal{D}%
_{2}^{\infty}\right\}  .
\end{equation}
Note that by II above (Proposition \ref{p:lemma_6.2}), $\operatorname{Ker}%
\mathcal{D}_{2}^{\infty}$ is of finite dimension. We shall define a certain
extension $s^{\infty}\in\mathcal{H}^{1}(M_{2}^{\infty};S)$ of $s$.

The reasoning would be easy, if we could extend $s$ to an eigensection of
$\mathcal{D}_{2}^{\infty}$ on all of $M_{2}^{\infty}$. Then it would follow at
once that the discrete part of the spectrum of $\mathcal{D}_{2}^{\infty}$ is
not empty, $\mu$ belongs to it, $\sqrt{\mu_{1}}$ is the smallest eigenvalue
$>0$, and hence we would have $\mu^{2}>\mu_{1}/2$ as desired. In general, such
a convenient extension of the given eigensection $s$ cannot be achieved. But
due to the spectral boundary condition satisfied by $s$ in the hypersurface
$\{-R\}\times\Sigma$, the eigensection $s$ over $M_{2}^{R} $ can be
continuously extended by a section over $(-\infty,-R]\times\Sigma$ on which
the Dirac operator vanishes. By construction, both the enlargement $\alpha$ of
the $L^{2}$ norm of $s$ by the chosen extension and the cosine, say $\beta$,
of the angle between $s^{\infty}$ and $\operatorname{Ker}\mathcal{D}%
_{2}^{\infty}$ can be estimated independently of the specific choice of $s$
and $\mu$. It turns out that they both decrease exponentially with growing $R
$.

Let $\{s_{1},\dots,s_{q}\}$ be an orthonormal basis of $\operatorname{Ker}%
\mathcal{D}_{2}^{\infty}$ and set
\[
\widetilde{s}:=s^{\infty}-\sum_{j=1}^{q}\left\langle s^{\infty};s_{j}%
\right\rangle s_{j}\,.
\]
Clearly, the section $\widetilde{s}$ belongs to $\mathcal{H}^{1}(M_{2}%
^{\infty};S)$ and is orthogonal to $\operatorname{Ker}\mathcal{D}_{2}^{\infty}
$. Hence, on the one hand,
\begin{equation}
\frac{\left|  \mathcal{D}_{2}^{\infty}\widetilde{s}\right|  ^{2}}{\left|
\widetilde{s}\right|  ^{2}}\geq\mu_{1}\,. \label{e:hat_estimate}%
\end{equation}
On the other hand, we have by construction
\[
\left|  \mathcal{D}_{2}^{\infty}\widetilde{s}\right|  ^{2}=\left|
\mathcal{D}_{2}^{\infty}s^{\infty}\right|  ^{2}=\left|  \mathcal{D}_{2}%
^{R}s\right|  _{M_{2}^{R}}^{2}=\mu^{2}\,.
\]
Finally, we shall prove that
\begin{equation}
\left|  \widetilde{s}\right|  \rightarrow1\text{ as $R\rightarrow\infty$}.
\label{e:6.20'}%
\end{equation}
Then the estimate
\begin{equation}
\mu^{2}>\frac{\mu_{1}}{2} \label{e:theorem_6.1}%
\end{equation}
follows for sufficiently large $R$. Since we have assumed that $\mu^{2}%
<\frac{\lambda_{1}^{2}}{2}$, we have also $\mu_{1}<\lambda_{1}^{2}$, hence
$\mu_{1}$ belongs to the discrete part of the spectrum of $\left(
\mathcal{D}_{2}^{\infty}\right)  ^{2}$ and, by the \emph{Min--Max Principle}
(e.g., see \cite{ReSi78}), must be its smallest eigenvalue $>0$.

\medskip Thus, to prove the Theorem \ref{t:6.1}\textbf{\ }we are left with the
task of first constructing a suitable extension $\widetilde{s}$ of $s$ and
then proving \eqref{e:6.20'}.

\medskip We expand $s|_{[-R,0]\times\Sigma}$ in terms of a spectral
resolution
\[
\{\varphi_{k},\lambda_{k};\sigma\varphi_{k},-\lambda_{k}\}_{k\in\mathbb{N}}%
\]
of $L^{2}(\Sigma;S)$ generated by $B$:
\[
s(u,y)=\sum_{k=1}^{\infty}f_{k}(u)\varphi_{k}(y)+g_{k}(u)\mathbb{\sigma
}\varphi_{k}(y)\,.
\]
Since $\mathbb{\sigma}(\partial_{u}+B)s|_{[-R,0]\times\Sigma}=0$, the
coefficients must satisfy the system of ordinary differential equations
\begin{equation}%
\begin{pmatrix}
f_{k}^{\prime}\\
g_{k}^{\prime}%
\end{pmatrix}
=\mathbf{A}_{k}
\begin{pmatrix}
f_{k}\\
g_{k}%
\end{pmatrix}
\qquad\text{with $\mathbf{A}_{k}:=
\begin{pmatrix}
-\lambda_{k} & \mu\\
-\mu & \lambda_{k}%
\end{pmatrix}
$}\,. \label{e:dgl_s}%
\end{equation}
Moreover, since $P_{>}s|_{\{-R\}\times\Sigma}=0$ we have
\begin{equation}
f_{k}(-R)=0\quad\text{ for any $k\geq1$}. \label{e:bc_s}%
\end{equation}
Thus, for each $k$ the pair $(f_{k},g_{k})$ is uniquely determined up to a
constant $a_{k}$. More explicitly, since the eigenvalues of $\mathbf{A}_{k}$
are $\pm(\lambda_{k}^{2}-\mu^{2})^{1/2}$, a suitable choice of the
eigenvectors of $\mathbf{A}_{k}$ gives
\begin{gather*}
f_{k}(u)=a_{k}\frac{\mu}{\sqrt{\lambda_{k}^{2}-\mu^{2}}}\sinh\sqrt{\lambda
_{k}^{2}-\mu^{2}}(R+u)\quad\text{and}\\
g_{k}(u)=a_{k}\left(  \cosh\sqrt{\lambda_{k}^{2}-\mu^{2}}(R+u)+\frac
{\lambda_{k}}{\sqrt{\lambda_{k}^{2}-\mu^{2}}}\sinh(\lambda_{k}^{2}-\mu
^{2})^{1/2}(R+u)\right)  .
\end{gather*}
We assume $\left|  s\right|  _{L^{2}}=1$. Then we have, with $v:=(\lambda
_{k}^{2}-\mu^{2})^{1/2}(R+u)$:
\begin{align*}
1  &  \geq\int_{\lbrack-R,0]\times\Sigma}\lvert s(u,y)\rvert^{2}%
\,dudy=\sum_{k=1}^{\infty}\int_{-R}^{0}\left(  \lvert f_{k}(u)\rvert
^{2}+\lvert g_{k}(u)\rvert^{2}\right)  du\\
&  =\sum_{k=1}^{\infty}\lvert a_{k}\rvert^{2}\frac{1}{(\lambda_{k}^{2}-\mu
^{2})^{1/2}}\int_{0}^{(\lambda_{k}^{2}-\mu^{2})^{1/2}R}\left(  \frac{\mu^{2}%
}{\lambda_{k}^{2}-\mu^{2}}\cdot\sinh^{2}v\right. \\
&  \qquad\qquad\left.  +\cosh^{2}v+2\frac{\lambda_{k}}{(\lambda_{k}^{2}%
-\mu^{2})^{1/2}}\cdot\cosh v\cdot\sinh v+\frac{\lambda_{k}^{2}}{\lambda
_{k}^{2}-\mu^{2}}\cdot\sinh^{2}v\right)  dv\\
&  =\sum_{k=1}^{\infty}\lvert a_{k}\rvert^{2}\left\{  -\frac{\lambda_{k}^{2}%
}{\lambda_{k}^{2}-\mu^{2}}\cdot R+(1/4)\cdot\frac{\mu^{2}}{(\lambda_{k}%
^{2}-\mu)^{3/2}}\cdot\sinh(2(\lambda_{k}^{2}-\mu^{2})^{1/2}R)\right. \\
&  \qquad+(1/4)\cdot\left(  1+\frac{\lambda_{k}^{2}}{\lambda_{k}^{2}-\mu^{2}%
}\right)  \cdot(\lambda_{k}^{2}-\mu^{2})^{1/2}\cdot\sinh(2(\lambda_{k}^{2}%
-\mu^{2})^{1/2}R)\\
&  \qquad\left.  +\frac{\lambda_{k}^{2}}{\lambda_{k}^{2}-\mu^{2}}\cdot
\cosh^{2}((\lambda_{k}^{2}-\mu^{2})^{1/2}R)\right\}  .
\end{align*}

Since $\lambda_{k}^{2}\geq\lambda_{1}^{2}>{2}\mu^{2}$ we have $2(\lambda
_{k}^{2}-\mu^{2})^{1/2}>\sqrt{2}\lambda_{k}>\lambda_{k}$. Moreover, we have
for all $k\geq1$
\begin{align*}
-\frac{\lambda_{k}^{2}}{\lambda_{k}^{2}-\mu^{2}}\cdot R  &  +\frac{\lambda
_{k}^{2}}{\lambda_{k}^{2}-\mu^{2}}\cdot\frac{(\lambda_{k}^{2}-\mu^{2})^{1/2}%
}{4}\cdot\sinh(2(\lambda_{k}^{2}-\mu^{2})^{1/2}R)\\
&  >-R+\frac{(\lambda_{1}^{2}-\mu^{2})^{1/2}}{4}\cdot\sinh(2(\lambda_{1}%
^{2}-\mu^{2})^{1/2}R)\\
&  >-R+\frac{\sqrt{2}}{8}\lambda_{1}\cdot\sinh(\sqrt{2}\lambda_{1}R)>0,
\end{align*}
if $R\geq R_{0}$ for some positive $R_{0}$ which depends only on $\lambda_{1}
$ and not on $\mu$, $s$ and $k$.

Thus, for any $k$ the sum in the braces can be estimated in the following
way:
\[
\left\{  \mathstrut\dots\right\}  >\frac{\lambda_{k}^{2}}{\lambda_{k}^{2}%
-\mu^{2}}\cdot\cosh^{2}((\lambda_{k}^{2}-\mu^{2})^{1/2}R)>\frac{1}%
{4}e^{2(\lambda_{k}^{2}-\mu^{2})^{1/2}R}>\frac{1}{4}e^{\lambda_{k}R}\,.
\]
Hence, we have
\begin{equation}
\sum_{k=1}^{\infty}\lvert a_{k}\rvert^{2}\cdot e^{\lambda_{k}R}\leq4\,.
\label{e:6.8}%
\end{equation}
Note that the preceding estimate does not depend on $R$ (provided that
$R>R_{0}$), $k$ or the specific choice of $s$, and that $R_{0}$ only depends
on $\lambda_{1}$.

According to \eqref{e:6.8} the absolute value of the coefficients $a_{k}$ is
rapidly decreasing in such a way that, in particular, we can extend the
eigensection $s$ of $\mathcal{D}_{2,P_{>}}^{R}$, given on $M_{2}^{R}$ to a
continuous section on $M_{2}^{\infty}$ by the formula
\[
s^{\infty}(x):=\left\{
\begin{array}
[c]{ll}%
s(x) & \text{ for $x\in M_{2}^{R}$}\\
\sum_{k=1}^{\infty}a_{k}e^{\lambda_{k}(R+u)}\mathbb{\sigma}\varphi_{k}(y) &
\text{ for $x=(u,y)\in(-\infty,-r]\times\Sigma$}.
\end{array}
\right.
\]
By construction, $s^{\infty}$ is smooth on $M_{2}^{\infty}\setminus
(\{-R\}\times\Sigma)$ and belongs to the Sobolev space $\mathcal{H}^{1}%
(M_{2}^{\infty};S)$. It follows from \eqref{e:6.8} that
\begin{align}
\left|  s^{\infty}\right|  _{L^{2}}^{2}  &  =\left|  s\right|  _{L^{2}}%
^{2}+\sum_{k=1}^{\infty}\lvert a_{k}\rvert^{2}\int_{-\infty}^{-R}%
e^{2\lambda_{k}(R+u)}\,du=1+\sum_{k=1}^{\infty}\lvert a_{k}\rvert^{2}\frac
{1}{2\lambda_{k}}\nonumber\\
&  \leq1+\frac{1}{2\lambda_{1}}\cdot\sum_{k=1}^{\infty}\lvert a_{k}\rvert
^{2}\leq1+\frac{1}{2\lambda_{1}}\cdot\left(  \sum\nolimits_{k=1}^{\infty
}\lvert a_{k}\rvert^{2}\cdot e^{\lambda_{k}R}\right)  \cdot e^{-\lambda_{1}%
R}\nonumber\\
&  \leq1+\frac{2}{\lambda_{1}}\cdot e^{-\lambda_{1}R}\,. \label{e:6.10}%
\end{align}

\medskip Next, let $\Psi\in\operatorname{Ker}\mathcal{D}_{2}^{\infty}$ and
assume that $\left|  \Psi\right|  =1$. By \eqref{e:gb_expansion}, on
$(-\infty,0]\times\Sigma$ the section $\Psi$ has the form
\begin{equation}
\Psi((u,y))=\sum_{k=1}^{\infty}b_{k}e^{\lambda_{k}u}G(y)\varphi_{k}(y)
\label{e:gb_expansion'}%
\end{equation}
with
\[
\sum_{k=1}^{\infty}\int_{-\infty}^{0}\lvert b_{k}\rvert^{2}e^{2\lambda_{k}%
u}\,du=\sum_{k=1}^{\infty}\frac{1}{2\lambda_{k}}\cdot\lvert b_{k}\rvert
^{2}<+\infty.
\]
Set $l:=\Psi|_{M_{2}^{R}}$. Then $l$ satisfies the equations
\[
\mathcal{D}_{2}^{R}l=0\text{\ and}\;P_{>}(l|_{\{-R\}\times\Sigma})=0.
\]
Hence, $l$ belongs to $\operatorname{Ker}\mathcal{D}_{2,P_{>}}^{R}$. This
implies the following equality:
\begin{multline*}
\int_{M_{2}^{R}}\left\langle s^{\infty}(x);\Psi(x)\right\rangle dx=\frac
{1}{\mu}\cdot\int_{M_{2}^{R}}\left\langle \mathcal{D}_{2}^{R}s^{\infty
}(x);l(x)\right\rangle dx\\
=\frac{1}{\mu}\cdot\int_{M_{2}^{R}}\left\langle s^{\infty}(x);\mathcal{D}%
_{2}^{R}l(x)\right\rangle dx-\frac{1}{\mu}\cdot\int_{\Sigma}\left\langle
\mathbb{\sigma}s^{\infty}(-R,y);l(-R,y)\right\rangle dx.
\end{multline*}
On the other hand,
\[
\int_{(-\infty,-r]\times\Sigma}\left\langle s^{\infty}(x);\Psi(x)\right\rangle
dx=\sum_{k>0}\frac{a_{k}\overline{b_{k}}}{2\mu}\cdot{e^{-\lambda_{k}R}}\leq
C_{1}e^{-\lambda_{1}R}\,.
\]
Therefore
\begin{equation}
\lvert\left\langle s^{\infty};\Psi\right\rangle \rvert\leq C_{1}%
e^{-\lambda_{1}R}\,. \label{e:6.14}%
\end{equation}

Hence, we have proved

\begin{lemma}
\label{l:eigenext}Any eigensection $s\in\mathcal{H}^{1}(M_{2}^{R};S)$ of
$\mathcal{D}_{2,P_{>}}^{R}$ with eigenvalue $\mu\in(-\lambda_{1}/\sqrt
{2},\lambda_{1}/\sqrt{2})$ can be extended to a continuous section $s^{\infty}
$ on $M_{2}^{\infty}$ which is smooth on $M_{2}^{\infty}\setminus\left(
\{-R\}\times\Sigma\right)  $ and belongs to the first Sobolev space
$\mathcal{H}^{1}(M_{2}^{\infty};S)$. Moreover, the enlargement of the norm of
$s$ by the extension and the cosine of the angle between $s^{\infty}$ and
$\operatorname{Ker}\mathcal{D}_{2}^{\infty}$ are exponentially decreasing by
formulae \eqref{e:6.10} and \eqref{e:6.14}.
\end{lemma}

The final step in proving the Theorem \ref{t:6.1}\textbf{\ }follows at once
from the preceding lemma. We recall: by definition of $\mu_{1}$ we have
$\left|  \mathcal{D}_{2}^{\infty}\widetilde{s}\right|  ^{2}/\left|
\widetilde{s}\right|  ^{2}\geq\mu_{1}$ and by construction of $\widetilde{s}$
we have $\left|  \mathcal{D}_{2}^{\infty}\widetilde{s}\right|  ^{2}=\mu^{2}$.
Thus, we have $\mu^{2}\geq\left|  \widetilde{s}\right|  ^{2}\cdot\mu_{1}$. To
get the desired bound $\mu^{2}>\mu_{1}/2$, it remains to show that $\left|
\widetilde{s}\right|  ^{2}>1/2$ for sufficiently large $R$.

Since $\widetilde{s}$ is the orthogonal projection of $s^{\infty}$ onto
$(\operatorname{Ker}\mathcal{D}_{2}^{\infty})^{\perp}$ and the basis
$\{s_{1},\dots,s_{q}\}$ of $\operatorname{Ker}\mathcal{D}_{2}^{\infty}$ is
orthonormal, we have
\[
\left|  \widetilde{s}\right|  ^{2}=\left|  s^{\infty}\right|  ^{2}-\sum
_{j=1}^{q}\lvert\left\langle s^{\infty};s_{j}\right\rangle \rvert^{2}%
\leq1+\frac{2}{\lambda_{1}}e^{-\lambda_{1}R}-\sum_{j=1}^{q}\lvert\left\langle
s^{\infty};s_{j}\right\rangle \rvert^{2}\,.
\]
Thus, Theorem \ref{t:6.1}\textbf{\ }follows from
\[
\lvert\left|  \widetilde{s}\right|  ^{2}-1\rvert\leq\frac{2}{\lambda_{1}%
}e^{-\lambda_{1}R}+qC_{1}e^{-2\lambda_{1}R}\leq C_{2}e^{-C_{3}R}.
\]

\bigskip We finish this section by proving the asymptotic estimate
\eqref{e:lemma.7.1}. Recall that $\mathcal{E}_{2}^{R}(t;x,x^{\prime})$ denotes
the kernel of the operator $\mathcal{D}_{2,P_{>}}^{R}e^{-t(\mathcal{D}%
_{2,P_{>}}^{R})^{2}}$, where $\mathcal{D}_{2,P_{>}}^{R}$ denotes the Dirac
operator over the manifold $M_{2}^{R}$ with the spectral boundary condition at
the boundary $\{-R\}\times\Sigma$. Then we have:

\begin{lemma}
\label{l:lemma.7.1}
\[
\frac{1}{\sqrt{\pi}}\int_{\sqrt{R}}^{\infty}\frac{dt}{\sqrt{t}}\int_{M_{2}%
^{R}}\operatorname{Tr}\mathcal{E}_{2}^{R}(t;x,x)\,dx=\operatorname{O}%
(e^{-cR}).
\]
\end{lemma}

\begin{proof}
For any eigenvalue $\mu\neq0$ of $\mathcal{D}_{2,P_{>}}^{R}$ and $R>0$
sufficiently large ($R\cdot c_{0}^{2}\geq1$, where $c_{0}$ denotes the lower
uniform bound for $\mu^{2}$ of Theorem \ref{t:6.1}) we have
\begin{multline}
\lvert\int_{\sqrt{R}}^{\infty}\frac{1}{\sqrt{t}}\mu e^{-t\mu^{2}}%
\,dt\rvert\leq\int_{\sqrt{R}}^{\infty}\frac{1}{\sqrt{t}}\lvert\mu\rvert
e^{-t\mu^{2}}\,dt=\int_{\lvert\mu\rvert R^{1/4}}^{\infty}e^{-\tau^{2}}%
\,d\tau\label{e:gromov?}\\
\leq\int_{\lvert\mu\rvert R^{1/4}}^{\infty}\tau e^{-\tau^{2}}\,d\tau= \left[
-\frac{1}{2}e^{-\tau^{2}}\right]  _{\lvert\mu\rvert R^{1/4}}^{\infty}\frac
{1}{2}e^{-\mu^{2}\sqrt{R}}\,,
\end{multline}
which gives
\begin{align*}
&  \lvert\int_{\sqrt{R}}^{\infty}\frac{dt}{\sqrt{t}}\int_{M_{2}^{R}%
}\operatorname{Tr}\left(  \mathcal{D}_{2,P_{>}}^{R}e^{-t(\mathcal{D}_{2,P_{>}%
}^{R})^{2}}\right)  \rvert\leq\int_{\sqrt{R}}^{\infty}\frac{dt}{\sqrt{t}}%
\int_{M_{2}^{R}}\sum_{\mu\neq0}\lvert\mu\rvert e^{-t\mu^{2}}\,dt\\
&  \qquad\qquad\leq\frac{1}{2}\cdot\sum_{\mu\neq0}e^{-\mu^{2}\sqrt{R}}%
=\frac{1}{2}\cdot\sum_{\mu\neq0}e^{-(\sqrt{R}-1)\mu^{2}}\cdot e^{-\mu^{2}}\\
&  \qquad\qquad\leq C_{1}\cdot e^{-\sqrt{R}\mu_{0}^{2}}\operatorname{Tr}%
\left(  (e^{-\left(  \mathcal{D}_{2,P_{>}}^{R}\right)  ^{2}}\right)  \leq
C_{2}\cdot e^{-\sqrt{R}\mu_{0}^{2}}\operatorname*{Vol}(M_{2}^{R})\\
&  \qquad\qquad\leq C_{3}\cdot e^{-\sqrt{R}\mu_{0}^{2}}\leq C_{3}\cdot
e^{-C_{4}\sqrt{R}}\,.
\end{align*}
Here we have exploited that the heat kernel $\operatorname{e}_{2}%
^{R}(t;x,x^{\prime}) $ of the operator $\mathcal{D}_{2,P_{>}}^{R}$ can be
estimated by
\[
\lvert\operatorname{e}_{2}^{R}(t;x,x^{\prime})\rvert\leq c_{1}\cdot
t^{-\frac{m}{2}}\cdot e^{c_{2}t}\cdot e^{-c_{3}\frac{d^{2}(x,x^{\prime})}{t}}%
\]
according to \eqref{e:4.3.a}. Thus,
\begin{equation}
\lvert\operatorname{Tr}\left(  e^{-\left(  \mathcal{D}_{2,P_{>}}^{R}\right)
^{2}}\right)  \rvert\leq\int_{M_{2}^{R}}\lvert\operatorname{Tr}%
\operatorname{e}_{2}^{R}(1;x,x)\rvert\,dx\leq c_{1}\cdot e^{c_{2}}\cdot
\int_{M_{2}^{R}}dx. \label{e:vol}%
\end{equation}
\bigskip
\end{proof}

\subsection{The Spectrum on the Closed Stretched Manifold\label{ss:Spectrum on
Closed Stretched Mfd}}

Thus far, we have proved the asymptotic equation
\[
\frac{1}{\sqrt{\pi}}\int_{0}^{\sqrt{R}}\frac{dt}{\sqrt{t}}\int_{M_{2}^{R}%
}\operatorname{Tr}\mathcal{E}^{R}(t;x,x)\,dx+\operatorname{O}(e^{-cR}%
)=\eta_{\mathcal{D}_{2,P_{>}}^{R}}(0)
\]
as $R\rightarrow\infty$. It follows that
\[
\lim_{R\rightarrow\infty}\eta_{R}=\lim_{R\rightarrow\infty}\left(
\eta_{\mathcal{D}_{1,P_{<}}^{R}}(0)+\eta_{\mathcal{D}_{2,P_{>}}^{R}%
}(0)\right)  ,
\]
where
\[
\eta_{R}:=\frac{1}{\sqrt{\pi}}\int_{0}^{\sqrt{R}}\frac{dt}{\sqrt{t}}%
\int_{M^{R}}\operatorname{Tr}\mathcal{E}^{R}(t;x,x)\,dx.
\]
To prove Theorem \ref{t:adiabatic_eta}, we still have to show \eqref{e:7.18},
i.e., that we can extend the integration from $\sqrt{R}$ to infinity:
\[
\frac{1}{\sqrt{\pi}}\int_{\sqrt{R}}^{\infty}\frac{dt}{\sqrt{t}}\int_{M_{2}%
^{R}}\operatorname{Tr}\mathcal{E}^{R}(t;x,x)\,dx=\operatorname{O}%
(e^{-cR})\text{ as $R\rightarrow\infty$.}%
\]
Recall that $\mathcal{E}^{R}(t;x,x^{\prime})$ denotes the kernel of the
operator $\mathcal{D}^{R}e^{-t(\mathcal{D}^{R})^{2}}$ on the closed stretched
manifold $M^{R}$.

Formally, our \emph{task} of proving the preceding estimate is reminiscent of
our previous task of proving the corresponding estimate for the kernel
$\mathcal{E}_{2}^{R}(t;x,x^{\prime})$ of the operator $\mathcal{D}_{2}%
^{R}e^{-t(\mathcal{D}_{2,P_{>}}^{R})^{2}}$ (see Lemma \ref{l:lemma.7.1}). Both
integrals are over the same prolonged compact manifold $M_{2}^{R}$ with
boundary $\{-R\}\times\Sigma$. However, the methods we can apply are
different: In the previous case, we had a uniform positive bound for the
absolute value of the smallest nonvanishing eigenvalue of the boundary value
problem $\mathcal{D}_{2,P_{>}}^{R}$ for sufficiently large $R$.

As mentioned above in Remark \ref{r:adiabatic}, such a bound does not exist
for the Dirac operator $\mathcal{D}^{R}$ on the closed stretched manifold
$M^{R}$. Moreover, for the spectral boundary condition we shall show
\[
\dim\operatorname{Ker}\mathcal{D}_{2,P_{>}}^{R}=\dim\operatorname{Ker}%
\mathcal{D}_{2,P_{>}}\;\text{and}\;\eta_{\mathcal{D}_{2,P_{>}}^{R}}%
(0)=\eta_{\mathcal{D}_{2,P_{>}}}(0)
\]
for any $R$ (see Proposition \ref{p:muller} below). For $\mathcal{D}^{R}$, on
the contrary, the dimension of the kernel can change and, thus, $\eta
_{\mathcal{D}^{R}}$ can admit an integer jump in value as $R\rightarrow\infty
$. This is due to the presence of `small' eigenvalues created by $L^{2}$
solutions of the operators $\mathcal{D}_{1}^{\infty}$ and $\mathcal{D}%
_{2}^{\infty}$ on the half-manifolds with cylindrical ends. We use a
straightforward analysis of small eigenvalues inspired by the proof of Theorem
\ref{t:6.1} to prove the following result

\begin{theorem}
\label{t:0.3} There exist $R_{0}>0$ and positive constants $a_{1},\,a_{2},$
and $a_{3}$, such that for any $R>R_{0}$, the eigenvalue $\mu$ of the operator
$\mathcal{D}^{R}$ is either bounded away from $0$ with $a_{1}<\lvert\mu\rvert
$, or is exponentially small $\lvert\mu\rvert<a_{2}e^{-a_{3}R}$. Let
$\mathcal{W}^{R}$ denote the subspace of $L^{2}(M^{R};S)$ spanned by the
eigensections of $\mathcal{D}^{R}$ corresponding to the exponentially small
eigenvalues. Then $\dim\mathcal{W}^{R}=q$, where $q=\dim(\operatorname{Ker}%
\mathcal{D}_{1}^{\infty})+\dim(\operatorname{Ker}\mathcal{D}_{2}^{\infty})$.
\end{theorem}

Recall from Proposition \ref{p:lemma_6.2} that the operator $\mathcal{D}%
_{j}^{\infty}$, acting on the first Sobolev space $\mathcal{H}^{1}%
(M_{j}^{\infty};S)$, is an (unbounded) self-adjoint Fredholm operator in
$L^{2}(M_{j}^{\infty};S)$ which has a discrete spectrum in the interval
$(-\lambda_{1},+\lambda_{1})$ where $\lambda_{1}$ denotes the smallest
positive eigenvalue of the tangential operator $B$. Thus, the space
$\operatorname{Ker}\mathcal{D}_{j}^{\infty}$ of $L^{2}$ solutions is of finite dimension.

\bigskip

To prove the theorem we first investigate the small eigenvalues of the
operator $\mathcal{D}^{R}$ and the pasting of $L^{2}$ solutions. Let $R>0$. We
reparametrize the normal coordinate $u$ such that $M_{1}^{R}=M_{1}\cup\left(
(-R,0]\times\Sigma\right)  $ and $M_{2}^{R}=\left(  [0,R)\times\Sigma\right)
\cup M_{2}$, and introduce the subspace $\mathcal{V}^{R}\subset L^{2}%
(M^{R};S)$ spanned by $L^{2}$ solutions of the operators $\mathcal{D}%
_{j}^{\infty}$. We choose an auxiliary smooth real function $f^{R}=f_{1}%
^{R}\cup f_{2}^{R}$ on $M^{R}$ with $f^{R}=1$ outside the cylinder
$[-R,R]\times\Sigma$, and where $f^{R}$ is a function of the normal variable
$u$ on the cylinder. Moreover, we assume $f^{R}(-u)=f^{R}(u)$ (i.e.,
$f_{1}^{R}(-u)=f_{2}^{R}(u)$), and that $f_{2}^{R}$ is an increasing function
of $u$ with
\[
f_{2}^{R}(u)=\left\{
\begin{array}
[c]{ll}%
0 & \text{for $0\leq u\leq\frac{R}{4}$}\\
1 & \text{for $\frac{R}{2}\leq u\leq R$}.
\end{array}
\right.  \,
\]
We also assume that there exists a constant $\gamma>0$ such that $\lvert
\frac{\partial^{p}f_{2}^{R}}{\partial u^{p}}(u)\rvert<\gamma R^{-p}$. If
$s_{j}\in C^{\infty}(M_{j}^{\infty};S)$, we define $s_{1}\cup_{f^{R}}s_{2}$ by
the formula
\[
\left(  s_{1}\cup_{f^{R}}s_{2}\right)  (x):=\left\{
\begin{array}
[c]{ll}%
f_{1}^{R}(x)s_{1}(x) & \text{for $x\in M_{1}^{R}$}\\
f_{2}^{R}(x)s_{2}(x) & \text{for $x\in M_{2}^{R}.$}%
\end{array}
\right.
\]
Clearly, we have
\begin{align}
s_{1}\cup_{f^{R}}s_{2}\,  &  =\,s_{1}\cup_{f^{R}}0\,+\,0\cup_{f^{R}}%
s_{2}\nonumber\\
\mathcal{D}^{R}(s_{1}\cup_{f^{R}}s_{2})\,  &  =\,(\mathcal{D}_{1}^{\infty
}s_{1})\cup_{f^{R}}(\mathcal{D}_{2}^{\infty}s_{2})\,+\,s_{1}\cup_{g^{R}}%
s_{2}\text{\ and}\label{e:1.3}\\
\left|  s_{1}\cup_{f^{R}}s_{2}\right|  ^{2}\,  &  =\,\left|  s_{1}\cup_{f^{R}%
}0\right|  ^{2}\,+\,\left|  \,0\cup_{f^{R}}s_{2}\right|  ^{2},\nonumber
\end{align}
where $g^{R}:=g_{1}^{R}\cup g_{2}^{R}$ with $g_{j}^{R}(u,y)=\sigma
(y)\frac{\partial f_{j}^{R}}{\partial u}(u,y)$ and $\left|  \cdot\right|  $
denotes the $L^{2}$ norm on the manifold $M^{R}$.

\begin{definition}
\label{d:VR}The subspace $\mathcal{V}^{R}\subset C^{\infty}(M^{R};S)$ is
defined by
\[
\mathcal{V}^{R}:=\operatorname{span}\{s_{1}\cup_{f^{R}}s_{2}\mid s_{j}%
\in\operatorname{Ker}\mathcal{D}_{j}^{\infty}\}.
\]
\end{definition}

Let $\{s_{1,1},\dots,s_{1,q_{1}}\}$ be a basis of $\operatorname{Ker}%
\mathcal{D}_{1}^{\infty}$ and $\{s_{2,1},\dots,s_{2,q_{2}}\}$ a basis of
$\operatorname{Ker}\mathcal{D}_{2}^{\infty}$. Then the $q=q_{1}+q_{2}$
sections $\{s_{1,\nu_{1}}\cup_{f^{R}}0\}\cup\{0\cup_{f^{R}}s_{2,\nu_{2}}\}$
form a basis of $\mathcal{V}^{R}$. We want to show that $\mathcal{V}^{R}$
approximates the space $\mathcal{W}^{R}$ of eigensections of $\mathcal{D}^{R}
$ corresponding to the `small' eigenvalues, for $R$ sufficiently large. We
begin with an elementary result:

\begin{lemma}
\label{l:R0}There exists $R_{0}$, such that for any $R>R_{0}$ and any
$s\in\mathcal{V}^{R}$, the following estimate holds
\[
\left|  \mathcal{D}^{R}s\right|  \leq e^{-\lambda_{1}R}\left|  s\right|  .
\]
\end{lemma}

\begin{proof}
It suffices to prove the estimate for basis sections of $\mathcal{V}^{R}$.
Thus, let $s=s_{1}\cup_{f^{R}}0$ with $s_{1}\in\operatorname{Ker}%
\mathcal{D}_{1}^{\infty}$. By \eqref{e:1.3} we have
\begin{equation}
\mathcal{D}^{R}s(x)=\left\{
\begin{array}
[c]{ll}%
0 & \text{for $x\in M_{1}\cup M_{2}$}\\
\sigma(y)\frac{\partial f_{1}^{R}}{\partial u}(u,y)\cdot s_{1}(u,y) &
\text{for $x=(u,y)\in\lbrack-R,R]\times\Sigma$}.
\end{array}
\right.  \label{e:1.3'}%
\end{equation}
Here $f_{1}^{R}$ is continued in a trivial way on the whole cylinder
$[-R,R]\times\Sigma$. Now, $s_{1}$ is a $L^{2}$ solution of $\mathcal{D}%
_{1}^{\infty}$, hence $s_{1}(u,y)=\sum_{k}c_{k}e^{-(R+u)\lambda_{k}}%
\varphi_{k}(y)$ on this cylinder where $\{\varphi_{k},\lambda_{k}%
;\sigma\varphi_{k},-\lambda_{k}\}_{k\in\mathbb{N}}$ is, as above, a spectral
resolution of $L^{2}(\Sigma;S)$ for $B$. We estimate the norm of
{$\mathcal{D}$}${^{R}s}$:
\begin{align*}
\left|  {\mathcal{D}^{R}s}\right|  ^{2}  &  =\left|  \frac{\partial f_{1}^{R}%
}{\partial u}\cdot{s}_{{1}}\right|  ^{2}\\
&  =\sum_{k}\int_{-\frac{R}{2}}^{-\frac{R}{4}}\int_{\Sigma}\left(
\frac{\partial f_{1}^{R}}{\partial u}\right)  ^{2}\cdot\lvert c_{k}\rvert
^{2}\cdot e^{-2(R+u)\lambda_{k}}(\varphi_{k}(y);\varphi_{k}(y))\,dy\,du\\
&  =\int_{-\frac{R}{2}}^{-\frac{R}{4}}\left(  \frac{\partial f_{1}^{R}%
}{\partial u}\right)  ^{2}\cdot\sum_{k}\lvert c_{k}\rvert^{2}\cdot
e^{-2(R+u)\lambda_{k}}\cdot1\cdot du\\
&  \leq\frac{\gamma^{2}}{R^{2}}\cdot\sum\nolimits_{k}\left(  \lvert
c_{k}\rvert^{2}\cdot\int_{-\frac{R}{2}}^{-\frac{R}{4}}e^{-2(R+u)\lambda_{k}%
}\,du\right) \\
&  =\frac{\gamma^{2}}{R^{2}}\cdot\sum\nolimits_{k}\left(  \lvert c_{k}%
\rvert^{2}\cdot\int_{R\lambda_{k}}^{\frac{3}{2}R\lambda_{k}}e^{-v}\,\frac
{dv}{2\lambda_{k}}\right) \\
&  \leq\frac{\gamma^{2}}{R^{2}}\cdot\sum\nolimits_{k}\lvert c_{k}\rvert
^{2}\cdot\frac{e^{-R\lambda_{k}}-e^{-\frac{3}{2}R\lambda_{k}}}{2\lambda_{k}}\\
&  \leq\frac{\gamma^{2}}{R^{2}}\cdot\sum\nolimits_{k}\frac{e^{-R\lambda_{k}}%
}{2\lambda_{k}}\lvert c_{k}\rvert^{2}\leq\frac{\gamma^{2}}{R^{2}}%
{e^{-R\lambda_{1}}}\cdot\sum_{k}\frac{\lvert c_{k}\rvert^{2}}{2\lambda_{k}}.
\end{align*}
On the other hand, we have the elementary inequality
\begin{multline*}
\left|  s\right|  ^{2}=\left|  s_{1}\cup_{f^{R}}0\right|  ^{2}\geq\int
_{-R}^{-R+1}\int_{\Sigma}\lvert s_{1}(u,y)\rvert^{2}\,dy\,du\\
=\sum\lvert c_{k}\rvert^{2}\cdot\frac{1-e^{-2\lambda_{k}}}{2\lambda_{k}}\geq
d\cdot\sum\frac{\lvert c_{k}\rvert^{2}}{2\lambda_{k}},
\end{multline*}
with $0<d\leq1-e^{-2\lambda_{1}}$. Thus, we have the following estimate for
any $s\in\mathcal{V}^{R}$ of the form $s_{1}\cup_{f^{R}}0$ and for
sufficiently large $R$
\[
\left|  {\mathcal{D}^{R}s}\right|  ^{2}\leq\frac{\gamma^{2}}{R^{2}%
d}{e^{-R\lambda_{1}}}\cdot d\cdot\sum_{k}\frac{\lvert c_{k}\rvert^{2}%
}{2\lambda_{k}}\leq\frac{\gamma^{2}}{R^{2}d}{e^{-R\lambda_{1}}}\cdot\left|
s\right|  ^{2}\leq e^{-R\lambda_{1}}\cdot\left|  s\right|  ^{2}\,.
\]
For $s=0\cup_{f^{R}}s_{2}$, we estimate the norm of $\mathcal{D}_{R}s$ in the
same way, in view of the fact that $s_{2}$ has the form $s_{2}(u,y)=\sum
_{k}d_{k}e^{(u+R)\lambda_{k}}\sigma(y)\varphi_{k}(y)$ on the cylinder.
\end{proof}

\bigskip

Let $\{\lambda_{k};\psi_{k}\}$ denote a spectral decomposition of the space
$L^{2}(M^{R};S)$ generated by the operator $\mathcal{D}^{R}$. For $a>0$, let
$P_{a}$ denote the orthogonal projection onto the space $\mathcal{H}%
_{a}:=\operatorname{span}\{\psi_{k}\mid\lvert\lambda_{k}\rvert>a\}$.

\begin{lemma}
\label{l:1.3} For sufficiently large $R$, we have the estimate
\[
\left|  \left(  \operatorname{I}-P_{e^{-R\lambda_{1}/4}}\right)  s\right|
\leq e^{-R\lambda_{1}/2}\cdot\left|  s\right|  \text{,\ for all }%
s\in\mathcal{V}^{R}.
\]
\end{lemma}

\begin{proof}
We represent $s$ as the series $s=\sum_{k}a_{k}\psi_{k}$\thinspace. We have
\begin{align*}
\left|  \left(  \operatorname{I}-P_{e^{-R\lambda_{1}/4}}\right)  s\right|   &
=\sum_{\lambda_{k}^{2}>e^{-R\lambda_{1}/2}}a_{k}^{2}\leq\sum_{\lambda_{k}%
^{2}>e^{-R\lambda_{1}/2}}e^{\frac{R\lambda_{1}}{2}}\cdot\lambda_{k}^{2}%
a_{k}^{2}\\
&  \leq\sum\nolimits_{k}e^{R\lambda_{1}/2}\cdot\lambda_{k}^{2}a_{k}%
^{2}=e^{R\lambda_{1}/2}\left|  \mathcal{D}^{R}s\right|  ^{2}\\
&  \leq e^{R\lambda_{1}/2}e^{-R\lambda_{1}}\left|  s\right|  ^{2}%
=e^{-R\lambda_{1}/2}\left|  s\right|  ^{2}.
\end{align*}
\end{proof}

\bigskip

\begin{proposition}
\label{p:specproj}The spectral projection $P_{e^{-R\lambda_{1}/4}}$ restricted
to the subspace $\mathcal{V}^{R}$ is an injection. In particular,
$\mathcal{D}^{R}$ has at least $q$ eigenvalues $\rho$ such that $\left|
\rho\right|  \leq e^{-R\lambda_{1}/4}$, where $q$ is the sum of the dimensions
of the spaces $\operatorname{Ker}\mathcal{D}_{j}^{\infty}$ of $L^{2}$
solutions of the operators $\mathcal{D}_{1}^{\infty}$ and $\mathcal{D}%
_{2}^{\infty}$.
\end{proposition}

\begin{proof}
Let $s\in\mathcal{V}^{R}$, and assume that $P_{e^{-R\lambda_{1}/4}}(s)=0$. We
have
\[
\left|  s\right|  =\left|  \left(  \operatorname{I}-P_{e^{-R\lambda_{1}/4}%
}\right)  s\right|  \leq e^{-\frac{R\lambda_{1}}{2}}\cdot\left|  s\right|
\leq\frac{1}{2}\left|  s\right|  ,
\]
for $R$ sufficiently large.
\end{proof}

\bigskip

The proposition shows that the operator $\mathcal{D}^{R}$ has at least $q$
exponentially small eigenvalues with corresponding eigensections, which we can
approximate by pasting together $L^{2}$ solutions. Now we will show that this
makes the list of eigenvalues approaching 0 as $R\rightarrow+\infty$ complete.

Let $\psi$ be an eigensection of $\mathcal{D}^{R}$ corresponding to an
eigenvalue $\mu$, where $\lvert\mu\rvert<\lambda_{1}$. As in the proof of
Theorem \ref{t:6.1}, we expand $\psi|_{[-R,R]\times\Sigma}$ in terms of a
spectral resolution
\[
\{\varphi_{k},\lambda_{k};\sigma\varphi_{k},-\lambda_{k}\}_{k\in\mathbb{N}}%
\]
of $L^{2}(\Sigma;S)$ generated by $B$:
\[
\psi(u,y)=\sum_{k=1}^{\infty}f_{k}(u)\varphi_{k}(y)+g_{k}(u)\sigma\varphi
_{k}\,,
\]
where the coefficients satisfy the system of ordinary differential equations
of \eqref{e:dgl_s}
\[%
\begin{pmatrix}
f_{k}^{\prime}\\
g_{k}^{\prime}%
\end{pmatrix}
=\mathbf{A}_{k}
\begin{pmatrix}
f_{k}\\
g_{k}%
\end{pmatrix}
\qquad\text{with $\mathbf{A}_{k}:=
\begin{pmatrix}
-\lambda_{k} & \mu\\
-\mu & \lambda_{k}%
\end{pmatrix}
$}\,.
\]
For the eigenvalues $\pm\sqrt{\lambda_{k}^{2}-\mu^{2}}$ of $\mathbf{A}_{k} $
and the eigenvectors
\[%
\begin{pmatrix}
\lambda_{k}+\sqrt{\lambda_{k}^{2}-\mu^{2}}\\
\mu
\end{pmatrix}
\text{\ and}\;
\begin{pmatrix}
\mu\\
\lambda_{k}+\sqrt{\lambda_{k}^{2}-\mu^{2}}%
\end{pmatrix}
,
\]
we get a natural splitting of $\psi(u,y)$ in the form $\psi(u,y)=\psi
_{+}(u,y)+\psi_{-}(u,y)$ with
\begin{align*}
\psi_{+}(u,y)  &  =\sum_{k}a_{k}e^{-\sqrt{\lambda_{k}^{2}-\mu^{2}}u}\left\{
\left(  \lambda_{k}+\sqrt{\lambda_{k}^{2}-\mu^{2}}\right)  \varphi_{k}%
(y)+\mu\sigma(y)\varphi_{k}(y)\right\}  ,\text{ and}\\
\psi_{-}(u,y)  &  =\sum_{k}b_{k}e^{\sqrt{\lambda_{k}^{2}-\mu^{2}}u}\left\{
\mu\varphi_{k}(y)+\left(  \lambda_{k}+\sqrt{\lambda_{k}^{2}-\mu^{2}}\right)
\sigma(y)\varphi_{k}(y)\right\}  .
\end{align*}

Then we have the following estimate of the $L^{2}$ norm of $\psi$ in the $y$
direction on the cylinder:

\begin{lemma}
\label{l:2.1} Assume that $\left|  \psi\right|  =1$. There exist positive
constants $c_{1},c_{2}$ such that $\left|  \psi_{|\left\{  u\right\}
\times\Sigma}\right|  \leq c_{1}e^{-c_{2}R}$ for $-\frac{3}{4}R\leq u\leq
\frac{3}{4}R$.
\end{lemma}

\begin{proof}
We have
\begin{align*}
&  \left|  \psi_{|\left\{  -R+r\right\}  \times\Sigma}\right|  ^{2}\\
&  \leq e^{-2r\sqrt{\lambda_{k}^{2}-\mu^{2}}}\cdot\left|  \sum\nolimits_{k}%
a_{k}e^{-R\sqrt{\lambda_{k}^{2}-\mu^{2}}}\left\{  \left(  \lambda_{k}%
+\sqrt{\lambda_{k}^{2}-\mu^{2}}\right)  f_{k}+\mu\sigma f_{k}\right\}
\right|  ^{2}\\
&  =e^{-2r\sqrt{\lambda_{k}^{2}-\mu^{2}}}\cdot\left|  \psi_{|\left\{
-R\right\}  \times\Sigma}\right|  ^{2}\,.
\end{align*}
In the same way we get
\[
\left|  \psi_{|\left\{  R-r\right\}  \times\Sigma}\right|  ^{2}\leq
e^{-2r\sqrt{\lambda_{k}^{2}-\mu^{2}}}\cdot\left|  \psi_{|\left\{  R\right\}
\times\Sigma}\right|  ^{2}.
\]
Let us observe that, in fact, the argument used here proves that
\begin{align*}
\left|  \psi_{+|\left\{  r\right\}  \times\Sigma}\right|   &  \leq
e^{-(r-s)\sqrt{\lambda_{k}^{2}-\mu^{2}}}\cdot\left|  \psi_{+|\left\{
s\right\}  \times\Sigma}\right|  ,\text{\ and}\\
\left|  \psi_{-|\left\{  s\right\}  \times\Sigma}\right|   &  \leq
e^{-(r-s)\sqrt{\lambda_{k}^{2}-\mu^{2}}}\cdot\left|  \psi_{-|\left\{
r\right\}  \times\Sigma}\right|  ,
\end{align*}
for any $-R<s<r<R$. We also have another elementary inequality
\[
\left|  \psi_{|\left\{  r\right\}  \times\Sigma}\right|  ^{2}\geq\left|
\psi_{+|\left\{  r\right\}  \times\Sigma}\right|  ^{2}-2\cdot\left|
\psi_{+|\left\{  r\right\}  \times\Sigma}\right|  \cdot\left|  \psi
_{-|\left\{  r\right\}  \times\Sigma}\right|  .
\]
This helps estimate the $L^{2}$ norm of $\psi_{\pm}$ in the $y$ direction. We
have
\begin{align*}
\left|  \psi\right|  ^{2}  &  \geq\int_{-R}^{-R+1}\left|  \psi_{|\left\{
u\right\}  \times\Sigma}\right|  ^{2}\,du\\
&  \geq\int_{-R}^{-R+1}\left(  \left|  \psi_{+|\left\{  u\right\}
\times\Sigma}\right|  ^{2}-2\left|  \psi_{+|\left\{  u\right\}  \times\Sigma
}\right|  \,\left|  \psi_{-|\left\{  u\right\}  \times\Sigma}\right|  \right)
\,du\\
&  \geq\left|  \psi_{+|\left\{  -R\right\}  \times\Sigma}\right|  ^{2}\\
&  \quad-2\int_{-R}^{-R+1}\left|  \psi_{+|\left\{  -R\right\}  \times\Sigma
}\right|  \,e^{-2R\sqrt{\lambda_{k}^{2}-\mu^{2}}}\left|  \psi_{-|\left\{
R\right\}  \times\Sigma}\right|  \,du\\
&  \geq\left|  \psi_{+|\left\{  -R\right\}  \times\Sigma}\right|
^{2}-2e^{-2R\sqrt{\lambda_{k}^{2}-\mu^{2}}}\left|  \psi_{+|\left\{
-R\right\}  \times\Sigma}\right|  \left|  \psi_{-|\left\{  R\right\}
\times\Sigma}\right|  .
\end{align*}
In the same way we obtain
\[
\left|  \psi\right|  ^{2}\geq\left|  \psi_{-|\left\{  R\right\}  \times\Sigma
}\right|  ^{2}-2e^{-2R\sqrt{\lambda_{k}^{2}-\mu^{2}}}\left|  \psi_{+|\left\{
-R\right\}  \times\Sigma}\right|  \left|  \psi_{-|\left\{  R\right\}
\times\Sigma}\right|  .
\]
We add the last two inequalities and use
\[
2\left|  \psi_{+|\left\{  -R\right\}  \times\Sigma}\right|  \text{ }\left|
\psi_{-|\left\{  R\right\}  \times\Sigma}\right|  \leq\left|  \psi_{+|\left\{
-R\right\}  \times\Sigma}\right|  ^{2}+\left|  \psi_{-|\left\{  R\right\}
\times\Sigma}\right|  ^{2}%
\]
to obtain
\[
2\left|  \psi\right|  ^{2}\geq\left(  1-e^{-2R\sqrt{\lambda_{k}^{2}-\mu^{2}}%
}\right)  \left(  \left|  \psi_{+|\left\{  -R\right\}  \times\Sigma}\right|
^{2}+\left|  \psi_{-|\left\{  R\right\}  \times\Sigma}\right|  ^{2}\right)  .
\]
This gives us the inequality we need, namely
\[
\left|  \psi_{\pm|\left\{  \mp R\right\}  \times\Sigma}\right|  ^{2}%
\leq4\left|  \psi\right|  ^{2}.
\]
Now we finish the proof of the lemma.
\begin{align*}
\left|  \psi_{|\left\{  u\right\}  \times\Sigma}\right|   &  =\left|
\psi_{+|\left\{  u\right\}  \times\Sigma}+\psi_{-|\left\{  u\right\}
\times\Sigma}\right| \\
&  \leq e^{-(u+R)\sqrt{\lambda_{k}^{2}-\mu^{2}}}\left|  \psi_{+|\left\{
-R\right\}  \times\Sigma}\right|  +e^{-(R-u)\sqrt{\lambda_{k}^{2}-\mu^{2}}%
}\left|  \psi_{-|\left\{  R\right\}  \times\Sigma}\right| \\
&  \leq2\left(  e^{-(u+R)\sqrt{\lambda_{k}^{2}-\mu^{2}}}+e^{-(R-u)\sqrt
{\lambda_{k}^{2}-\mu^{2}}}\right)  \left|  \psi\right|  \leq c_{1}e^{-c_{2}%
R}\,,
\end{align*}
for certain positive constants $c_{1},c_{2}$ when $-\frac{3}{4}R\leq
u\leq\frac{3}{4}R$.
\end{proof}

\medskip We are ready to state the technical main result of this section.

\begin{theorem}
\label{t:2.2} Let $\psi$ denote an eigensection of the operator $\mathcal{D}%
^{R}$ corresponding to an eigenvalue $\mu$, where $\lvert\mu\rvert<\lambda
_{1}$\thinspace. Assume that $\psi$ is orthogonal to the subspace
$P_{e^{-R\lambda_{1}/4}}\mathcal{V}^{R}\subset L^{2}(M^{R};S)$. Then there
exists a positive constant $c$, such that $\lvert\mu\rvert>c$.
\end{theorem}

To prove the theorem we may assume that $\left|  \psi\right|  =1$. We begin
with an elementary consequence of Lemma \ref{l:1.3}.

\begin{lemma}
\label{l:2.3} For any $s\in\mathcal{V}^{R}$ we have
\[
\lvert\left\langle \psi;s\right\rangle \rvert\leq e^{-R\lambda_{1}/2}\left|
s\right|  .
\]
\end{lemma}

\begin{proof}
We have
\begin{align*}
\lvert\left\langle \psi;s\right\rangle \rvert &  =\lvert\left\langle
\psi;P_{e^{-R\lambda_{1}/4}}(s)+\left(  s-P_{e^{-R\lambda_{1}/4}}(s)\right)
\right\rangle \rvert=\lvert\left\langle \psi;P_{e^{-R\lambda_{1}/4}%
}(s)\right\rangle \rvert\\
&  \leq\left|  \psi\right|  \,\left|  P_{e^{-R\lambda_{1}/4}}(s)\right|  \leq
e^{-R\lambda_{1}/2}\left|  s\right|  .
\end{align*}
\end{proof}

We want to compare $\psi$ with the eigensections on the corresponding
manifolds with cylindrical ends. We use $\psi$ to construct a suitable section
on $M_{2}^{\infty}=\left(  (-\infty,R]\times\Sigma\right)  \cup M_{2}$ (Note
the reparametrization compared with the convention chosen in the beginning of
this chapter). Let $h:M_{2}^{\infty}\rightarrow\mathbb{R}$ be a smooth
increasing function such that $h$ is equal to $1$ on $M_{2}$ and $h$ is a
function of the normal variable on the cylinder, equal to 0 for $u\leq\frac
{1}{2}R$, and equal to 1 for $\frac{3}{4}R\leq u$. We also assume, as usual,
that $\lvert\frac{\partial^{p}h}{\partial u^{p}}\rvert\leq\gamma R^{-p}$ for a
certain constant $\gamma>0$. We define
\[
\psi_{2}^{\infty}(x):=\left\{
\begin{array}
[c]{ll}%
h(x)\psi(x) & \text{for $x\in M_{2}^{R}$}\\
0 & \text{for $x\in(-\infty,0]\times\Sigma$.}%
\end{array}
\right.
\]

\begin{proposition}
\label{p:2.4} There exist positive constants $c_{1},c_{2}$, such that
\[
\lvert\left\langle \psi_{2}^{\infty};s\right\rangle \rvert\leq c_{1}%
e^{-c_{2}R}\left|  s\right|
\]
for any $s\in\operatorname{Ker}\mathcal{D}_{2}^{\infty}$.
\end{proposition}

\begin{proof}
For a suitable cutoff function $f_{2}^{R}$ we have
\begin{multline*}
\lvert\left\langle \psi_{2}^{\infty};s\right\rangle \rvert=\left|  \int
_{M_{2}^{\infty}}\left(  \psi_{2}^{\infty}(x);s(x)\right)  dx\right|  =\left|
\int_{M_{2}^{R}}\left(  h(x)\psi(x);f_{2}^{R}(x)s(x)\right)  dx\right| \\
\leq\left|  \int_{M_{2}^{R}}\left(  \psi(x);f_{2}^{R}(x)s(x)\right)
dx\right|  +\left|  \int_{M_{2}^{R}}\left(  (1-h(x))\psi(x);f_{2}%
^{R}(x)s(x)\right)  dx\right|  .
\end{multline*}
We use Lemma \ref{l:2.3} to estimate the first summand:
\begin{align*}
\left|  \int_{M_{2}^{R}}\left(  \psi(x);f_{2}^{R}(x)s(x)\right)  dx\right|
&  =\left|  \int_{M_{2}^{R}}\left(  \psi(x);\left(  0\cup_{f^{R}}s\right)
(x)\right)  dx\right| \\
&  =\lvert\left\langle \psi;0\cup_{f^{R}}s\right\rangle \rvert\leq
e^{-R\lambda_{1}/2}\,\left|  s\right|  .
\end{align*}
We use Lemma \ref{l:2.1} to estimate the second summand:
\begin{align*}
&  \left|  \int_{M_{2}^{R}}\left(  (1-h(x))\psi(x);f_{2}^{R}(x)s(x)\right)
\,dx\right| \\
&  \leq\int_{M_{2}^{R}}\left|  \left(  (1-h(x))\psi(x);f_{2}^{R}%
(x)s(x)\right)  \right|  \text{\thinspace}dx\\
&  \leq\int_{M_{2}^{R}}\left|  \left(  (1-h(x))\psi(x)\right)  \right|
\ \left|  f_{2}^{R}(x)s(x)\right|  \text{\thinspace}dx\\
&  \leq\left(  \int_{M_{2}^{R}}\left|  \left(  (1-h(x))\psi(x)\right)
\right|  ^{2}\,dx\right)  ^{\frac{1}{2}}\left|  s\right| \\
&  \leq\left(  \int_{0}^{\frac{3}{4}R}\left|  \psi_{|\left\{  u\right\}
\times\Sigma}\right|  ^{2}\right)  ^{\frac{1}{2}}\left|  s\right|  \,du\\
&  \leq\left(  c_{1}^{2}e^{-2c_{2}R}\frac{3}{4}R\right)  ^{\frac{1}{2}}\left|
s\right|  \leq c_{3}e^{-c_{4}R}\left|  s\right|  .
\end{align*}
\bigskip
\end{proof}

\begin{proof}
[Proof of Theorem \ref{t:2.2}]Now we estimate $\mu^{2}$ from below by
following the proof of Theorem \ref{t:6.1}. We choose $\left\{  s_{k}\right\}
_{k=1}^{q_{2}}$an orthonormal basis of the kernel of the operator
$\mathcal{D}_{2}^{\infty}$. Let us define
\[
\widetilde{\psi}:=\psi_{2}^{\infty}-\sum_{k=1}^{q_{2}}\left\langle \psi
_{2}^{\infty},s_{k}\right\rangle s_{k}.
\]
Then $\widetilde{\psi}$ is orthogonal to $\operatorname{Ker}\mathcal{D}%
_{2}^{\infty}$, and it follows from Proposition \ref{p:2.4} that
\[
\left|  \widetilde{\psi}\right|  \geq\tfrac{1}{3}\left|  \psi_{2}^{\infty
}\right|  >\kappa>0,
\]
for $R$ large enough, where $\kappa$ is independent of $R$, of the specific
choice of the eigensection $\psi$, and of the cutoff function $h$. Let
$\mu_{1}^{2}$ denote the smallest nonzero eigenvalue of the operator $\left(
\mathcal{D}_{2}^{\infty}\right)  ^{2}$. Once again, it follows from the
Min-Max Principle, that $\left\langle \left(  \mathcal{D}_{2}^{\infty}\right)
^{2}\widetilde{\psi};\widetilde{\psi}\right\rangle \geq\mu^{2}\kappa^{2}$. We
have
\begin{align*}
\mu^{2}  &  =\left\langle \left(  \mathcal{D}^{R}\right)  ^{2}\psi
;\psi\right\rangle \geq\int_{M^{R}}\left|  \mathcal{D}^{R}\psi\left(
x\right)  \right|  ^{2}\,dx\\
&  =\int_{M^{R}}\left|  \mathcal{D}^{R}\left(  h(x)\psi\left(  x\right)
+\left(  1-h(x)\right)  \psi\left(  x\right)  )\right)  \right|  ^{2}\,dx\\
&  \geq\int_{M^{R}}\left|  \mathcal{D}^{R}h(x)\psi\left(  x\right)  \right|
^{2}dx-\int_{M^{R}}\left|  \mathcal{D}^{R}(\left(  1-h(x)\right)  \psi\left(
x\right)  )\right|  ^{2}\,dx
\end{align*}
It is not difficult to estimate the first term from below. We have
\[
\int_{M_{2}^{\infty}}\left|  \left(  \mathcal{D}_{2}^{\infty}\psi_{2}^{\infty
}\right)  \left(  x\right)  \right|  ^{2}\,dx=\left\langle \left(
\mathcal{D}_{2}^{\infty}\right)  ^{2}\psi_{2}^{\infty};\psi_{2}^{\infty
}\right\rangle =\left\langle \left(  \mathcal{D}_{2}^{\infty}\right)
^{2}\widetilde{\psi};\widetilde{\psi}\right\rangle \geq\mu_{1}^{2}\kappa^{2}%
\]
We estimate the second term as follows:
\begin{align*}
&  \int_{M^{R}}\left|  \mathcal{D}^{R}(1-h(x)\psi\left(  x\right)  )\right|
^{2}\,dx=\int_{M_{2}^{R}}\left|  \left(  1-h(x)\right)  \left(  \mathcal{D}%
^{R}\psi\right)  \left(  x\right)  )-\sigma\left(  x\right)  \tfrac{\partial
h}{\partial u}(x)\psi\left(  x\right)  \right|  ^{2}\,dx\\
&  \leq\int_{M_{2}^{R}}\left(  \left|  \mu\left(  1-h(x)\right)  \psi\left(
x\right)  )\right|  ^{2}+2\left|  \mu\left(  1-h(x)\right)  \psi\left(
x\right)  )\right|  +\left|  \sigma\left(  x\right)  \tfrac{\partial
h}{\partial u}(x)\psi\left(  x\right)  \right|  ^{2}\right)  \,dx
\end{align*}
Now we use Lemma \ref{l:2.1} successively to estimate each summand on the
right side by $c_{1}e^{-c_{2}R}$. This gives us
\[
\int_{M_{2}^{R}}\left|  \mathcal{D}^{R}\left(  \left(  1-h\right)
\psi\right)  (x)\right|  ^{2}\,dx\leq c_{3}e^{-c_{4}R},
\]
and finally we have $\mu^{2}\geq\mu_{1}^{2}\kappa^{2}-c_{3}e^{-c_{4}R}%
\geq\tfrac{1}{2}\mu_{1}^{2}\kappa^{2}$ for $R$ large enough.
\end{proof}

Theorem \ref{t:0.3} is an easy consequence of Theorem \ref{t:2.2}.

\subsection{The Additivity for Spectral Boundary Conditions}

\label{s:adiabatic_limit}

We finish the proof of Theorem \ref{t:adiabatic_eta}. We still have to show
equation \eqref{e:7.18}, i.e.

\begin{lemma}
\label{l:etaR}We have $\eta^{R}=\operatorname{O}(e^{-cR})$ $\operatorname{mod}%
\mathbb{Z}$ where
\[
\eta^{R}:=\frac{1}{\sqrt{\pi}}\int_{\sqrt{R}}^{\infty}\frac{1}{\sqrt{t}%
}\operatorname{Tr}\left(  \mathcal{D}^{R}e^{-t\mathcal{D}_{R}^{2}}\right)
\,dt.
\]
\end{lemma}

\begin{proof}
It follows from Theorem \ref{t:0.3} that we have `exponentially small'
eigenvalues corresponding to the eigensections from the subspace
$\mathcal{W}^{R}$ and the eigenvalues $\mu$ bounded away from 0, with
$\lvert\mu\rvert\geq a_{1}$\thinspace, corresponding to the eigensections from
the orthogonal complement of $\mathcal{W}^{R}$\thinspace. First we show that
we can neglect the contribution due to the eigenvalues that are bounded away
from 0. We are precisely in the same situation as with the large $t$
asymptotic of the corresponding integral for the Atiyah--Patodi--Singer
boundary problem on the half manifold with the cylinder attached. Literally,
we can repeat the proof of Lemma \ref{l:lemma.7.1} by replacing $\mathcal{D}%
_{2,P_{>}}^{R}$ by $\mathcal{D}^{R}$ and the uniform bound for the smallest
positive eigenvalue of $\mathcal{D}_{2,P_{>}}^{R}$ by our present bound
$a_{1}$. Thus, we have
\begin{align*}
&  \left|  \int_{\sqrt{R}}^{\infty}\frac{1}{\sqrt{t}}\operatorname{Tr}\left(
\mathcal{D}^{R}e^{-t\mathcal{D}_{R}^{2}}|_{(\mathcal{W}^{R})^{\perp}}\right)
\,dt\right|  \leq\int_{\sqrt{R}}^{\infty}\frac{1}{\sqrt{t}}\left\{
\sum\nolimits_{\lvert\mu\rvert\geq a_{1}}\lvert\mu\rvert e^{-t\mu^{2}%
}\right\}  \,dt\\
&  \leq\int_{\sqrt{R}}^{\infty}\frac{1}{\sqrt{t}}\left\{  \sum
\nolimits_{\lvert\mu\rvert\geq a_{1}}e^{-(t-1)\mu^{2}}\right\}  \,dt\leq
\int_{\sqrt{R}}^{\infty}\frac{1}{\sqrt{t}}\left\{  \sum\nolimits_{\lvert
\mu\rvert\geq a_{1}}e^{-\mu^{2}}\right\}  e^{-(t-2)a_{1}^{2}}\,dt\\
&  \leq e^{2a_{1}^{2}}\operatorname{Tr}\left(  e^{-t\mathcal{D}_{R}^{2}%
}\right)  \int_{\sqrt{R}}^{\infty}\frac{1}{\sqrt{t}}e^{-ta_{1}^{2}%
}\,dt=e^{2a_{1}^{2}}\operatorname{Tr}\left(  e^{-t\mathcal{D}_{R}^{2}}\right)
\frac{1}{a_{1}}\int_{\sqrt{R}}^{\infty}\frac{1}{\sqrt{t}}e^{-ta_{1}^{2}%
}\,a_{1}\,dt\\
&  \leq\frac{e^{2a_{1}^{2}}}{2a_{1}}\operatorname{Tr}\left(  e^{-t\mathcal{D}%
_{R}^{2}}\right)  e^{-a_{1}^{2}\sqrt{R}}\,.
\end{align*}
For the last inequality see \eqref{e:gromov?}. A standard estimate on the heat
kernel of the operator $\mathcal{D}^{R}$ gives (as in \eqref{e:vol}) the
inequality $\operatorname{Tr}\left(  e^{-t\mathcal{D}_{R}^{2}}\right)  \leq
b_{3}\cdot\operatorname*{Vol}(M^{R})\leq b_{4}R$, which implies that
\begin{equation}
\left|  \int_{\sqrt{R}}^{\infty}\frac{1}{\sqrt{t}}\operatorname{Tr}\left(
\mathcal{D}^{R}e^{-t\mathcal{D}_{R}^{2}}|_{(\mathcal{W}^{R})^{\perp}}\right)
\,dt\right|  \leq b_{5}e^{-b_{6}\sqrt{R}}\,. \label{e:3.10}%
\end{equation}
That proves that the contribution from the large eigenvalues disappears as
$R\rightarrow\infty$. The essential part of $\eta^{R}$ comes from the subspace
$\mathcal{W}^{R}$:
\begin{align}
\frac{1}{\sqrt{\pi}}\int_{\sqrt{R}}^{\infty}\frac{1}{\sqrt{t}}%
\operatorname{Tr}\left(  \mathcal{D}^{R}e^{-t\mathcal{D}_{R}^{2}%
}|_{\mathcal{W}^{R}}\right)  \,dt  &  =\nonumber\\
\sum_{\lvert\mu\rvert<a_{1}}\frac{1}{\sqrt{\pi}}\int_{\sqrt{R}}^{\infty}%
\frac{1}{\sqrt{t}}\mu e^{-t\mu^{2}}\,dt  &  =\sum_{\lvert\mu\rvert<a_{1}%
}\operatorname*{sign}(\mu)\frac{2}{\sqrt{\pi}}\int_{\left|  \mu\right|
R^{1/4}}^{\infty}e^{-v^{2}}\,dv. \label{e:3.11}%
\end{align}
It follows from Theorem \ref{t:0.3} that $\lim_{R\rightarrow\infty}\left|
\mu\right|  R^{1/4}=0$. Thus, the right side of (\ref{e:3.11}) is equal to
\begin{equation}
\operatorname*{sign}\nolimits_{R}(\mathcal{D}):=\sum_{\lvert\mu\rvert<a_{1}%
}\operatorname*{sign}(\mu) \label{d:signRD}%
\end{equation}
plus the smooth error term which is rapidly decreasing as $R\rightarrow\infty$.
\end{proof}

\bigskip

Thus we have proved Theorem \ref{t:adiabatic_eta}. In particular, we have
proved
\[
\lim_{R\rightarrow\infty}\eta_{\mathcal{D}^{R}}(0)\equiv\lim_{R\rightarrow
\infty}\left\{  \eta_{\mathcal{D}_{1,P_{<}}^{R}}(0)+\eta_{\mathcal{D}%
_{2,P_{>}}^{R}}(0)\right\}  \operatorname{mod}\mathbb{Z}.
\]
To establish the true additivity assertion of Corollary \ref{c:0.1}, we show
that the preceding $\eta$ invariants do not depend on $R$ modulo integers.

\begin{proposition}
\emph{(W. M\"{u}ller.)}\label{p:muller} The eta invariant $\eta_{\mathcal{D}%
_{2,P_{>}}^{R}}(0)\in\mathbb{R}/\mathbb{Z}$ is independent of the cylinder
length $R$.
\end{proposition}

\begin{proof}
Near to the boundary of $M_{2}^{R}$ we parametrize the normal coordinate
$u\in\lbrack-R,1)$ with the boundary at $u=-R$. First we show that
$\dim\operatorname{Ker}\mathcal{D}_{2,P_{>}}^{R}$ is independent of $R$. Let
$s\in\operatorname{Ker}\mathcal{D}_{2,P_{>}}^{R}$\thinspace, namely
\begin{equation}
s\in C^{\infty}(M_{2}^{R};S),\quad\mathcal{D}_{2}^{R}s=0,\text{ and}%
\;P_{>}(s|_{\{-R\}\times\Sigma})=0. \label{e:2.17}%
\end{equation}
As in equation \eqref{e:gb_expansion} (and in equation \eqref{e:gb_expansion'}
of the proof of Theorem \ref{t:6.1}) we may expand $s|_{[-R,0]\times\Sigma}$
in terms of the eigensections of the tangential operator $B$:
\[
s(u,y)=\sum_{k=1}^{\infty}e^{\lambda_{k}u}\sigma(y)\varphi_{k}(y).
\]
Let $R^{\prime}>R$. Then $s$ can be extended in the obvious way to
$\widetilde{s}\in\operatorname{Ker}\mathcal{D}_{2,P_{>}}^{R^{\prime}}$, and
the map $s\mapsto\widetilde{s}$ defines an isomorphism of $\operatorname{Ker}%
\mathcal{D}_{2,P_{>}}^{R}$ onto $\operatorname{Ker}\mathcal{D}_{2,P_{>}%
}^{R^{\prime}}$\thinspace. Next, observe that there exists a smooth family of
diffeomorphisms $f_{R}:[0,1)\rightarrow\lbrack-R,1)$ which have the following
cutoff properties
\[
f_{R}(u)=\left\{
\begin{array}
[c]{ll}%
u & \text{for $\frac{2}{3}<u<1$}\\
u+R & \text{for $0\leq u<\frac{1}{3}$}.
\end{array}
\right.
\]
Let $\psi_{R}:[0,1)\times\Sigma\rightarrow\lbrack-R,1)\times\Sigma$ be defined
by $\psi_{R}(u,y):=(f_{R}(u),y)$, and extend $\psi_{R}$ to a diffeomorphism
$\psi_{R}:M_{2}\rightarrow M_{2}^{R}$ in the canonical way, i.e., $\psi_{R}$
becomes the identity on $M_{2}\setminus\left(  (0,1)\times\Sigma\right)  $.
There is also a bundle isomorphism which covers $\psi_{R}$. This induces an
isomorphism $\psi_{R}^{\ast}:C^{\infty}(M_{2}^{R};S)\rightarrow C^{\infty
}(M_{2};S)$. Let $\widetilde{\mathcal{D}_{2}^{R}}:=\psi_{R}^{\ast}%
\circ\mathcal{D}_{2}^{R}\circ\left(  \psi_{R}^{\ast}\right)  ^{-1}$\thinspace.
Then $\{\widetilde{\mathcal{D}_{2}^{R}}\}_{R}$ is a family of Dirac operators
on $M_{2}$, and $\widetilde{\mathcal{D}_{2}^{R}}=\sigma(\partial_{u}+B)$ near
$\Sigma$. We pick the self-adjoint $L^{2}$ extension defined by
$\operatorname{Dom}\widetilde{\mathcal{D}_{2,P_{>}}^{R}}:=\psi_{R}^{\ast
}\left(  \operatorname{Dom}\mathcal{D}_{2,P_{>}}^{R}\right)  $. Hence,
\[
\eta_{\widetilde{\mathcal{D}_{2,P_{>}}^{R}}}(s)=\eta_{\mathcal{D}_{2,P_{>}%
}^{R}}(s)\text{\ \ and}\;\;\operatorname{Ker}\widetilde{\mathcal{D}_{2,P_{>}%
}^{R}}=\psi_{R}^{\ast}\operatorname{Ker}\mathcal{D}_{2,P_{>}}^{R}\,.
\]
In particular, $\dim\widetilde{\mathcal{D}_{2,P_{>}}^{R}}$ is constant, and we
apply variational calculus to get
\[
\frac{d}{dR}\left(  \eta_{\mathcal{D}_{2,P_{>}}^{R}}(0)\right)  =-\frac
{2}{\sqrt{\pi}}c_{m}(R),
\]
where $c_{m}(R)$ is the coefficient of $t^{-1/2}$ in the asymptotic expansion
of
\[
\operatorname{Tr}\dot{A}_{R}e^{-tA_{R}^{2}}\sim\sum_{j=0}^{\infty}%
c_{j}(R)t^{(j-m-1)/2}%
\]
with $A_{R}:=\widetilde{\mathcal{D}_{2,P_{>}}^{R}}$ and $m:=\dim M$. Now let
$S_{R}^{1}$ denote the circle of radius $2R$. We lift the Clifford bundle from
$\Sigma$ to the torus $\mathbf{T}_{R}:=S_{R}^{1}\times\Sigma$. We define the
action of $\widehat{D^{R}}:C^{\infty}(\mathbf{T}_{R},S)\rightarrow C^{\infty
}(\mathbf{T}_{R},S)$ by $\widehat{D^{R}}=\sigma(\partial_{u}+B)$. Since
$c_{m}(R)$ is locally computable, it follows in the same way as above that
\[
\frac{d}{dR}\left(  \eta_{\widehat{D^{R}}}(0)\right)  =-\frac{2}{\sqrt{\pi}%
}c_{m}(R).
\]
But a direct computation shows that the spectrum of $\widehat{D^{R}}$ is
symmetric. Hence $\eta_{\widehat{D^{R}}}(s)=0$ and, therefore, $c_{m}(R)=0$.
\end{proof}

\medskip In the same way we show that $\eta_{\mathcal{D}^{R}}$ is independent
of $R$. This proves the additivity assertion of Corollary \ref{c:0.1}. In
fact, we have proved a little bit more:

\begin{theorem}
\label{t:0.5} The following formula holds for $R$ large enough
\[
\eta_{\mathcal{D}}(0)=\eta_{\mathcal{D}_{1,P_{<}}}(0)+\eta_{\mathcal{D}%
_{2,P_{>}}}(0)+\operatorname*{sign}\nolimits_{R}(\mathcal{D}),
\]
where $\operatorname*{sign}\nolimits_{R}(\mathcal{D}):=\sum_{\lvert\mu
\rvert<a_{1}}\operatorname*{sign}(\mu)$; see (\ref{d:signRD}).
\end{theorem}

Theorem \ref{t:0.5} has an immediate corollary which describes the case in
which our additivity formula holds in $\mathbb{R}$, not just in $\mathbb{R}%
/\mathbb{Z}$.

\begin{corollary}
If $\operatorname{Ker}\mathcal{D}_{1}^{\infty}=\{0\}=\operatorname{Ker}%
\mathcal{D}_{2}^{\infty}$, then
\[
\eta_{\mathcal{D}}(0)=\eta_{\mathcal{D}_{1,P_{<}}}(0)+\eta_{\mathcal{D}%
_{2,P_{>}}}(0).
\]
\end{corollary}

%
%
%
%
%
%
%
%
%
%
%
%

\bigskip
\end{document}